\documentstyle[amscd,verbatim,leqno,12pt]{amsart}

\numberwithin{equation}{section}

\theoremstyle{plain}
\newtheorem{thm}{Theorem}[section]
\newtheorem{defn}[thm]{Definition}
\newtheorem{prop}[thm]{Proposition}
\newtheorem{lem}[thm]{Lemma}
\newtheorem{cor}[thm]{Corollary}
\newtheorem{conj}[thm]{Conjecture}

\theoremstyle{definition}
\newtheorem{rem}[thm]{Remark}

\renewcommand{\b}{\bullet}
\newcommand{\beast}{\begin{eqnarray*}}
\newcommand{\east}{\end{eqnarray*}}

\newcommand{\N}{{\Bbb N}}
\newcommand{\Z}{{\Bbb Z}}
\newcommand{\Q}{{\Bbb Q}}

\newcommand{\F}{{\Bbb F}}

\renewcommand{\Pr}{{\Bbb P}}

\newcommand{\A}{{\bold A}}
\renewcommand{\a}{{\bold a}}

\newcommand{\Spec}{{\mathrm{Spec}}\,}

\newcommand{\Spf}{{\mathrm{Spf}}\,}
\newcommand{\Spm}{{\mathrm{Spm}}\,}
\newcommand{\Zar}{{\mathrm{Zar}}}
\newcommand{\lra}{\longrightarrow}
\newcommand{\ra}{\rightarrow}
\newcommand{\hra}{\hookrightarrow}

\newcommand{\lla}{\longleftarrow}

\newcommand{\isom}{\overset{\sim}{=}}
\newcommand{\ti}[1]{\widetilde{#1}}
\newcommand{\wt}[1]{\widetilde{#1}}

\newcommand{\ul}[1]{\underline{#1}}
\newcommand{\ol}[1]{\overline{#1}}
\newcommand{\os}{\overset}

\newcommand{\DR}{{\mathrm{DR}}}
\newcommand{\et}{{\mathrm{et}}}

\newcommand{\crys}{{\mathrm{crys}}}

\newcommand{\conv}{{\mathrm{conv}}}

\newcommand{\Hom}{{\mathrm{Hom}}}

\newcommand{\Ker}{{\mathrm{Ker}}}

\newcommand{\triv}{{\mathrm{triv}}}

\newcommand{\id}{{\mathrm{id}}}

\newcommand{\res}{{\mathrm{res}}}

\newcommand{\Coh}{{\mathrm{Coh}}}
\newcommand{\Strat}{{\mathrm{Str}}}

\newcommand{\red}{{\mathrm{red}}}

\newcommand{\HPDI}{{\mathrm{HPDI}}}

\newcommand{\an}{{\mathrm{an}}}

\newcommand{\cA}{{\cal A}}
\newcommand{\cB}{{\cal B}}
\newcommand{\cC}{{\cal C}}
\newcommand{\cD}{{\cal D}}
\newcommand{\cE}{{\cal E}}
\newcommand{\cF}{{\cal F}}

\newcommand{\cI}{{\cal I}}
\newcommand{\cJ}{{\cal J}}
\newcommand{\cO}{{\cal O}}
\newcommand{\cP}{{\cal P}}
\newcommand{\cQ}{{\cal Q}}

\newcommand{\cS}{{\cal S}}
\newcommand{\cT}{{\cal T}}

\newcommand{\cV}{{\cal V}}
\newcommand{\cX}{{\cal X}}
\newcommand{\cY}{{\cal Y}}
\newcommand{\cZ}{{\cal Z}}

\newcommand{\fP}{{\frak P}}

\newcommand{\rk}{{\mathrm{rk}}\,}
\newcommand{\dti}{\breve}

\renewcommand{\sp}{{\mathrm{sp}}}
\newcommand{\rig}{{\mathrm{rig}}}
\newcommand{\logrig}{{\mathrm{log}\text{-}\rig}}
\newcommand{\ex}{{\mathrm{ex}}}

\newcommand{\cosk}{{\mathrm{cosk}}}
\newcommand{\vpl}{\varprojlim}

\newcommand{\fl}{{\mathrm{fl}}}

\renewcommand{\wt}{\widetilde}

\newcommand{\LF}{{({\mathrm{LFS}}/\cB)}}
\newcommand{\pLF}{{(p{\mathrm{LFS}}/\cB)}}
\renewcommand{\L}{{({\mathrm{LS}}/\cB)}}
\newcommand{\LB}{{({\mathrm{LS}}/B)}}
\renewcommand{\res}{{\mathrm{res}}}
\renewcommand{\d}{\dagger}

\newcommand{\lam}{\lambda}
\newcommand{\lamm}{\lambda\text{-}}
\newcommand{\mum}{\mu\text{-}}
\renewcommand{\lor}{\langle 0 \rangle}
\newcommand{\lr}{\langle \b \rangle}
\newcommand{\lmr}{\langle m \rangle}

\newcommand{\dsucc}{\succ \succ}
\newcommand{\dint}{\displaystyle\int}
\newcommand{\dvil}{\displaystyle\varinjlim}
\newcommand{\dvpl}{\displaystyle\varprojlim}

\newcommand{\FI}{F\text{-}I}

\begin{document}

\title[Relative Log Convergent Cohomology II]
{Relative Log Convergent Cohomology and Relative Rigid Cohomology II}

\author{Atsushi Shiho}
\address{Graduate School of Mathematical Sciences, 
University of Tokyo, 3-8-1 Komaba, 
Meguro-ku, Tokyo 153-8914, JAPAN.}
                          

\begin{abstract}
In this paper, we develop the theory of relative log convergent 
cohomology of radius $\lambda$ ($0 < \lambda \leq 1$), which is 
a generalization of the notion of relative log convergent cohomology 
in the previous paper. By comparing this cohomology with 
relative log crystalline cohomology, relative rigid cohomology 
and its variants and by 
using some technique of hypercovering, we prove a version of 
Berthelot's conjecture on the overconvergence of 
relative rigid cohomology for proper smooth families. 
\end{abstract}

\maketitle

\tableofcontents


\section*{Introduction}

Let $k$ be a perfect field of characteristic $p>0$, 
let $V$ be a complete discrete valuation ring of mixed characteristic 
with residue field $k$ and let $X$ be a scheme separated of finite type
over $k$. Then it is expected that the correct $p$-adic analogue on $X$ 
of the notion of local systems in complex analytic case or that of 
smooth $l$-adic sheaves in $l$-adic case is the notion of overconvergent 
$F$-isocrystals on $X$, which is defined by Berthelot (\cite{berthelot1}, 
\cite{berthelot2}). So, when we are given a proper smooth morphism 
$f:X \lra Y$ and an overconvergent $F$-isocrystal $\cE$ on $X$, 
the higher direct image of $\cE$ ($=$ relative rigid cohomology of 
$X/Y$ with coefficient $\cE$, \cite{berthelot1}, \cite{chts}) should 
have the canonical structure of an overconvergent $F$-isocrystal. 
This is the content of Berthelot's conjecture on the coherence and the 
overconvergence of relative rigid cohomology. 
(As we noted in \cite[\S 5]{shiho3} and we will repeat again in this 
introduction, one can consider several versions of this conjecture. 
See also \cite[(4.3)]{berthelot1}, \cite[\S 4]{tsuzuki3}.) 
The purpose of this paper is to prove a version of Berthelot's 
conjecture. \par 
Let me introduce some terminologies and give some explanation on 
Berthelot's conjecture. In this introduction, 
a pair $(X,\ol{X})$ is a pair of schemes 
$X,\ol{X}$ separated of finite type over $k$ endowed with 
an open immersion $X \hra \ol{X}$. 
A triple $(X,\ol{X},\cX)$ is a pair $(X,\ol{X})$ endowed with 
a $p$-adic formal scheme $\cX$ separated, topologically 
of finite type over $\Spf V$ and a closed immersion 
$\ol{X} \hra \cX$ over $\Spf V$. 
For a triple of the form $(S,S,\cS)$ and a pair $(X,\ol{X})$ over 
$(S,S)$, an $(X,\ol{X})$-triple over $(S,S,\cS)$ is a triple 
$(Y,\ol{Y},\cY)$ over $(S,S,\cS)$ endowed with a morphism 
of pairs $(Y,\ol{Y}) \lra (X,\ol{X})$ over $(S,S)$. 
Then Berthelot's conjecture is described as follows 
(\cite[5.3]{shiho3}): 

\begin{conj}\label{bconj1}
Let us take a triple of the form 
$(S,S,\cS)$ $($endowed with a lift of 
Frobenius on $\cS)$
and assume we are given a diagram of pairs 
$$ (X,\ol{X}) \os{f}{\lra} (Y,\ol{Y}) \os{g}{\lra} (S,S) $$ 
such that $f:\ol{X} \lra \ol{Y}$ is proper, $f^{-1}(Y)=X$ 
and that $f|_X:X \lra Y$ is smooth. Then, for an overconvergent 
$(F$-$)$isocrystal $\cE$ on $(X,\ol{X})/\cS_K$ and $q \geq 0$, 
there exists uniquely an overconvergent $(F$-$)$isocrystal $\cF$ 
$($which is called the $q$-th rigid cohomology overconvergent 
isocrystal$)$ on 
$(Y,\ol{Y})/\cS_K$ satisfying the following condition$:$ 
For any $(Y,\ol{Y})$-triple $(Z,\ol{Z},\cZ)$ over $(S,S,\cS)$ 
with $\cZ$ formally smooth over $\cS$ on a neighborhood of $Z$, 
the restriction of $\cF$ to $I^{\d}((Z,\ol{Z})/\cS_K, \cZ)$ is given 
functorially by 
$(R^qf_{(X \times_Y Z,\ol{X} \times_{\ol{Y}} \ol{Z})/\cZ, \rig *}\cE, 
\epsilon)$, where $\epsilon$ is given by 
\begin{align*}
p_2^*R^qf_{(X \times_Y Z,\ol{X} \times_{\ol{Y}} \ol{Z})/\cZ, \rig *}\cE 
& \os{\simeq}{\rightarrow} 
R^qf_{(X \times_Y Z,\ol{X} \times_{\ol{Y}} \ol{Z})/\cZ \times_{\cS} \cZ, 
\rig *}\cE \\ 
& \os{\simeq}{\leftarrow} 
p_1^*R^qf_{(X \times_Y Z,\ol{X} \times_{\ol{Y}} \ol{Z})/\cZ, \rig *}\cE. 
\end{align*} 
$($Here $I^{\d}((Z,\ol{Z})/\cS_K, \cZ)$ denotes the category of 
overconvergent isocrystals on $(Z,\ol{Z})\allowbreak / \allowbreak \cS_K$ 
over $\cZ$, 
$R^qf_{(X \times_Y Z,\ol{X} \times_{\ol{Y}} \ol{Z})/\cZ, \rig *}\cE, 
R^qf_{(X \times_Y Z,\ol{X} \times_{\ol{Y}} \ol{Z})/\cZ \times_{\cS} \cZ, 
\rig *}\cE$ are relative rigid cohomologies and 
$p_i$ is the morphism 
$]\ol{Z}[_{\cZ \times_{\cS} \cZ} \lra \,]\ol{Z}[_{\cZ}$ induced by the 
$i$-th projection.$)$ 
\end{conj} 

We would like to 
consider the following version of Berthelot's conjecture 
(\cite[5.5]{shiho3}), which is 
slightly weaker than Conjecture \ref{bconj1} but strong enough to assure 
the unique existence of `the $q$-th rigid cohomology overconvergent 
isocrystal': 

\begin{conj}\label{bconj2}
Let us take a triple of the form 
$(S,S,\cS)$ 
$($endowed with a lift of 
Frobenius on $\cS)$
and assume we are given a diagram of pairs 
$$ (X,\ol{X}) \os{f}{\lra} (Y,\ol{Y}) \os{g}{\lra} (S,S) $$ 
such that $f:\ol{X} \lra \ol{Y}$ is proper, $f^{-1}(Y)=X$ 
and that $f|_X:X \lra Y$ is smooth. Then 
there exists a subcategory $\cC$ of the category of 
$(Y,\ol{Y})$-triples over $(S,S,\cS)$ such that, 
for an overconvergent 
$(F$-$)$isocrystal $\cE$ on $(X,\ol{X})/\cS_K$ and $q \geq 0$, 
there exists uniquely an overconvergent $(F$-$)$isocrystal $\cF$ 
$($which is called the $q$-th rigid cohomology overconvergent 
isocrystal$)$ on 
$(Y,\ol{Y})/\cS_K$ satisfying the following condition$:$ 
For any $(Z,\ol{Z},\cZ) \in \cC$ 
with $\cZ$ formally smooth over $\cS$ on a neighborhood of $Z$, 
the restriction of $\cF$ to $I^{\d}((Z,\ol{Z})/\cS_K, \cZ)$ is given 
functorially by 
$(R^qf_{(X \times_Y Z,\ol{X} \times_{\ol{Y}} \ol{Z})/\cZ, \rig *}\cE, 
\epsilon)$, where $\epsilon$ is given by 
\begin{align*}
p_2^*R^qf_{(X \times_Y Z,\ol{X} \times_{\ol{Y}} \ol{Z})/\cZ, \rig *}\cE 
& \os{\simeq}{\rightarrow} 
R^qf_{(X \times_Y Z,\ol{X} \times_{\ol{Y}} \ol{Z})/\cZ \times_{\cS} \cZ, 
\rig *}\cE \\ 
& \os{\simeq}{\leftarrow} 
p_1^*R^qf_{(X \times_Y Z,\ol{X} \times_{\ol{Y}} \ol{Z})/\cZ, \rig *}\cE. 
\end{align*}
$($Here the notations are as in Conjecture \ref{bconj1}.$)$
\end{conj} 

The difference of Conjecture \ref{bconj2} compared to Conjecture \ref{bconj1} 
is the introduction of the category $\cC$: Conjecture \ref{bconj2} is weaker 
than Conjecture \ref{bconj1} in the sense that $\cC$ need not be equal to 
the category of all $(Y,\ol{Y})$-triples $(Z,\ol{Z},\cZ)$ over 
$(S,S,\cS)$. However, 
Conjecture \ref{bconj2} is strong enough to require the unique existence 
of $\cF$. As long as we know, there is no non-trivial example for which 
Conjecture \ref{bconj1} is known. On the other hand, Conjecture \ref{bconj2} 
is known in liftable case (\cite[Thm 5]{berthelot1}, 
\cite[4.1.4]{tsuzuki3}, see also \cite[5.7]{shiho3}) and the case where 
the morphism $f:\ol{X} \lra \ol{Y}$ admits a `nice' log structure and 
$\cE$ comes from a locally free isocrystal on log convergent site 
$(\ol{X}/\cS)^{\log}_{\conv}$. (For precise statement, see 
\cite[5.14, 5.15]{shiho3}.) The main result in this paper is the 
following theorems ($=$ Theorems \ref{theorem2}, 
\ref{theorem5} in this paper). 

\begin{thm}\label{thm1}
In the situation of Conjecture \ref{bconj2} 
$($without a lift of Frobenius on $\cS)$, assume moreover that 
$\cS = \Spf V$ and $Y$ is smooth over $k$. 
Then Conjecture \ref{bconj2} 
$($the version without Frobenius structure$)$ is true. 
\end{thm} 

\begin{thm}\label{thm2}
In the situation of Conjecture \ref{bconj2} 
$($with a lift of Frobenius on $\cS)$, 
assume moreover that 
$\cS = \Spf V$. Then Conjecture \ref{bconj2} 
$($the version with Frobenius structure$)$ is true. 
\end{thm} 

In the course of the proof of or as variants of the above theorems, 
we also prove the following results: \\
\quad \\
(a)\,\,\, If $f: \ol{X} \lra \ol{Y}$ admits a nice log structure and 
$\cE$ is locally free, then Conjecture \ref{bconj2} 
(the version without Frobenius structure) is true, without any assumption 
on $Y,\cS$ (Theorem \ref{thm00}). Moreover, if $S=\Spf V$, 
the version with Frobenius structure is also true 
(Theorem \ref{theorem4}). \\
(b)\,\,\, Even if we do not assume the smoothness of $f$, 
the analogue of Conjecture \ref{bconj2} (the version without 
Frobenius structure) is true for locally free $\cE$, if we shrink 
$Y$ (Theorem \ref{theorem1}). Moreover, if $S=\Spf V$, 
the version with Frobenius structure is also true 
(Theorem \ref{theorem4}). \\
(c)\,\,\, Theorem \ref{thm1} is true without the assumption of smoothness on 
$Y$ if we assume the full-faithfulness conjecture of Tsuzuki 
(\cite[1.2.1]{tsuzuki}, see also Conjecture \ref{tconj}) 
(Theorem \ref{theorem3}). 
\\

Now we give an outline of the proof of Theorems \ref{thm1}, \ref{thm2} and 
explain the content of this paper. The key to the proof is the comparison 
of the following six cohomology theories: 
\begin{enumerate} 
\item 
Relative log crystalline cohomology \,\,(\cite{kkato}, see also 
\cite[\S 1]{shiho3}). 
\item 
Relative log convergent cohomology of radius $\lam$ \,\,($\lam \in (0,1]$).
\item 
Relative log analytic cohomology of radius $\lam$ \,\,($\lam \in (0,1]$).
\item 
$\lam$-restriction of relative log rigid cohomology \,\,($\lam \in (0,1)$). 
\item 
Relative log rigid cohomology. 
\item 
Relative rigid cohomology. 
\end{enumerate} 
First, in Section 1, we introduce the notion of relative log convergent 
cohomology of radius $\lam$ and prove the comparison theorem with 
relative log crystalline cohomology when $\lam$ is sufficiently close 
to $1$. In the case $\lam=1$, relative log convergent cohomology of 
radius $\lam$ is the same as relative log convergent cohomology which was 
introduced in \cite{shiho3}, \cite{nakkshiho}. In Section 2, 
by using the results in Section 1 and the technique developped in 
\cite[\S 3]{shiho3}, we prove the finiteness and the base change property 
of relative log convergent cohomology of radius $\lam$ for 
proper log smooth integral morphism having log smooth parameter 
(for definition, see Definition \ref{parameter} or 
\cite[3.4]{shiho3}) when $\lam$ is sufficiently close to $1$. 
In Section 3, we introduce the notion of relative log analytic cohomology 
of radius $\lam$ and prove a relation of it with relative log convergent 
cohomology of radius $\lam$. Using this, we prove a kind of 
the coherence and the convergence property of relative log analytic cohomology 
of radius $\lam$ for 
proper log smooth integral morphism having log smooth parameter. 
In Section 4, first we introduce the notion of a log pair and 
an overconvergent isocrystal on it. Then, for a morphism of log pairs and 
an overconvergent isocrystal on the source, we define the notion of 
relative log rigid cohomology. Then we relate it to the relative rigid 
cohomology (defined in \cite{berthelot1}, \cite{chts}) in certain case. 
After that, we introduce the notion of $\lam$-restriction of relative log 
rigid cohomology and compare it with relative rigid cohomology when 
$\lam$ is sufficiently close to $1$. In Section 5, we prove the 
assertion (a) (the version without Frobenius structure) 
above: By the results proved up to Section 4, 
relative log analytic cohomology of radius $\lam$ has certain 
coherence property and relative rigid cohomology is related to the 
$\lam$-restriction of relative log rigid cohomology. So we compare 
the relative log analytic cohomology of radius $\lam$ and 
the $\lam$-restriction of relative log rigid cohomology. 
To deduce the assertion (a), the base change theorem of relative 
rigid cohomology proved by Tsuzuki \cite[2.3.1]{tsuzuki3} is also 
important. In Section 6, we prove a result of altering a given 
proper morphism of schemes over $k$ to a certain simplicial morphism 
of schemes whose components admit nice log structures, by 
using results of de Jong \cite{dejong}, \cite{dejong2}. This result allows us 
to reduce the proof of the main theorems to the case of 
proper log smooth integral morphisms having log smooth parameter. 
Indeed, in Section 7, first we prove the assertion (b) 
(the version without Frobenius structure) by combining 
the results in Sections 5 and 6. Then, by combining this with 
(a special case of) the main result of \cite{shiho3} and 
a result of Kedlaya in \cite{kedlaya2} on the extendability of 
overconvergent isocrystals, we prove Theorem \ref{thm1}. 
Finally, by working slightly more, we prove the assertion (c), 
the assertions (a), (b) (the version with Frobenius structure) and 
Theorem \ref{thm2}. \par 
The author would like to thank to 
Nobuo Tsuzuki and Makoto Matsumoto 
for inviting him to give a talk at a conference 
held at Hiroshima University, and to Takeshi Tsuji for 
encouraging him to give a talk at a conference held at Tambara 
Institute for Mathematical Sciences. 
The author would like to thank to Kiran Kedlaya 
for kindly sending the preprint \cite{kedlaya-val}.
The author is partly supported by Grant-in-Aid for Young Scientists (B), 
the Ministry of Education, Culture, Sports, Science and Technology of 
Japan and JSPS Core-to-Core program 18005 whose representative is 
Makoto Matsumoto.

\section*{Convention}
\hspace{-16pt}
(1) \,\, Throughout this paper, $k$ is a field of characteristic 
$p>0$, $W$ is a fixed Cohen ring of $k$ and $K$ is the fraction field 
of $W$. As we will remark in the text, we will assume that $k$ is perfect 
in Sections 6 and 7. 
We fix a $p$-adic fine log formal scheme 
$(\cB,M_{\cB})$ separated and topologically of finite type 
over $\Spf W$ as a base log formal scheme. (Note that $\cB$ is 
Noetherian.) 
We denote the 
scheme $(\cB,M_{\cB}) \otimes_{\Z_p} \Z/p\Z$ by $(B,M_B)$. 
We denote the 
category of fine log (not necessarily $p$-adic) formal schemes 
which are separated and 
topologically of finite type over $(\cB,M_{\cB})$ by $\LF$ and denote 
the full subcategory of $\LF$ which consists of $p$-adic ones by $\pLF$. 
We denote the full subcategory of $\LF$ which consists of usual fine 
log schemes by $\L$ and the full subcategory which consists of 
fine log schemes over $(B,M_B)$ by $\LB$. 
We call an object in $\LF$ 
a fine log formal $\cB$-scheme, an object in $\pLF$ a $p$-adic fine 
log formal $\cB$-scheme, an object in $\L$ a fine log $\cB$-scheme and 
an object in $\LB$ a fine log $B$-scheme. (When log structure is trivial, 
we omit the term `fine log'.) Note that, in Sections 4--7, 
we will impose more assumptions on $\cB$. \\
(2) \,\, For a formal $\cB$-scheme $T$, we denote the rigid analytic space 
assocated to $T$ by $T_K$. \\
(3) \,\, In this paper, we freely use terminologies concerning log structures 
defined in \cite{kkato}, \cite{shiho1} and \cite{shiho2}. 
For a fine log (formal) scheme $(X,M_X)$, $(X,M_X)_{\triv}$ denotes the 
maximal open sub (formal) scheme of $X$ on which the log structure $M_X$ is 
trivial. A morphism $f:(X,M_X) \lra (Y,M_Y)$ is said to be strict if 
$f^*M_Y = M_X$ holds. \\ 
(4) \,\, For a site ${\cal S}$, we will denote the topos associated to
${\cal S}$ by ${\cal S}^{\sim}$. \\
(5) \,\, Fiber products of log formal schemes are completed unless 
otherwise stated. On the other hand, the completed tensor pruduct of 
topological modules are denoted by $\hat{\otimes}$ to distinguish 
with the usual tensor product $\otimes$. \\


\section{Relative log convergent cohomology of radius $\lambda$ (I)} 

In this section, first we introduce the notion of 
relative log convergent site of radius $\lambda$ ($0 < \lambda \leq 1$), 
which is a generalization of the relative log convergent site 
introduced in \cite[2.4]{shiho3}. Then we prove some basic properties 
of it and we prove the log convergent Poincar\'{e} lemma for this site. 
After that, we prove the comparison theorem between relative log convergent 
cohomology of radius $\lambda$ with relative log crystalline cohomology 
when $\lambda$ is sufficiently close to $1$. 
The proof is similar to that given in \cite[\S 2]{shiho3} but 
slightly more complicated. \par 
First let us recall the definition of pre-widening, widening and enlargement, 
which was defined in \cite[Definition 2.1, 2.2]{shiho3}: 

\begin{defn}\label{quaddef}
Let $f:(\cX,M_{\cX}) \lra (\cY,M_{\cY})$ be a morphism in $\pLF$. 
Define the category 
${\cal Q}((\cX/\cY)^{\log})$ 
of quadruples on
$(\cX,M_{\cX})/(\cY, \allowbreak 
M_{\cY})$ as follows$:$ The objects are 
the data 
 $((\cZ,M_{\cZ}),(Z,M_Z),i,z),$ where $(\cZ,M_{\cZ})$ is an object in $\LF$
over $(\cY,M_{\cY})$, 
$(Z,M_Z)$ is an object in $\L$ over $(\cY,M_{\cY})$, 
$i$ is a closed immersion $(Z,M_Z) \hra (\cZ,M_{\cZ})$ over
 $(\cY,M_{\cY})$ and $z$ is a morphism $(Z,M_Z) \lra (\cX,M_{\cX})$ 
 in $\LF$ over $(\cY,M_{\cY})$. 
We define a morphism of quadruples on 
$(\cX,M_{\cX})/(\cY,M_{\cY})$ 
in an obvious way. $($See 
\cite[2.1.8]{shiho2}.$)$ 
\end{defn}

\begin{defn}\label{widedef}
Let $f:(\cX,M_{\cX}) \lra (\cY,M_{\cY})$ be as above. \\
$(1)$ \,
A quadruple $(({\cZ},M_{\cZ}),(Z,M_Z),i,z)$ on $(\cX,M_{\cX})/(\cY,M_{\cY})$ 
is called a pre-wid\-ening on $(\cX,M_{\cX})/(\cY,M_{\cY})$ 
if $\cZ$ is in $\pLF$. \\
$(2)$ \,
A quadruple $((\cZ,M_{\cZ}),(Z,M_Z),i,z)$ 
on $(\cX,M_{\cX})/(\cY,M_{\cY})$ 
is called a widening on $(\cX,M_{\cX})/(\cY,M_{\cY})$ 
if $i$ is a homeomorphic closed immersion $($that is, 
if $Z$ is a scheme of definition of $\cZ$ via $i)$. \\
$(3)$ \,
A $($pre-$)$widening $((\cZ,M_{\cZ}),(Z,M_Z),i,z)$ 
on $(\cX,M_{\cX})/(\cY,M_{\cY})$ 
is said to be exact
if $i$ is exact. It is said to be affine if $\cZ, z$ and the structure 
morphism $\cZ \lra \cY$ are affine. \\
$(4)$ \, A quadruple $((\cZ,M_{\cZ}),(Z,M_Z),i,z)$ on 
$(\cX,M_{\cX})/(\cY,M_{\cY})$ is called an 
enlargement if it is both pre-widening and widening, it is exact and 
$\cZ$ is flat over $\Spf W$. 
\end{defn}

We introduce the notion of the radius of an enlargement, which is 
one of the key notion in this paper. 

\begin{defn} 
\begin{enumerate} 
\item 
For a sheaf of rings $\cA$ on a site and a sheaf of ideals 
$\cI$ of $\cA$, we define $r(\cA,\cI) \in [0,\infty]$ by 
$r(\cA,\cI):=\sup \{b/a \,\vert\, \cI^a \subseteq p^b\cA\}$. 
\item 
Let $f:(\cX,M_{\cX}) \lra (\cY,M_{\cY})$ be a morphism in $\pLF$, let 
$\cZ := ((\cZ,M_{\cZ}),\allowbreak (Z, \allowbreak M_Z),
\allowbreak i,z)$ be an enlargement of 
$(\cX,M_{\cX})/(\cY,M_{\cY})$ and put $\cI := \Ker(\cO_{\cZ} \allowbreak \lra \cO_Z)$. 
Then we define $r(\cZ)$ by $r(\cZ) := r(\cO_{\cZ},\cI)$ and call it 
the logarithmic radius of $\cZ$. We define $\lam(\cZ)$ by 
$$ \lam(\cZ) := p^{-r(\cZ)} = 
\inf \{p^{-b/a} \,\vert\, \cI^a \subseteq p^b\cO_{\cZ}\} $$ 
and call it the radius of $\cZ$. By definition, we have 
$r(\cZ) \in (0,\infty)$ and $\lam(\cZ) \in (0,1)$. 
\end{enumerate}
\end{defn} 

We prove basic properties of $r(\cA,\cI)$ which we will use later. 

\begin{lem} 
Let $\cA$ be as above. 
\begin{enumerate} 
\item 
For sheaves of ideals $\cI,\cJ$ of $\cA$, we have 
$(r(\cA,\cI)^{-1}+r(\cA,\cJ)^{-1})^{-1} \leq r(\cA,\cI+\cJ)$. 
Moreover, if $r(\cA,\cJ)=\infty$ $($this is the case if, for example, 
$\cJ$ is nilpotent$)$, we have $r(\cA,\cI)=r(\cA,\cI+\cJ)$. 
\item 
For sheaves $\cJ \subseteq \cI$ of $\cA$ such that $\cJ$ is nilpotent, 
we have $r(\cA/\cJ,\cI/\cJ)=r(\cA,\cI)$. 
\end{enumerate}
\end{lem}

\begin{pf} 
If we have $\cI^a \subseteq p^b\cA$ and 
$\cJ^{a'} \subseteq p^{b'}\cA$, we have 
$$ (\cI+\cJ)^{ab'+a'b} \subseteq \cI^{ab'}+\cJ^{a'b} \subseteq p^{bb'}\cA. $$
The former assertion of (1) follows from this. the latter assertion 
follows from the inequality 
$$ r(\cA,\cI) = (r(\cA,\cI)^{-1}+\infty^{-1})^{-1} \leq 
r(\cA,\cI+\cJ) \leq r(\cA,\cI). $$ \par 
Next let us prove (2). 
If we have $\cI^a \subseteq p^b\cA$, 
we have $(\cI/\cJ)^a \subseteq p^b(\cA/\cJ)$. On the other hand, 
assume we have $(\cI/\cJ)^a \subseteq p^b(\cA/\cJ)$. Then we have 
$\cI^a \subseteq p^b\cA + \cJ$. 
If we take 
$N\in \N$ satisfying $\cJ^{N+1}=0$, we have, for any $n \geq N$, the 
inclusion 
$$ 
\cI^{an} \subseteq (p^b\cA + \cJ)^n \subseteq p^{b(n-N)}\cA. 
$$ 
From this, we can deduce the equality $r(\cA/\cJ,\cI/\cJ)=r(\cA,\cI)$. 
\end{pf} 

Now we give the definition of relative log convergent site of radius 
$\lam$ and the category of isocrystals on it, as follows: 

\begin{defn}
Let $\tau$ be one of the words $\{ \Zar(={\mathrm{Zariski}}), 
\et(={\mathrm{etale}})\}$.  
For a morphism $(\cX,M_{\cX}) \lra (\cY,M_{\cY})$ in $\pLF$ and 
$\lam \in (0,1]$, 
we define the log convergent site $(\cX/\cY)^{\log}_{\lamm\conv,\tau}$ 
of 
$(\cX,M_{\cX})/(\cY,M_{\cY})$ of radius $\lam$ 
with respect to $\tau$-topology as follows$:$ 
The objects are the enlargements $\cZ$ on $(\cX,M_{\cX})/(\cY,M_{\cY})$ 
 with $\lam(\cZ) < \lam$ and the morphisms are the morphism of enlargements. 
A family of morphisms 
$$ \{ 
((\cZ_{\alpha}, M_{\cZ_{\alpha}}), 
 (Z_{\alpha}, M_{Z_{\alpha}}), i_{\alpha}, z_{\alpha}) \lra 
 ((\cZ,M_{\cZ}),(Z,M_Z),i,z) \}_{\alpha \in I} $$
is a covering if the morphisms $(\cZ_{\alpha}, M_{\cZ_{\alpha}})
\ra (\cZ,M_{\cZ})$ are strict, form a covering of $\cZ$ with respect to 
$\tau$-topology and $(Z_{\alpha}, M_{Z_{\alpha}})$ is canonically 
isomorphic to 
$(\cZ_{\alpha}, M_{\cZ_{\alpha}}) \times_{(\cZ,M_{\cZ}),i} 
\allowbreak (Z,M_Z)$. 
When the log structures are 
 trivial, we omit the superscript $^{\log}$ in 
$(\cX/\cY)^{\log}_{\lamm\conv,\tau}$. 
We denote 
the right derived functor 
$($resp. the 
$q$-th right derived functor$)$ of the functor 
$$ 
(\cX/\cY)^{\log,\sim}_{\lamm\conv,\tau} \lra \cY^{\sim}_{\Zar}; 
\,\, \cE \mapsto (U \mapsto 
\Gamma((\cX \times_{\cY} U/U)^{\log}_{\lamm\conv,\tau},\cE)) 
$$ 
by $Rf_{\cX/\cY,\lamm\conv *}\cE$ 
$($resp. $R^qf_{\cX/\cY,\lamm\conv *}\cE)$. We call 
$R^qf_{\cX/\cY,\lamm\conv *}\cE$ 
the $q$-th relative log convergent cohomology of $(\cX,M_{\cX})/
(\cY,M_{\cY})$ with coefficient $\cE$. 
\end{defn}

\begin{rem}\label{radrem}
\begin{enumerate}
\item 
When $\lam$ is equal to $1$, 
$(\cX/\cY)^{\log}_{\lamm\conv,\tau}$ is nothing but 
the relative log convergent site $(\cX/\cY)^{\log}_{\conv,\tau}$ defined 
in \cite[2.4]{shiho3}. 
\item 
For a morphism of enlargements 
$g:\cZ' := ((\cZ', M_{\cZ'}),(Z', M_{Z'}), i', z') 
\lra ((\cZ,M_{\cZ}),(Z,M_Z),i,z)=:\cZ$ 
satisfying 
$(\cZ', M_{\cZ'}) \times_{(\cZ,M_{\cZ}),i} (Z,M_Z) = (Z', \allowbreak 
M_{Z'})$,  
we have $\lam(\cZ') \leq \lam(\cZ)$. So the category of coverings 
of an object $\cZ$ in $(\cX/\cY)^{\log}_{\lamm\conv,\tau}$ is equal to 
the category of coverings of $\cZ$ when considered in 
$(\cX/\cY)^{\log}_{\conv,\tau}$. 
\end{enumerate}
\end{rem}

\begin{defn}
Let the notations be as above. An isocrystal on the log 
convergent site of radius $\lam$ 
$(\cX/\cY)^{\log}_{\lamm\conv,\tau}$ is a sheaf 
${\cal E}$ on $(\cX/\cY)^{\log}_{\lamm\conv,\tau}$ satisfying 
the following conditions$:$ 
\begin{enumerate}
\item
For any object $\cZ$ in $(\cX/\cY)^{\log}_{\lamm\conv,\tau}$, 
the sheaf ${\cal E}_{\cZ}$ on $\cZ$ 
induced by ${\cal E}$ is an isocoherent sheaf. $($That is, 
$\cE_{\cZ}$ has the form $\Q \otimes_{\cZ} \cF$ with $\cF$ a coherent 
$\cO_{\cZ}$-module.$)$ 
\item 
For any morphism $f:\cZ' \lra \cZ$ in $(\cX/\cY)^{\log}_{\lamm\conv,\tau}$, 
the homomorphism $f^*{\cal E}_{\cZ} \allowbreak \lra {\cal E}_{\cZ'}$ of 
sheaves on $\cZ'$ induced by ${\cal E}$ is an 
isomorphism. 
\end{enumerate}
We denote the category of isocrystals on 
$(\cX/\cY)^{\log}_{\lamm\conv,\tau}$ 
by 
$I_{\lamm\conv,\tau}((\cX/\cY)^{\log})$. 
When the log structures are trivial, we omit the superscript 
$^{\log}$. 
\end{defn}

\begin{defn} 
Let the notations be as above. Then an isocrystal 
$\cE$ is said to be locally free if, 
for any enlargement $\cZ$, the sheaf ${\cal E}_{\cZ}$ on $\cZ$ 
induced by ${\cal E}$ is a locally free $\Q \otimes_{\Z} 
{\cal O}_{\cZ}$-module in the sense of \cite[1.9]{shiho3}. 
$($Note that it does not mean that $\cE$ is a free 
$\Q \otimes_{\Z} \cO_{\cZ}$-module locally on $\cZ_{\Zar}.)$ 

\end{defn}

For a morphism $(\cX,M_{\cX}) \lra (\cY,M_{\cY})$ in $\pLF$ and 
a (pre-)widening $\cZ$ on $(\cX,M_{\cX})/(\cY,M_{\cY})$, we can define 
the notion of localized log convergent site of radius $\lam$ 
$(\cX/\cY)^{\log}_{\lamm\conv,\tau} \vert_{\cZ}$ 
and the category of isocrystals 
$I_{\lamm\conv,\tau}((\cX \allowbreak /\cY)^{\log}|_{\cZ})$ on 
it. Also, 
for a simplicial object $(\cX^{(\b)},M_{\cX^{(\b)}})$ in $\pLF$ and a 
morphism $(\cX^{(\b)},\allowbreak M_{\cX^{(\b)}}) \allowbreak 
\lra (\cY,M_{\cY})$ in $\pLF$, 
we can define the relative log convergent topos 
$(\cX^{(\b)}/\cY)^{\log,\sim}_{\lamm\conv,\tau}$
of radius $\lam$. As in \cite[2.7]{shiho3} and \cite[2.1.20]{shiho2}, 
we have the following: 

\begin{prop}\label{cohdes}
Let $(\cX,M_{\cX}) \lra (\cY,M_{\cY})$ be a morphism in $\pLF$ and 
let $\tau$ be one of the words 
$\{ \Zar(={\mathrm{Zariski}}), \et(={\mathrm{etale}})\}$.
Let $g^{(\b)}: (\cX^{(\b)},M_{\cX^{(\b)}}) \lra 
(\cX,M_{\cX})$ be a strict $\tau$-hypercovering of $\cX$. Let 
$ \theta:= (\theta_*,\theta^{-1}): 
({\cX}^{(\bullet)}/\cY)^{\log,\sim}_{\lamm\conv,\tau} \lra 
(\cX/\cY)^{\log,\sim}_{\lamm\conv,\tau} $
be the morphism of topoi characterized by 
$\theta^{-1}(E)^{(i)} := g^{(i),-1}(E)$. 
Then, for any abelian sheaf $E$ on $(\cX/\cY)^{\log}_{\lamm\conv,\tau}$, 
the canonical homomorphism 
$E \lra R\theta_* \theta^{-1} E$
is a quasi-isomorphism. 
\end{prop}

We also need a similar but slightly different descent property, 
which we will recall now. 
Let $(\cX,M_{\cX}) \lra (\cY,M_{\cY})$ be as above and assume we are 
given an open covering 
$\cX = \bigcup_{j \in J}\cX_j$ by finite number of sub formal schemes. 
For non-empty subset $L \subseteq J$, put 
$\cX_L := \bigcap_{j \in L}\cX_j, M_{\cX_L} := M_{\cX_j} |_{\cX_L}$ and 
for $m \in \N$, let $(\cX^{(m)}, M_{\cX^{(m)}})$ be the disjoint union 
of $(\cX_L,M_{\cX_L})$'s for $|L|=m+1$. ($\cX^{(m)}$ is empty 
for $m \geq |J|=:N+1$.) 
 Let $\Delta^+_N$ be the category 
such that the 
objects are the sets $[m]:=\{0,1,2,\cdots m\}$ with $m \leq N$ and that 
$\Hom_{\Delta^+_N}([m],[m'])$ is the set of strictly increasing 
maps $[m] \ra [m']$. Then, if we fix a total order on $J$, 
we can regard $(\cX^{(\b)},M_{\cX^{(\b)}}) := 
\{(\cX^{(m)},M_{\cX^{(m)}})\}_{0 \leq m \leq N}$ naturally as a diagram 
indexed by $\Delta^+_N$ and we have a morphism 
$g^{(\b)}: (\cX^{(\b)},M_{\cX^{(\b)}}) \lra (\cX,M_{\cX})$. 
(cf. \cite[\S 1]{shiho3}.) Then we have the following proposition, 
whose proof is the same as Proposition \ref{cohdes} (we omit it): 

\begin{prop}\label{cohdes2}
With the above notation, let 
$ \theta:= (\theta_*,\theta^{-1}): 
({\cX}^{(\bullet)}/\cY)^{\log,\sim}_{\lamm\conv,\tau} \lra 
(\cX/\cY)^{\log,\sim}_{\lamm\conv,\tau} $
be the morphism of topoi characterized by 
$\theta^{-1}(E)^{(i)} := g^{(i),-1}(E)$. 
Then, for any abelian sheaf $E$ on $(\cX/\cY)^{\log}_{\lamm\conv,\tau}$, 
the canonical homomorphism 
$E \lra R\theta_* \theta^{-1} E$
is a quasi-isomorphism. 
\end{prop}

We have the following as in \cite[2.9]{shiho3}, \cite[2.1.21]{shiho2}: 

\begin{prop}\label{zid}
Let $(\cX,M_{\cX}) \lra (\cY,M_{\cY})$ 
be as above and let $\cZ$ be a 
$($pre-$)$widening of $(\cX,M_{\cX})/(\cY,M_{\cY})$. Let us denote the 
canonical morphism of topoi 
$$ (\cX/\cY)^{\log,\sim}_{\lamm\conv,\et} \lra 
(\cX/\cY)^{\log,\sim}_{\lamm\conv,\Zar} \quad \text{$($resp.} \,\,\,
(\cX/\cY)^{\log,\sim}_{\lamm\conv,\et} \vert_{\cZ} 
\lra (\cX/\cY)^{\log,\sim}_{\lamm\conv,\Zar} \vert_{\cZ} \,\,\, \text{$)$} $$ 
by $\epsilon$. Then$:$ 
\begin{enumerate}
\item
For any $\cE \in I_{\lamm\conv,\et}((\cX/\cY)^{\log})$ 
$($resp. $\cE \in I_{\lamm\conv,\et}((\cX/\cY)^{\log} \vert_{\cZ}))$, we have 
$R\epsilon_* \cE = \epsilon_* \cE$. 
\item 
The functor $\cE \mapsto \epsilon_* \cE$ induces the equivalence 
of categories 
\begin{align*}
I_{\lamm\conv,\et}((\cX/\cY)^{\log}) & \overset{\sim}{\lra} 
I_{\lamm\conv,\Zar}((\cX/\cY)^{\log}) \\ 
\text{$($resp.} \,\,\, 
I_{\lamm\conv,\et}((\cX/\cY)^{\log} \vert_{\cZ}) & \overset{\sim}{\lra} 
I_{\lamm\conv,\Zar}((\cX/\cY)^{\log} \vert_{\cZ}) 
\,\,\, \text{$).$} 
\end{align*}
\end{enumerate}
\end{prop}

In the rest of this paper, we denote the category 
$I_{\lamm\conv,\et}((\cX/\cY)^{\log})=
I_{\lamm\conv,\Zar}((\cX \allowbreak /\cY \allowbreak )^{\log})$ (resp. 
$I_{\lamm\conv,\et}((\cX/\cY)^{\log} \vert_{\cZ})= 
I_{\lamm\conv,\Zar}((\cX/\cY)^{\log} \vert_{\cZ})$) simply by 
$I_{\lamm\conv} \allowbreak ((\cX/\allowbreak \cY \allowbreak 
)^{\log})$ (resp. 
$I_{\lamm\conv}((\cX/\cY)^{\log} \vert_{\cZ})$) and call it an isocrystal 
on $(\cX/\cY)^{\log}_{\lamm\conv}$ 
(resp. $(\cX/\allowbreak \cY)^{\log}_{\lamm\conv} |_{\cZ}$), 
by abuse of terminology. \par 
Next we introduce the notion of the system of 
universal enlargements of radius $\lambda$, which is a 
generalization of the notion of the system of 
universal enlargements defined in \cite[2.12]{shiho3}. To this end, first 
we define the oriented set which will be the index category of the system 
of universal enlargements: 

\begin{defn} 
Let $\A$ be the set $\{(a,b)\,\vert\,a,b \in \N_{\geq 1}\}$. For 
$(a,b),(a',b') \in \A$, 
we write $(a,b) \dsucc (a',b')$ if $(a',b')$ is equal to one of 
the elements $(a+1,b), (a,b-1), (na,nb)\,(\text{for some 
$n \in \N_{\geq 1}$})$. 
We write $(a,b) \succ (a',b')$ if there exists a sequence of elements in 
$\A$ of the form 
$$ (a,b) = (a_0,b_0) \dsucc (a_1,b_1) \dsucc \cdots \dsucc 
(a_n,b_n) = (a',b'). $$
For $r>0$, let $\A_r$ be the subset $\{(a,b) \in \A \,\vert\, b/a > r\}$. 
\end{defn}

\begin{lem} 
The relation $\succ$ defines a partial order on $\A$ 
and for any $r>0$, 
the set $\A_r$ is oriented with respect to the order $\succ$. 
$($That is, for any $(a,b), (a',b')\in \A$, there exists an element 
$(a'',b'')$ in $\A$ with $(a,b) \succ (a'',b''), (a',b') \succ (a'',b'').)$ 
\end{lem} 

\begin{pf} 
If we have $(a,b) \dsucc (a',b')$, we have $b/a > b'/a'$ or 
$b/a = b'/a', a \leq a'$. So the same is true if we have 
$(a,b) \succ (a',b')$. Hence, if we have both $(a,b) \succ (a',b')$ and 
$(a',b') \succ (a,b)$, we have $b/a=b'/a', a=a'$. So we have 
$(a,b)=(a',b')$ and so $\succ$ is a partial order on $\A$. \par 
If we have $(a,b), (a',b') \in \A_r$, we have $(a,b) \succ (a'',b''), 
(a',b') \succ (a'',b'')$ if we put 
$(a'',b''):= (aa', \min\{ab',a'b\})$. So $\A_r$ is oriented. 
\end{pf} 

Sometimes it is convenient to take a sequence which is cofinal 
with $\A_r$: 

\begin{defn}
A sequence $\a := \{(a_n,b_n)\}$ of elements of $\A_r$ is said to be 
a good sequence if the following conditions are satisfied$:$ 
\begin{enumerate}
\item 
We have $(a_1,b_1) \succ (a_2,b_2) \succ \cdots$, 
$b_n/a_n > r \,(\forall n)$ and $\lim_{n \to \infty} (b_n/a_n) = r$. 
\item 
$a_{n+1}$ is divisible by $a_n$ for any $n \in \N$ and 
for any $m \in \N_{\geq 1}$, there exists $n \in \N$ 
such that $a_n$ is divisible by $m$. 
\end{enumerate}
\end{defn} 

\begin{lem} 
For any $r>0$, there exists a good sequence $\a$ in $\A_r$ and 
any good sequence $\a$ in $\A_r$ 
is cofinal with respect to $\succ$ in $\A_r$. 
\end{lem} 

\begin{pf} 
For $n \in \N$, let $r_n := b_n/n!$ be the least rational number 
strictly bigger than $r$ with denominator $n!$. Then one can check that 
the sequence $\{(n!,b_n)\}_{n}$ is a good sequence in $\A_r$. \par 
Let us prove the latter statement. 
For any good sequence $\a:=\{(a_n,b_n)\}$ in $\A_r$ and 
$(a,b) \in \A_r$, there exists $n \in \N$ satisfying 
$b/a > b_n/a_n$ and that $a_n$ is divisible by $a$. Then we have 
$(a,b) \succ (a_n,b_n)$. So $\a$ is cofinal in $\A_r$. 
\end{pf} 

Now we define the notion of the system of universal enlargements 
of radius $\lam$. 
Let $\cZ := ((\cZ,M_{\cZ}),(Z,M_Z),i,z)$ be a $($pre-$)$widening, 
let $((\cZ^{\ex}, \allowbreak M_{\cZ^{\ex}}), \allowbreak 
(Z,M_Z),i^{\ex},z)$ be 
the exactification of $\cZ$ and let ${\cal I}$ be the 
ideal $\Ker ({\cal O}_{\cZ^{\ex}} \ra {\cal O}_{Z})$. For 
$(a,b) \in \A$, let 
$B_{a,b}(\cZ)$ be the 
formal blow-up of $\cZ^{\ex}$ with respect 
to the ideal $ {\cal I}^a + p^b {\cal O}_{\cZ^{\ex}}$, 
let 
$T'_{a,b}(\cZ)$ be the open sub formal scheme
$$ \{x \in B_{a,b}(\cZ) \,\vert\, ({\cal I}^a + p^b 
{\cal O}_{\cZ^{\ex}}) \cdot {\cal O}_{B_{n}(\cZ),x} = 
p^b{\cal O}_{B_{n}(\cZ),x} \} $$ 
and 
let $T_{a,b}(\cZ)$ be $T'_{a,b}(\cZ)_{\fl}$, where, for a 
$p$-adic formal $\cB$-scheme $T$, we denote by $T_{\fl}$ the closed 
subscheme of $T$ defined by the ideal 
$\{x \in {\cal O}_{T}\,\vert\, \exists n, \allowbreak p^nx=0\}$. 
Let $\lambda_{a,b}:T_{a,b}(\cZ) \lra \cZ$ be the canonical 
morphism and put $Z_{a,b} := \lambda_n^{-1}(Z)$. 
Then the quadruple 
{\tiny{ $$ T_{a,b}(\cZ) := ((T_{a,b}(\cZ),M_{\cZ^{\ex}} 
|_{T_{a,b}(\cZ)}),(Z_{a,b}, M_{Z}|_{Z_{a,b}}), 
Z_{a,b} \hra T_{a,b}(\cZ), (Z_{a,b},M_{Z_{a,b}}) \overset{\lambda_{a,b}}{\lra} 
(Z,M_Z) \overset{z}{\lra} (X,M))  $$ }}
is an enlargement for each $(a,b) \in \A$. 

\begin{lem} 
Let $T_{a,b}(\cZ)$ be the enlargement defined above. Then$:$ 
\begin{enumerate}
\item  
We have $\lam(T_{a,b}(\cZ)) \leq p^{-b/a}$. 
\item 
The enlargements $\{T_{a,b}(\cZ)\}_{(a,b)}$ forms an inductive system 
of enlargements indexed by $\A$. 
\item 
Let us take $\lambda \in (0,1]$ and put $\lambda = p^{-r}$. 
For $T=\cZ$ or $T=T_{a,b}(\cZ)$, denote the sheaf represented by $T$ on 
$(\cX/\cY)^{\log}_{\lamm\conv,\tau}$ by $h_T$. Then
we have 
$h_{\cZ} = \varinjlim_{(a,b) \in \A_r} h_{T_{a,b}(\cZ)}$. 
$($That is, $h_{\cZ}$ is ind-represented by the system 
$\{T_{a,b}(\cZ)\})$. 
\end{enumerate}
\end{lem} 

\begin{pf} 
(1) is obvious by definition. Let us prove (2). Take $(a,b) \in \A$. Then 
we have ${\cal I}^a \subseteq p^b\cO_{T_{a,b}(\cZ)}$. This implies 
the inclusions $\cI^{a+1} \subseteq p^b\cO_{T_{a,b}(\cZ)}, 
\cI^a \subseteq p^{b-1}\cO_{T_{a,b}(\cZ)}$ and 
$\cI^{na} \subseteq p^{nb}\cO_{T_{a,b}(\cZ)}\,(n \in \N_{\geq 1})$. 
Then, by the universality 
of blow-up, we have the canonical mmorphisms 
$$ T_{a,b}(\cZ) \lra T_{a+1,b}(\cZ), \,\,\, T_{a,b}(\cZ) \lra 
T_{a,b-1}(\cZ), \,\,\, T_{a,b}(\cZ) \lra T_{na,nb}(\cZ), $$
that is, we have $T_{a,b}(\cZ) \lra T_{a',b'}(\cZ)$ for 
$(a',b')$'s with $(a,b) \dsucc (a',b')$. By composition, we obtain the 
morphisms $T_{a,b}(\cZ) \lra T_{a',b'}(\cZ)$ for $(a',b')$'s with 
$(a,b) \succ (a',b')$. This morphism defines the structure of inductive 
system on $\{T_{a,b}(\cZ)\}_{(a,b)\in\A}$ indexed by $\A$. (Transitivity 
follows again from the universality of blow-up.) \par 
Let us prove (3). For an enlargement 
$\cZ' := ((\cZ',M_{\cZ'}),(Z',M_{Z'}),i',z')$ with $\lam(\cZ) < \lam$ with 
a morphism of quadruples $g: \cZ' \lra \cZ$, there exists $(a,b) \in \A_r$ 
such that $\Ker(\cO_{\cZ'} \ra \cO_{Z'})^a \subseteq p^b\cO_{\cZ'}$. 
Then, by the universality of blow-up, the morphism $g$ factors uniquely as 
$\cZ' \lra T_{a,b}(\cZ) \lra \cZ$. From this fact, we obtain the 
assertion (3). 
\end{pf} 

We have the following explicit description of $T_{a,b}(\cZ)$ when 
$\cZ^{\ex}$ is affine: If we put 
$A := \Gamma(\cZ^{\ex},\cO_{\cZ^{\ex}})$, 
$I := \Gamma(\cZ^{\ex},\Ker(\cO_{\cZ^{\ex}} \ra \cO_{Z}))$ and 
write $I^a = (f_1, \cdots, f_k)$, we have 
$T_{a,b}(\cZ) = \Spf A_{a,b}$, where 
\begin{equation}\label{explicit}
A_{a,b} = (A[t_1,\cdots,t_k]/(p^bt_i-f_i\,(1 \leq i \leq k)))^{\wedge}
/(\text{$p$-torsion}). 
\end{equation}
(${}^{\wedge}$ denotes the $p$-adic completion.) 
Then we have the following, as in \cite[0.2.2]{ogus2}:

\begin{lem}\label{dense}
The image of the homomorphism $\Q \otimes_{\Z} A \lra \Q \otimes_{\Z} A_{a,b}$ 
is dense in $p$-adic topology. 
\end{lem} 

\begin{pf}
If we put $B := 
(A[t_1,\cdots,t_k]/(p^bt_i-f_i\,(1 \leq i \leq k)))
/(\text{$p$-torsion})$, the above homomorphism factors as 
$\Q \otimes_{\Z} A {\lra} \Q \otimes_{\Z} B \lra 
\Q \otimes_{\Z} A_{a,b}$ and the first homomorphism is an isomorphism. 
Moreover, we have $A_{a,b}=\widehat{B}$, where $\widehat{\phantom{B}}$ 
denotes the $p$-adic completion. If $x$ is an element in 
$\Q \otimes_{\Z} A_{a,b}$, $p^nx$ is in $A_{a,b}$ for some $n$. 
Then, for any $m \in \N$, there exists some $y \in B$ and 
$z \in A_{a,b}$ such that $p^nx = y + p^{m+n}z$ holds. Then we have 
$x = p^{-n}y + p^mz$ and $p^{-n}y \in \Q \otimes_{\Z} B = 
\Q \otimes_{Z} A$. So we are done. 
\end{pf} 

\begin{cor} 
For a $($pre-$)$widining $\cZ$ and $(a,b) \in \A$, the map 
$h_{T_{a,b}(\cZ)} \lra h_{\cZ}$ is injective. 
\end{cor} 

\begin{pf} 
The proof is same as \cite[0.2.2]{ogus2}. 
We may assume that $\cZ = \Spf A, T_{a,b}(\cZ) = \Spf A_{a,b}$ as above and 
it suffices to show that, for any enlargement $\cZ'$ with 
$\cZ' = \Spf B$ and any two homomorphisms $A_{a,b} \lra B$ with the same 
composition $A \lra A_{a,b} \lra B$ are equal. We can check it 
after applying $\Q \otimes_{\Z}$ (because $B$ is flat over $W$), and in 
this case, the claim follows from the separatedness of $B$ with respect to 
the $p$-adic topology and Lemma \ref{dense}. 
\end{pf} 

We have the following lemmas as in \cite[2.15, 2.16]{shiho3}, 
\cite[2.1.25, 2.1.26]{shiho2}. 
(Since we can prove them in the same way by using the explicit 
description \eqref{explicit}, we omit the proof. )

\begin{lem} 
Let $\cZ$ be a pre-widening and let $\widehat{\cZ}$ be the associated 
widening. Then we have $T_{a,b}(\cZ) \cong T_{a,b}(\widehat{\cZ})$ for 
any $(a,b) \in \A$. 
\end{lem} 

\begin{lem}\label{2.1.26}
Let 
$g: ((\cZ',M_{\cZ'}),(Z',M_{Z'}),i',z') \lra ((\cZ,M_{\cZ}),(Z,M_Z),i,z)$ 
be a morphism of $($pre-$)$widenings such that 
$(Z,M_Z) \times_{(\cZ,M_{\cZ})} (\cZ',M_{\cZ'}) = 
(Z',M_{Z'})$ holds naturally and that $\cZ' \lra \cZ$ is flat. 
Then $g$ induces the natural isomorphism of enlargements 
$T_{a,b}(\cZ') \os{=}{\lra} T_{a,b}(\cZ) \times_{\cZ} \cZ'$. 
\end{lem}

The following is a variant of \cite[2.1.27]{shiho2}, \cite[2.17]{shiho3}. 

\begin{prop}\label{2.1.27}
Let $(\cX,M_{\cX}) \lra (\cY,M_{\cY})$ be a morphism in 
$\pLF$. 
Let 
$$f:\cZ_1 := ((\cZ,M_{\cZ}),(Z,M_Z),i,z) \lra 
\cZ_2 := ((\cZ,M_{\cZ}),(Z',M_{Z'}),i',z')$$ 
be a morphism of $($pre-$)$widenings of $(\cX,M_{\cX})/(\cY,M_{\cY})$ 
such that $f$ is identity on $\cZ$ and 
$(Z,M_Z) \lra (Z',M_{Z'})$ is an exact closed immersion. 
For $(a,b) \in \A$, denote the map of enlargements 
$T_{a,b}(\cZ_1) \lra T_{a,b}(\cZ_2)$ induced by $f$ by $f_{(a,b)}$ 
and for $(a,b),(a',b') \in \A$ with $(a,b)\succ (a',b')$, 
denote the transition 
morphism $T_{a,b}(\cZ_i) \lra T_{a',b'}(\cZ_i)$ by 
$g_{i,(a,b)(a',b')}$. 
Let $\cI, \cI'$ be the defining ideal of $Z,Z'$ in $\cZ$ respectively and 
assume that there exists an ideal sheaf $\cJ$ of $\cO_{\cZ}$ with 
$\cJ^{m+1}=0$ satisfying $\cI \subseteq \cI' + \cJ$. 
Then$:$ 
\begin{enumerate}
\item 
For any $(\alpha,\beta) \in \A$ with $\alpha>m$, 
there exists a morphism of $p$-adic fine log formal $\cB$-schemes 
$h_{(\alpha,\beta)}: 
T_{\alpha-m,\beta}(\cZ_2) \lra T_{\alpha,\beta}(\cZ_1)$, satisfying the 
following condition$:$ For any 
$(a,b),(\alpha,\beta) \in \A$ with $\alpha >m, (\alpha,\beta)\succ (a,b)$, 
the composite 
$$ 
T_{\alpha-m,\beta}(\cZ_2) \os{h_{(\alpha,\beta)}}{\lra} 
T_{\alpha,\beta}(\cZ_1) \os{g_{1,(\alpha,\beta)(a,b)}}{\lra} 
T_{a,b}(\cZ_1) \os{f_{(a,b)}}{\lra} T_{a,b}(\cZ_2) 
$$ 
is equal to $g_{2,(\alpha-m,\beta)(a,b)}$ and for any 
$(a,b),(\alpha,\beta) \in \A$ with $\alpha > m, 
(a,b) \succ (\alpha-m,\beta)$, the composite 
$$ 
T_{a,b}(\cZ_1) \os{f_{(a,b)}}{\lra} T_{a,b}(\cZ_2) 
\os{g_{2,(a,b)(\alpha-m,\beta)}}{\lra} T_{\alpha-m,\beta}(\cZ_2) 
\os{h_{(\alpha,\beta)}}{\lra} T_{\alpha,\beta}(\cZ_1) 
$$ 
is equal to $g_{1,(a,b)(\alpha,\beta)}$. 
\item 
Let $\lam \in (0,1]$ and put $\lam = p^{-r}$. Then, 
for an isocrystal $\cE$ on $(\cX/\cY)^{\log}_{\lamm\conv,\tau}$ 
and any $(\alpha,\beta) \in \A_r$ with $\alpha>m$, 
we have the natural isomorphism 
$$ \varphi_{(\alpha,\beta)}: 
h_{(\alpha,\beta)}^*\cE_{T_{(\alpha,\beta)}(\cZ_1)} \os{\sim}{\lra} 
\cE_{T_{(\alpha-m,\beta)}(\cZ_2)} $$
which satisfies the following condition$:$ 
For any 
$(a,b),(\alpha,\beta) \in \A_r$ with 
$\alpha >m$ and $(\alpha,\beta)\succ (a,b)$, 
the composite 
\begin{align*}
(h_{(\alpha,\beta)} \circ g_{1,(\alpha,\beta)(a,b)} \circ 
f_{(a,b)})^*\cE_{T_{a,b}(\cZ_2)}
& \os{h_{(\alpha,\beta)}^*
\cE(g_{1,(\alpha,\beta)(a,b)} \circ f_{(a,b)})}{\lra} 
h_{(\alpha,\beta)}^*\cE_{T_{\alpha,\beta}(\cZ_1)} \\ & 
\os{\varphi_{(\alpha,\beta)}}{\lra} \cE_{T_{\alpha-m,\beta}(\cZ_2)}
\end{align*}
is equal to the homomorphism induced by $g_{2,(\alpha-m,\beta)(a,b)}$ and 
for any 
$(a,b),(\alpha, \allowbreak \beta) \allowbreak \in \A_r$ with $\alpha > m$ 
and 
$(a,b) \succ (\alpha-m,\beta)$, the composite
{\tiny{ \begin{align*}
(f_{(a,b)} \circ g_{2,(a,b)(\alpha-m,\beta)} \circ h_{(\alpha,\beta)})^*
\cE_{T_{\alpha,\beta}(\cZ_1)} 
& \os{(f_{(a,b)} \circ 
g_{2,(a,b)(\alpha-m,\beta)})^*\varphi_{\alpha,\beta}}{\lra} 
(f_{(a,b)} \circ g_{2,(a,b)(\alpha-m,\beta)})^*\cE_{T_{\alpha-m,\beta}(\cZ_2)} 
\\ & 
\os{\cE(f_{(a,b)} \circ g_{2,(a,b)(\alpha-m,\beta)})}{\lra} 
\cE_{T_{a,b}(\cZ_1)}
\end{align*}}}
is equal to the homomorphism induced by 
$g_{1,(a,b)(\alpha,\beta)}$. 
\end{enumerate}
\end{prop} 

\begin{pf}
We may replace $\cZ$ by $\cZ^{\ex}$. So we may assume that 
$\cZ_1, \cZ_2$ are exact widenings. \par 
For $(\alpha,\beta) \in \A$ with $\alpha>m$, we have 
$\cI^{\alpha} \subseteq (\cI'+\cJ)^{\alpha} \subseteq {\cI'}^{\alpha-m}$. 
So, since we have 
${\cI'}^{\alpha-m} \subseteq 
p^{\beta}\cO_{T_{(\alpha,\beta)}(\cZ_2)}$, we have 
$\cI^{\alpha} \subseteq p^{\beta}\cO_{T_{(\alpha,\beta)}(\cZ_2)}$. 
So, by the universality of blow-up, we have the morphism 
$h_{(\alpha,\beta)}$. \par 
Let us denote the pull-back of $Z$ to 
$T_{(\alpha,\beta)}(\cZ_1)$ by $Z_{(\alpha,\beta)}$ and the pull-back of 
$Z'$ to $T_{(\alpha-m,\beta)}(\cZ_2)$ by $Z'_{(\alpha-m,\beta)}$. 
On the other hand, denote the pull-back of $Z$ to 
$T_{(\alpha-m,\beta)}(\cZ_2)$ by $Z_{(\alpha-m,\beta)}$. Then we have 
the diagram of enlargements 
\begin{align*}
T_{(\alpha-m,\beta)}(\cZ_2) = 
(T_{(\alpha-m,\beta)}(\cZ_2),Z'_{(\alpha-m,\beta)}) 
& \os{j}{\lla} 
(T_{(\alpha-m,\beta)}(\cZ_2),Z_{(\alpha-m,\beta)}) \\
& \os{h_{(\alpha,\beta)}}{\lra} 
(T_{(\alpha,\beta)}(\cZ_1),Z_{(\alpha,\beta)}) =: T_{(\alpha,\beta)}(\cZ_1), 
\end{align*}
where $j$ is the morphism of enlargements induced by 
the identity map on $T_{(\alpha-m,\beta)}(\cZ_2)$. 
Moreover, we have $\lam(T_{(\alpha-m,\beta)}(\cZ_2)) 
\leq p^{-\beta/(\alpha-m)} < \lam, 
\lam(T_{(\alpha,\beta)}(\cZ_1)) \leq p^{-\beta/\alpha} < \lam$ by definition 
and  
$\lam((T_{(\alpha-m,\beta)}(\cZ_2),Z_{(\alpha-m,\beta)})) \leq 
\lam(T_{(\alpha,\beta)}(\cZ_1)) < \lam$ by Remark \ref{radrem} (2). 
So we have the isomorphisms 
$$ 
h_{(\alpha,\beta)}^*E_{T_{(\alpha,\beta)}(\cZ_1)} 
\os{\sim}{\lra} 
E_{(T_{(\alpha-m,\beta)}(\cZ_2),Z_{(\alpha-m,\beta)})} 
\os{\sim}{\lla} 
j^*E_{T_{(\alpha-m,\beta)}(\cZ_2)} = E_{T_{(\alpha-m,\beta)}(\cZ_2)}. 
$$ 
We define the isomorphism
 $\varphi_{(\alpha,\beta)}$ by the morphism induced by 
the above diagram. Then it is easy to check that 
$h_{(\alpha,\beta)}$ and $\varphi_{(\alpha,\beta)}$ 
satisfy the desired properties. 
\end{pf} 

Now we introduce the site $\vec{\cZ}_{r,\a}$, 
which is a generalization of the site $\vec{\cZ}$ introduced in 
\cite{ogus2}, \cite{shiho2}, \cite{shiho3}: 

\begin{defn} 
Let $r>0$, let $\a$ be a good sequence in $\A_r$ and let 
$\cZ$ be a $($pre-$)$widening. Then we define the site $\vec{\cZ}_{r,\a}$ 
as follows$:$ Objects are the open sets of $T_{a,b}(\cZ)$ for some 
$(a,b) \in \a$. For open sets $U \subseteq T_{a,b}(\cZ)$ and 
$V \subseteq T_{a',b'}(\cZ) \,((a,b),(a',b') \allowbreak \in \a)$, 
$\Hom_{\vec{\cZ}_{r,\a}}(U,V)$ is the set of morphisms 
$U \lra V$ commutes with the transition morphism $T_{a,b}(\cZ) \lra 
T_{a',b'}(\cZ)$ when $(a,b) \succ (a',b')$ holds and is empty otherwise. 
The coverings are defined by Zariski topology for each object. 
We define the strucrure sheaf $\cO_{\vec{\cZ}_{r,\a}}$ by 
$\cO_{\vec{\cZ}_{r,\a}}(U) := \Gamma(U,\cO_U)$. 
\end{defn} 

\begin{defn} 
A sheaf of $\Q \otimes_{\Z} \cO_{\vec{\cZ}_{r,\a}}$-module $E$ is called 
crystalline if, for any transition map 
$\psi: T_{a,b}(\cZ) \lra T_{a',b'}(\cZ)$, the transition map of sheaves 
$\psi^{-1}E_{a',b'} \lra E_{a,b}$ induces the isomorphism 
$\cO_{T_{a,b}(\cZ)} \otimes_{\psi^{-1}(\cO_{T_{a',b'}(\cZ)})} 
\psi^{-1}E_{a',b'} \os{\sim}{\lra} E_{a,b}, $
where $E_{a,b}, \allowbreak E_{a',b'}$ denotes the Zariski sheaf on 
$T_{a,b}(\cZ),T_{a',b'}(\cZ)$ induced by $E$. 
\end{defn}

We define the morphism of ringed 
topoi $\gamma: \vec{\cZ}^{\sim}_{r,\a} \lra 
\cZ^{\sim}_{\Zar}$ by $\gamma_*E := 
\varprojlim_{(a,b) \in \a}\alpha_{a,b *} \allowbreak E_{a,b}$ and 
$\gamma^*(E)|_{T_{a,b}(\cZ)} := \alpha_{a,b}^*E$. 
(Here $E_{a,b}$ is as above and $\alpha_{a,b}$ is the morphism 
$T_{a,b}(\cZ) \allowbreak \lra \cZ$.) 
Then we have the following as in \cite[0.3.7]{ogus2}: 

\begin{prop}\label{gamma-acyclic}
Let $\cZ$ be an affine widening, 
let $\a$ be a good sequence in $\A_r$ and let $E$ be a crystalline 
sheaf of $\Q \otimes_{\Z} \cO_{\vec{\cZ}_{r,\a}}$-modules. Then we have 
$H^q(\vec{\cZ}_{r,\a},E)=0$ for $q>0$. 
\end{prop} 

\begin{pf} 
Since the proof is similar to \cite[0.3.7]{ogus2}, we only give a sketch 
of the proof. 
Let us put $\cZ^{\ex} = \Spf A$, $T_{a,b}(\cZ) = \Spf A_{a,b} 
\, ((a,b) \in \a)$, $\a := \{(a_n,b_n)\}_{n \in \N}$ and 
put $A_n := A_{a_n,b_n}, E_n := E_{a_n,b_n}$ for simplicity. 
Then, exactly in the same way as the proof of \cite[0.3.7]{ogus2}, 
we can reduce the proof of the proposition to the following claim: \\
\quad \\
{\bf claim.} \,\,\, 
Let $\{F_n\}_n$ be a family of finitely generated 
$\Q \otimes_{\Z} A_n$-modules endowed with transition maps 
$\psi_n: F_{n+1} \lra F_n$ 
such that each $\psi_n$ induces 
the isomorphism $A_{n} \otimes_{A_{n+1}} F_{n+1} \os{\simeq}{\lra} F_{n}$. 
Then we have $R^q\varprojlim \{F_n\}_n = 0$ for $q>0$. \\
\quad \\ 
We can prove this claim in the same way as \cite[0.3.8]{ogus2}: 
It suffices to prove $R^1\varprojlim \allowbreak \{F_n\}_n \allowbreak 
= 0$, and to show this, 
we may assume the existence of lattices $L_n \subseteq F_n \,(n \in \N)$ 
satisfying $\psi_n(L_{n+1}) \subseteq L_n$. 
Since the image of the map $\psi_n$ is dense by Lemma \ref{dense}, 
for any $(f_n)_n \in \prod_n F_n$, we have $(g_n)_n \in \prod_n F_n$ and 
$(f'_n)_n \in \prod p^nL_n$ such that 
$f_n + g_n = \psi(g_{n+1}) + f'_n$ holds for any $n \in \N$. 
Then $(f_n)_n$ is cohomologous to $(f'_n)_n$. Moreover, if 
we put $h_n := -\sum_{m=n}^{\infty}f'_m \in F_n$ (the infinite sum converges 
by definition of $(f'_n)_n$), we have 
$f'_n + h_n = \psi_n(h_{n+1})$. So $(f'_n)_n$ is cohomologous to $0$. 
So we are done. 
\end{pf} 

We define the notion of the log tubular neighborhood of 
radius $\lam$, as follows. 

\begin{defn} 
Let $\lam \in (0,1], \lam=p^{-r}$ and 
let $(Z,M_Z) \hra (\cZ,M_{\cZ})$ be a closed immersion of 
a fine log $\cB$-scheme to a fine log formal $\cB$-scheme. 
$($Then we can regard $((\cZ,M_{\cZ}), (Z,M_Z))$ as a widening on 
$(\cB,M_{\cB})/(\cB,M_{\cB})$.$)$ Then we define the log tubular 
neighborhood $]Z[^{\log}_{\cZ,\lam}$ of radius $\lam$ of 
$(Z,M_Z)$ in $(\cZ,M_{\cZ})$ as the rigid analytic space 
$\bigcup_{(a,b) \in \A_r} T_{a,b}(\cZ)_K$. We define the 
specialization map $\sp: ]Z[^{\log}_{\cZ,\lam} \lra \cZ$ as 
the union of composites 
$T_{a,b}(\cZ)_K \lra T_{a,b}(\cZ) \lra \cZ$, where the first map 
is the specialization map for $T_{a,b}(\cZ)$. 
\end{defn} 

Next, we show 
the relation between the category of isocrystals and 
certain categories of stratifications and 
give the definition of the log de Rham complex on tubular neighborhood 
associated to an isocrystal. (This is also a `$\lam$-version' of 
what was shown in \cite[2.21, 2.22]{shiho3}.) 
Let us consider the situation
\begin{equation}
\begin{CD}
(X,M_{X}) @>i>> (\cP,M_{\cP}) \\
@VfVV @VgVV \\
(\cY,M_{\cY}) @= (\cY,M_{\cY}), 
\end{CD}
\end{equation}
where $(X,M_X)$ is an object in $\L$, 
$f$ is a morphism in $\pLF$, $i$ is a closed immersion in 
$\pLF$ and $g$ is a formally log smooth morphism 
in $\pLF$. 
For $n \in \N$, 
let $(\cP(n),M_{\cP(n)})$ be the $(n+1)$-fold 
fiber product of $(\cP,M_{\cP})$ over $(\cY,M_{\cY})$, let 
$i(n): (X,M_{X}) \hra (\cP(n),M_{\cP(n)})$ be the closed immersion 
induced by $i$ and denote the pre-widening 
$((\cP(n),M_{\cP(n)}),(X,M_X),i(n),\id)$ of $(X,M_{X})/(\cY,M_{\cY})$ 
simply by $\cP(n)$. Then we have the system 
$\{T_{a,b}(\cP(n))\}_{(a,b) \in \A}$, 
the projections 
$$ p_{i,(a,b)}: T_{a,b}(\cP(1)) \lra T_{a,b}(\cP) 
\,\,\,\, (i=1,2), $$ 
$$ p_{ij,(a,b)}: T_{a,b}(\cP(2)) \lra T_{a,b}(\cP(1)) 
\,\,\,\, (1 \leq i < j \leq 3) $$
and the diagonal map 
$\Delta_{(a,b)}: T_{a,b}(\cP) \lra T_{a,b}(\cP(1)).$  
For $\lam \in (0,1]$ with $\lam=p^{-r}$, 
let $\Strat'_{\lam}((X \hra \cP/\cY)^{\log})$ 
be the category 
of compatible family of isocoherent sheaves $E_{a,b}$ on 
$T_{a,b}(\cP) \, ((a,b) \in \A_r)$ endowed with compatible isomorphisms 
$\epsilon_{a,b}: p^*_{2,(a,b)}E_{a,b} \overset{\sim}{\lra} 
p^*_{1,(a,b)}E_{a,b}$ satisfying 
$\Delta_{(a,b)}^*(\epsilon_{a,b}) = \id$, 
$p_{12,(a,b)}^*(\epsilon_{a,b}) \circ 
p_{23,(a,b)}^*(\epsilon_{a,b}) = p_{13,(a,b)}^*(\epsilon_{a,b}).$ \par 
On the other hand, from the closed immersion $i(n)$, we can form the 
tubular neighbouhood $]X[^{\log}_{\cP(n),\lam}$ and 
we have the projections 
$$ p_i: \, ]X[^{\log}_{\cP(1),\lam} \lra ]X[^{\log}_{\cP,\lam} \,\,\,\, 
(i=1,2), \quad   
p_{ij}: \, ]X[^{\log}_{\cP(2),\lam} \lra ]X[^{\log}_{\cP(1),\lam} \,\,\,\, 
(1 \leq i < j \leq 3) $$
and the diagonal map 
$\Delta: \, ]X[^{\log}_{\cP,\lam} \lra ]X[^{\log}_{\cP(1),\lam}.$
Let $\Strat''_{\lam}((X \hra \cP/\cY)^{\log})$ be the 
category of pairs $(E,\epsilon)$, where $E$ is a coherent 
${\cal O}_{]\cX[^{\log}_{\cP,\lam}}$-module and $\epsilon$ is an 
${\cal O}_{]\cX[^{\log}_{\cP(1),\lam}}$-linear isomorphism 
$p_2^*E \overset{\sim}{\lra} p_1^*E$ satisfying 
$\Delta^*(\epsilon) = \id$, 
$p_{12}^*(\epsilon) \circ p_{23}^*(\epsilon) = 
p_{13}^*(\epsilon)$. Then we have the following proposition: 

\begin{prop}
With the above notation, we have the functorial equivalence of categories 
\begin{equation*}
I_{\conv}(\cX/\cY)^{\log}) \simeq 
\Strat'((\cX \hra \cP/\cY)^{\log}) \simeq 
\Strat''((\cX \hra \cP/\cY)^{\log}). 
\end{equation*}
\end{prop}

\begin{pf} 
We can prove the proposition 
in the same way as \cite[2.2.7]{shiho2}, once we note the 
following fact: If we have an enlargement $\cZ:=((\cZ,M_{\cZ}),(Z,M_Z),j,z)$ 
with $\lam(\cZ) < \lam = p^{-r}$ 
which fits into the commutative diagram 
\begin{equation*}
\begin{CD}
(Z,M_Z) @>j>> (\cZ,M_{\cZ}) \\ 
@VzVV @VVV \\ 
(X,M_X) @>i>> (\cP,M_{\cP}), 
\end{CD}
\end{equation*}
there exists $(a,b) \in \A_r$ such that the map $(\cZ,M_{\cZ}) \lra 
(\cP,M_{\cP})$ factors through $(T_{a,b}(\cP), M_{T_{a,b}(\cP)})$. 
\end{pf} 

Let $\cE$ be an isocrystal on $(X/\cY)^{\log}_{\lamm\conv}$ and let 
$(E,\epsilon)$ be the associated object in 
$\Strat''_{\lam}((X \hra \cP/\cY)^{\log})$. 
If we denote the first log infinitesimal neighborhood of 
$(\cP,M_{\cP})$ into $(\cP(1),M_{\cP(1)})$ by $(\cP^1, M_{\cP^1})$, then 
$\epsilon$ induces the $\cO_{]X[^{\log}_{\cP^1,\lam}}$-linear isomorphism 
$$ \epsilon_1: 
\cO_{]X[^{\log}_{\cP^1,\lam}} \otimes_{\cO_{]X[^{\log}_{\cP,\lam}}} E \lra E 
\otimes_{\cO_{]X[^{\log}_{\cP,\lam}}} \cO_{]X[^{\log}_{\cP^1,\lam}} $$ 
and it induces the log connection 
$$ \nabla: E \lra E \otimes_{{\cal O}_{]X[^{\log}_{\cP,\lam}}} 
\omega^1_{]X[^{\log}_{\cP,\lam}/\cY_K} $$
by $\nabla(e) := \epsilon_1(1 \otimes e) - e \otimes 1$. 
Then we have the following lemma (cf.\cite[2.22]{shiho3}). 

\begin{lem}
The log connection $\nabla$ above is integrable. 
\end{lem}

\begin{pf} 
The proof is similar to that of \cite[2.22]{shiho3}: When there exists a 
chart of the log scheme $(\cP,M_{\cP})$, we can prove the lemma in 
the same way as \cite[1.2.7, 1.2.8]{shiho2}. So it suffices to prove that 
we may work etale locally on $\cP$. \par 
Let $\{(E_{a,b},\epsilon_{a,b})\}_m$ be the object in 
$\Strat'_{\lam}((X\hra\cP/\cY)^{\log})$ associated to 
$\cE$, let 
$({T_{a,b}(\cP)^n}',M_{{T_{a,b}(\cP)^n}'})$ 
be the $n$-th log infinitesimal neighborhood 
of $(T_{a,b}(\cP),M_{T_{a,b}(\cP)})$ in 
$(T_{a,b}(\cP),M_{T_{a,b}(\cP)}) \allowbreak 
\times_{(\cY,M_{\cY})} (T_{a,b}(\cP),M_{T_{a,b}(\cP)})$ and 
let $(T_{a,b}(\cP)^n,M_{T_{a,b}(\cP)^n}) \hra \allowbreak 
(T_{a,b}\allowbreak {(\cP)^n}',M_{{T_{a,b}(\cP)^n}'})$ be the exact closed immersion 
defined by the ideal $\{ x \in \cO_{{T_{a,b}(\cP)^n}'}\,\vert\, \allowbreak 
\exists m,\,
p^mx=0\}$. 
(Then, if we denote $T_{a,b}(\cP):=((T_{a,b}(\cP),\allowbreak 
M_{T_{a,b}(\cP)}), \allowbreak 
(X_{a,b},M_{X_{a,b}}))$ as enlargement, 
$((T_{a,b}(\cP)^n,M_{T_{a,b}(\cP)^n}),(X_{a,b},M_{X_{a,b}}))$ 
is regarded also as an enlargement, which we denote by $T_{a,b}(\cP)^n$.) 
Let us consider the morphism of pre-widenings 
\begin{equation}\label{map1}
T_{a,b}(\cP)^n \lra T_{a,b}(\cP) \times T_{a,b}(\cP) \lra \cP(1). 
\end{equation} 
If we put $\cI := \Ker (\cO_{T_{a,b}(\cP)^n} \ra \cO_{X_{a,b}})$ and 
$\cJ := \Ker (\cO_{T_{a,b}(\cP)^n} \ra \cO_{T_{a,b}(\cP)})$, 
we have $\cJ^{n+1}=0, (\cI/\cJ)^a \subseteq p^b(\cO_{T_{a,b}(\cP)^n}/\cJ)$. 
Hence, for any $N \geq n$, we have 
$$ 
\cI^{Na} \subseteq (p^b\cO_{T_{a,b}(\cP)^n}+\cJ)^N 
\subseteq p^{b(N-n)}\cO_{T_{a,b}(\cP)^n} + \cJ^{n+1} = 
p^{b(N-n)}\cO_{T_{a,b}(\cP)^n}. 
$$ 
So, by the universality of blow-up, the morphism 
$T_{a,b}(\cP)^n \lra \cP(1)$ in \eqref{map1} factors through 
$T_{Na,(N-n)b}(\cP)$. If we take $N$ sufficiently large, 
we have $(N-n)b/(Na) > r$. So, by pulling back the isomorphism 
$\epsilon_{Na,(N-n)b}$, we obtain the isomorphism 
$$ \epsilon'_{a,b,n}: \cO_{T_{a,b}(\cP)^n} \otimes_{\cO_{T_{a,b}(\cP)}} 
E_{a,b} \os{\sim}{\lra} E_{a,b} \otimes_{\cO_{T_{a,b}(\cP)}} 
\cO_{T_{a,b}(\cP)^n} $$
and $\epsilon'_{a,b,1}$ induces the connection 
$$ \nabla_{a,b}: E_{a,b} \lra E_{a,b} \otimes_{\cO_{T_{a,b}(\cP)}} 
\omega^1_{T_{a,b}(\cP)/\cY}. $$
Then it is easy to see that the compatible family 
$\{(E_{a,b},\nabla_{a,b})\}$ induces $(E,\nabla)$ via the equivalence 
\begin{equation*}
\left( 
\begin{aligned}
& \text{compatible family of} \\
& \text{isocoherent sheaves on} \\
& \{ T_{a,b}(\cP) \}_{(a,b) \in \A_r} 
\end{aligned}
\right)
\overset{\sim}{\lra} 
\left( 
\begin{aligned}
& \text{coherent} \\
& \text{${\cal O}_{]X[^{\log}_{\cP,\lam}}$-module} \\
\end{aligned}
\right).
\end{equation*}
So it suffices to prove the integrability of $(E_{a,b},\nabla_{a,b})$ to prove 
the lemma and we may work etale locally on $\cP$ to check it. 
So we are done. 
\end{pf}

By the above lemma, 
we can define the log de Rham complex 
{\tiny{
$$ \DR(]X[^{\log}_{\cP,\lam}/\cY_K,{\cal E}) := [0 @>>> E @>{\nabla}>> 
E \otimes_{{\cal O}_{]X[^{\log}_{\cP,\lam}}} 
\omega^1_{]X[^{\log}_{\cP,\lam}/\cY_K} 
@>{\nabla}>> \cdots @>{\nabla}>> 
E \otimes_{{\cal O}_{]X[^{\log}_{\cP,\lam}}} 
\omega^q_{]X[^{\log}_{\cP,\lam}/\cY_K} @>{\nabla}>> \cdots] $$ }}
on $]X[^{\log}_{\cP,\lam}/\cY_K$ 
associated to the isocrystal ${\cal E}$ 
in standard way.

Next we give a proof of the log convergent Poincar\'{e} lemma for 
relative log convergent cohomology of radius $\lam$. 

Let $\lam \in (0,1], \lam=p^{-r}$ and take a good sequence 
$\a$ in $\A_r$. 
Let $(X,M_X)$ be an object in $\L$, let 
$(X,M_{X}) \lra (\cY,M_{\cY})$ be a morphism in 
$\pLF$ and let $\cZ$ be a widening on $(X,M_{X})/(\cY,M_{\cY})$. 
Then we have morphisms of topoi 
$$ u: (X/\cY)^{\log,\sim}_{\lamm\conv,\Zar} \lra X^{\sim}_{\Zar}, \quad 
j_{\cZ}: (X/\cY)^{\log,\sim}_{\lamm\conv,\Zar} |_{\cZ} \lra 
(X/\cY)^{\log,\sim}_{\lamm\conv,\Zar} $$
as in \cite[pp.90--91]{shiho2}, \cite[\S 2]{shiho3}. 
Also, we have 
a morphism of ringed topoi 
$$ \phi_{\cZ}: 
(X/\cY)^{\log,\sim}_{\conv,\Zar} |_{\cZ} \lra \cZ^{\sim}_{\Zar} $$
defined by 
$(\phi_{\cZ,*}E)(U) := \varprojlim_{(a,b) \in \A_r}E(T_{a,b}(U)) = 
\varprojlim_{(a,b) \in \a}E(T_{a,b}(U)), 
\phi_{\cZ}^*E(g:T \ra \cZ) := g^*E(T)$, and a functor 
$$\phi_{\vec{\cZ}_{r,\a},*}: 
(\cX/\cY)^{\log,\sim}_{\lamm\conv,\Zar}\vert_{\cZ} 
  \lra \vec{\cZ}^{\sim}_{r,\a} $$
defined by $\phi_{\vec{\cZ}_{r,\a},*}E(U) := E(U \os{\subset}{\ra} 
T_{a,b}(\cZ) \ra \cZ)$.  
Then, by the same agrument as in 
\cite[2.28]{shiho3} and \cite[2.3.4]{shiho2}, we have 
the following (we omit the proof): 

\begin{prop}\label{acyclic}
Let $\cZ$ be an affine widening and $\cE$ an isocrystal on
$(X/\cY)^{\log,\sim}_{\lamm\conv}\vert_{\cZ}$. Then $j_{\cZ,*}\cE$ is 
$u_*$-acyclic.
\end{prop}

\begin{rem} 
In the course of the proof of Proposition 
\ref{acyclic}, we use Proposition \ref{gamma-acyclic}. 
\end{rem} 

Now let us consider the situation
\begin{equation}\label{embeddable}
\begin{CD}
(X,M_{X}) @>i>> (\cP,M_{\cP}) \\
@VfVV @VgVV \\
(\cY,M_{\cY}) @= (\cY,M_{\cY}), 
\end{CD}
\end{equation}
where $(X,M_X)$ is an object in $\L$, 
$f$ is a morphism in $\pLF$, $i$ is a closed immersion in 
$\pLF$ and $g$ is a formally log smooth morphism 
in $\pLF$. 
Then $\cP:=((\cP,M_{\cP}),(X,M_X),i,\id)$ is a pre-widening 
of $(X,M_{X})/(\cY,M_{\cY})$. Denote the associated widening of 
$\cP$ by $\hat{\cP}$ and, 
for an isocrystal $\cE$ on 
$(X/\cY)^{\log}_{\lamm\conv}$, put 
$$ \omega^i_{\hat{\cP}}(\cE) := 
j_{\hat{\cP},*}(j^*_{\hat{\cP}}\cE \otimes_{\cO_{\cX/\cY}} 
\phi^*_{\hat{\cP}}(\omega^i_{\cP/\cY} \vert_{\hat{P}})). $$
Then we have the following theorem 
(`radius $\lam$ version' of \cite[2.29]{shiho3}): 

\begin{thm}\label{lcpl1}
Let the notations be as above. Then there exists 
a canonical structure of complex on 
$\omega^{\b}_{\hat{\cP}}(\cE)$ and 
the 
adjoint homomorphism $\cE \lra j_{\hat{\cP},*}j^*_{\hat{\cP}}\cE 
\allowbreak = \omega^0_{\hat{\cP}}(\cE)$ induces the quasi-isomor\-phism 
$\cE \overset{\simeq}{\lra} \omega^{\b}_{\hat{\cP}}(\cE).$ 
\end{thm}

\begin{pf} 
The proof is similar to that of \cite[2.29]{shiho3}. First, one can check that, for an enlargement $\cZ$ with $\lam(\cZ)<\lam$, 
the sheaf $\omega^i_{\hat{\cP}}(\cE)_{\cZ}$ on $\cZ_{\Zar}$ induced by 
$\omega^i_{\hat{\cP}}$ is given by 
$$ \omega^i_{\hat{\cP}}(\cE)_{\cZ} = 
\vpl_{(a,b) \in \A_r} \pi_{\cZ,a,b,*}(\cE_{T_{a,b}(\cZ \times \cP)} 
\otimes_{\cO_{T_{a,b}(\cZ \times \cP)}} \pi_{\cP,a,b}^*\omega^i_{\cP/\cY}), $$
where $r:=-\log_p \lam$, $\cZ \times \cP$ is the product of $\cZ$ and $\cP$ 
taken in the category of pre-widenings, 
$\pi_{\cP,a,b}$ is the canonical map $T_{a,b}(\cZ \times \cP) \lra \cP$ and 
$\pi_{\cZ,a,b}$ is the canonical map $T_{a,b}(\cZ \times \cP) \lra \cZ$. 
To define a canonical structure of a complex on $\omega^{\b}_{\hat{\cP}}(\cE)$, it suffices to construct 
a canonical, functorial structure of a complex on 
$\omega^{\b}_{\hat{\cP}}(\cE)_{\cZ}$ for affine enlargements 
$\cZ := ((\cZ,M_{\cZ}),(Z,M_Z),i,z)$ with $\lam(\cZ)<\lam$. 
Let $(\cP^m,M_{\cP^m})$ 
be the $m$-th log 
infinitesimal neighborhood of $(\cP,M_{\cP})$ in $(\cP,M_{\cP}) 
\times_{(\cY,M_{\cY})} (\cP,M_{\cP})$, let $p_{i,m}:(\cP^m, M_{\cP^m}) \lra 
(\cP,M_{\cP}) \, (i=1,2)$ be the morphism induced by the $i$-th 
projection and put $X_i^m := p_{i,m}^{-1}(X), 
M_{X_i^m} := M_X |_{X_i^m}$. Then we have pre-widenings 
$$ \cP^m := ((\cP^m,M_{\cP^m}), (X,M_X)), \quad 
\cP_i^m := ((\cP^m,M_{\cP^m}), (X_i^m,M_{X_i^m})), $$
and the following diagram of pre-widenings (for $i=1,2$): 
\begin{equation*}
\begin{CD}
\cP^m @>>> \cP_i^m \\
@VVV @V{p_{i,m}}VV \\
\cP @= \cP. 
\end{CD}
\end{equation*}
It induces the diagram of pre-widenings
\begin{equation}\label{zzzz}
\begin{CD}
\cP^m \times \cZ @>{r_i}>> \cP^m_i \times \cZ \\
@VVV @V{q_i}VV \\
\cP \times \cZ @= \cP \times \cZ 
\end{CD}
\end{equation}
for $i=1,2$. Now let us consider the following diagram ($i=1,2$)
{\tiny{ \begin{equation*}
\begin{CD} 
(\cP,M_{\cP}) \times_{(\cY,M_{\cY})} (\cZ,M_{\cZ}) 
@>{\Delta}>> 
(\cP^m,M_{\cP^m}) \times_{(\cY,M_{\cY})} (\cZ,M_{\cZ}) 
@>{q_i}>> 
(\cP,M_{\cP}) \times_{(\cY,M_{\cY})} (\cZ,M_{\cZ}) \\ 
@AAA @AAA @AAA \\ 
(Z,M_Z) @>{\Delta'}>> (X^m_i,M_{X^m_i}) \times_{(X,M_X)} (Z,M_Z) 
@>{q'_i}>> (Z,M_Z), 
\end{CD}
\end{equation*}}}
where $q_i,q'_i$ are 
the morphism induced by the morphism of pre-widenings $q_i$ 
in \eqref{zzzz} and $\Delta,\Delta'$ are the morphisms induced by the 
diagonal $(\cP,M_{\cP}) \lra (\cP^m,M_{\cP^m})$. Then it is easy to see 
that the right square and the large rectangle are Cartesian. Hence 
the left square is also Cartesian. Hence, by Proposition \ref{2.1.27}, 
$r_i \,(i=1,2)$ induce 
the isomorphisms of inductive systems 
\begin{equation}\label{compa1}
\{T_{a,b}(\cP^m \times \cZ)\}_{(a,b) \in \A_r} 
\isom \{T_n(\cP_i^m \times \cZ) \}_{(a,b) \in \A_r}. 
\end{equation}
On the other hand, by Lemma \ref{2.1.26}, we see that 
$q_i \, (i=1,2)$ induce the isomorphisms 
\begin{equation}\label{str-eq}
T_{a,b}(\cP_1^m \times \cZ) 
\isom T_{a,b}(\cP \times \cZ) \times_{\cP} \cP^m, \quad 
T_{a,b}(\cP_2^m \times \cZ) 
\isom \cP^m \times_{\cP} T_{a,b}(\cP \times \cZ).
\end{equation}

By evaluating $\cE$ on $\{T_{a,b}(\cP \times \cZ)\}_{(a,b) \in \A_r}$, 
we see that $\cE$ naturally induces a coherent sheaf $E$ on 
$]Z[^{\log}_{\cP^m \times \cZ,\lam} = 
\bigcup_{(a,b) \in \A_r} T_{a,b}(\cP \times \cZ)_K$, and the isomorphisms 
\eqref{compa1}, \eqref{str-eq} induce the isomorphisms 
$$ \theta_m: \pi_{\cP}^*\cO_{\cP^m} 
\otimes_{\cO_{]Z[^{\log}_{\cP \times \cZ,\lam}}} E 
\os{\simeq}{\lra} 
E \otimes_{\cO_{]Z[^{\log}_{\cP \times \cZ,\lam}}} \pi_{\cP}^*\cO_{\cP^m} 
\quad (m \in \N). $$
(Here $\pi_{\cP}^*$ denotes the functor 
$\Coh(\cO_{\cP}) \lra \Coh(\cO_{]Z[^{\log}_{\cP \times \cZ,\lam}})$ 
defined as composite 
$$\Coh(\cO_{\cP}) \lra \Coh(\Q \otimes \cO_{\cP}) \simeq 
\Coh(\cO_{\cP_K}) \lra \Coh(\cO_{]Z[^{\log}_{\cP \times \cZ,\lam}}),$$ 
where the last arrow is the functor induced by the map 
$]Z[^{\log}_{\cP \times \cZ,\lam} \lra \cP_K$.) Then we can define the map 
$\ti{d}: E \lra E \otimes_{\cO_{]Z[^{\log}_{\cP \times \cZ,\lam}}} 
\pi_{\cP}^*\omega^1_{\cP/\cY}$ by 
$\wt{d}(e) := \theta_1(1\otimes e) - e \otimes 1$, and extend it 
to the diagram 
$$ E @>{\wt{d}}>> E \otimes_{\cO_{]Z[^{\log}_{\cP \times \cZ,\lam}}} 
\pi_{\cP}^*\omega^1_{\cP/\cY} @>{\wt{d}}>> \cdots @>{\wt{d}}>> 
E \otimes_{\cO_{]Z[^{\log}_{\cP \times \cZ,\lam}}} 
\pi_{\cP}^*\omega^q_{\cP/\cY} @>{\wt{d}}>> \cdots $$ 
in standard way. By applying the direct image by the map 
$\pi_{\cZ}: ]Z[^{\log}_{\cP \times \cZ,\lam} \lra \cZ_K \lra \cZ$ 
to the above diagram, we obtain the diagram 
$$\omega^{\b}_{\hat{\cP}}(\cE)_{\cZ} := 
[\omega^0_{\hat{\cP}}(\cE)_{\cZ} @>d>> \omega^1_{\hat{\cP}}(\cE)_{\cZ} 
@>d>> \cdots @>d>> \omega^q_{\hat{\cP}}(\cE)_{\cZ} @>d>> \cdots ], 
$$
since we have 
$$ 
\pi_{\cZ,*} (E \otimes_{\cO_{]Z[^{\log}_{\cP \times \cZ,\lam}}} 
\pi_{\cP}^*\omega^q_{\cP/\cY}) = 
\vpl_{(a,b)\in \A_r} \pi_{\cZ,a,b,*}(\cE_{T_{a,b}(\cZ \times \cP)} 
\otimes_{\cO_{T_{a,b}(\cZ \times \cP)}} \pi_{\cP,a,b}^*\omega^q_{\cP/\cY})
$$ 
by definition. 
This construction is functorial with respect to affine enlargement 
$\cZ$ and so it induces the diagram 
$$ \omega^{\b}_{\hat{\cP}}(\cE) := 
[\omega^0_{\hat{\cP}}(\cE) @>d>> \omega^1_{\hat{\cP}}(\cE) 
@>d>> \cdots @>d>> \omega^q_{\hat{\cP}}(\cE) @>d>> \cdots ]. $$
We should prove that the diagram $\omega^{\b}_{\hat{\cP}}(\cE)$ 
forms a complex and that the adjoint map 
$\cE \allowbreak \lra \allowbreak 
j_{\hat{\cP},*}j^*_{\hat{\cP}}\cE = \omega^0_{\hat{\cP}}(\cE)$ 
induces the quasi-isomorphism 
$ \cE \overset{\sim}{\lra} \omega^{\b}_{\hat{\cP}}(\cE)$. 
We can prove this in the same way as \cite[2.29]{shiho3}, once we establish 
Lemmas \ref{add1}, \ref{add2} below. So we are done. 
\end{pf} 

We give (a sketch of) the 
proof of the following two lemmas which is used in the 
above proof (cf. \cite[2.31, 2.32]{shiho3}): 

\begin{lem}\label{add1} 
Let $(X,M_X)$ be a fine log $\cB$-scheme and assume we are given the 
commutative diagram 
\begin{equation}\label{add11}
\begin{CD}
(X,M_X) @>>> (\cP,M_{\cP}) \\ 
@\vert @VfVV \\ 
(X,M_X) @>>> (\cQ,M_{\cQ}), 
\end{CD}
\end{equation}
where $(\cP,M_{\cP}), (\cQ,M_{\cQ})$ are $p$-adic fine log formal 
$\cB$-scheme, horizontal arrows are closed immersions and 
$f$ is formally log smooth. Then, Zariski locally on $\cQ$, 
we have the isomorphism 
$]X[^{\log}_{\cP,\lam} \cong 
]X[^{\log}_{\cQ,\lam} \times D_{\lam}^r$ for some $r$ 
$($where $D_{\lam}$ denotes the open disc of radius $\lam)$ such that 
the morphism $]X[^{\log}_{\cP,\lam} \lra ]X[^{\log}_{\cQ,\lam}$ 
induced by $f$ is 
identified with the first projection 
$]X[^{\log}_{\cQ,\lam} \times D_{\lam}^r \lra ]X[^{\log}_{\cQ,\lam}$ 
via this isomorphism. 
\end{lem}

\begin{pf} 
By localizing enough, we have 
$\cQ^{\ex} =: \Spf A, \cP^{\ex} =: \Spf B$ with 
$B \cong A[[x_1, \cdots, \allowbreak x_r]]$, 
by \cite[2.31]{shiho3}. 
Then one can see easily that 
$]X[^{\log}_{\cP,\lam} = ]X[^{\log}_{\cQ,\lam} \times D_{\lam}^r$ holds. 
\end{pf} 

\begin{lem}\label{add2}
Let $\cZ := \Spf A_0$ be a $p$-adic affine formal $\cB$-scheme and 
$\pi:\cZ_K \times D_{\lam}^r \lra \cZ_K$ be the projection. 
$($Here $D_{\lam}^r$ is as in the previous lemma.$)$
Then the complex 
\begin{equation}\label{add21}
\sp_*\cO_{\cZ_K} \lra 
(\sp \circ \pi)_*\Omega^{\b}_{\cZ_K \times D_{\lam}^r/\cZ_K} 
\end{equation}
is homotopic to zero. 
\end{lem} 

\begin{pf} 
The proof is the same as \cite[2.32]{shiho3}. 
The crucial point is that, 
for $f(x_1,x_2,\cdots, \allowbreak x_r) \in \Gamma(\cZ_K \times D_{\lam}^r, 
\cO_{\cZ_K \times D_{\lam}^r})$, the integral 
$\dint_0^{x_j}f(x_1, \cdots, \underset{\widehat{j}}{t}, \cdots x_r)dt$ 
is again contained in $\Gamma(\cZ_K \times D_{\lam}^r, \allowbreak 
\cO_{\cZ_K \times D_{\lam}^r})$. 
\end{pf}

Let $\ti{u}: (X/\cY)^{\log,\sim}_{\lamm\conv,\et} 
\lra X^{\sim}_{\Zar}$ be the 
composite 
$(X/\cY)^{\log,\sim}_{\lamm\conv,\et} \os{\epsilon}{\lra} 
(X/\cY)^{\log,\sim}_{\lamm\conv,\Zar} \allowbreak \os{u}{\lra}
X^{\sim}_{\Zar}$. Then we can prove the following propositions in 
the same way as \cite[2.3.8]{shiho2}, \cite[2.33]{shiho3}: 

\begin{cor}\label{lcpl3}
Assume that the diagram \eqref{embeddable} is given and 
let $\cE$ be an 
isocrystal on $(X/\cY)^{\log}_{\lamm\conv,\et}$.  
Then we have the canonical quasi-isomorphism 
$$ R\ti{u}_*\cE = \sp_* \DR(]X[^{\log}_{\cP,\lam}/\cY_K, \cE). $$
\end{cor} 

\begin{pf}
We may assume that $((\cP,M_{\cP}),(X,M_X),i,\id)$ is an affine 
widening. 
By Theorem \ref{lcpl1}, 
we have the quasi-isomorphism 
$$ R\ti{u}_*\cE \os{\sim}{\lra} 
Ru_* \omega^{\b}_{\hat{\cP}}(\epsilon_*\cE). $$ 
Next, by Proposition \ref{acyclic}, we have the quasi-isomorphism 
\begin{align*}
Ru_*\omega^{\b}_{\hat{\cP}}(\epsilon_*\cE) & = 
Ru_* j_{\hat{\cP},*}(\phi^*_{\hat{\cP}}(\omega^{\b}_{\cP/\cY} 
\vert_{\hat{\cP}}) \otimes j^*_{\hat{\cP}} \epsilon_*\cE) \\
& = 
u_* j_{\hat{\cP},*}(\phi^*_{\hat{\cP}}(\omega^{\b}_{\cP/\cY} 
\vert_{\hat{\cP}}) \otimes j^*_{\hat{\cP}} \epsilon_*\cE) 
 = u_*\omega^{\b}_{\hat{\cP}}(\epsilon_*\cE). 
\end{align*}
So, it suffices to prove the isomorphism 
\begin{equation}\label{usp}
\ti{u}_* \omega^{\b}_{\hat{\cP}}(\epsilon_*\cE) 
\isom \sp_* \DR(]X[^{\log}_{\cP,\lam}/\cY_K, \cE). 
\end{equation}
We can check that the each component of both hand sides in \eqref{usp} 
are equal, as in \cite[p.106 line ($-$5) -- p.107 line 7]{shiho2}. 
So it suffices to prove that 
the differentials are compatible. Let $\cP^1, \cP^1_i\,(i=1,2)$ are as in 
the proof of Theorem \ref{lcpl1} and let us consider the following diagram: 
\begin{equation*}
\begin{CD}
\cO_{\cP^1} \otimes_{\cO_{\cP}} \cE_{T_{a,b}(\cP)} @>>> \cE_{T_{a,b}(\cP^1)} 
@<<< \cE_{T_{a,b}(\cP)} \otimes_{\cO_{\cP}} \cO_{\cP^1} \\ 
@VVV @\vert @VVV \\ 
\cE_{T_{a,b}(\cP^1_2)} @>>> \cE_{T_{a,b}(\cP^1)} @<<< \cE_{T_{a,b}(\cP^1_1)}. 
\end{CD}
\end{equation*}
Note that the vertical arrows are isomorphisms, and the lower horizontal 
arrows become isomorphic after we take the inverse limit on 
$(a,b) \in \A_r$. So all the morphisms are isomorphic 
after we take the inverse limit on $(a,b) \in \A_r$. 
The differential on 
$\sp_*\DR(]X[^{\log}_{\cP}/\cY_K, \cE)$ is, by definition, induced by 
$\varprojlim_{(a,b) \in \A_r}$ of the upper horizontal arrows. 
On the other hand, $\varprojlim_{(a,b) \in \A_r}$ of the lower 
horizontal arrows is compatible with the 
isomorphism $\pi_{\cZ,*}\theta_1$ (where $\pi_{\cZ},\theta_1$ are as in the 
proof of Theorem \ref{lcpl1}), 
which induces the differential on 
$\omega^{\b}_{\hat{\cP}}(\epsilon_*\cE)_{\cZ}$ 
(for any affine enlargement $\cZ$ with $\lam(\cZ)<\lam$). 
So the differential on $\sp_*\DR(]X[^{\log}_{\cP}/\cY_K, \cE)$ is 
compatible with the differential on 
$\ti{u}_* \omega^{\b}_{\hat{\cP}}(\epsilon_*\cE)$. So 
we are done. 
\end{pf}

We can prove the following corollary in the same way as 
\cite[2.3.9]{shiho2},\cite[2.34]{shiho3} (we omit the proof). 

\begin{cor}\label{lcpl4}
Let $(X,M_X)$ be an object in the category $\L$, 
let $f:(X,M_{X}) \allowbreak \lra (\cY,M_{\cY})$ be a morphism in $\pLF$ and 
let $\cE$ be an isocrystal on 
$(X/\cY)^{\log}_{\lamm\conv}$. If we take an embedding system 
$$ (X,M_{X}) \overset{g}{\lla} (X^{(\bullet)}, M_{X^{(\bullet)}}) 
\overset{i}{\hra} (\cP^{(\bullet)},M_{\cP^{(\bullet)}}), $$
we have the isomorphism 
$$ R^nf_{X/\cY,\lamm\conv *}\cE \isom 
R^n(f \circ g)_* \sp_*\DR(]X^{(\b)}[^{\log}_{\cP^{(\b)},\lam}/\cY_K,\cE). $$
\end{cor}

Now we would like to compare relative log convergent cohomology 
of radius $\lam$ and relative log crystalline cohomology. 
Assume we are given the diagram
\begin{equation}\label{diag-p}
(X,M_X) \os{f}{\lra} (Y,M_Y) \os{\iota}{\hra} (\cY,M_{\cY}), 
\end{equation}
where $f$ is a log smooth morphism in $\LB$, 
 $(\cY,M_{\cY})$ is an object in $\pLF$ and $\iota$ is the exact closed 
immersion defined by $p\cO_{\cY}$. Denote the canonical PD-structure 
on $p\cO_{\cY}$ by $\gamma$. \par 
First we construct 
a functor from the category of isocrystals on log convergent site of 
radius $\lam$ to that on log crystalline site when $\lam$ is 
close to $1$: 

\begin{prop}\label{funct}
With the above notation, let $d := \max_{x \in X}\rk(\omega^1_{X/Y})$ 
and $\lam_0 := p^{-1/(26d(p-1)+1)}$. Then, for $\lam \in (\lam_0,1]$, 
we have the canonical functor 
$$ \Phi: I_{\lamm\conv}((X/\cY)^{\log}) \lra I_{\crys}((X/\cY)^{\log}) $$
such that $\Phi$ sends locally free isocrystals on 
$(X/\cY)^{\log}_{\lamm\conv}$ 
to locally free isocrystals on $(X/\cY)^{\log}_{\crys}$. 
\end{prop}

\begin{pf} 
The proof is similar to \cite[2.35]{shiho3}, 
but it is a little more complicated. 
First, let us consider the situation that $(X,M_X)$ admits a closed immersion 
$(X,M_X) \hra (\cP,M_{\cP})$ into a 
$p$-adic fine log formal $\cB$-scheme 
which is formally log smooth over $(\cY,M_{\cY})$. 
Denote the $(i+1)$-fold fiber product of $(\cP,M_{\cP})$ over 
$(\cY,M_{\cY})$ by $(\cP(i), \allowbreak M_{\cP(i)})$ and let 
$(\cP(i)^{\ex},M_{\cP(i)^{\ex}})$ be the exactification of the closed 
immersion $(X,M_X) \allowbreak \hra (\cP(i),M_{\cP(i)})$. 
Let $\cI(i):= \Ker (\cO_{\cP(i)}^{\ex} \lra \cO_X)$, 
let $B_{a,b}(i)$ be the formal blow-up of $\cP(i)^{\ex}$ with respect to 
the ideal $\cI(i)^a + p^b\cO_{\cP^{\ex}}$, let $T'_{a,b}(i)$ be the open 
sub formal scheme of $B_{a,b}(i)$ defined as the set of points 
$x \in B_{a,b}(i)$ satisfying 
$(\cI(i)^a+p^b\cO_{\cP(i)})\cO_{B_{a,b}(i),x} = p^b\cO_{B_{a,b}(i),x}$ 
and let 
$T_{a,b}(i)$ be the closed sub formal scheme of $T_{a,b}(i)$ 
defined by the ideal 
$\{x \in \cO_{T_{a,b}(i)} \,\vert\, p^nx=0 \,\text{for some $n>0$}\}$.
On the other hand, let $D(i)$ be 
the $p$-adically completed log PD-envelope of $(X,M_X)$ in 
$(\cP(i),M_{\cP(i)})$. Let 
$m=m(X,\cP)$ be the minimal integer such that, for $i=0,1,2$, 
$\Ker(\cO_{\cP^{\ex}(i)} \lra \cO_X)$ 
is generated locally by $p$ and some other $m$ elements. 
Then there exist canonical diagrams 
$D(i) \lra T'_{n,1}(i) \hookleftarrow T_{n,1}(i)$ 
for $n=(p-1)m+1$ and $i=0,1,2$, 
where the second map is the canonical closed immersion. Since we have 
the canonical equivalences of categories 
$\Coh(\Q \otimes \cO_{T_{n,1}(i)}) \cong 
\Coh(\Q \otimes \cO_{T'_{n,1}(i)})$, 
we can define, by the `pull-back by $D(i) \lra T'_{n,1}(i)$', 
the functor 
$$
I_{\lamm\conv}((X/\cY)^{\log}) \os{\sim}{\lra}
\Strat''((X \hra \cP/\cY)^{\log}) \lra 
\HPDI((X \hra \cP/\cY)^{\log}) $$
if we have $\lam > p^{-1/n} = p^{-1/((p-1)m(X,\cP)+1)}$. 
(Here $\HPDI((X \hra \cP/\cY)^{\log})$ is 
the category of HPD-isostratifications defined in \cite[1.17]{shiho3}.) 
\par 
Next, let us consider the general situation and take an embedding system 
\begin{equation}\label{convcrysembsys}
(X,M_{X}) \overset{g^{(\b)}}{\lla} (X^{(\bullet)}, M_{X^{(\bullet)}}) 
\overset{i}{\hra} (\cP^{(\bullet)},M_{\cP^{(\bullet)}})
\end{equation}
over $(\cY,M_{\cY})$ such that 
$g$ is a Zariski \v{C}ech hypercovering and that 
$(\cP^{(\b)},M_{\cP^{(\b)}}) = 
\cosk_0^{(\cY,M_{\cY})}(\cP^{(0)},M_{\cP^{(0)}})$ holds. 
Let 
$\Delta_{\lamm\conv}$ be the category of descent data on 
$I_{\lamm\conv} \allowbreak 
((X^{(n)}/ \allowbreak \cY)^{\log}) \,(n=0,1,2)$ 
(that is, the category of 
objects in $I_{\lamm\conv}((X^{(0)}/\cY)^{\log})$ endowed with 
isomorphism of glueing in $I_{\lamm\conv}((X^{(1)}/\cY)^{\log})$ satisfying 
the cocycle condition in $I_{\lamm\conv}(X^{(2)}/\cY)^{\log})).$ 
Similarly, let $\Delta_{\crys}$ be the 
category of descent data on $I_{\crys}((X^{(n)}/\cY)^{\log}) \, (n=0,1,2)$ 
and let $\Delta_{\HPDI}$ be the 
category of descent data on 
$\HPDI((X^{(n)} \hra \cP^{(n)}/\cY)^{\log}) \, (n=0,1,2)$. 
Let $m(X^{(\b)},\cP^{(\b)})$ be the integer 
$\max_{0 \leq i \leq 2}m( \allowbreak X^{(i)},\allowbreak \cP^{(i)})$. 
Then, for $\lam \in (p^{-1/((p-1)m(X^{(\b)},\cP^{(\b)})+1)},1]$, 
we have the diagram 
$$ \Delta_{\lamm\conv} \os{\Phi}{\lra} \Delta_{\HPDI} 
\os{\Lambda}{\lla} \Delta_{\crys},$$
where $\Lambda$ is the functor from the category of isocrystals on 
log crystalline site to the category of HPD-isostratifications defined 
in \cite[\S 1]{shiho3}. \par 
We know that the functor $\Lambda$ is fully faithful, and that, 
if we have 
\begin{equation}\label{liftable}
(\cP^{(0)},M_{\cP^{(0)}}) 
\times_{(\cY,M_{\cY})} (Y,M_Y) = (X^{(0)}, M_{X^{(0)}}), 
\end{equation} 
the functor $\Lambda$ is an equivalence of categories. 
(Note that it is possible to take an embedding system 
\eqref{convcrysembsys} satisfying the condition \eqref{liftable}.) 
So, under the condition \eqref{liftable}, we have the functor 
$$ \Phi: I_{\lamm\conv}((X/\cY)^{\log}) \os{\sim}{\lra} \Delta_{\conv} 
\os{\Psi}{\lra} \Delta_{\HPDI} \os{\Lambda^{-1}}{\lra} \Delta_{\crys} 
\os{\sim}{\lra} I_{\crys}((X/\cY)^{\log}) $$
for $\lam \in (p^{-1/((p-1)m(X^{(\b)},\cP^{(\b)})+1)},1]$. 
Moreover, under the condition \eqref{liftable}, we can give a bound for 
$m(X^{(\b)},\cP^{(\b)})$: Under \eqref{liftable}, 
$m(X^{(i)},\cP^{(i)})\,(i=0,1,2)$ is the minimal 
number of (local) generators of the closed immersion 
$$ 
X^{(i)} \hra {(\cP^{(i)}(2),M_{\cP^{(i)}(2)})}^{\ex}
\otimes_{\Z_p} \F_p  = 
((X^{(0)},M_{X^{(0)}}) 
\times_{(Y,M_Y)} \cdots \times_{(Y,M_Y)} (X^{(0)},M_{X^{(0)}}))^{\ex}, $$
where the fiber product is taken $3(i+1)$-times. So 
$m(X^{(i)},\cP^{(i)})$ is bounded by 
$(3i+2)d$. So we have $m(X^{(\b)},\cP^{(\b)}) \leq 8d$. Therefore, 
if we have $\lam > p^{-1/(8d(p-1)+1)}$, we have the functor 
$\Phi: I_{\lamm\conv}((X/\cY)^{\log}) \lra I_{\crys}((X/\cY)^{\log}).$ \par 
In order to prove the well-definedness of $\Phi$, we need to make 
$\lam$ closer to $1$: Assume we have another embedding system 
\begin{equation}\label{emb2}
(X,M_{X}) \overset{g^{(\b)}}{\lla} ({X'}^{(\bullet)}, M_{{X'}^{(\bullet)}}) 
\overset{i}{\hra} ({\cP'}^{(\bullet)},M_{{\cP'}^{(\bullet)}}) 
\end{equation}
as above (satisfying the analogue of \eqref{liftable}) and let us define 
\begin{equation}\label{emb3}
(X,M_{X}) \overset{g^{(\b)}}{\lla} ({X''}^{(\bullet)}, 
M_{{X''}^{(\bullet)}}) 
\overset{i}{\hra} ({\cP''}^{(\bullet)},M_{{\cP''}^{(\bullet)}})
\end{equation}
by 
$$({X''}^{(\b)},M_{{X''}^{(\b)}}) := 
(X^{(\b)},M_{X^{(\b)}}) \times_{(X,M_X)} ({X'}^{(\b)},M_{{X'}^{(\b)}}), $$
$$ ({\cP''}^{(\b)}, M_{{\cP''}^{(\b)}}) := 
(\cP^{(\b)},M_{\cP^{(\b)}}) \times_{(\cY,M_{\cY})} 
({\cP'}^{(\b)},M_{{\cP'}^{(\b)}}). $$ 
Then 
$m({X''}^{(i)},{\cP''}^{(i)})\,(i=0,1,2)$ is the minimal 
number of (local) generators of the closed immersion 
{\tiny{ $$ 
{X''}^{(i)} \hra {({\cP''}^{(i)}(2),M_{{\cP''}^{(i)}(2)})}^{\ex}
 \otimes_{\Z_p} \F_p  = 
(({X''}^{(0)},M_{{X''}^{(0)}}) 
\times_{(Y,M_Y)} \cdots \times_{(Y,M_Y)} 
({X''}^{(0)},M_{{X''}^{(0)}}))^{\ex}, $$ }}
where the fiber product is taken $6(i+1)$-times. 
So we have 
$m({X''}^{(\b)},{\cP''}^{(\b)}) \leq 17d$. Hence, for 
$\lam \in (p^{-1/(17d(p-1)+1)},1]$, 
we have the commutative diagram of functors 
\begin{equation}\label{welldef}
\begin{CD}
I_{\lamm\conv}((X/\cY)^{\log}) @>{\sim}>> 
\Delta_{\conv} @>{\Psi}>> \Delta_{\HPDI} @<{\Lambda}<< \Delta_{\crys} 
@>{\sim}>> I_{\crys}((X/\cY)^{\log}) \\ 
@\vert @VVV @VVV @VVV @\vert \\ 
I_{\lamm\conv}((X/\cY)^{\log}) @>{\sim}>> 
\Delta''_{\conv} @>{\Psi''}>> \Delta''_{\HPDI} @<{\Lambda''}<< 
\Delta''_{\crys} @>{\sim}>> I_{\crys}((X/\cY)^{\log}) \\ 
@\vert @AAA @AAA @AAA @\vert \\
I_{\lamm\conv}((X/\cY)^{\log}) @>{\sim}>> 
\Delta'_{\conv} @>{\Psi'}>> \Delta'_{\HPDI} @<{\Lambda'}<< 
\Delta'_{\crys} @>{\sim}>> I_{\crys}((X/\cY)^{\log}), 
\end{CD}
\end{equation}
where $\Delta'_{\conv}, ...$ are $\Delta_{\conv}, ...$ for the 
embedding system \eqref{emb2} and $\Delta''_{\conv}, ...$ 
are $\Delta_{\conv}, ...$ for the embedding system \eqref{emb3}. 
Then $\Lambda, \Lambda'$ are equivalences and $\Lambda''$ is fully 
faithful. From this, we see that the functor 
$\Phi: I_{\lamm\conv}((X/\cY)^{\log}) \lra  I_{\crys}((X/\cY)^{\log})$ 
induced by the top horizontal arrows in \eqref{welldef} is 
isomorphic to the functor 
$\Phi': I_{\lamm\conv}((X/\cY)^{\log}) \lra  I_{\crys}((X/\cY)^{\log})$ 
induced by the bottom horizontal arrows. 
We can prove, in the same way, that the isomorphism $\Phi \cong \Phi'$ 
we constructed is canonical (in the sense that it satisfies the 
cocycle condition) when we have $\lam \in (p^{-1/(26d(p-1)+1)},1]$. 
(We should consider the fiber product of three embedding systems.) 
So the functor $\Phi$ is canonically defined when we have 
$\lam \in (p^{-1/(26d(p-1)+1)},1]$. 
So the proof is finished. 
\end{pf} 

Now we prove the comparison theorem between relative log convergent 
cohomology of radius $\lam$ and relative log crystalline cohomology. 
Before the statement of the theorem, we introduce some notations. 

\begin{defn}\label{al}
Let us assume given the diagram as \eqref{diag-p}. Then we define 
the integers $a := a(X)$ and $l:=l(X,Y,\cY)$ as follows$:$ 
$a$ is the least positive integer such that $X$ admits a Zariski 
covering $X = \bigcup_{j=1}^aX_j$ by $a$ affine open subschemes. 
$l$ is the least positive integer such that $X$ admits 
a Zariski covering $X = \bigcup_{j=1}^lX_j$ by $l$ open subschemes, 
satisfying the following condition $(*):$ \\
$(*)$ \,\,\, There exists an exact closed immersion 
$(X_j,M_{X_j}) \hra (\cP_j,M_{\cP_j})$ over $(\cY,M_{\cY})$ such that 
$(\cP_j,M_{\cP_j})$ is formally log smooth over $(\cY,M_{\cY})$ and that 
$(\cP_j,M_{\cP_j}) \times_{(\cY,M_{\cY})} (Y,M_Y) \cong (X,M_X)$ holds. 
\end{defn} 

By \cite[(3.14)]{kkato}, we have $l(X,Y,\cY) \leq a(X)$. So they are finite.

\begin{thm}\label{conv-crys} 
Assume we are given the diagram \eqref{diag-p} and put 
$$\lam_0 := p^{-1/(\max(26,3l(X,Y,\cY)-1)d(p-1)+1)}.$$ 
Then, for 
$\lam \in (\lam_0,1]$ and $\cE \in I_{\lamm\conv}((X/\cY)^{\log})$, 
we have the quasi-isomorphism 
$$ Rf_{X/\cY,\conv *}\cE \isom Rf_{X/\cY,\crys *}\Phi(\cE). $$ 
\end{thm}

\begin{pf} 
Put $l:=l(X,Y,\cY)$ and 
let us take a Zariski covering $X = \bigcup_{j=1}^lX_j$ and an exact 
closed immersion 
$i_j: (X_j,M_{X_j}) \hra (\cP_j,M_{\cP_j})$ over $(\cY,M_{\cY})$ satisfying 
the condition $(*)$ in Definition \ref{al}. Put 
$J := \{1,2, \cdots, l\}$ (endowed with usual total order). 
For non-empty subset $L \subseteq J$, let 
$(X_L,M_{X_L})$ (resp. $(\cP_L,M_{\cP_L})$) be the fiber product of 
$(X_j,M_{X_j})$'s (resp. $(\cP_j,M_{\cP_j})$'s) for $j \in L$ over 
$(X,M_X)$ (resp. $(\cY,M_{\cY})$) and let $i_L: (X_L, M_{X_L}) \hra 
(\cP_L,M_{\cP_L})$ be the closed immersion induced by $i_j$'s ($j \in L$). 
Then, for $m \in \N$, let $(X^{(m)},M_{X^{(m)}})$ (resp. 
$(\cP^{(m)},M_{\cP^{(m)}})$) be the disjoint union of 
$(X_L,M_{X_L})$'s (resp. $(\cP_L,M_{\cP_L})$'s) for $L \subseteq J$ 
with $|L|=m+1$ and 
let $i^{(m)}:(X^{(m)},M_{X^{(m)}}) \hra (\cP^{(m)},M_{\cP^{(m)}})$ be the 
disjoint union of $i_L$'s for $L \subseteq L$ with $|L|=m+1$. 
Then $(X^{(\b)},M_{X^{(\b)}}) := \{(X^{(m)},M_{X^{(m)}})\}_{m \leq l-1}$ 
and $(\cP^{(\b)},M_{(\b)}) := \{(\cP^{(m)},M_{\cP^{(m)}})\}_{m \leq l-1}$ 
naturally forms a diagram indexed by the category $\Delta^+_{l-1}$ 
(where $\Delta^+_{l-1}$ is as in the paragraph before Proposition 
\ref{cohdes2}) and we have the diagram 
$$ 
(X,M_X) @<g<< (X^{(\b)},M_{X^{(\b)}}) @>{i^{(\b)}}>> 
(\cP^{(\b)},M_{\cP^{(\b)}}).  
$$ 
Let us denote the $p$-adically completed 
PD-envelope of the closed immersion $i^{(\b)}$ by 
$(D^{(\b)},M_{D^{(\b)}})$. \par 
For $\lam \in (p^{-1/(26d(p-1)+1)},1]$, we have the functor 
$\Phi:I_{\lamm\conv}(X/\cY) \lra I_{\crys}(X/\cY)$ by Proposition \ref{funct}. 
On the other hand, if we define $m(X^{(j)},\cP^{(j)})$ as in the proof of 
Proposition \ref{funct}, it is the minimal 
number of (local) generators of the closed immersion 
$$ 
X^{(j)} \hra {(\cP^{(j)}(2),M_{\cP^{(j)}(2)})}^{\ex} \otimes_{\Z_p} \F_p  = 
((X^{(0)},M_{X^{(0)}}) 
\times_{(Y,M_Y)} \cdots \times_{(Y,M_Y)} (X^{(0)},M_{X^{(0)}}))^{\ex}, $$
where the fiber product is taken $3(j+1)$-times. So 
$m(X^{(j)},\cP^{(j)})$ is bounded by $(3j+2)d$. 
So if we put $n:=(3l-1)d(p-1)+1$, 
we have the diagram 
\begin{equation}\label{dtt}
D^{(\b)}(i) \lra T'_{n,1}(\cP^{(\b)}(i)) 
\hookleftarrow T_{n,1}(\cP^{(\b)}(i))
\end{equation} 
(where the notation are as in the proof of Proposition \ref{funct}), 
and by using this diagram, we can define the functor 
$I_{\lamm\conv}(X/\cY) \lra I_{\crys}(X/\cY)$ for 
$\lam \in (p^{-1/(3l-1)d(p-1)+1}\allowbreak , \allowbreak 1]$, 
which is identical with $\Phi$ when 
$\lam \in (\lam_0,1]$. \par 
By the diagram \eqref{dtt}, we can define the map of complexes 
\begin{equation}\label{drmap}
\sp_*\DR(]X^{(\b)}[^{\log}_{\cP^{(\b)},\lam}/\cY, \cE) 
\lra 
\DR(\cD^{(\b)}/\cY,\Phi(\cE))
\end{equation}
as in \cite[\S 2]{shiho3}, which induces the map 
$$ Rf_{X/\cY,\lamm\conv *}\cE \lra Rf_{X/\cY,\crys *}\Phi(\cE). $$
 (Here $\DR(\cD^{(\b)}/\cY,\Phi(\cE))$ denotes the 
log de Rham complex on $\cD^{(\b)}$ associated to $\Phi(\cE)$. 
See \cite[\S 1]{shiho3}.) So it suffices to prove that the map 
\eqref{drmap} is a quasi-isomorphism. To see this, first we may replace 
$\b$ by $j$, and then, we may replace $\cP^{(j)}$ by the open subscheme 
of $\cP^{(0)}$ which is naturally homeomorphic to $X^{(j)}$, because 
the both hand sides are unchanged in derived category. 
Then we can reduce to the claim in the case $j=0$, 
and in this case, both hand sides coincide. So we are done. 
\end{pf} 

As a corollary, we have the following: 

\begin{cor}\label{perf1}
Assume we are given the diagram 
\begin{equation*}
(X,M_X) \os{f}{\lra} (Y,M_Y) \os{\iota}{\hra} (\cY,M_{\cY}), 
\end{equation*}
where $f$ is a proper log smooth integral morphism in $\LB$, 
 $(\cY,M_{\cY})$ is an object in $\pLF$ and $\iota$ is the exact closed 
immersion defined by the ideal sheaf $p{\cal O}_{\cY}$. 
Then, for $\lam \in (\lam_0,1]$ $(\lam_0$ is as in Theorem \ref{conv-crys}$)$
 and 
a locally free isocrystal $\cE$ on $(X/\cY)^{\log}_{\lamm\conv}$, 
$Rf_{X/\cY,\lamm\conv *}\cE$ is a perfect complex of 
$\Q \otimes_{\Z} \cO_{\cY}$-modules on $\cY_{\Zar}$, 
that is, it is quasi-isomorphic to 
a bounded complex of locally free $\Q \otimes_{\Z} \cO_{\cY}$-modules in 
the sense of \cite[1.9]{shiho3} Zariski locally on $\cY$. 
\end{cor}

\begin{pf}
This is the immediate consequence of Theorem \ref{conv-crys} and 
\cite[1.16]{shiho3}. 
\end{pf} 

\begin{cor}\label{bccrys}
Assume we are given a diagram 
\begin{equation}\label{eq1}
\begin{CD}
(X',M_{X'}) @>>> (Y',M_{Y'}) @>>> (\cY',M_{\cY'}) \\
@VVV @VVV @V{\varphi}VV \\
(X,M_X) @>f>> (Y,M_Y) @>{\iota}>> (\cY,M_{\cY}), 
\end{CD}
\end{equation}
where 
$f$ is proper log smooth integral, 
$\iota$ is the exact closed immersion defined by 
the ideal sheaf $p\cO_{\cY}$ and the squares are Cartesian. 
Then, for $\lam \in (\lam_0,1]$ $(\lam_0$ is as in Theorem \ref{conv-crys}$)$
 and a 
locally free isocrystal $\cE$ on $I_{\lamm\conv}((X/\cY)^{\log})$, 
we have the quasi-isomorphism 
$$ 
L\varphi^*
Rf_{X/\cY,\lamm\conv *}\cE \os{\sim}{\lra} 
Rf_{X'/\cY',\lamm\conv *} \varphi^*\cE. $$
\end{cor}

\begin{pf}
This is the immediate consequence of Theorem \ref{conv-crys} and 
\cite[1.19]{shiho3}. 
\end{pf} 


\section{Relative log convergent cohomology of radius $\lam$ (II)}

Assume we are given a diagram 
\begin{equation}\label{diag-s2}
(X,M_X) \os{f}{\lra} (Y,M_Y) \os{\iota}{\hra} (\cY,M_{\cY}), 
\end{equation}
where $f$ is a proper log smooth integral morphism in $\LB$, 
 $(\cY,M_{\cY})$ is an object in $\pLF$ and $\iota$ is a homeomorphic 
exact closed immersion in $\pLF$. 
In this section, we prove the coherence and the base change property 
for relative log convergent cohomology of radius $\lam$ of 
$(X,M_X)/(\cY,M_{\cY})$ 
when $f$ has log smooth parameter (\cite[3.4]{shiho3}, see also 
Definition \ref{parameter}) and $\lam$ is sufficiently close to $1$. 
The strategy of proof is the same as \cite[\S 3]{shiho3}, but we need 
more subtle argument to know for which $\lam$ the proof works. \par 
The first proposition we need is 
the topological invariance of the category of isocrystals on 
relative log convergent site with radius and that of 
relative log convergent cohomology with radius. 
Assume that we are given the diagram 
\begin{equation}\label{topinv}
\begin{CD}
(X,M_X) @>i>> (X_1,M_{X_1}) \\ 
@VfVV @V{f_1}VV  \\ 
(Y,M_Y) @>{\iota'}>> (Y_1,M_{Y_1}) @>{\iota_1}>> (\cY,M_{\cY}), 
\end{CD}
\end{equation} 
where $f,f_1$ are morphisms in $\L$, 
$i,\iota'$ are homeomorphic exact closed immersions in $\L$, 
$(\cY,M_{\cY})$ is an object in $\pLF$, $\iota_1$ is 
a homeomorphic exact closed immersion in $\pLF$ and the 
square is Cartesian. Put $\iota:= \iota_1 \circ \iota'$, 
$\cI_0:=\Ker(\iota^*\cO_{\cY} \ra \cO_{Y})$ and 
$r_0:=r(\cO_{\cY},\cI_0)$. Then we have, for any 
$\lam \in (0,1]$, the restriction functor 
$$ i^*: I_{\lamm\conv}((X_1/\cY)^{\log}) \lra I_{\lamm\conv}((X/\cY)^{\log}) $$
and for any $\lam,\lam' \in (0,1]$ with $\lam' > \lam$, we have another 
kind of restriction functors 
$$ 
I_{\lam'\text{-}\conv}((X/\cY)^{\log}) \lra 
I_{\lamm\conv}((X/\cY)^{\log}), 
\,\,\,\, 
I_{\lamm\conv}((X_1/\cY)^{\log}) \lra 
I_{\lam'\text{-}\text{-}\conv}((X_1/\cY)^{\log}), 
$$ 
which we denote by $\res_{\lam',\lam}$. Then we have the following 
propositions: 

\begin{prop}\label{21}
Let the notations be as above. Then, for any $\lam = p^{-r} \in (0,1]$ and 
$\lam' \in (p^{-(r^{-1}+r_0^{-1})^{-1}},1]$, 
there exists a functor 
$$ i_*: I_{\lam'\text{-}\conv}((X/\cY)^{\log}) \lra 
I_{\lamm\conv}((X_1/\cY)^{\log}) $$
satisfying the following conditions$:$ With $\lam,\lam'$ as above, 
the composite 
$$ I_{\lam'\text{-}\conv}((X/\cY)^{\log}) \os{i_*}{\lra} 
I_{\lamm\conv}((X_1/\cY)^{\log}) \os{i^*}{\lra} 
I_{\lamm\conv}((X/\cY)^{\log}) 
$$ 
is equal to $\res_{\lam',\lam}$ and the composite 
$$ 
I_{\lam'\text{-}\conv}((X_1/\cY)^{\log}) \os{i^*}{\lra} 
I_{\lam'\text{-}\conv}((X/\cY)^{\log}) \os{i_*}{\lra} 
I_{\lamm\conv}((X_1/\cY)^{\log}) 
$$ 
is equal to $\res_{\lam',\lam}$. 
\end{prop} 

\begin{pf} 
For $\cE \in I_{\lam'\text{-}\conv}((X/\cY)^{\log})$ and an enlargement 
$\cZ := ((\cZ,M_{\cZ}),(Z,M_Z),j,z)$ of $(X_1,M_{X_1})/(\cY,M_{\cY})$ with 
$\lam(\cZ) < \lam$, put 
$(Z',M_{Z'}) := (Z,M_Z) \times_{(X_1,M_{X_1})} (X,M_X)$ and let 
$j':(Z',M_{Z'}) \hra (\cZ,M_{\cZ}), z':(Z',M_{Z'}) \lra (X,M_X)$ be the 
morphisms naturally induced by $j,z$. Then 
$\cZ' := ((\cZ,M_{\cZ}),(Z',M_{Z'}),j',z')$ is an enlargement of 
$(X,M_X)/(\cY,M_{\cY})$. If we put $\cI:=\Ker(\cO_{\cZ} \ra \cO_{Z})$ and 
$\cJ:= \Ker(\cO_{\cZ} \ra \cO_{Z'})$, we have 
$\cJ \subseteq \cI + \cI_0\cO_{\cZ}$. So we have 
$$
r(\cO_{\cZ},\cJ) \geq 
(r(\cO_{\cZ},\cI)^{-1} + r(\cO_{\cZ},\cI_0\cO_{\cZ})^{-1})^{-1} \geq 
(r^{-1} + r_0^{-1})^{-1},  
$$ 
that is, $\lam(\cZ') < \lam'$. So we can define $i_*\cE(\cZ)$ by 
$i_*\cE(\cZ) := \cE(\cZ')$. So we have the functor $i_*$, and 
we can check the required property easily from this definition. 
\end{pf} 

\begin{prop}\label{22}
Assume given the diagram \eqref{topinv} and assume that $f_1$ is log smooth. 
Then, for $\lam \in (p^{-r_0},1]$ and $\cE \in 
I_{\lamm\conv}((X_1/\cY)^{\log})$, we have the quasi-isomorphism 
$$ 
Rf_{X_1/\cY,\lamm\conv *}\cE \os{\sim}{\lra} 
Rf_{X/\cY,\lamm\conv *}i^*\cE. 
$$ 
\end{prop}

\begin{pf}
We may work locally on $X_1$. So we may assume that there exists 
an exact closed immersion $(X_1,M_{X_1}) \hra (\cP,M_{\cP})$ over 
$(\cY,M_{\cY})$ such that $(\cP,M_{\cP})$ is formally log smooth 
over $(\cY,M_{\cY})$ and that we have the canonical isomorphism 
$(\cP,M_{\cP}) \times_{(\cY,M_{\cY})} (Y_1,M_{Y_1}) = (X_1,M_{X_1})$. 
Then we have 
$(\cP,M_{\cP}) \times_{(\cY,M_{\cY})} (Y,M_{Y}) = (X,M_{X})$, and from this, 
one can see the isomorphisms 
$]X[^{\log}_{\cP,\lam} \os{\cong}{\lra} ]X_1[^{\log}_{\cP,\lam} = \cP_K$
for $\lam \in (p^{-r_0},1]$. Then one can see easily that 
the canonical map of complexes 
$$ 
\sp_*\DR(]X_1[^{\log}_{\cP,\lam}/\cY_K,\cE) \lra 
\sp_*\DR(]X[^{\log}_{\cP,\lam}/\cY_K, i^*\cE)
$$ 
is an isomorphism. So we are done. 
\end{pf} 

We also have the following invariance property with respect to 
the radius: 

\begin{prop}\label{23}
Let us assume given the diagram 
$$
(X,M_X) \os{f}{\lra} (Y,M_Y) \os{\iota}{\hra} (\cY,M_{\cY}), 
$$ 
where $f$ is a log smooth morphism in $\L$, 
$(\cY,M_{\cY})$ is an object in $\pLF$ and $\iota$ is a homeomorphic 
exact closed immersion in $\pLF$. Put $r_0:=r(\cO_{\cY},
\Ker(\iota^*:\cO_{\cY} \ra \cO_Y))$ and let us take 
$\lam,\lam' \in (p^{-r_0},1]$ 
with $\lam'>\lam$. Then, for $\cE \in 
I_{\lam'\text{-}\conv}((X \allowbreak / \allowbreak \cY)^{\log})$, 
the restriction map 
$$ Rf_{X/\cY,\lam'\text{-}\conv *}\cE \lra 
Rf_{X/\cY,\lamm\conv *}\res_{\lam',\lam}(\cE) $$ 
is a quasi-isomorphism. 
\end{prop} 

\begin{pf} 
We can reduce to the case that there exists 
an exact closed immersion $(X,M_{X}) \allowbreak \hra (\cP,M_{\cP})$ over 
$(\cY,M_{\cY})$ such that $(\cP,M_{\cP})$ is formally log smooth 
over $(\cY,M_{\cY})$ and that we have the canonical isomorphism 
$(\cP,M_{\cP}) \times_{(\cY,M_{\cY})} (Y,M_{Y}) = (X,M_{X})$. 
Then we have the isomorphism 
$]X[^{\log}_{\cP,\lam} \os{\cong}{\lra} ]X[^{\log}_{\cP,\lam'} = \cP_K$
and the proposition follows from this as Proposition \ref{22}. 
\end{pf} 

Now let us assume that we are given a diagram \eqref{diag-s2} (with $r_0$ as 
above). Then we define the rational number $r(X,Y,\cY)$ by 
\begin{equation}\label{r}
r(X,Y,\cY) := ((\max(26,3a(X)-1)d(p-1)+1)+r_0^{-1})^{-1}, 
\end{equation}
where $d = \max_{x \in X}(\rk \omega^1_{X/Y})$. Then, as a consequence 
of Propositions \ref{21}, \ref{22} and \ref{23}, we have the following: 

\begin{cor}\label{perf2}
Assume given the diagram \eqref{diag-s2}, put 
$(Y_1, M_{Y_1}) := (\cY,M_{\cY}) \otimes_{\Z_p} \Z/p\Z$ and 
assume that there exists a proper 
log smooth integral morphism $(X_1,M_{X_1}) \allowbreak 
\lra (Y_1, M_{Y_1})$ 
satisfying $(X_1,M_{X_1}) \times_{(Y_1,M_{Y_1})} (Y,M_Y) = (X,M_X)$. Then, 
for $\lam \in \allowbreak (p^{-r(X,Y,\cY)}, \allowbreak 1]$ 
and a locally free isocrystal $\cE$ on 
$(X/\cY)^{\log}_{\lamm\conv}$, 
$Rf_{X/\cY,\conv *}\cE$ is a perfect complex 
of $\Q \otimes_{\Z} \cO_{\cY}$-modules. 
\end{cor}

\begin{pf} 
Denote the exact closed immersion $(X,M_X) \hra (X_1,M_{X_1})$ by $i$. 
Then we have $i_*\cE \in I_{\lam'\text{-}\conv}((X_1/\cY)^{\log})$, where 
$\lam' := p^{-1/(\log_p(\lam^{-1})+r_0^{-1})} >
p^{-1/(\max(26,3a(X)-1)d(p-1)+1)}$, and by Propositions \ref{21}, 
\ref{22} and \ref{23}, we have the isomorphism 
$$ 
Rf_{X_1/\cY, \lam'\text{-}\conv *}i_*\cE
= 
Rf_{X/\cY, \lam'\text{-}\conv *}i^*i_*\cE 
= 
Rf_{X/\cY,\lam'\text{-}\conv *}\res_{\lam,\lam'}\cE
= 
Rf_{X/\cY,\lamm\conv *}\cE. 
$$ 
So it suffices to see that $Rf_{X_1/\cY, \lam'\text{-}\conv *}i_*\cE$ is a 
perfect complex of $\Q \otimes_{\Z} \cO_{\cY}$-modules. 
To see this, let us note the inequality $l(X_1,\cY) \leq a(X)$: Indeed, 
if $X = \bigcup_{j=1}^{a(X)} X^{(j)}$ is a Zariski covering of $X$ by 
affine open subschemes and if we denote the open subscheme of $X_1$ 
which is homeomorphic to $X^{(j)}$ by $X_1^{(j)}$, then 
$X_1=\bigcup_{j=1}^{a(X)} X_1^{(j)}$ is a Zariski covering 
of $X_1$ by affine open subschemes. So we have 
$a(X) \geq a(X_1) \geq l(X_1,\cY)$. Then the perfectness of 
$Rf_{X_1/\cY, \lam'\text{-}\conv *}i_*\cE$ follows from Corollary 
\ref{perf1}. 
\end{pf} 

\begin{cor}\label{bc1}
Assume we are given a diagram 
\begin{equation}
\begin{CD}
(X',M_{X'}) @>>> (Y',M_{Y'}) @>>> (\cY',M_{\cY'}) \\
@VVV @VVV @V{\varphi}VV \\
(X,M_X) @>>> (Y,M_Y) @>>> (\cY,M_{\cY}), 
\end{CD}
\end{equation}
where 
the horizontal lines are as in the diagram \eqref{diag-s2} and 
the left square is Cartesian. Let 
$(Y_1,M_{Y_1}) := (\cY,M_{\cY}) \otimes_{\Z_p} \Z/p\Z$ 
and assume that there exists a proper 
log smooth integral morphism $(X_1,M_{X_1}) \lra (Y_1, M_{Y_1})$ 
saisfying $(X_1,M_{X_1}) \times_{(Y_1,M_{Y_1})} (Y,M_Y) = (X,M_X)$. 
Then, 
for $\lam \in (p^{-r(X,Y,\cY)},1]$ and 
a locally free isocrystal $\cE$ on $I_{\lamm\conv}((X/\cY)^{\log})$, 
we have the canonical quasi-isomorphism 
$$ L\varphi^*Rf_{X/\cY,\lamm\conv *}\cE \os{\sim}{\lra} 
Rf_{X'/\cY',\lamm\conv *} \varphi^*\cE. $$
\end{cor}

\begin{pf} 
Denote the base change of 
$(X_1,M_{X_1}) \lra (Y_1, M_{Y_1}) \hra (\cY,M_{\cY})$ by 
$(\cY',M_{\cY'}) \allowbreak \lra (\cY,M_{\cY})$ by 
$(X'_1,M_{X'_1}) \lra (Y'_1, M_{Y'_1}) \hra (\cY',M_{\cY'})$. Then 
we have $l(X'_1,\cY') \leq l(X_1,\cY) \leq a(X)$. Noting this, 
we can reduce the claim to Corollary \ref{bccrys} in the same way as 
the proof of Corollary \ref{perf2}. 
\end{pf} 

Now let us recall the notion of `having log smooth parameter', which 
is defined in \cite[3.4]{shiho3}: 

\begin{defn}\label{parameter}
We say that a proper log smooth integral 
morphism $f: (X,M_X) \lra (Y,M_Y)$ of fine log $B$-schemes has 
log smooth parameter in strong sense $($over $(B,M_B))$, 
if there exists a diagram of 
fine log formal $B$-schemes
\begin{equation}\label{star}
\begin{CD}
(X,M_X) @<<< (X',M_{X'}) @>>> (X_0,M_{X_0}) \\
@VfVV @V{f'}VV @V{f_0}VV \\
(Y,M_Y) @<g<< (Y',M_{Y'}) @>{g'}>> (Y_0,M_{Y_0}), 
\end{CD}
\end{equation}
where two squares are Cartesian, $g$ is strict etale and surjective, 
$f_0$ is proper log smooth integral and $(Y_0,M_{Y_0})$ is log smooth 
over $(B,M_B)$. We say a proper log smooth integral 
morphism $f: (X,M_X) \lra (Y,M_Y)$ of fine log $B$-schemes has 
log smooth parameter if there exists a decomposition 
$X = \coprod_i X_i$ into open and closed subschemes such that 
the composite 
$$ (X_i,M_X|_{X_i}) \hra (X,M_X) \os{f}{\lra} (Y,M_Y) $$ 
has log smooth parameter in strong sense for each $i$. 
\end{defn}

Let us assume given the diagram \eqref{diag-s2} such that $f$ has log 
smooth parameter in strong sense 
and take a diagram \eqref{star}. Let $(\cY',M_{\cY'})$ 
be the fine log formal $\cB$-scheme which is strict formally etale over 
$(\cY,M_{\cY})$ satisfying $(\cY',M_{\cY'}) \times_{(\cY,M_{\cY})} (Y,M_Y) 
= (Y',M_{Y'})$ and put $(Y'_1,M_{Y'_1}) := (\cY',M_{\cY'}) 
\otimes_{\Z_p} \F_p$. Then, by \cite[3.7]{shiho3}, 
we may assume the following condition $(\star)$ by replacing $Y'$ by 
its etale covering if necessary: \\
\quad \\
$(\star)$: \,\,\, The morphism $Y' \lra Y$ is affine and 
there exists a proper log smooth integral morphism 
$(X'_1,M_{X'_1}) \lra (Y'_1,M_{Y'_1})$ satisfying 
$(X'_1,M_{X'_1}) \times_{(Y'_1,M_{Y'_1})} (Y,M_Y) = (X,M_X)$. \\
\quad \\
Now let us prove one of the main results in this section: 

\begin{thm}\label{perfmain}
Assume we are given the diagram \eqref{diag-s2} such that 
$f$ is a proper log smooth integral morphism having log smooth parameter. 
Then, for $\lam \in (p^{-r(X,Y,\cY)},1]$ and 
a locally free isocrystal $\cE$ on 
$(X/\cY)^{\log}_{\lamm\conv}$, 
$Rf_{X/\cY,\lamm\conv *}\cE$ is a perfect complex 
of $\Q \otimes_{\Z} \cO_{\cY}$-modules. 
\end{thm}

\begin{pf} 
The proof is similar to that in \cite[3.6]{shiho3}. 
We may assume that $f$ has log smooth parameter in strong sense. 
Take a diagram \eqref{star} satisfying 
the condition $(\star)$ above. Let 
$\epsilon: 
(\cY^{(\b)},M_{\cY^{(\b)}}) \lra (\cY,M_{\cY})$ be the strict 
formally etale \v{C}ech hypercovering 
associated to the morphism $(\cY',M_{\cY'}) \lra (\cY,M_{\cY})$. 
We denote 
the base change of the diagram \eqref{diag-s2} by 
$(\cY^{(\b)},M_{\cY^{(\b)}}) \lra (\cY,M_{\cY})$ by 
$$ 
(X^{(\b)},M_{X^{(\b)}}) \lra 
(Y^{(\b)},M_{Y^{(\b)}}) \hra (\cY^{(\b)},M_{\cY^{(\b)}}). $$
For $n \in \N$, let $\pi: (\cY^{(n)},M_{\cY^{(n)}}) \lra 
(\cY^{(0)},M_{\cY^{(0)}}) = (\cY',M_{\cY'})$ be one of the projections and 
denote by 
$(X_1^{(n)},M_{X_1^{(n)}}) \lra (Y_1^{(n)},M_{Y_1^{(n)}}) 
\hra (\cY^{(n)},M_{\cY^{(n)}})$ 
the base change of the diagram 
$(X'_1,M_{X'_1}) \lra (Y'_1,M_{Y'_1}) \hra (\cY,M_{\cY})$ by $\pi$. 
Then we have 
$(X_1^{(n)},M_{X_1^{(n)}}) \times_{(Y_1^{(n)},M_{Y_1^{(n)}})}
(Y^{(n)},M_{Y^{(n)}}) = (X^{(n)},M_{X^{(n)}})$ for each $n$. 
(Note that $(X_1^{(n)},M_{X_1^{(n)}})$, $(Y_1^{(n)},M_{Y_1^{(n)}})$ 
does not form a simplicial log scheme.) 
Since $Y^{(n)} \lra Y$ is affine, 
we have $a(X^{(n)}) \leq a(X)$ and so we have 
$r(X^{(n)},Y^{(n)},\cY^{(n)}) \geq r(X,Y,\cY)$. So, 
if we denote the restriction of $\cE$ to 
 $I_{\lamm\conv}((X^{(\b)}/\cY^{(\b)})^{\log})$ by $\cE^{(\b)}$, 
$Rf_{X^{(n)}/\cY^{(n)},\lamm\conv *}\cE^{(n)}$ is perfect for each $n$, 
by Corollary \ref{perf2}. \par 
Next we prove the quasi-isomorphism 
\begin{equation}\label{aux1}
Rf_{X/\cY,\lamm\conv *}\cE = 
R\epsilon_*Rf_{X^{(\b)}/\cY^{(\b)}, \lamm\conv *}\cE^{(\b)}.
\end{equation}
To see this, we can assume that $(X,M_X)$ admits an exact closed immersion 
$(X,M_X) \allowbreak \hra (\cP,M_{\cP})$ to a fine log formal $\cB$-scheme 
$(\cP,M_{\cP})$ 
which is formally log smooth over $(\cY,M_{\cY})$ and satisfies 
$(\cP,M_{\cP}) \times_{(\cY,M_{\cY})} (Y,M_Y) 
= (X,M_X)$. In this case, we have the strict formally 
etale \v{C}ech hypercovering 
$\delta: (\cP^{(\b)},M_{\cP^{(\b)}}) \lra (\cP,M_{\cP})$ 
satisfying $(\cP^{(\b)},M_{\cP^{(\b)}}) \times_{(\cY,M_{\cY})} (Y,M_Y) 
= (X^{(\b)},M_{X^{(\b)}})$ since we have $\cP_{\et} \simeq X_{\et}$. 
If we denote the morphism 
$\cP \lra \cY$ by $h$, we have 
$$ Rf_{X/\cY,\lamm\conv *}\cE = Rh_* \sp_* \DR(\cP_K/\cY_K, \cE), $$
$$ R\epsilon_*Rf_{X^{(\b)}/\cY^{(\b)}, \lamm\conv *}\cE^{(\b)}
=  Rh_* R\delta_* \sp^{(\b)}_* \DR(\cP^{(\b)}_K/\cY^{(\b)}_K, 
\cE^{(\b)}). $$ 
(Note that we have $]X[^{\log}_{\cP,\lam}=\cP_K, 
]X^{(\b)}[^{\log}_{\cP^{(\b)},\lam}=\cP^{(\b)}_K$ since we have 
$\lam > p^{-r_0}$.) Using this, we can prove the 
quasi-isomorphism \eqref{aux1} in the same way as \cite[3.6]{shiho3}. \par 
By \eqref{aux1}, we have the spectral sequence 
$$ 
E_2^{s,t} = R^s\epsilon_*R^tf_{X^{(\b)}/\cY^{(\b)},\lamm\conv *}\cE^{(\b)} 
\,\Longrightarrow\, R^{s+t}f_{X/\cY,\lamm\conv *}\cE. $$
Note that 
$(R^tf_{X^{(\b)}/\cY^{(\b)},\lamm\conv *}\cE^{(\b)})_n = 
R^tf_{X^{(n)}/\cY^{(n)},\lamm\conv *}\cE^{(n)}$ is known to be isocoherent and 
that they are compatible with respect to $n$ by Corollary \ref{bc1}. 
(Since we have $r(X^{(n)},Y^{(n)},\cY^{(n)}) \geq r(X,Y,\cY)$, we can 
use Corollary \ref{bc1}.) 
So there exists (by etale descent of isocoherent sheaves) 
an isocoherent sheaf $\cF^t$ on $\cY$ such that 
$\cF^t \otimes_{\cO_{\cY}} \cO_{\cY^{(\b)}} = 
R^tf_{X^{(\b)}/\cY^{(\b)},\lamm\conv *}\cE^{(\b)}$ holds. 
Then we have 
$$ 
R^s\epsilon_*R^tf_{X^{(\b)}/\cY^{(\b)},\lamm\conv *}\cE^{(\b)} = 
\cF^t \,\, (s=0), \,\,\, 0 \,\, (s>0). $$
So we have $R^qf_{X/\cY,\lamm\conv *}\cE = \cF^q$ and it is isocoherent. 
Moreover, $Rf_{X/\cY,\lamm\conv *}\cE$ is bounded and it 
has finite tor-dimension Zariski locally 
because so does $Rf_{X^{(0)}/\cY^{(0)},}\allowbreak{}_{\lamm\conv *}\cE^{(0)}
= Rf_{X^{(0)}/\cY^{(0)},}\allowbreak{}_{\crys *} \allowbreak \Phi(\cE^{(0)})$. 
Therefore $Rf_{X/\cY,\lamm\conv *}\cE$ is a perfect complex
of $\Q \otimes_{\Z} \cO_{\cY}$-modules. 
\end{pf} 

\begin{rem}\label{remrem}
As in \cite[3.8]{shiho3}, we have the following: 
With the notation in the proof of Theorem \ref{perfmain}, 
we have the isomorphism 
$$ R^qf_{X/\cY,\lamm\conv *}\cE \otimes_{\cO_{\cY}} \cO_{\cY^{(n)}} = 
R^qf_{X^{(n)}/\cY^{(n)},\lamm\conv *}\cE^{(n)}. $$
\end{rem} 

Next we prove the base change property, which is the second main 
result in this section. 

\begin{thm}\label{bc2}
Assume we are given a diagram 
\begin{equation}\label{eqeqeqeq1}
\begin{CD}
(X'',M_{X''}) @>>> (Y'',M_{Y''}) @>>> (\cY'',M_{\cY''}) \\
@VVV @VVV @V{\varphi}VV \\
(X,M_X) @>f>> (Y,M_Y) @>>> (\cY,M_{\cY}), 
\end{CD}
\end{equation}
where 
the horizontal lines are as in Theorem \ref{perfmain} and the 
left square is Cartesian. 
Then, for $\lam \in (p^{-r(X,Y,\cY)},1]$ and 
a locally free isocrystal $\cE$ on $I_{\lamm\conv}((X/\cY)^{\log})$, 
we have the quasi-isomorphism 
$$ 
L\varphi^*
Rf_{X/\cY,\conv *}\cE \os{\sim}{\lra} 
Rf_{X''/\cY'',\conv *} \varphi^*\cE. $$
\end{thm}

\begin{pf}
The proof is the same as \cite[3.9]{shiho3}. So we only give a sketch. 
We may assume that $f$ has log smooth parameter in strong sense. 
If we take the diagram \eqref{star} satisfying the condition $(\star)$, 
we see that it suffices to prove the theorem after pulling-back the 
diagram \eqref{eqeqeqeq1} by $(\cY',M_{\cY'}) \lra (\cY,M_{\cY})$. 
(We use Remark \ref{remrem} here.) Then the theorem is reduced to 
Corollary \ref{bc1}. 
\end{pf} 

Recall that a morphism $(\cY',M_{\cY'}) \lra (\cY,M_{\cY})$ 
in $\pLF$ is said to be analytically flat if the induced functor 
$\Coh(\Q \otimes \cO_{\cY}) \lra \Coh(\Q \otimes \cO_{\cY'})$ is 
exact. Then, as in \cite[3.10]{shiho3}, we have the following result, 
which is a corollary of Corollary \ref{bc1} and Theorem 
\ref{bc2}: 

\begin{cor}\label{bc3}
Assume we are given a diagram in $\pLF$ 
\begin{equation*}
\begin{CD}
(X',M_{X'}) @>>> (Y',M_{Y'}) @>>> (\cY',M_{\cY'}) \\
@VVV @VVV @V{\varphi}VV \\
(X,M_X) @>f>> (Y,M_Y) @>{\iota}>> (\cY,M_{\cY}), 
\end{CD}
\end{equation*}
where $f$ is a proper log smooth integral morphism, $\iota$ is a 
homeomorphic 
exact closed immersion and the left square is Cartesian. 
Let $\lam \in (p^{-r(X,Y,\cY)},1]$ and 
let $\cE$ be a locally free isocrystal on $(X/\cY)^{\log}_{\lamm\conv}$. 
Assume one of the conditions $(1)$, $(1)'$ below and 
one of the conditions $(2)$, $(2)'$ below are true$:$ \\
$(1)$ \,\, $\varphi$ is analytically flat. \\
$(1)'$ \,\, $R^qf_{X/\cY, \lamm\conv *}\cE$ is a locally free $\Q \otimes_{\Z}
\cO_{\cY}$-module $($in the sense of \cite[1.9]{shiho3}$)$ for any $q$. \\
$(2)$ \,\, $f$ has log smooth parameter. \\
$(2)'$ \,\, If we put $(Y_1,M_{Y_1}) := (\cY,M_{\cY}) \otimes_{\Z_p} \Z/p\Z$, 
there exists a proper log smooth integral morphism 
$(X_1,M_{X_1}) \lra (Y_1,M_{Y_1})$ satisfying 
$(X_1,M_{X_1}) \times_{(Y_1,M_{Y_1})} (Y,M_Y) = (X,M_X)$. \\ 
Then we have the canonical isomorphism 
$$ \varphi^*R^qf_{X/\cY,\lamm\conv *}\cE \os{\cong}{\lra} 
R^qf_{X'/\cY',\lamm\conv *}\varphi^*\cE \quad (q \in \N). $$
\end{cor} 

In the case where the condition (1) is satisfied, we call the result 
of Corollary \ref{bc3} `analytically flat base change theorem'. 


\section{Relative log analytic cohomology of radius $\lam$} 

In this section, first we introduce `radius $\lam$ version' of 
the relative log analytic cohomology introduced in \cite[\S 4]{shiho3}. 
Then we prove 
the relation between relative log convergent cohomology of radius 
$\lam$ and the 
relative log analytic cohomology of radius $\lam$ 
for proper log smooth integral morphisms 
having log smooth parameter when $\lam$ is sufficiently close to $1$. 
This implies the coherence of 
relative log analytic cohomology under the same assumption on 
$\lam$. After that, we 
prove the existence of a canonical structure of an isocrystal on 
relative log analytic cohomology of radius $\lam$. 
The proofs are again similar to \cite[\S 4]{shiho3}, but we need 
more subtle argument. 
\par 
First we give a definition of relative log analytic 
cohomology of radius $\lam$. 

\begin{defn}\label{defrellogancoh}
Assume we are given a diagram 
\begin{equation}\label{diag3-1}
(X,M_X) \os{f}{\lra} (Y,M_Y) \os{\iota}{\hra} (\cY,M_{\cY}), 
\end{equation}
where $f$ is a morphism in $\LB$ and $\iota$ is a closed immersion 
in $\pLF$. Let $\lam \in (0,1]$ 
and let $\cE$ be an isocrystal on $(X/\cY)^{\log}_{\lamm\conv}$. 
Take an embedding system 
\begin{equation}\label{diag3-2}
 (X,M_{X}) \overset{g}{\lla} (X^{(\bullet)}, M_{X^{(\bullet)}}) 
\overset{i}{\hra} (\cP^{(\bullet)},M_{\cP^{(\bullet)}}), 
\end{equation}
let $\cE^{(\b)}$ be the restriction of $\cE$ to 
$(X^{(\b)},M_{X^{(\b)}})$, 
let $\DR(]X^{(\b)}[^{\log}_{\cP^{(\b)},\lam}/\cY_K,\cE^{(\b)})$ be the 
log de Rham complex associated to $\cE^{(\b)}$ and let $h$ be the morphism 
$]X^{(\b)}[^{\log}_{\cP^{(\b)},\lam} \lra ]Y[^{\log}_{\cY,\lam}$. 
Then we define 
$Rf_{X/\cY, \lamm\an *}\cE, R^qf_{X/\cY,\lamm\an *}\cE$ by 
$$ 
Rf_{X/\cY,\lamm\an *}\cE := 
Rh_*\DR(]X^{(\b)}[^{\log}_{\cP^{(\b)},\lam}/\cY_K,\cE^{(\b)}), $$
$$ 
R^qf_{X/\cY,\lamm\an *}\cE := 
R^qh_*\DR(]X^{(\b)}[^{\log}_{\cP^{(\b)},\lam}/\cY_K,\cE^{(\b)}) $$
and we call $R^qf_{X/\cY,\lamm\an *}\cE$ the $q$-th relative log analytic 
cohomology of $(X,M_X)/(\cY,\allowbreak 
M_{\cY})$ of radius $\lam$ with coefficient $\cE$. 
It is a sheaf of $\cO_{]Y[^{\log}_{\cY,\lam}}$-modules. 
\end{defn}

\begin{rem}\label{rem3-1}
We have the remark similar to \cite[4.2]{shiho3}: 
With the above notation, put $\lam=p^{-r}$. Let 
$\{(\cY_{a,b},M_{\cY_{a,b}})\}_{(a,b) \in \A}$ 
be the system of universal enlargements of 
$((\cY,M_{\cY}), (Y,M_Y), \iota,\id)$ 
and denote the base change of 
the diagrams \eqref{diag3-1} by 
$(\cY_{a,b},M_{\cY_{a,b}}) \lra (\cY,M_{\cY})$ by 
\begin{equation*}
(X_{a,b},M_{X_{a,b}}) {\lra} 
(Y_{a,b},M_{Y_{a,b}}) \os{\iota}{\hra} (\cY_{a,b},M_{\cY_{a,b}}). 
\end{equation*}
Then we have $]Y[^{\log}_{\cY,\lam}= \bigcup_{(a,b)\in\A_r}  \cY_{a,b,K}$ 
and if we denote the restriction of $\cE$ to 
$I_{\lamm\conv}((X_{a,b}/\cY_{a,b})^{\log})$ by $\cE_{a,b}$,  
$R^qf_{X_{a,b}/\cY_{a,b},\lamm\an *}\cE_{a,b}$ 
is nothing but the restriction of $R^qf_{X/\cY,\lamm\an *}\cE$ to 
$\cY_{a,b,K}$. 
\end{rem}

As in \cite[4.3]{shiho3}, we should prove the following proposition: 

\begin{prop}\label{wdprop}
Let the notations be as in Definition \ref{defrellogancoh}. Then 
the definition of the relative log analytic cohomology 
$R^qf_{X/\cY,\lamm\an *}\cE$ of radius $\lam$ 
is independent of the choice of the embedding system. 
\end{prop}

The proof of Proposition \ref{wdprop} is the same as that of 
\cite[4.3]{shiho3}, once we establish Lemmas \ref{wdlem1}, \ref{logsminv} 
below (which correspond to \cite[4.4, 4.5]{shiho3}). So we only give 
a proof of the following two key lemmas: 

\begin{lem}\label{wdlem1}
Let $\lam \in (0,1]$ and 
assume we are given the Cartesian diagram 
\begin{equation}\label{wd1}
\begin{CD}
(X^{(\b)},M_{X^{(\b)}}) @>>> (\cP^{(\b)},M_{\cP^{(\b)}}) \\ 
@VVV @VgVV \\
(X,M_X) @>{\iota}>> (\cP,M_{\cP}), 
\end{CD}
\end{equation}
where $(X,M_X)$ is an object in $\LB$, $(\cP,M_{\cP})$ is an object in 
$\LF$, $\iota$ is a closed immersion 
and $g$ is a strict formally etale hypercovering. Let 
$g_K:]X^{(\b)}[^{\log}_{\cP^{(\b)},\lam} \lra ]X[^{\log}_{\cP,\lam}$ be the 
morphism induced by $g$ and let $\cE$ be a coherent 
$\cO_{]X[^{\log}_{\cP,\lam}}$-module. Then we have the isomorphism 
$$ \cE \overset{=}{\lra} Rg_{K,*}g_K^*\cE. $$
\end{lem} 

\begin{pf}
Let $\{(\cP_{a,b},M_{\cP_{a,b}})\}_{(a,b) \in \A}$ 
be the system of universal enlargements 
of the pre-widening $((\cP,M_{\cP}),(X,M_X),\iota,i)$ and let 
\begin{equation*}
\begin{CD}
(X_{a,b}^{(\b)},M_{X_{a,b}^{(\b)}}) @>>> (\cP_{a,b}^{(\b)},
M_{\cP_{a,b}^{(\b)}}) \\ 
@VVV @V{g_{a,b}}VV \\
(X_{a,b},M_{X_{a,b}}) @>{\iota}>> (\cP_{a,b},M_{\cP_{a,b}}) 
\end{CD}
\end{equation*}
be the diagram obtained by applying $\times_{(\cP,M_{\cP})} (\cP_{a,b},
M_{\cP_{a,b}})$ 
to the diagram \eqref{wd1}. Let $g_{a,b,K}:P^{(\b)}_{a,b,K} \lra P_{a,b,K}$ be the 
map induced by $g_{a,b}$. Then, if we put $r:=-\log_p \lam$, 
 we have 
{\small{
$$ ]X[^{\log}_{\cP,\lam} = \bigcup_{(a,b) \in \A_r} \cP_{a,b,K}, \quad 
]X^{(\b)}[^{\log}_{\cP^{(\b)},\lam} = \bigcup_{(a,b) \in \A_r} 
\cP^{(\b)}_{a,b,K}, $$ 
$$ g_K^{-1}(\cP_{a,b,K}) = \cP^{(\b)}_{a,b,K}, \quad 
g_{a,b,K} = g_K \vert_{\cP^{(\b)}_{a,b,K}}.$$ }}
So we have 
$(Rg_{K,*}g_K^*\cE)\vert_{\cP_{a,b,K}} = 
Rg_{a,b,K,*}g^*_{a,b,K}(\cE \vert_{\cP_{a,b,K}})$ for 
$(a,b) \in \A_r$ 
and so it suffices to prove the 
isomorphism $\cE \vert_{\cP_{a,b,K}} \os{=}{\lra} 
Rg_{a,b,K,*}g^*_{a,b,K}(\cE \vert_{\cP_{a,b,K}})$. 
This is proven in \cite[7.3.3]{chts}. 
\end{pf} 

\begin{lem}\label{logsminv}
Let $\lam \in (0,1]$, 
let $(X,M_X) \lra (Y,M_Y) \hra (\cY,M_{\cY})$ be as in Definition 
\ref{defrellogancoh} and assume we are given a commutative diagram 
over $(\cY,M_{\cY})$ 
\begin{equation}\label{wd5}
\begin{CD}
(X,M_X) @>{\iota_1}>> (\cP_1,M_{\cP_1}) \\ 
@\vert @V{\varphi}VV \\ 
(X,M_X) @>{\iota_2}>> (\cP_2,M_{\cP_2}), 
\end{CD}
\end{equation}
where $(\cP_j,M_{\cP_j}) \, (j=1,2)$ 
is an object in $\pLF$ which is formally log 
smooth over $(\cY,\allowbreak M_{\cY})$, 
$\iota_j \, (j=1,2)$ is a closed immersion and $\varphi$ is a fornally 
log smooth morphism. Let 
$\varphi_K: ]X[^{\log}_{\cP_1,\lam} \lra ]X[^{\log}_{\cP_2,\lam}$ 
be a morphism 
of log tubular neighborhoods induced by $\varphi$ and let $\cE$ be 
an isocrystal on $(X/\cY)^{\log}_{\lamm\conv}$. 
Then we have a quasi-isomorphism 
$$ \DR(]X[^{\log}_{\cP_2,\lam}/\cY,\cE) \os{=}{\lra} 
R\varphi_{K,*}\DR(]X[^{\log}_{\cP_1,\lam}/\cY,\cE). $$
\end{lem}

\begin{pf} 
The proof of similar to that of \cite[4.5]{shiho3}. 
For $j=1,2$, let $\cP^{\ex}_j$ be the exactification of the closed immersion 
$\iota_j$. Then it suffices to prove the lemma Zariski locally on 
$\cP^{\ex}_2$. So we may assume that $\cP^{\ex}_2 = \Spf A$ is affine, 
$\cP^{\ex}_1 = \Spf A[[t_1, \cdots, t_r]]$ for some $r$ and the morphism 
$\varphi^{\ex}:\cP^{\ex}_1 \lra \cP^{\ex}_2$ is induced by the canonical 
inclusion $A \hra A[[t_1, \cdots, t_r]]$. Then we have the isomorphism 
$]X[^{\log}_{\cP_1,\lam} \cong ]X[^{\log}_{\cP_2,\lam} \times D^r_{\lam}$ 
(where $D^r_{\lam}$ is the $r$-dimensional open polydisc of radius $\lam$) and 
the morphism $\varphi_K$ is equal to the projection 
$]X[^{\log}_{\cP_2,\lam} \times D^r_{\lam} \lra 
]X[^{\log}_{\cP_2,\lam}$. 
Let 
$\Omega^{\b}$ be the complex 
{\tiny{ $$ 
[\cO_{]X[^{\log}_{\cP_1,\lam}} \lra 
\bigoplus_{i=1}^r\cO_{]X[^{\log}_{\cP_1,\lam}} dt_i \lra 
\bigoplus_{1 \leq i_1 < i_2 \leq r} \cO_{]X[^{\log}_{\cP_1,\lam}} 
dt_{i_1} \wedge dt_{i_2} \lra \cdots \lra 
\cO_{]X[^{\log}_{\cP_1,\lam}} dt_1 \wedge \cdots \wedge dt_r]. $$ }}
Then $\DR(]X[^{\log}_{\cP_1,\lam}/\cY_K,\cE)$ is equal to the total complex 
associated to the double complex 
$\varphi_K^*\DR(]X[^{\log}_{\cP_2,\lam}/\cY_K,\cE) 
\otimes_{\cO_{]X[^{\log}_{\cP_1,\lam}}} \Omega^{\b}$. 
So, to prove the lemma, it suffices to prove the quasi-isomorphism 
$$ E \os{\simeq}{\lra} R\varphi_{K,*}(\varphi_K^*E \otimes \Omega^{\b})
= \varphi_{K,*}(\varphi_K^*E \otimes \Omega^{\b}) = 
E \otimes_{\cO_{]X[^{\log}_{\cP_2,\lam}}} \varphi_{K,*}\Omega^{\b} $$ 
for a coherent $\cO_{]X[^{\log}_{\cP_2,\lam}}$-module $E$. This is reduced to 
showing that the complex 
$\cO_{]X[^{\log}_{\cP_2,\lam}} \lra \varphi_{K,*}\Omega^{\b}$
 is homotopic to zero. We can prove it in the same way as the proof of 
 Lemma \ref{add2} (see also \cite[2.32]{shiho3}). So we are done. 
\end{pf}

We have the following independence result with respect to radius: 

\begin{prop}\label{33}
Assume we are given a diagram 
\begin{equation}\label{diag-33}
(X,M_X) \os{f}{\lra} (Y,M_Y) \os{\iota}{\hra} (\cY,M_{\cY}), 
\end{equation}
where $f$ is a log smooth morphism in $\LB$ and $\iota$ is a closed immersion 
in $\pLF$. Let $\lam,\lam' \in (0,1]$ with $\lam'>\lam$. Then, for 
an isocrystal $\cE$ on $(X/\cY)^{\log}_{\lam'\text{-}\conv}$, 
we have the quasi-isomorphism 
$$ (Rf_{X/\cY,\lam'\text{-}\an *}\cE)|_{]Y[^{\log}_{\cY,\lam}} \os{=}{\lra} 
Rf_{X/\cY,\lamm\an *}\res_{\lam',\lam}\cE. $$
\end{prop}

\begin{pf} 
Put $\lam = p^{-r}$. 
For $(a,b) \in \A_r$, let 
\begin{equation}\label{diag-33-2}
(X_{a,b},M_{X_{a,b}}) \lra (Y_{a,b},M_{Y_{a,b}}) 
\hra (\cY_{a,b},M_{\cY_{a,b}})
\end{equation}
be as in Remark \ref{rem3-1}. 
To prove the quasi-isomorphism, we may check it on $\cY_{a,b,K}$ for 
$(a,b) \in \A_r$. So, by Remark \ref{rem3-1}, we may replace 
the diagram \eqref{diag-33} by the diagram \eqref{diag-33-2} to prove 
the proposition: So, if we put 
$r_0:=r(\cO_{\cY},\Ker(\iota^*:\cO_{\cY} \ra \cO_Y))$, we may 
assume $r_0 < r$. Then we can prove the proposition in the same way as 
Proposition \ref{23}. 
\end{pf} 

The next theorem establishes the relation of relative log convergent 
cohomology of radius $\lam$ 
and relative log analytic cohomology of radius $\lam$: 

\begin{thm}\label{coherence0}
Assume we are given a diagram 
\begin{equation}
(X,M_X) \os{f}{\lra} (Y,M_Y) \os{\iota}{\hra} (\cY,M_{\cY}), 
\end{equation}
where $f$ is a proper log smooth integral morphism having log smooth parameter 
in $\LB$ and $\iota$ is a homeomorphic exact closed immersion 
in $\pLF$. Then, for $\lam \in (p^{-r(X,Y,\cY)},1]$, 
a locally free isocrystal $\cE$ on 
$(X/\cY)^{\log}_{\lamm\conv}$ and $q \geq 0$, 
$R^qf_{X/\cY,\lamm\an *}\cE$ is a coherent 
sheaf on $]Y[^{\log}_{\cY,\lam} = \cY_K$ and we have the isomorphism 
$\sp_*R^qf_{X/\cY,\lamm\an *}\cE = R^qf_{X/\cY,\lamm\conv *}\cE$. 
\end{thm}

\begin{pf}
We can prove the theorem exactly in the same way as 
\cite[4.6]{shiho3}, using Theorem \ref{perfmain} and Corollary \ref{bc3}. 
\end{pf} 

Next we prove the existence of a structure of an isocrystal on 
relative log convergent cohomology with radius. 
For a proper log smooth integral morphism 
$(X,M_X) \allowbreak \lra (Y,M_Y)$ in $\LB$ and $\mu \in (0,1]$, we define 
$r(X,Y;\mu) \in (0,\infty]$ by 
\begin{equation}\label{rr}
r(X,Y;\mu) := ((\max(26,3a(X)-1)d(p-1)+1)-(\log_p\mu)^{-1})^{-1}, 
\end{equation}
where $d:=\max_{x \in X}\rk(\omega^1_{X/Y})$. Then we have the 
following: 

\begin{thm}\label{iso1}
Assume we are given a diagram 
$$ (X,M_X) \os{f}{\lra} (Y,M_Y) \os{g}{\lra} (\cS,M_{\cS}), $$
where $f$ is a proper log smooth integral morphism having log smooth parameter 
in $\LB$ and $g$ is a morphism in $\pLF$. 
Let $\mu = p^{-r} \in (0,1]$ and fix $\lam \in (p^{r(X,Y;\mu)},1]$. 
Then, for a locally free isocrystal 
$\cE$ on $(X/\cS)^{\log}_{\lamm\conv}$ and a non-negative integer $q$, 
there exists a unique isocrystal $\cF$ on $(Y/\cS)^{\log}_{\mum\conv}$ 
satisfying the following condition$:$ 
For any pre-widening $\cZ := ((\cZ,M_{\cZ}),(Z,M_Z),i,z)$ such that 
$z$ is a strict morphism and that $(\cZ,M_{\cZ})$ is formally log smooth 
over $(\cS,M_{\cS})$, the restriction of $\cF$ to 
$I_{\mum\conv}((Z/\cS)^{\log}) \cong 
\Strat'_{\mu}((Z \hra \cZ/\cS)^{\log})$ is 
functorially given by 
$\{(R^qf_{X \times_Y Z_{a,b}/T_{a,b}(\cZ),\lamm\conv *}\cE, \allowbreak 
\epsilon_{a,b})\}_{(a,b) \in \A_r}$, 
where 
$$ \{T_{a,b}(\cZ)\}_{(a,b)} := 
\{((T_{a,b}(\cZ),M_{T_{a,b}(\cZ)}),(Z_{a,b},M_{Z_{a,b}}))\}_{(a,b)}$$ is 
the system of universal enlargement of $\cZ$ and $\epsilon_{a,b}$ is the 
isomorphism 
{\allowdisplaybreaks{
\begin{align*}
p^*_{2,(a,b)}Rf_{X \times_Y Z_{a,b}/T_{a,b}(\cZ), \lamm\conv *}\cE 
& \os{\simeq}{\rightarrow} 
Rf_{X \times_Y Z(1)_{a,b}/T_{a,b}(\cZ(1)), \lamm\conv *}
\cE \\ & \os{\simeq}{\leftarrow} 
p^*_{1,(a,b)}Rf_{X \times_Y Z_{a,b}/T_{a,b}(\cZ), \lamm\conv *}\cE.
\end{align*}}}
$($Here $(\cZ(1),M_{\cZ(1)}):=(\cZ,M_{\cZ}) \times_{(\cS,M_{\cS})} 
(\cZ,M_{\cZ})$, 
$$\{T_{a,b}(\cZ(1))\}_{(a,b)} := \{((T_{a,b}(\cZ(1)),M_{T_{a,b}(\cZ(1))}), 
(Z(1)_{a,b}, M_{Z(1)_{a,b}}))\}_{(a,b)}$$ 
is the system of universal enlargements of 
$((\cZ(1),M_{\cZ(1)}),(Z,M_Z))$ and $p_{i,(a,b)}$ is the morphism 
$(T_{a,b}(\cZ(1)),M_{T_{a,b}(\cZ(1))}) \lra 
(T_{a,b}(\cZ),M_{T_{a,b}(\cZ)})$ induced by the $i$-th projection.$)$ 
\end{thm}

\begin{pf} 
The proof is similar to that of \cite[4.8]{shiho3}. 
First let us define $\cF$ in the case where there exists a closed immersion 
$(Y,M_Y) \os{\iota}{\hra} (\cP,M_{\cP})$ in $\pLF$ over $(\cS,M_{\cS})$ 
such that $(\cP,M_{\cP})$ is formally log smooth over $(\cS,M_{\cS})$. 
Let $(\cP(i),M_{\cP(i)}) \,(i=0,1,2)$ be the $(i+1)$-fold fiber product 
of $(\cP,M_{\cP})$ over $(\cS,M_{\cS})$ and denote the canonical 
closed immersion $(Y,M_Y) \hra (\cP(i),M_{\cP(i)})$ induced by $\iota$ by 
$\iota(i)$. 
Denote the system of universal 
enlargements of $((\cP(i),M_{\cP(i)}),(Y,M_Y),
\iota(i),\allowbreak \id)$ by 
$\{T_{a,b}(\cP(i)) := ((T_{a,b}(\cP(i)),M_{T_{a,b}(\cP(i))}),
(Z_{a,b}(i),M_{Z_{a,b}(i)}))\}_{(a,b) \in \A}$ and 
denote the pull-back of $f$ by $(Z_{a,b}(i),M_{Z_{a,b}(i)}) \lra (Y,M_Y)$ by 
$(X_{a,b}(i),M_{X_{a,b}(i)}) \lra (Z_{a,b}(i),M_{Z_{a,b}(i)})$. 
Then, since the morphism $(Z_{a,b}(i),M_{Z_{a,b}(i)}) \lra (Y,M_Y)$ is 
affine, we have $a(X_{a,b}(i)) \leq a(X)$. So, for $(a,b) \in \A_r$, 
we have $r(X_{a,b}(i),Z_{a,b}(i),T_{a,b}(\cP(i))) \leq r(X,Y;\mu)$. So 
if we put $\cF_{a,b}(i) \allowbreak := 
R^qf_{X_{a,b}(i)/T_{a,b}(\cP(i)),\lamm\conv *}\cE$, 
it is isocoherent and it satisfies the analytically flat base change property. 
Note first that the transition morphisms 
$(T_{a,b}(\cP(0)),M_{T_{a,b}(\cP(0))}) \allowbreak 
\lra (T_{a',b'}(\cP(0)),M_{T_{a',b'}(\cP(0))}) \, ((a,b) \succ (a',b'))$ are 
analytically flat (\cite[3.11]{shiho3}). So, by 
analytically flat base change theorem, the family 
$\{\cF_{a,b}(0)\}_{(a,b) \in \A_r}$ 
defines a compatible family of isocoherent sheaves 
on $\{T_{a,b}(\cP(0))\}\allowbreak{}_{(a,b) \in \A_r}$. Let us note next that 
the projections 
$(T_{a,b}(\cP(i+1)),M_{T_{a,b}(\cP(i+1))}) \allowbreak 
\lra (T_{a,b}(\cP(i)), \allowbreak M_{T_{a,b}(\cP(i))})$ 
are analytically flat. (This can be shown exactly by the same proof as 
\cite[3.12]{shiho3}.) So, again by 
analytically flat base change theorem, 
$\{\cF_{a,b}(1)\}_{(a,b) \in \A_r}, \{\cF_{a,b}(2)\}_{(a,b) \in \A_r}$ 
induce a structure of compatible family of stratification 
on $\{\cF_{a,b}(0)\}_{(a,b)\in \A_r}$. In this way, 
$\{\cF_{a,b}(i)\}_{(a,b),i}$ induces an object $\cF$ in 
$\Strat'_{\mu}((Y \allowbreak 
\hra \cP/\cS)^{\log}) = I_{\mum\conv}((Y/\cS)^{\log}).$ \par 
Next we define the isocrystal $\cF$ in general case. 
Take an embedding system 
\begin{equation*}
(Y,M_{Y}) \overset{g}{\lla} (Y^{(\bullet)}, M_{Y^{(\bullet)}}) 
{\hra} (\cP^{(\bullet)},M_{\cP^{(\bullet)}})
\end{equation*}
such that $g$ is a strict formally etale \v{C}ech hypercovering 
with $Y^{(0)} \lra Y$ affine 
and that 
$(\cP^{(i)},M_{\cP^{(i)}})$ is the $(i+1)$-fold fiber product 
of $(\cP^{(0)},M_{\cP^{(0)}})$ over $(\cS,M_{\cS})$. 
Let us denote the pull-back of $f$ by $(Y^{(i)},M_{Y^{(i)}}) \lra (Y,M_Y)$ 
by $(X^{(i)},M_{X^{(i)}}) \lra (Y^{(i)},M_{Y^{(i)}})$. 
Then, since $Y^{(i)} \lra Y$ is affine, 
we have $r(X^{(i)},Y^{(i)};\mu) \leq r(X,Y,\mu)$. 
Hence, 
by the construction in the previous paragraph for 
$(Y^{(i)},M_{Y^{(i)}}) \hra (\cP^{(i)},M_{\cP^{(i)}})$, we have 
$\cF^{(i)} \in I_{\conv}((Y^{(i)}/\cS)^{\log}) \, (i=0,1,2)$ and they are 
compatible thanks to 
analytically flat base change theorem. 
So, by etale descent of isocrystals on relative log convergent site of 
radius $\mu$, 
$\{\cF^{(i)}\}_{i=0,1,2}$ descents to an isocrystal $\cF$ on 
$(Y/\cS)^{\log}_{\mum\conv}$. \par 
Note that the isocrystal $\cF$ constructed in the previous paragraph 
satisfies the required condition for 
$\cZ = ((\cP^{(i)},M_{\cP^{(i)}}), (Y^{(i)},M_{Y^{(i)}})) \, (i=0,1,2)$. 
By etale descent for isocrystals on relative log convergent site 
of radius $\mu$, this property characterizes $\cF$. So we have the 
uniqueness of $\cF$. \par 
We can check that the isocrystal $\cF$ satisfies the required property in 
the same way as \cite[4.8]{shiho3}. So we are done. 
\end{pf} 

As a corollary, we have the existence of 
a structure of isocrystal on relative log analytic 
cohomology with radius. 

\begin{cor}\label{iso2}
Assume we are given a diagram 
$$ (X,M_X) \os{f}{\lra} (Y,M_Y) \os{g}{\lra} (\cS,M_{\cS}), $$
where $f$ is a proper log smooth integral 
morphism having log smooth parameter 
in $\LB$ and $g$ is a morphism in $\pLF$. Let 
$\mu \in (0,1]$. 
Then, for $\lam \in (p^{-r(X,Y;\mu)},1]$, 
a locally free isocrystal 
$\cE$ on $(X/\cS)^{\log}_{\lamm\conv}$ and a non-negative integer $q$, 
there exists a unique isocrystal $\cF$ on $(Y/\cS)^{\log}_{\mum\conv}$ 
such that, for any pre-widening 
$\cZ:=((\cZ,M_{\cZ}),(Z,M_Z),i,z)$ of $(Y,M_Y)/(\cS,M_{\cS})$
 such that $z$ is strict and $(\cZ,M_{\cZ})$ is formally 
 log smooth over 
$(\cS,M_{\cS})$, $\cF$ induces, via the functor 
$$ I_{\mum\conv}((Y/\cS)^{\log}) \lra 
I_{\mum\conv}((Z/\cS)^{\log}) \simeq 
\Strat''_{\mu}((Z \hra \cZ/\cS)^{\log}), $$
an object of the form 
$(R^qf_{X \times_Y Z/\cZ,\mum\an *}\cE, \epsilon)$, where 
$\epsilon$ is the canonical isomorphism 
$$ p_2^*R^qf_{X \times_Y Z/\cZ,\mum\an *}\cE  \os{\sim}{\ra} 
R^qf_{X \times_Y Z/\cZ(1),\mum\an *}\cE \os{\sim}{\leftarrow}
p_1^*R^qf_{X \times_Y Z/\cZ,\mum\an *}\cE. $$ 
$($Here $(\cZ(1),M_{\cZ(1)}):=(\cZ,M_{\cZ}) \times_{(\cS,M_{\cS})} 
(\cZ,M_{\cZ})$ and 
$p_i$ denotes the $i$-th projection $]Z[^{\log}_{\cZ(1),\mu} \lra 
]Z[^{\log}_{\cZ,\mu}.)$ 
\end{cor}

\begin{pf}
It is immediate from Theorems \ref{iso1}, 
\ref{coherence0} and Proposition \ref{33}. 
Detail is left to the reader. 
\end{pf}


\section{Relative log rigid cohomology and its $\lam$-restriction}

In this section, first we introduce the notion of relative log rigid 
cohomology of log pairs, which is a log version of relative rigid cohomology. 
Then we prove a relation between relative rigid cohomology and 
relative log rigid cohomology in certain case. We also introduce the notion 
of $\lambda$-restriction of relative log rigid cohomology for $\lam \in 
(0,1)$ and compare it with relative log rigid cohomology. The results in this 
section is used in the next section, where we prove the coherence and the 
overconvergence of relative rigid cohomology in the case where a given 
morphism of pairs admits a nice log structure. \par 
The treatment of relative log rigid cohomology here is not 
complete. Here we give a part of the theory of relative log rigid 
cohomology which is related to the purpose of this paper. We hope to 
develop a more general theory of relative log rigid cohomology in a future 
paper. \par 
{\it From now on, we assume that the log structure $M_{\cB}$ on 
the base log formal scheme $\cB$ is trivial}. 
All the pairs (resp. all the triples)
(for definition, see \cite[2.1, 2.3]{chts}, \cite[\S 5]{shiho3}) 
are assumed to be 
pairs over $(B,B)$ (resp. triples over $(B,B,\cB)$) 
whose structure morphism is separated of finite 
type and all the morphisms of pairs (resp. triples) 
are assumed to be a separated 
morphism of finite type over $(B,B)$ (resp. $(B,B,\cB)$). 
As for notation on overconvergent isocrystals and 
relative rigid cohomologies, we follow 
the notation in \cite[\S 5]{shiho3}. (Note that some notations are different 
from that in \cite{chts}.) \par
A log pair is a pair of fine log $B$-schemes 
$((X,M_X),(\ol{X},M_{\ol{X}}))$ endowed with a strict open immersion 
$(X,M_X) \hra (\ol{X},M_{\ol{X}})$ over $B$. A log triple is a triple 
$((X,M_X), (\ol{X},M_{\ol{X}}), (\cP,M_{\cP}))$, where 
$((X,M_X), (\ol{X},M_{\ol{X}}))$ is a log pair, $(\cP,M_{\cP})$ is a fine log 
formal $\cB$-scheme endowed with a closed immersion 
$(\ol{X},M_{\ol{X}}) \hra (\cP,M_{\cP})$. A morphism of 
log pairs or log triples is defined in natural way. For a log triple 
$((X,M_X),(\ol{X},M_{\ol{X}}),(\cP,M_{\cP}))$, a strict neighborhood 
of $]X[^{\log}_{\cP}$ in $]\ol{X}[^{\log}_{\cP}$ is defined to be 
an admissible open subset $V$ such that 
$\{V,]\ol{X}[^{\log}_P-]X[^{\log}_{\cP}\}$ forms an admissible covering 
of $]\ol{X}[^{\log}_{\cP}$. 
For a sheaf $E$ on a strict neighborhood $V$, 
we define the sheaf of overconvergent sections $j_V^{\d}E$ by 
$j_V^{\d}E := \varinjlim_U \alpha_{U]\ol{X}[^{\log}_X *}\alpha^{-1}_{UV}E$, 
where $U$ runs through the strict neighborhoods contained in $V$ and 
$\alpha_{TS}$ denotes the admissible open immersion $T \hra S$. 
In the case $V=]\ol{X}[^{\log}_{\cP}$, we denote $j^{\d}_VE$ by 
$j^{\d}E$ or $j^{\d}_XE$, by abuse of notation. 
As in the case without log structures, any 
coherent $j_X^{\d}\cO_{]\ol{X}[^{\log}_{\cP}}$-module has the form 
$j^{\d}_VE$ for some strict neighborhood $V$ and a coherent $\cO_{V}$-module 
$E$. (See also the paragraph before Proposition \ref{pp}.) \par 
Now let us assume given a diagram 
\begin{equation}\label{logrig-diag}
(\ol{X},M_{\ol{X}}) \os{f}{\lra} (\ol{Y},M_{\ol{Y}}) 
\os{\iota}{\hra} (\cY,M_{\cY})
\end{equation}
and a log pair $((X,M_X),(\ol{X},M_{\ol{X}}))$, where $f$ is a morphism 
in $\LB$, $(\cY,M_{\cY})$ is an object in $\pLF$ and 
$\iota$ is a closed immersion. In this situation, we define the notion 
of overconvergent isocrystals on 
$((X,M_X),(\ol{X},M_{\ol{X}}))/(\cY,M_{\cY})$. \par 
First let us assume that we have a commutative diagram 
\begin{equation}\label{logrig-emb}
\begin{CD}
(\ol{X},M_{\ol{X}}) @>i>> (\cP,M_{\cP}) \\ 
@VfVV @VgVV \\ 
(\ol{Y},M_{\ol{Y}}) @>{\iota}>> (\cY,M_{\cY}), 
\end{CD}
\end{equation}
where $i$ is a closed immersion and $g$ is formally log smooth. 
For $n \in \N$, we denote by $(\cP(n),M_{\cP(n)})$ the $(n+1)$-fold 
fiber product of $(\cP,M_{\cP})$ over $(\cY,M_{\cY})$. Then we have 
the projections 
$$ 
p_i:\,]\ol{X}[^{\log}_{\cP(1)} \lra ]\ol{X}[^{\log}_{\cP} \,\,\,(i=1,2), 
\qquad 
p_{ij}:\,]\ol{X}[^{\log}_{\cP(2)} \lra ]\ol{X}[^{\log}_{\cP(1)}\,\,\, 
(1 \leq i < j \leq 3) 
$$
and the diagonal $\Delta:\,]\ol{X}[^{\log}_{\cP} \hra 
]\ol{X}[^{\log}_{\cP(1)}$. Then: 

\begin{defn}\label{defocl}
With the above notation, 
the category $I^{\d}(((X,\ol{X})/\cY_K,\cP)^{\log})$ of realization of 
overconvergent isocrystals on $((X,M_X),(\ol{X},M_{\ol{X}}))/\cY_K$ over 
$(\cP,M_{\cP})$ is defined to be the category of pairs $(E,\epsilon)$, where 
$E$ is a coherent $j^{\d}\cO_{]\ol{X}[^{\log}_{\cP}}$-module and 
$\epsilon$ is a $j^{\d}\cO_{]\ol{X}[^{\log}_{\cP(1)}}$-linear 
isomorphism $p_2^*E \os{\sim}{\lra} p_1^*E$ satisfying 
$\Delta^*(\epsilon)= \id, p_{12}^*(\epsilon) \circ p_{23}^*(\epsilon) = 
p_{13}^*(\epsilon)$. 
$(E,\epsilon)$ is called locally free if $E$ is a locally free 
$j^{\d}\cO_{]\ol{X}[^{\log}_{\cP}}$-module. 
\end{defn} 

We give a slight generalization of the above definition as follows: 
Assume given a commutative diagram 
\begin{equation}\label{logrig-diag-frak}
\begin{CD}
(\ol{X},M_{\ol{X}}) @>i>> (\fP,M_{\fP}) \\ 
@VfVV @VgVV \\ 
(\ol{Y},M_{\ol{Y}}) @>{\iota}>> (\cY,M_{\cY}), 
\end{CD}
\end{equation}
where $(\fP,M_{\fP})$ is (not necessarily $p$-adic) fine log formal 
$\cB$-scheme, $i$ is a closed immersion, and 
$g$ is a morphism in $\LF$ satisfying the following condition $(*)$: \\
\quad \\
$(*)$ \,\,\, Zariski locally on $\fP$, there exists a diagram 
\begin{equation}
\begin{CD}
(\ol{X},M_{\ol{X}}) @>i>> (\fP,M_{\fP}) 
@>{i'}>> (\ti{\fP},M_{\ti{\fP}}) \\ 
@VfVV @VgVV @V{g'}VV \\ 
(\ol{Y},M_{\ol{Y}}) @>{\iota}>> (\cY,M_{\cY}) @= (\cY,M_{\cY}), 
\end{CD}
\end{equation}
where $g'$ is a formally log smooth morphism in $\pLF$, 
$i'$ is a morphism in $\LF$ such that $i' \circ i$ is again a closed 
immersion such that $i'$ induces the isomorphism 
$(\fP^{\ex},M_{\fP^{\ex}}) \os{=}{\lra} (\ti{\fP}^{\ex},M_{\ti{\fP}^{\ex}})$, 
where ${}^{\ex}$ denotes the exactification of the closed immersion from 
$(\ol{X},M_{\ol{X}})$. \\
\quad \\
For $n \in \N$, we denote by $(\fP(n),M_{\fP(n)})$ the $(n+1)$-fold 
fiber product of $(\fP,M_{\fP})$ over $(\cY,M_{\cY})$ and we denote 
$(\fP(n)^{\ex},M_{\fP(n)^{\ex}})$ the exactification of 
$(\ol{X},M_{\ol{X}}) \hra (\fP(n),M_{\fP(n)})$. 
(So we have $\fP(0)^{\ex}=\fP^{\ex}$.) 
Let ${\fP(n)^{\ex}}'$ be 
the open sub formal scheme of $\fP(n)^{\ex}$ which is homeomorphic to 
$X$. Then, starting from the admissible open 
immersion ${\fP(n)^{\ex}}'_K \hra 
\fP(n)^{\ex}_K$, we can define the structure sheaf of overconvergent 
sections $j^{\d}\cO_{\fP(n)^{\ex}_K}$. We also have the projections 
$p_i,p_{ij}$ and the diagonal $\Delta$ between $\fP(n)^{\ex}$'s. 
Then we have the following generalization of Definition \ref{defocl}: 

\begin{defn}\label{defocl2}
With the above notation, 
the category $I^{\d}(((X,\ol{X})/\cY_K,\fP)^{\log})$ of realization of 
overconvergent isocrystals on $((X,M_X),(\ol{X},M_{\ol{X}}))/\cY_K$ over 
$(\fP,M_{\fP})$ is defined to be the category of pairs $(E,\epsilon)$, where 
$E$ is a coherent $j^{\d}\cO_{\fP^{\ex}_K}$-module and 
$\epsilon$ is a $j^{\d}\cO_{\fP(1)^{\ex}_K}$-linear 
isomorphism $p_2^*E \os{\sim}{\lra} p_1^*E$ satisfying 
$\Delta^*(\epsilon)= \id, p_{12}^*(\epsilon) \circ p_{23}^*(\epsilon) = 
p_{13}^*(\epsilon)$. 
\end{defn} 

Note that, in the notation in Definition \ref{defocl2}, 
we have the canonical equivalence of categories 
$$ I^{\d}(((X,\ol{X})/\cY,\fP)^{\log}) \cong 
I^{\d}(((X,\ol{X})/\cY,\fP^{\ex})^{\log}). $$

\begin{lem}\label{logrig-des1}
With the notation in Definition \ref{defocl2}, the category 
$I^{\d}(((X,\ol{X})/\cY,\fP)^{\log})$ satisfies the descent property 
for formally etale covering of $\fP$. 
\end{lem} 

\begin{pf} 
We can prove this lemma in the same way as etale descent of 
the category of overconvergent isocrystals proved in \cite[5.1]{shiho3}. 
The detail is left to the reader. 
\end{pf} 

\begin{lem}\label{logrig-des2}
With the notation in Definition \ref{defocl}, assume moreover that 
we are given the commutative diagram 
\begin{equation}
\begin{CD}
(\ol{X},M_{\ol{X}}) @>>> (\cQ,M_{\cQ}) \\ 
@\vert @VhVV \\
(\ol{X},M_{\ol{X}}) @>i>> (\cP,M_{\cP}), 
\end{CD}
\end{equation} 
where the top horizontal arrow is a closed immersion and $h$ is a 
formally log smooth morphism in $\pLF$. Then the restriction functor 
$$ h^*: I^{\d}(((X,\ol{X})/\cY,\cP)^{\log}) \lra 
I^{\d}(((X,\ol{X})/\cY,\cQ)^{\log}) $$ 
is an equivalence of categories. 
\end{lem} 

\begin{pf} 
We may replace $(\cP,M_{\cP}), (\cQ,M_{\cQ})$ by 
$(\cP^{\ex},M_{\cP^{\ex}}), (\cQ^{\ex},M_{\cQ^{\ex}})$ to prove the lemma, 
and we may work Zariski locally on $\cP^{\ex}$. So we may assume that 
$\cP^{\ex}=\Spf A$ is affine and the map $\cQ^{\ex} \lra \cP^{\ex}$ is 
given by the natural morphism 
$\cQ^{\ex} = \Spf A[[t_1, \cdots, t_r]] \lra \Spf A = \cP^{\ex}$. Let 
$s: \cP^{\ex} \hra \cQ^{\ex}$ be the map given by 
the closed immersion of $\Spf A$ into the origin of 
$\Spf A[[t_1, \cdots, t_r]]$. Then we have the functor 
$s^*: I^{\d}(((X,\ol{X})/\cY,\cQ)^{\log}) \lra 
I^{\d}(((X,\ol{X})/\cY,\cP)^{\log})$ with $s^* \circ h^* = \id$. 
On the other hand, for any $(E,\epsilon) \in 
I^{\d}(((X,\ol{X})/\cY,\cQ)^{\log})$, we have the isomorphism 
$\varphi_E: (h^* \circ s^*)(E) \lra E$ given by pulling back the isomorphism 
$\epsilon: p_2^*E \os{\sim}{\lra} p_1^*E$ by 
$(h \circ s, \id)_K: \cQ^{\ex}_K \lra \cP^{\ex}(1)^{\ex}_K$. 
Since one can check easily the commutativity 
$\epsilon \circ p_2^*\varphi_E = p_1^*\varphi_E \circ 
(h^* \circ s^*)\epsilon$ and the functoriality of $\varphi_E$ with respect 
to $(E,\epsilon)$, we see that $\varphi_E$'s induce the isomorphism 
of functors $h^* \circ s^* \cong \id$. So $h^*$ is an equivalence 
of categories. 
\end{pf} 

Now, let us go back to the situation \eqref{logrig-diag} and 
let us define the category of overconvergent isocrystals on 
$((X,M_X),(\ol{X},M_{\ol{X}}))/(\cY,M_{\cY})$ in general case, that is, 
the case where we do not necessarily have the diagram 
\eqref{logrig-emb}. In this case, we have an embedding system 
\begin{equation}\label{logrig-emb2}
(\ol{X},M_{\ol{X}}) \os{g}{\lla} (\ol{X}^{(\b)},M_{\ol{X}^{(\b)}})
\os{i^{(\b)}}{\hra} (\cP^{(\b)},M_{\cP^{(\b)}}) 
\end{equation}
over $(\cY,M_{\cY})$. Let us put 
$(X^{(\b)},M_{X^{(\b)}}) := (X,M_X) \times_{(\ol{X},M_{\ol{X}})} 
(\ol{X}^{(\b)},M_{\ol{X}^{(\b)}})$. 
Then, we denote by 
$I^{\d}(((X^{(\b)},\ol{X}^{(\b)})/\cY, \cP^{(\b)})^{\log})$ 
the category of descent data with respect to 
$I^{\d}(((X^{(n)},\ol{X}^{(n)})/\cY, \cP^{(n)})^{\log}) \,(n=0,1,2)$. 
Then we have the following: 

\begin{lem}\label{logrig-indep}
With the above notation, the category 
$I^{\d}(((X^{(\b)},\ol{X}^{(\b)})/\cY, \cP^{(\b)})^{\log})$ 
is independent of the choice of the embedding system. 
\end{lem} 

\begin{pf} 
By standard argument, we are reduced to showing the following claim: 
Assume we are given the diagram 
\begin{equation*}
\begin{CD}
(\ol{X}^{(\b)},M_{\ol{X}^{(\b)}}) @>{i^{(\b)}}>> 
(\cP^{(\b)},M_{\cP^{(\b)}}) \\ 
@VgVV @VhVV \\ 
(\ol{X},M_{\ol{X}}) @>i>> (\cP,M_{\cP}), 
\end{CD}
\end{equation*}
where $i^{(\b)},g$ are as in \eqref{logrig-emb2}, 
$h$ is a formally log smooth morphism in $\pLF$, 
$i$ is a closed immersion and 
$(\cP,M_{\cP})$ is a $p$-adic fine log formal $\cB$-scheme 
formally log smooth over $(\cY,M_{\cY})$. Then the restriction functor 
$$ h^*: I^{\d}(((X,\ol{X})/\cY_K,\cP)^{\log}) \lra 
I^{\d}(((X^{(\b)},\ol{X}^{(\b)})/\cY_K,\cP^{(\b)})^{\log}) $$ 
gives an equivalence of categories. (Here 
$I^{\d}(((X^{(\b)},\ol{X}^{(\b)})/\cY_K,\cP^{(\b)})^{\log})$ denotes the 
category of descent data with respect to  
$I^{\d}(((X^{(n)},\ol{X}^{(n)})/\cY_K,\cP^{(n)})^{\log}) \, (n=0,1,2)$.) 
We prove this claim. 
Let $(\cP^{\ex},M_{\cP^{\ex}})$ be 
the exactification of $i$ and for $n \in \N$, 
let $(\cP^{\ex,(n)},M_{\cP^{\ex,(n)}})$ be the unique fine log formal 
$\cB$-scheme strict formally etale over $(\cP^{\ex},M_{\cP^{\ex}})$ 
satisfying 
$(\ol{X}^{(n)},M_{\ol{X}^{(n)}}) = 
(\ol{X},M_{\ol{X}}) \times_{(\cP^{\ex},M_{\cP^{\ex}})} 
(\cP^{\ex,(n)},M_{\cP^{\ex,(n)}})$. Then, by Lemma \ref{logrig-des1}, 
we have the equivalence 
$$I^{\d}(((X,\ol{X})/\cY,\cP)^{\log}) \cong 
I^{\d}(((X^{(\b)},\ol{X}^{(\b)})/\cY,\cP^{\ex,(\b)})^{\log}).$$
(Here $I^{\d}(((X^{(\b)},\ol{X}^{(\b)})/\cY,\cP^{\ex,(\b)})^{\log})$ 
denotes the category of descent data with respect to  
$I^{\d}(((X^{(n)},\ol{X}^{(n)})/\cY,\cP^{\ex,(n)})^{\log}) \, (n=0,1,2)$.) 
On the other hand, let $(\cP^{(n),\ex},\allowbreak M_{\cP^{(n),\ex}})$ be the 
exactification of $i^{(n)}$. Then we have the morphism 
$$ h^{(n)}: 
(\cP^{(n),\ex},\allowbreak M_{\cP^{(n),\ex}}) \allowbreak \lra 
\allowbreak (\cP^{\ex,(n)},M_{\cP^{\ex,(n)}}) $$
 induced by $h$. Then, it suffices to prove that the restriction functor 
$$ h^{(n),*}: 
I^{\d}(((X^{(n)},\ol{X}^{(n)})/\cY,\cP^{\ex,(n)})^{\log})
\lra 
I^{\d}(((X^{(n)},\ol{X}^{(n)})/\cY,\cP^{(n),\ex})^{\log})
$$ 
gives an equivalence of categories. 
Let 
$((\cP^{\ex,(n)},M_{\cP^{\ex,(n)}}) 
\times_{(\cY,M_{\cY})} (\cP^{(n),\ex},M_{\cP^{(n),\ex}}))^{\ex}$ be the 
exactification of the closed immersion 
$$(\ol{X}^{(n)},M_{\ol{X}^{(n)}}) \hra 
(\cP^{\ex,(n)},M_{\cP^{\ex,(n)}}) 
\times_{(\cY,M_{\cY})} (\cP^{(n),\ex},M_{\cP^{(n),\ex}}).$$ 
Let $$p_i: 
((\cP^{\ex,(n)},M_{\cP^{\ex,(n)}}) 
\times_{(\cY,M_{\cY})} (\cP^{(n),\ex},M_{\cP^{(n),\ex}}))^{\ex} \lra 
(\cP^{\ex,(n)},M_{\cP^{\ex,(n)}}) \,\,\,\, (i=1,2)$$
be the morphism induced by projections 
and let 
$$\gamma: (\cP^{(n),\ex},M_{\cP^{(n),\ex}}) \lra 
((\cP^{\ex,(n)},M_{\cP^{\ex,(n)}}) 
\times_{(\cY,M_{\cY})} (\cP^{(n),\ex},M_{\cP^{(n),\ex}}))^{\ex}$$ 
be the morphism induced by the graph of $h^{(n)}$. Then 
the restriction functor $h^{(n),*}$ is written as the composite 
$\gamma^* \circ p_1^*$. Note now that 
the fine log formal $\cB$-schemes 
$(\cP^{\ex,(n)},M_{\cP^{\ex,(n)}}), (\cP^{(n),\ex},M_{\cP^{(n),\ex}})$ 
are obtained Zariski locally as the exactification of some 
$p$-adic fine log formal $\cB$-schemes which are formally log smooth 
over $(\cY,M_{\cY})$. (It is clear in the case of 
$(\cP^{(n),\ex},M_{\cP^{(n),\ex}})$ and in the case 
$(\cP^{\ex,(n)},M_{\cP^{\ex,(n)}})$, this follows from 
\cite[claim in p.81]{shiho2}.) So the morphisms 
$p_i\,(i=1,2)$ are obtained, via taking exactification, from 
formally log smooth morphisms (projections) between $p$-adic fine 
log formal $\cB$-schemes which are formally log smooth over $(\cY,M_{\cY})$. 
So, by Lemma \ref{logrig-des2}, the restriction functors 
$p_i^*\,(i=1,2)$ are equivalences. On the other hand, since 
$\gamma$ is a section of $p_2$, we see that the functors 
$p_2^*$ and $\gamma^*$ are the inverse of each other. (The proof is the same 
as the proof in Lemma \ref{logrig-des2}.) So the functor $\gamma^*$ is 
also an equivalence. So $h^{(n),*}=\gamma^* \circ p_1^*$ is an equivalence 
of categories. So we are done. 
\end{pf} 

So we can define the category of overconvergent isocrystals as follows: 

\begin{defn}
Let us assume given the diagram \eqref{logrig-diag} and take an embedding 
system \eqref{logrig-emb2}. Then we define the category 
$I^{\d}(((X,\ol{X})/\cY)^{\log})$ 
of overconvergent 
isocrystals on $((X,M_X),(\ol{X},M_{\ol{X}}))/(\cY,M_{\cY})$ by 
$I^{\d}(((X,\ol{X})/\cY)^{\log}) := 
I^{\d}(((X^{(\b)},\ol{X}^{(\b)})/\cY, \allowbreak \cP^{(\b)} \allowbreak
)^{\log})$. 
$($By Lemma \ref{logrig-indep}, this is independent of the choice of 
the embedding system.$)$ 
\end{defn} 

Let us assume given a diagram \eqref{logrig-diag} and a log pair 
$((X,M_X),(\ol{X},M_{\ol{X}}))$. Assume moreover for the moment that 
we have the commutative diagram \eqref{logrig-emb}. Then, for an object 
$\cE := (E,\epsilon)$ in 
$I^{\d}(((X,\ol{X})/\cY_K)^{\log})=I^{\d}(((X,\ol{X})/\cY_K,\cP)^{\log})$, 
we can naturally define the log de Rham complex of the form 
$$ E \lra E \otimes_{j^{\d}\cO_{]\ol{X}[^{\log}_{\cP}}} 
j^{\d}\omega^1_{]\ol{X}[^{\log}_{\cP}/\cY_K}
 \lra E \otimes_{j^{\d}\cO_{]\ol{X}[^{\log}_{\cP}}} 
j^{\d}\omega^2_{]\ol{X}[^{\log}_{\cP}/\cY_K}
 \lra \cdots 
$$ 
(where 
$\omega^i_{]\ol{X}[^{\log}_{\cP}/\cY_K}$ 
is the restriction of $\Q \otimes \omega^i_{\cP/\cY} \in 
\Coh(\Q \otimes \cO_{\cP}) \cong \Coh(\cP_K)$ to 
$]\ol{X}[^{\log}_{\cP}$) on $]\ol{X}[^{\log}_{\cP}$, 
which we denote by 
$\DR^{\d}(]\ol{X}[^{\log}_{\cP}/\cY_K, \cE)$. \par 
In the case where we have the diagram \eqref{logrig-diag-frak} 
(satisfying the condition $(*)$) instead of the diagram \eqref{logrig-emb}, 
we also have a similar log de Rham complex on $\fP^{\ex}_K$, 
which we denote by $\DR^{\d}(\fP^{\ex}_K/\cY_K, \cE)$. \par

Now let us assume given the diagram \eqref{logrig-diag} and 
log pairs $((X,M_X),(\ol{X},M_{\ol{X}})), \allowbreak ((Y, 
\allowbreak M_Y),(\ol{Y},M_{\ol{Y}}))$ 
satisfying $f(X) \subseteq Y$ and an overconvegent isocrystal $\cE$ on 
$((X,M_X),(\ol{X},M_{\ol{X}}))/ \allowbreak (\cY, 
\allowbreak M_{\cY})$. Let us take the embedding system 
\eqref{logrig-emb2} and denote the morphism 
$]\ol{X}^{(\b)}[^{\log}_{\cP^{(\b)}} \lra ]\ol{Y}[^{\log}_{\cY}$ by $h$. 
Then we define 
$Rf_{(X,\ol{X})/\cY, \logrig *}\cE, R^qf_{(X,\ol{X})/\cY,}\allowbreak
{}_{\logrig *}\cE$ by 
\begin{align*}
& Rf_{(X,\ol{X})/\cY, \logrig *}\cE := 
Rh_*\DR^{\d}(]\ol{X}^{(\b)}[^{\log}_{\cP^{(\b)}}/\cY_K, \cE), \\ 
& R^qf_{(X,\ol{X})/\cY, \logrig *}\cE := 
R^qh_*\DR^{\d}(]\ol{X}^{(\b)}[^{\log}_{\cP^{(\b)}}/\cY_K, \cE) 
\end{align*}
and we call $R^qf_{(X,\ol{X})/\cY, \logrig *}\cE$ the $q$-th 
relative log rigid cohomology. 
(It is an $j^{\d}\cO_{]\ol{Y}[^{\log}_{\cY}}$-module.)
Then we have the following: 

\begin{prop}\label{logrig-wd}
The above definition is independent of the choice of 
the embedding system. 
\end{prop}

Before the proof of Proposition \ref{logrig-wd}, we prove a lemma: 

\begin{lem}\label{acyc0}
Let $\cP$ be a $($not necessarily $p$-adic$)$ affine 
fine log formal $\cB$-scheme 
and let $g_1, \cdots, g_c$ be elements of $\Gamma(\cP,\cO_{\cP})$. 
For $\nu \in p^{\Q_{<0}}$, let $U_{\nu}$ be the admissible open set 
$\bigcup_{i=1}^c \{|g_i| \geq \nu\}$ of $\cP_K$. 
Then, for any affinoid admissible open set $W$ of $\cP_K$ admits 
an affinoid admissible 
open covering $W = \bigcup_{i=0}^cW_i$ such that 
$W_i \cap U_{\nu}$ is affinoid for any $\nu$ sufficiently close to $1$ and for 
any $0 \leq i \leq c$. 
\end{lem} 

\begin{pf} 
For $\nu \in p^{\Q_{<0}}$ and $1 \leq i \leq c$, 
define $U_{\nu,i}$ to be the admissible open set 
$\{|g_i| \geq \nu, |g_i| \geq |g_j|\,(\forall j \not= i)\}$ of $\cP_K$ and 
define $U_{\nu,0}$ to be the admissible open set 
$\{|g_j| \leq \nu \,(\forall j)\}$. Then 
$\cP_K = \bigcup_{i=0}^c U_{\nu,i}$ is an admissible covering. 
If we fix $\nu_0 \in p^{\Q_{<0}}$ and put $W_i := W \cap U_{\nu_0,i}$, 
we obtain the admissible covering $W = \bigcup_{i=0}^c W_i$ by 
affinoid subsets. Moreover, for $\nu > \nu_0$, we have 
$W_0 \cap U_{\nu} = \emptyset$ and for $1 \leq i \leq c$, 
$W_i \cap U_{\nu} = W \cap U_{\nu,i}$ is affinoid. So we are done. 
\end{pf} 

\begin{pf*}{Proof of Proposition \ref{logrig-wd}} 
The proof is similar to \cite[4.3]{shiho3} and the case without log 
structures (\cite[8.3.5]{chts}, see also \cite[1.4]{berthelot3}). 
By looking at carefully the proof of \cite[4.3]{shiho3}, we see that it 
suffices to prove the following two claims (cf. \cite[4.4, 4.5]{shiho3}): \\
\quad \\
{\bf claim 1.} \,\,\, Assume we are given a Cartesian diagram 
\begin{equation}
\begin{CD}
(\ol{X}^{(\b)},M_{\ol{X}^{(\b)}}) @>>> (\cP^{(\b)},M_{\cP^{(\b)}}) \\ 
@VVV @VgVV \\ 
(\ol{X},M_{\ol{X}}) @>{\iota}>> (\cP,M_{\cP}), 
\end{CD}
\end{equation}
where $(\ol{X},M_{\ol{X}})$ is as above, $(\cP,M_{\cP})$ is an object in 
$\LF$, $\iota$ is a closed immersion and $g$ is a strict etale hypercovering. 
Let $g_K:]\ol{X}^{(\b)}[^{\log}_{\cP^{(\b)}} \lra ]\ol{X}[^{\log}_{\cP}$ 
be the morphism induced by $g$ and let $\cE$ be a coherent 
$j^{\d}\cO_{]\ol{X}[^{\log}_{\cP}}$-module. Then we have the isomorphism 
$$ \cE \os{=}{\lra} Rg_{K,*}g^*_K \cE. $$ 

\quad \\
{\bf claim 2.} \,\,\, Assume we are given a diagram over $(\cY,M_{\cY})$ 
\begin{equation}
\begin{CD}
(\ol{X},M_{\ol{X}}) @>{\iota_1}>> (\cQ,M_{\cQ}) \\ 
@\vert @VgVV \\ 
(\ol{X},M_{\ol{X}}) @>{\iota_2}>> (\cP,M_{\cP}), 
\end{CD}
\end{equation}
where $(\ol{X},M_{\ol{X}})$ is as above, 
$(\cP,M_{\cP}),(\cQ,M_{\cQ})$ are objects in 
$\pLF$, $\iota_1,\iota_2$ are closed immersions and $g$ is a formally log 
smooth 
morphism. Let $g_K:]\ol{X}[^{\log}_{\cQ} \lra ]\ol{X}[^{\log}_{\cP}$ be the 
morphism induced by $g$. Then we have a quasi-isomorphism 
$$ \DR^{\d}(]\ol{X}[^{\log}_{\cP}/\cY_K, \cE) 
\os{=}{\lra} Rg_{K,*}\DR^{\d}(]\ol{X}[^{\log}_{\cQ}/\cY_K,\cE).$$ 

\quad \\
The proof of claim 1 can be reduced to the case where $\iota$ is a 
homeomorphic exact closed immersion into a $p$-adic fine log formal 
$\cB$-scheme 
(see the proof of 
\cite[4.4]{shiho3}) and 
in this case, the claim 
is a special case of the etale cohomological 
descent of Chiarellotto-Tsuzuki (\cite[7.1.2]{chts}). \par 
Let us prove the claim 2. 
We may work Zariski locally on $\cP$. 
Let $(\cQ^{\ex},M_{\cQ^{\ex}}), 
(\cP^{\ex},\allowbreak M_{\cP^{\ex}})$ be the exactification of 
$\iota_1, \iota_2$, respectively. Since we may work Zariski locally on 
$\cP^{\ex}$, we may assume that $\cQ^{\ex}_{K}$ is isomorphic to 
$\cP^{\ex}_{K} \times D^r$ and the morphism $g_K$ is equal to the 
projection $\cP^{\ex}_{K} \times D^r \lra \cP^{\ex}_{K}$. 
Let $(t_1, \cdots, t_r)$ be the coordinate of $D^r$ and let 
$\Omega^{\b}$ be the complex 
$$ 
[\cO_{\cQ^{\ex}_{K}} \lra 
\bigoplus_{i=1}^r\cO_{\cQ^{\ex}_{K}} dt_i \lra 
\bigoplus_{1 \leq i_1 < i_2 \leq r} \cO_{\cQ^{\ex}_{K}} dt_{i_1} \wedge 
dt_{i_2} \lra \cdots \lra 
\cO_{\cQ^{\ex}_{K}} dt_1 \wedge \cdots \wedge dt_r]. $$ 
Then $\DR^{\d}(]X[^{\log}_{\cQ}/\cY,\cE)$ is equal to the total complex 
associated to the double complex 
$g_K^*\DR^{\d}(]X[^{\log}_{\cP}/\cY,\cE) 
\otimes_{\cO_{\cQ^{\ex}_{K}}} \Omega^{\b}$. 
So, to prove the lemma, it suffices to prove the quasi-isomorphism 
\begin{equation}\label{logrig-proj}
E \os{\simeq}{\lra} Rg_{K,*}(g_K^*E \otimes \Omega^{\b})
\end{equation}
for any coherent $j_X^{\d}\cO_{\cP^{\ex}_{K}}$-module $E$. 
Moreover, to prove the quasi-isomorphism \eqref{logrig-proj}, 
we can reduce to the case $r=1$ by induction. So we may assume $r=1$. 
Let 
$\alpha_1, \cdots, \alpha_c \in 
\Gamma(\cP^{\ex},\cO_{\cP^{\ex}})$ 
be a lift of generators of $\Ker (\Gamma(\ol{X},\cO_{\ol{X}}) \ra 
\Gamma(\ol{X}-X,\cO_{\ol{X}-X})$ 
and for $\nu \in p^{\Q_{<0}}$, let 
$U_{\nu}$ be the admissible open subset $\bigcup_{i=1}^c\{|\alpha_i| 
\geq \nu\}$ in $\cP^{\ex}_K$. 
To prove the quasi-isomorphism \eqref{logrig-proj}, it suffices to check it 
on sufficiently small affinoid admissible open subset 
of $\cP^{\ex}_K$. So it suffices to check on an affinoid admissible open 
subset $W$ 
of $\cP^{\ex}_K$ such that $W_{\nu} := W \cap U_{\nu}$ is also affinoid 
for any $\nu$ sufficiently close to $1$, by Lemma \ref{acyc0}. 
For $\lam, \nu \in p^{\Q_{<0}}$, let 
$g_K^{-1}(W)_{\lam,\nu}$ be the affinoid admissible open set 
$\{|t_1|\leq \lam\} \cap g_K^{-1}(W_{\nu})$ of $g_K^{-1}(W)$. 
When $\nu$ is sufficiently close to $1$,  
$E$ is obtained from a coherent $\cO_{\cP^{\ex}_K}$-module 
on $W_{\nu}$, which we denote by $E_{\nu}$. 
Then we have 
$$ \Gamma(W,E)= 
\dvil_{\nu\to1}\Gamma(W_{\nu},E_{\nu}), $$ 
$$ R\Gamma(g_K^{-1}(W),g_K^{*}E \otimes \Omega^{\b})= 
R\dvpl_{\lam\to1}\dvil_{\nu\to1}
\Gamma(g_K^{-1}(W)_{\lam,\nu},g_K^*E_{\nu} \otimes \Omega^{\b}). $$ 
So, to prove the quasi-isomorphism \eqref{logrig-proj}, it suffices 
to prove the following claim: \\
\quad \\
{\bf claim 3.}\,\,\, Under the above situation, \\
(1)\,\, The complex 
\begin{align*}
0 \ra \dvil_{\nu\to 1}\Gamma(W_{\nu},E_{\nu}) & \ra 
\dvpl_{\lam\to 1}\dvil_{\nu\to 1}
\Gamma(g_K^{-1}(W)_{\lam,\nu},g_K^*E_{\nu}) \\ & 
\os{d}{\ra} 
\dvpl_{\lam\to 1}\dvil_{\nu\to 1}
\Gamma(g_K^{-1}(W)_{\lam,\nu},g_K^*E_{\nu} \otimes \Omega^1) \ra 0
\end{align*}
is exact. \\
(2)\,\, The map 
$$ R^1\dvpl_{\lam\to 1}\dvil_{\nu\to 1}
\Gamma(g_K^{-1}(W)_{\lam,\nu},g_K^{*}E_{\nu}) \lra 
R^1\dvpl_{\lam\to 1}\dvil_{\nu\to 1}
\Gamma(g_K^{-1}(W)_{\lam,\nu},g_K^{*}E_{\nu} \otimes \Omega^1)$$ 
is an isomorphism. \\
\quad \\
To prove claim 3, note first that we have 
$$ 
\Gamma(g_K^{-1}(W)_{\lam,\nu},g_K^{*}E_{\nu})= \{\sum_{n=0}^{\infty}a_nt_1^n
\,\vert\, a_n \in \Gamma(W_{\nu},E_{\nu}), |a_n|\lam^n\to 0\}, 
$$ 
where $|-|$ denotes the Banach norm on $\Gamma(W_{\nu},E_{\nu})$, 
and similar description holds also 
for $\Gamma(g_K^{-1}(W)_{\lam,\nu},g_K^{*}E_{\nu} \otimes \Omega^1)$. 
Then, 
$\sum_{n=0}^{\infty}a_nt_1^n \in 
\Gamma(g_K^{-1}(W)_{\lam,\nu},g_K^{*}E_{\nu})$ is sent to zero by $d$ 
if and only if 
$a_n=0$ for $n>0$, and for any 
$\sum_{n=0}^{\infty}a_nt_1^ndt_1 \in 
\Gamma(g_K^{-1}(W)_{\lam,\nu},g_K^{*}E_{\nu} \otimes \Omega^1)$, 
$\sum_{n=0}^{\infty}(a_n/(n+1))t_1^{n+1}$ is contained in 
$\Gamma(g_K^{-1}(W)_{\lam',\nu},g_K^{*}E_{\nu})$ for any $\lam'<\lam$. 
From these facts, we see claim 3 (1). Claim 3 (2) can be shown 
as \cite[Thm 1.4]{berthelot3}, using claim 3 (1). So we are done. 
\end{pf*}

\begin{rem} 
The final step of the proof of 
\cite[Prop 8.3.5]{chts} is incomplete, because claim 3 (2) is omitted in 
their proof. 
\end{rem} 

By Proposition \ref{logrig-wd}, we see that 
$Rf_{(X,\ol{X})/\cY,\logrig *}\cE$, 
$R^qf_{(X,\ol{X})/\cY,\logrig *}\cE$ are well-defined. \par 
We introduce another construction of relative log rigid cohomology in 
log smooth case. 
Assume given the diagram \eqref{logrig-diag} and 
log pairs $((X,M_X),(\ol{X},M_{\ol{X}})), ((Y,M_Y),\allowbreak 
(\ol{Y}, \allowbreak M_{\ol{Y}}))$ 
Assume moreover that $f$ is log 
smooth. Then we have an open covering 
$\ol{X} = \bigcup_{j \in J}\ol{X}_j$ by finite number of affine subschemes and 
exact closed immersions $i_j: 
(\ol{X}_j,M_{\ol{X}_j}) := 
(\ol{X}_j,M_{\ol{X}} |_{X_j}) \hra (\cP_j,M_{\cP_j}) \,(j \in J)$ 
into a fine log formal $\cB$-scheme $(\cP_j,M_{\cP_j})$ formally log 
smooth over $(\cY,M_{\cY})$ such that each $\cP_j$ is affine and 
$(\cP_j,M_{\cP_j}) \times_{(\cY,M_{\cY})} (\ol{Y},M_{\ol{Y}}) = 
(\ol{X}_j,M_{\ol{X}_j})$ 
 holds. For any 
 non-empty subset $L \subseteq J$, let 
$(\ol{X}_L,M_{\ol{X}_L})$ (resp. $(\cP_L, M_{\cP_L})$) 
be the 
fiber product of $(\ol{X}_j,M_{\ol{X}_j})$'s (resp. 
$(\cP_j,M_{\cP_j})$'s) for $j \in L$ over 
$(\ol{X},M_{\ol{X}})$ (resp. $(\cY,M_{\cY})$), and for 
$m \in \N$, 
let $(\ol{X}^{(m)}, M_{\ol{X}^{(m)}})$ (resp. $(\cP^{(m)},M_{\cP^{(m)}})$) 
be the disjoint union of $(\ol{X}_L,M_{\ol{X}_L})$'s (resp. 
$(\cP_L,M_{\cP_L})$'s) for $L \subseteq J$ with $|L|=m+1$. 
Note that if $m$ is equal to or greater than 
$|J|=:N+1$, $\ol{X}^{(m)}$ is empty. Now let $\Delta^+_N$ be the category 
such that the 
objects are the sets $[m]:=\{0,1,2,\cdots m\}$ with $m \leq N$ and that 
$\Hom_{\Delta^+_N}([m],[m'])$ is the set of strictly increasing 
maps $[m] \ra [m']$. Then, if we fix a total order on the set $J$, 
$(\ol{X}^{(m)},M_{\ol{X}^{(m)}})$'s 
and $(\cP^{(m)},M_{\cP^{(m)}})$'s $(m \leq N)$ 
naturally form diagrams 
$(\ol{X}^{(\b)},M_{\ol{X}^{(\b)}}), (\cP^{(\b)},M_{\cP^{(\b)}})$ 
indexed by the category $\Delta^+_N$. (It is not a truncated simplicial 
object since it does not have degeneracy maps.) Then we have a 
canonical diagram 
\begin{equation}\label{logrig-finemb}
(\ol{X},M_{\ol{X}}) \os{g^{(\b)}}{\lla} (\ol{X}^{(\b)},M_{\ol{X}^{(\b)}}) 
\hra (\cP^{(\b)},M_{\cP^{(\b)}}) 
\end{equation} 
over $(\cY,M_{\cY})$. Then, for an overconvergent isocrystal $\cE$ 
on $((X,M_X),(\ol{X},M_{\ol{X}}))/(\cY,\allowbreak 
M_{\cY})$, we have the log de Rham 
complex $\DR^{\d}(]\ol{X}^{(\b)}[^{\log}_{\cP^{(\b)}}/\cY_K, \cE)$. 
Let us denote the morphism 
$]\ol{X}^{(\b)}[^{\log}_{\cP^{(\b)}} \lra ]\ol{Y}[^{\log}_{\cY}$ by $h$. 
Then we have the following: 

\begin{prop}\label{logrig-fin}
With the above notation, we have the quasi-isomorphism 
$$ Rh_*\DR^{\d}(]\ol{X}^{(\b)}[^{\log}_{\cP^{(\b)}}/\cY_K, \cE) 
= Rf_{(X,\ol{X})/\cY, \logrig *}\cE. $$
\end{prop}

\begin{pf} 
For $m \in \N$, let 
$(\ti{X}^{(m)},M_{\ti{X}^{(m)}})$ 
(resp. $(\ti{\cP}^{(m)},M_{\ti{\cP}^{(m)}})$) be the $(m+1)$-fold fiber 
product of $(\ol{X}^{(0)},M_{\ol{X}^{(0)}})$ 
(resp. $(\cP^{(0)},M_{\cP^{(0)}})$) over $(\ol{X},M_{\ol{X}})$ (resp. 
$(\cY,M_{\cY})$). Then we have the embedding system 
$$ 
(\ol{X},M_{\ol{X}}) \lla (\ti{X}^{(\b)},M_{\ti{X}^{(\b)}}) \hra 
(\ti{\cP}^{(\b)},M_{\ti{\cP}^{(\b)}}) $$ 
and 
each $(\ol{X}^{(m)},M_{\ol{X}^{(m)}}) \hra (\cP^{(m)},M_{\cP^{(m)}})$ 
is contained as a direct summand in 
$(\ti{X}^{(m)},M_{\ti{X}^{(m)}}) \hra (\ti{\cP}^{(m)},M_{\ti{\cP}^{(m)}})$
 in a compatible way with respect to $m$ (indexed by $\Delta^+_N$). 
So, if we denote the morphism 
$]\ti{X}^{(\b)}[^{\log}_{\cP^{(\b)}} \lra ]\ol{Y}[^{\log}_{\cY}$ by 
$\ti{h}$, 
we have the homomorphism 
$$ 
Rf_{(X,\ol{X})/\cY, \logrig *}\cE = 
R\ti{h}_*\DR^{\d}(]\ti{X}^{(\b)}[^{\log}_{\ti{\cP}^{(\b)}}/\cY_K, \cE) 
\lra 
Rh_*\DR^{\d}(]\ol{X}^{(\b)}[^{\log}_{\cP^{(\b)}}/\cY, \cE). 
$$ 
We will prove that it is a quasi-isomorphism. \par 
Note that the above homorphism gives the homomorphism 
\begin{align*}
& \phantom{\lra} 
[(R^t\ti{h}_*\DR^{\d}(]\ti{X}^{(s)}[^{\log}_{\ti{\cP}^{(s)}}/\cY_K, \cE) 
\,\Rightarrow\,
R^{s+t}\ti{h}_*\DR^{\d}(]\ti{X}^{(\b)}[^{\log}_{\ti{\cP}^{(\b)}}/\cY_K, \cE)] \\ & \lra 
[(R^th_*\DR^{\d}(]\ol{X}^{(s)}[^{\log}_{\cP^{(s)}}/\cY, \cE) 
\,\Rightarrow\, 
Rh_*\DR^{\d}(]\ol{X}^{(\b)}[^{\log}_{\cP^{(\b)}}/\cY, \cE)]
\end{align*}
of spectral sequences. 
So it suffices to prove that the homomorphism 
$$ 
\alpha: R^t\ti{h}_*\DR^{\d}(]\ti{X}^{(\b)}[^{\log}_{\cP^{(\b)}}/\cY_K, \cE) 
\lra 
R^th_*\DR^{\d}(]\ol{X}^{(\b)}[^{\log}_{\cP^{(\b)}}/\cY, \cE)
$$ 
between the complex of $E_1^{\b,t}$-terms is a quasi-isomorphism. 
To prove this, it suffices to prove that the homomorphism 
$$ \Gamma(U,\alpha): 
\Gamma(U,R^t\ti{h}_*\DR^{\d}(]\ti{X}^{(\b)}[^{\log}_{\cP^{(\b)}}/\cY_K, \cE)) 
\lra 
\Gamma(U,R^th_*\DR^{\d}(]\ol{X}^{(\b)}[^{\log}_{\cP^{(\b)}}/\cY, \cE))
$$ 
is a quasi-isomorphism for any admissible open set 
$U \subseteq\, ]\ol{Y}[^{\log}_{\cY}$. Let $T_J$ be the topological space 
whose underlying set is the set of non-empty subsets of $J$ and 
whose open sets are the sets of the form $O(L):= \{L' \subseteq J \,\vert\, 
L' \supseteq L\} \,(\emptyset \not= L\subseteq J)$. Let $\cF$ be the 
presheaf on $T_J$ given by 
$$ \cF(O(L)) := 
\Gamma(U,R^th_L*\DR^{\d}(]\ol{X}_L[^{\log}_{\cP_L}/\cY_K,\cE)), $$ 
(where $h_L$ is the morphism $]\ol{X}_L[^{\log}_{\cP_L} \lra \cY_K$), 
and let us denote the open covering $T_J = \bigcup_{j \in J}O(\{j\})$ 
by $\cT$. Then 
$\Gamma(U,R^t\ti{h}_*\DR^{\d}(]\ti{X}^{(\b)}[^{\log}_{\cP^{(\b)}}/\cY_K, 
\cE))$ 
is nothing but the \v{C}ech complex $C^{\b}(\cT,\cF)$ of $\cF$ 
associated to $\cT$, 
$\Gamma(U,R^th_*\DR^{\d}(]\ol{X}^{(\b)}[^{\log}_{\cP^{(\b)}}/\cY, \cE))$ 
is nothing but the \v{C}ech complex consisting of alternating cochains 
$\ol{C}^{\b}(\cT,\cF)$ of $\cF$ and $\Gamma(U,\alpha)$ is the canonical 
homomorphism $C^{\b}(\cT,\cF) \lra \ol{C}^{\b}(\cT,\cF)$. 
It is well-known that this homomorphism is a quasi-isomorphism. 
So we are done. 
\end{pf} 

Next we prove a relation between relative log rigid cohomology and 
relative rigid cohomology in special case. 
(Part of the content here seems to be essentially 
appeared in \cite[\S 5]{shiho3}. 
But here we give it in a more systematic, simple way.) \par 
To do this, first we give a preliminary result concerning 
coherent modules on rigid analytic spaces 
and overconvergent isocrystals. Let us assume given morphisms 
$$ j_X: X \hra \ol{X}, \,\,\, i: (\ol{X},M_{\ol{X}}) \hra (\cP,M_{\cP}), $$
where $j_X$ is an open immersion of $B$-schemes and 
$i$ is a closed immersion from a fine log $B$-scheme to a $p$-adic 
fine log formal $\cB$-scheme such that $X \subseteq (\cP,M_{\cP})_{\triv}$ 
holds. Under this condition, let us first define the category 
$\wt{\Coh}(\cO_{]\ol{X}[_{\cP}})$ as follows: 
An object in this category is a pair $(U,E)$, where $U$ is a strict 
neighborhood of $]X[_{\cP}$ in $]\ol{X}[_{\cP}$ and $E$ is a coherent 
$\cO_{U}$-module. For objects $(U,E),(U',E')$, the set of morphisms 
is defined by 
$\Hom((U,E),(U',E')) := \varinjlim_V\Hom_{\cO_V}(E|_V,E'|_V)$, where 
$V$ runs through the set of strict neighborhoods of 
$]X[_{\cP}$ in $]\ol{X}[_{\cP}$ contained in $U \cap U'$. 
By \cite[(2.1.10)]{berthelot2}, we have 
the equivalence of categories 
$$ \wt{\Coh}(\cO_{]\ol{X}[_{\cP}}) \os{\sim}{\lra} 
\Coh(j_X^{\d}\cO_{]\ol{X}[_{\cP}}) $$
which is given by $(U,E) \mapsto j_U^{\d}E := 
\varinjlim_{V}\alpha_{V]\ol{X}[_{\cP},*}\alpha_{VU}^{-1}E$, where 
$V$ runs through the set of strict neighborhoods of 
$]X[_{\cP}$ in $]\ol{X}[_{\cP}$ contained in $U$ and $\alpha_{V?}$ denotes 
the admissible open immersion $V \hra ?$. 
On the other hand, let us define the category 
$\wt{\Coh}(\cO_{]\ol{X}[^{\log}_{\cP}})$ 
as follows: 
An object in this category is a pair $(U,E)$, where $U$ is a strict 
neighborhood of $]X[^{\log}_{\cP}=\,]X[_{\cP}$ in $]\ol{X}[^{\log}_{\cP}$ and $E$ is a coherent 
$\cO_{U}$-module. For objects $(U,E),(U',E')$, the set of morphisms 
is defined by 
$\Hom((U,E),(U',E')) := \varinjlim_V\Hom_{\cO_V}(E|_V,E'|_V)$, where 
$V$ runs through the set of strict neighborhoods of 
$]X[_{\cP}$ in $]\ol{X}[^{\log}_{\cP}$ contained in $U \cap U'$. 
Then, by the same argument as \cite[(2.1.10)]{berthelot2}, 
we see the equivalence of categories 
$$ \wt{\Coh}(\cO_{]\ol{X}[^{\log}_{\cP}}) \os{\sim}{\lra} 
\Coh(j_X^{\d}\cO_{]\ol{X}[^{\log}_{\cP}}). $$
Let us denote the morphism $]\ol{X}[^{\log}_{\cP} \lra 
]\ol{X}[_{\cP}$ by $\varphi_X$. 
Then we have the following: 

\begin{prop}\label{pp}
With the above notation, 
\begin{enumerate} 
\item 
The restriction functor 
$\varphi_X^*: \Coh(j_X^{\d}\cO_{]\ol{X}[_{\cP}}) \lra 
\Coh(j_X^{\d}\cO_{]\ol{X}[^{\log}_{\cP}})$ is an equivalence of 
categories. 
\item 
We have $R\varphi_{X,*}\varphi_X^*E=E$ for any $E \in 
\Coh(j^{\d}_X\cO_{]\ol{X}[_{\cP}})$. 
$($Consequently, the functor $\varphi_{X,*}$ sends 
$\Coh(j^{\d}_X\cO_{]\ol{X}[^{\log}_{\cP}})$ to 
$\Coh(j^{\d}_X\cO_{]\ol{X}[_{\cP}})$ and the resulting functor 
$$\varphi_{X,*}: \Coh(j^{\d}_X\cO_{]\ol{X}[^{\log}_{\cP}}) \lra 
\Coh(j^{\d}_X\cO_{]\ol{X}[_{\cP}})$$ 
is the exact functor which is the inverse of the functor $\varphi_X^*)$.  
\end{enumerate}
\end{prop}

\begin{pf} 
First we prove the assertion (1). To prove it, it suffices to prove 
the restriction functor 
$\varphi_X^*: \wt{\Coh}(\cO_{]\ol{X}[_{\cP}}) \lra 
\wt{\Coh}(\cO_{]\ol{X}[^{\log}_{\cP}})$ is an equivalence of 
categories. Let us take a Cartesian diagram 
\begin{equation}\label{dddd}
\begin{CD}
(\ol{X}^{(\b)},M_{\ol{X}^{(\b)}}) @>{i^{(\b)}}>> 
(\cP^{(\b)},M_{\cP^{(\b)}}) \\ 
@VVV @VgVV \\ 
(\ol{X},M_{\ol{X}}) @>i>> (\cP,M_{\cP}) 
\end{CD}
\end{equation}
such that $g$ is a strict formally etale \v{C}ech hypercovering and 
for each $n$, $i^{(n)}$ admits a factorization 
$$ (\ol{X}^{(n)},M_{\ol{X}^{(n)}}) \hra 
({\cP'}^{(n)},M_{{\cP'}^{(n)}}) \lra (\cP^{(n)},M_{\cP^{(n)}}), $$
where the first map is an exact closed immersion and the second map 
is formally log etale. Put $X^{(n)} := X \times_{\ol{X}} \ol{X}^{(n)}$, 
and let $\wt{\Coh}(\cO_{]X^{(\b)}[_{\cP^{(\b)}}}), 
\wt{\Coh}(\cO_{]X^{(\b)}[^{\log}_{\cP^{(\b)}}})$ be the category 
of descent data with respect to 
$\wt{\Coh}(\cO_{]X^{(n)}[_{\cP^{(n)}}}) \,(n=0,1,2), 
\wt{\Coh}(\cO_{]X^{(n)}[^{\log}_{\cP^{(n)}}})\,(n \allowbreak = 0,1,2)$, 
respectively. \par 
Then, by \cite[5.8]{shiho3}, the natural map 
$]\ol{X}^{(n)}[^{\log}_{\cP^{(n)}} \lra ]\ol{X}^{(n)}[_{\cP^{(n)}}$ induces 
an isomorphism between some strict neighborhood of $]X^{(n)}[_{\cP^{(n)}}$ in 
$]\ol{X}^{(n)}[^{\log}_{\cP^{(n)}}$ and some strict neighborhood 
of $]X^{(n)}[_{\cP^{(n)}}$ in $]\ol{X}^{(n)}[_{\cP^{(n)}}$ (for each $n$). 
Therefore, the restriction functor 
$\wt{\Coh}(\cO_{]X^{(\b)}[_{\cP^{(\b)}}}) \lra 
\wt{\Coh}(\cO_{]X^{(\b)}[^{\log}_{\cP^{(\b)}}})$ is an equivalence of 
categories. 
On the other hand, 
for any family of strict neighborhoods $\{U^{(n)}\}_{n=0,1,2}$ of 
$]X^{(n)}[_{\cP}$ in $]\ol{X}^{(n)}[_{\cP}$ 
(resp. $]\ol{X}^{(n)}[^{\log}_{\cP}$), there exists a strict 
neighborhood $U$ of $]X[_{\cP}$ in $]\ol{X}[_{\cP}$ 
(resp. $]\ol{X}[^{\log}_{\cP}$) such that the pull-back of $U$ in 
$]\ol{X}^{(n)}[_{\cP^{(n)}}$ (resp. $]\ol{X}^{(n)}[^{\log}_{\cP^{(n)}}$)
 is contained in $U^{(n)}$ (for $n=0,1,2$). 
 (This follows from the fact that the diagram \eqref{dddd} is Cartesian 
and the definition of $X^{(\b)}$.) So we can shrink $\{U^{(n)}\}_{n=0,1,2}$ 
so that they are the pull-back of $U$. Then the morphism $U^{(\b)} \lra U$ 
satisfies the following property: For any affinoid admissible open set 
$V =\Spm (\Q \otimes_{\Z} A)\subseteq U$, the 
morphism $V^{(\b)}:=V \times_U U^{(\b)} \lra V$ is induced from a 
etale \v{C}ech covering of $\Spf A$. So, by rigid analytic faifully 
flat descent, the restriction functor 
$\Coh(\cO_U) \lra \Coh(\cO_{U^{(\b)}})$ is an equivalence of categories. 
This implies that the restriction functors 
$$ 
\wt{\Coh}(\cO_{]\ol{X}[_{\cP}}) \lra 
\wt{\Coh}(\cO_{]X^{(\b)}[_{\cP^{(\b)}}}), \,\,\, 
\wt{\Coh}(\cO_{]\ol{X}[^{\log}_{\cP}}) \lra 
\wt{\Coh}(\cO_{]X^{(\b)}[^{\log}_{\cP^{(\b)}}})
$$ 
are equivalence of categories. Since $\varphi_X^*$ can be factorized as 
$$ 
\wt{\Coh}(\cO_{]\ol{X}[_{\cP}}) \os{\simeq}{\lra} 
\wt{\Coh}(\cO_{]X^{(\b)}[_{\cP^{(\b)}}})
\os{\simeq}{\lra} 
\wt{\Coh}(\cO_{]X^{(\b)}[^{\log}_{\cP^{(\b)}}}) 
\os{\simeq}{\lla}
\wt{\Coh}(\cO_{]\ol{X}[^{\log}_{\cP}}) 
$$
(where each arrows are defined by restriction), 
we can conclude that it is an equivalence of categories. \par 
Next we prove the assertion (2). Let us denote the morphisms 
$]\ol{X}^{(\b)}[_{\cP^{(\b)}} \lra ]\ol{X}[_{\cP}, 
\allowbreak ]\ol{X}^{(\b)}[^{\log}_{\cP^{(\b)}} \allowbreak 
\lra ]\ol{X}[^{\log}_{\cP}, 
]\ol{X}^{(\b)}[^{\log}_{\cP^{(\b)}} \lra ]\ol{X}^{(\b)}[_{\cP^{(\b)}}$ by 
$\pi, \pi^{\log}, \varphi_{X^{(\b)}}$ respectively. First we claim that, 
for a coherent module $E$ on $]\ol{X}[_{\cP}, ]\ol{X}[^{\log}_{\cP}, 
]\ol{X}^{(\b)}[_{\cP^{(\b)}}$ respectively, we have 
$R\pi_*\pi^*E = E, R\pi^{\log}_*\pi^{\log,*}E =E, 
R\varphi_{X^{(\b)}.*}\varphi_{X^{(\b)}}^*E = E$ respectively: 
In the case of $\pi$, it follows from \cite[7.3.1]{chts}. 
In the case of $\varphi_{X^{(\b)}}$, it follows from \cite[8.3.5]{chts}, 
because $]\ol{X}^{(n)}[^{\log}_{\cP^{(n)}} = ]\ol{X}^{(n)}[_{{\cP'}^{(n)}}$ 
holds and ${\cP'}^{(n)} \lra \cP^{(n)}$ is isomorphic on a neighborhood of 
$X^{(n)}$. In the case of $\pi^{\log}$, let us denote the system 
of universal enlargements of $((\cP,M_{\cP}), (\ol{X},M_{\ol{X}}))$ 
by $\{(T_{a,b}(\cP),M_{T_{a,b}(\cP)}), 
(\ol{X}_{a,b},M_{\ol{X}_{a,b}}))\}_{(a,b) \in \A}$. Then we have 
$]\ol{X}[^{\log}_{\cP} = \bigcup T_{a,b}(\cP)_K = 
\bigcup ]\ol{X}_{a,b}[_{T_{a,b}(\cP)}$. So it suffices to check the 
claim after pulling back the diagram \eqref{dddd} by 
$(T_{a,b}(\cP),M_{T_{a,b}(\cP)}) \lra (\cP,M_{\cP})$, and in this case, 
the claim follows from \cite[7.3.1]{chts} because 
we do not need to consider the log structure (since 
$(\ol{X}_{a,b},M_{\ol{X}_{a,b}}) 
\hra (T_{a,b}(\cP),M_{T_{a,b}(\cP)})$ is now an 
exact closed immersion). \par 
By the claim in the previous paragraph, 
we have the following equality for a 
coherent $j^{\d}\cO_{]\ol{X}[_{\cP}}$-module $E$: 
$$ 
R\varphi_{X,*}\varphi_X^*E  = 
R\varphi_{X,*}R\pi^{\log}_*  \pi^{\log,*}\varphi_X^*E = 
R\pi_*R\varphi_{X^{(\b)},*}\varphi_{X^{(\b)}}^*\pi^*E = E. 
$$ 
So we are done. 
\end{pf} 

Next let us assume given a diagram \eqref{logrig-diag} and open 
immersions 
$j_X: X \hra \ol{X}, j_Y:Y \hra \ol{Y}$ such that $f(X) \subseteq Y$, 
$X \subseteq (\ol{X},M_{\ol{X}})_{\triv}$ and
$Y \subseteq (\cY,M_{\cY})_{\triv}$ 
holds. Then there exists an embedding system 
$$ (\ol{X},M_{\ol{X}}) \lla (\ol{X}^{(\b)},M_{\ol{X}^{(\b)}}) \hra 
(\cP^{(\b)},M_{\cP^{(\b)}}) $$ 
such that, if we put $X^{(\b)}:= X \times_{\ol{X}} \ol{X}^{(\b)}$, 
$X^{(\b)} \subseteq (\cP^{(\b)},M_{\cP^{(\b)}})_{\triv}$ holds. 
(Hence $\cP^{(n)}$ is formally smooth over $\cY$ on a neighborhood of 
$X^{(n)}$.) Then we have the category of overconvergent 
isocrystals $I^{\d}((X,\ol{X})/\cY_K) = 
I^{\d}((X^{(\b)},\ol{X}^{(\b)})/\cY_K, \cP^{(\b)})$ 
on usual pairs defined in \cite{chts} and 
the category of overconvergent isocrystals 
$I^{\d}(((X,\ol{X})/\cY_K)^{\log}) = 
I^{\d}(((X^{(\b)},\ol{X}^{(\b)})/\cY_K, \cP^{(\b)})^{\log})$ on log pairs 
defined in this section. Let us denote the morphism 
$]\ol{X}^{(\b)}[^{\log}_{\cP^{(\b)}} \lra ]\ol{X}^{(\b)}[_{\cP^{(\b)}}, 
]\ol{Y}[^{\log}_{\cY} \lra ]\ol{Y}[_{\cY}, 
]\ol{X}^{(\b)}[_{\cP^{(\b)}} \lra ]\ol{Y}[_{\cY}, 
]\ol{X}^{(\b)}[^{\log}_{\cP^{(\b)}} \lra ]\ol{Y}[^{\log}_{\cY}$ by 
$\varphi_X, \varphi_Y, h, h^{\log}$, respectively. 
Then we have the following: 

\begin{prop}\label{pppp}
Under the above assumption, 
\begin{enumerate} 
\item 
The restriction functor 
{\tiny{ 
$$ I^{\d}((X,\ol{X})/\cY_K) = 
I^{\d}((X^{(\b)},\ol{X}^{(\b)})/\cY_K, \cP^{(\b)}) \os{\varphi_X^*}{\lra} 
I^{\d}(((X^{(\b)},\ol{X}^{(\b)})/\cY_K, \cP^{(\b)})^{\log})
= I^{\d}(((X,\ol{X})/\cY_K)^{\log})
$$ }}
is an equivalence of categories. 
\item 
For $\cE \in I^{\d}((X,\ol{X})/\cY_K)=I^{\d}(((X,\ol{X})/\cY_K)^{\log})$, 
we have the quasi-isomorph\-isms 
$$ 
\varphi_X^*\DR^{\d}(]\ol{X}^{(\b)}[_{\cP^{(\b)}}/\cY_K,\cE) = 
\DR^{\d}(]\ol{X}^{(\b)}[^{\log}_{\cP^{(\b)}}/\cY_K,\cE), $$ 
$$ 
\DR^{\d}(]\ol{X}^{(\b)}[_{\cP^{(\b)}}/\cY_K,\cE) = 
\varphi_{X,*}\DR^{\d}(]\ol{X}^{(\b)}[^{\log}_{\cP^{(\b)}}/\cY_K,\cE) = 
R\varphi_{X,*}\DR^{\d}(]\ol{X}^{(\b)}[^{\log}_{\cP^{(\b)}}/\cY_K,\cE). 
$$ 
\end{enumerate}
\end{prop} 

\begin{pf} 
The assertion (1) follows easily from Proposition \ref{pp} (1) and the 
assertion (2) follows from Proposition \ref{pp} (2). 
\end{pf} 

As a corollary of Proposition \ref{pppp}, we have the quasi-isomorphism 
\begin{align*}
Rf_{(X,\ol{X})/\cY, \rig *}\cE & = 
Rh_*\DR^{\d}(]\ol{X}^{(\b)}[_{\cP^{(\b)}}/\cY_K,\cE) \\
& = R\varphi_{Y,*}Rh^{\log}_*
\DR^{\d}(]\ol{X}^{(\b)}[^{\log}_{\cP^{(\b)}}/\cY_K, \cE) \\ 
& = R\varphi_{Y,*}Rf_{(X,\ol{X})/\cY_K, \logrig *}\cE. 
\end{align*}

So we have the following theorem, which we will use in the following section: 

\begin{thm}\label{logrig-rig}
Let the notations be as above. Then, if the relative log rigid 
cohomologies $R^qf_{(X,\ol{X})/\cY, \logrig *} \cE$ are 
coherent $j^{\d}_Y\cO_{]\ol{Y}[^{\log}_{\cY}}$-modules for any 
$q\geq 0$, the relative rigid cohomologies 
$R^qf_{(X,\ol{X})/\cY, \rig *} \cE$ are 
coherent $j^{\d}_Y\cO_{]\ol{Y}[_{\cY}}$-modules for any 
$q\geq 0$ and we have the isomorphism 
$R^qf_{(X,\ol{X})/\cY, \rig *} \cE = 
\varphi_{Y,*}R^qf_{(X,\ol{X})/\cY, \logrig *} \cE$. 
\end{thm} 

Next we introduce the notion of $\lam$-restriction of relative log rigid 
cohomology and compare it to relative log rigid cohomology. 
Let us fix $\lam \in (0,1)$, 
let us assume given the diagram \eqref{logrig-diag} and 
log pairs $((X,M_X),(\ol{X},M_{\ol{X}})), ((Y,M_Y),(\ol{Y},M_{\ol{Y}}))$ 
satisfying $f(X) \subseteq Y$, and let us assume given 
an overconvegent isocrystal $\cE$ on 
$((X,M_X),(\ol{X},M_{\ol{X}}))/\cY$. First, assume for the moment that 
we have the diagram \eqref{logrig-emb}. Then we have the log de Rham complex 
$\DR^{\d}(]\ol{X}[^{\log}_{\cP}/\cY_K, \cE)$. Then we can see the 
following in the same way as \cite[(2.2.3)]{berthelot2}: For any 
sufficiently small strict neighborhood $V$ of $]X[^{\log}_{\cP}$ in 
$]\ol{X}[^{\log}_{\cP}$, there exists a log de Rham complex 
$\DR(V/\cY_K,\cE)$ consisting of coherent $\cO_V$-modules which is 
compatible with respect to $V$ and satisfies 
$$\DR^{\d}(]\ol{X}[^{\log}_{\cP}/\cY_K,\cE) 
= \varinjlim_V\alpha_{V,*}\DR(V/\cY_K,\cE),$$ 
where $\alpha_V$ is the admissible open immersion 
$V \hra ]\ol{X}[^{\log}_{\cP}$. Let us put $V_{\lam}:=V \cap \,
]\ol{X}[^{\log}_{\cP,\lam}$, 
let us denote the composite $V_{\lam} \hra V \os{\alpha_V}{\hra} 
]\ol{X}[^{\log}_{\cP}$ by $j_{V,\lam}$ and we define the 
complex $\DR^{\d}_{\lam}(]\ol{X}[^{\log}_{\cP}/\cY_K, \cE)$ by 
$$ \DR^{\d}_{\lam}(]\ol{X}[^{\log}_{\cP}/\cY_K, \cE) := 
\varinjlim_{V} j_{V,\lam *}(\DR(V/\cY_K,\cE)|_{V_{\lam}}). $$ 
Note that, when we have the diagram \eqref{logrig-diag-frak}
(satisfying $(*)$) instead of \eqref{logrig-emb}, 
we can define the complex 
$\DR^{\d}_{\lam}(\fP^{\ex}_K/\cY_K,\cE)$ in the same way. \par 
Now let us consider the situation where there does not necessarily exist 
the diagram \eqref{logrig-emb}. Let us take an embedding system 
\eqref{logrig-emb2} and let $h$ be the morphism 
$]\ol{X}^{(\b)}[^{\log}_{\cP^{(\b)}} \lra ]\ol{Y}[^{\log}_{\cY}$. 
Then we have the complex 
$\DR^{\d}_{\lam}(]\ol{X}^{(\b)}[^{\log}_{\cP^{(\b)}}/\cY_K,\cE)$ on 
$]\ol{X}^{(\b)}[^{\log}_{\cP^{(\b)}}$. 
We define $Rf_{(X,\ol{X})/\cY, \logrig, \lam *}\cE, 
R^qf_{(X,\ol{X})/\cY, \logrig, \lam *}\cE$ by 
$$ 
Rf_{(X,\ol{X})/\cY, \logrig, \lam *}\cE := 
Rh_*\DR^{\d}_{\lam}(]\ol{X}^{(\b)}[^{\log}_{\cP^{(\b)}}/\cY_K,\cE), $$
$$
R^qf_{(X,\ol{X})/\cY, \logrig, \lam *}\cE := 
R^qh_*\DR^{\d}_{\lam}(]\ol{X}^{(\b)}[^{\log}_{\cP^{(\b)}}/\cY_K,\cE), 
$$ 
and we call $R^qf_{(X,\ol{X})/\cY, \logrig, \lam *}\cE$ the 
$\lam$-restriction of the $q$-th relative log rigid cohomology. It is 
a $j_Y^{\d}\cO_{]\ol{Y}[^{\log}_{\cY}}$-module. Then we have the 
following: 

\begin{prop}\label{logrig-lam-wd}
The above definition is independent of the choice of 
the embedding system. 
\end{prop}

Before the proof, we prove the following lemma: 

\begin{lem}\label{acyc1}
Let $((X,M_X),(\ol{X},M_{\ol{X}}),(\cP,M_{\cP}))$ be a log triple. 
Fix $\lam \in (0,1)$ and for a strict neighborhood $V$ of 
$]X[^{\log}_{\cP}$ in $]\ol{X}[^{\log}_{\cP}$, let us put 
$V_{\lam} := V \,\cap \,]\ol{X}[^{\log}_{\cP,\lam}$ and let us denote the 
admissible open immersion $V_{\lam} \hra ]\ol{X}[^{\log}_{\cP}$ by 
$j_{V,\lam}$. Then there exists an oriented set 
$\cV$ consisting of 
strict neighborhoods of $]X[^{\log}_{\cP}$ in 
$]\ol{X}[^{\log}_{\cP}$ which is cofinal in the set of 
all the strict neighborhoods such that the functor 
$j_{V,\lam *}$ is acyclic for coherent $\cO_V$-modules for any $V \in \cV$. 
\end{lem} 

\begin{pf} 
Since we may work Zariski locally on $\cP$, we may assume 
it is affine. Let us take a generator 
$g_1,g_2,\cdots, g_c$ of 
$\Ker(\Gamma(\cP,\cO_{\cP}) \ra \Gamma(\ol{X}-X, 
\cO_{\ol{X}-X}))$ and for $\nu \in p^{\Q_{<0}}$, 
let $U_{\nu,i} \,(1 \leq i \leq c), U_{\nu,0}, 
U_{\nu}$ be the admissible open set 
$\{|g_i| \geq \nu, |g_i| \geq |g_j|\,(\forall j\not=i)\}, 
\{|g_j| \leq \nu\,(\forall i)\}, \bigcup_{i=1}^cU_{\nu,i}$ of 
$]\ol{X}[^{\log}_{\cP}$, respectively. Then, 
any strict neighborhood $\ti{V}$ contains a strict neighborhood 
$V$ satisfying $V_{\lam} = 
U_{\nu} \,\cap\, ]\ol{X}[^{\log}_{\cP,\lam}$ for some $\nu$. 
(Here the assumption $\lam<1$ is crucial.) So, if we define $\cV$ to be 
the category of strict neighborhoods $V$ satisfying this condition, 
$\cV$ is cofinal in the 
set of all the strict neighborhoods. \par
Now we prove that, for $V \in \cV$, the functor $j_{V,\lam *}$ 
is acyclic for any coherent $\cO_V$-module $E$. 
Take $\nu_0 \in p^{\Q_{<0}}$ satisfying $V_{\lam} = 
U_{\nu_0} \,\cap\, ]\ol{X}[^{\log}_{\cP,\lam}$ and take $\nu \in (0,\nu_0)$. 
Let us take any affinoid admissible open set $W$ of 
$]\ol{X}[^{\log}_{\cP}$. If we put 
$W_i := W \cap U_{\nu,i}\,(0 \leq i \leq c)$, then $W = \bigcup_{i=0}^cW_i$ 
is an admissible covering of $W$ by affinoid admissible open sets. 
Then we have $j_{V,\lam}^{-1}(W_i) = 
W_i \cap U_{\nu_0} \cap ]\ol{X}[^{\log}_{\cP,\lam}=\emptyset$ for $i=0$ and 
$=W \cap U_{\nu_0,i} \cap ]\ol{X}[^{\log}_{\cP,\lam}$ for $i\not= 0$. 
If we take a good sequence $\a$ in $\A_r$ (where $r:=-\log_p \lam$), 
the right hand side for $i\not= 0$ 
is the union of the affinoid admissible open sets 
$W \cap U_{\nu_0,i} \cap T_{a,b}(\cP)_K\,((a,b) \in \a)$. So 
it is quasi-Stein. Hence we have 
$H^q(j_{V,\lam}^{-1}(W_i),E) = 0$ for any $i$ and $q>0$. 
This implies that the 
sheaf associated to the presheaf $W \mapsto H^q(j_{V,\lam}^{-1}(W),E)$ 
is zero for $q>0$, that is, $R^qj_{V,\lam *}E = 0$ for $q>0$. So we are 
done. 
\end{pf}

\begin{pf*}{Proof of Propsition \ref{logrig-lam-wd}} 
The proof is similar to that of Proposition \ref{logrig-wd}. First, we 
are reduced to proving the following claims: \\
\quad \\
{\bf claim 1.} \,\,\, Assume we are given the Cartesian diagram 
\begin{equation}
\begin{CD}
(\ol{X}^{(\b)},M_{\ol{X}^{(\b)}}) @>>> (\cP^{(\b)},M_{\cP^{(\b)}}) \\ 
@VVV @VgVV \\ 
(\ol{X},M_{\ol{X}}) @>{\iota}>> (\cP,M_{\cP})
\end{CD}
\end{equation}
and a log pair $((X,M_X),(\ol{X},M_{\ol{X}}))$, 
where 
$(\cP,M_{\cP})$ is an object in 
$\LF$, $\iota$ is a closed immersion and $g$ is a strict etale hypercovering. 
Let $V_0$ be a strict neighborhood of $]X[^{\log}_{\cP}$ in 
$]\ol{X}[^{\log}_{\cP}$ and let $E$ be a coherent module on $V_0$. 
For a strict neighborhood $V$ contained in $V_0$, put 
$V_{\lam}:=V \,\cap\, ]\ol{X}[^{\log}_{\cP}$, 
$E_{V,\lam} := E|_{V_{\lam}}$ and 
denote the admissible open immersion $V_{\lam} \hra ]\ol{X}[^{\log}_{\cP}$ 
by $j_{V,\lam}$. 
On the other hand, let $V_0^{(\b)}:=g^{-1}(V_0)$ and for 
a strict neighborhood $V$ contained in $V_0$, put 
$V^{(\b)}_{\lam}:=g^{-1}(V_{\lam}), 
E_{V^{(\b)},\lam} := E|_{V^{(\b)}_{\lam}}$ and denote the 
the admissible open immersion $V^{(\b)}_{\lam} \hra 
]\ol{X}^{(\b)}[^{\log}_{\cP^{(\b)}}$ by $j_{V^{(\b)},\lam}$. Then 
we have the quasi-isomorphism 
$$ \varinjlim_V j_{V,\lam *} E_{V,\lam} \os{=}{\lra}
 Rg_{K,*}(\varinjlim_{V} j_{V^{(\b)},\lam *}E_{V^{(\b)},\lam}). $$

\quad \\
{\bf claim 2.} \,\,\, Assume we are given a diagram over $(\cY,M_{\cY})$ 
\begin{equation}
\begin{CD}
(\ol{X},M_{\ol{X}}) @>{\iota_1}>> (\cQ,M_{\cQ}) \\ 
@\vert @VgVV \\ 
(\ol{X},M_{\ol{X}}) @>{\iota_2}>> (\cP,M_{\cP}) 
\end{CD}
\end{equation}
and a log pair $((X,M_X),(\ol{X},M_{\ol{X}}))$, 
where $(\cP,M_{\cP}),(\cQ,M_{\cQ})$ are objects in 
$\pLF$, $\iota_1,\iota_2$ are closed immersions and $g$ is a formally log 
smooth 
morphism. Let $g_K:]\ol{X}[^{\log}_{\cQ} \lra ]\ol{X}[^{\log}_{\cP}$ be the 
morphism induced by $g$. 
Then, for an overconvergent isocrystal $\cE$ on 
$((X,M_X),(\ol{X},M_{\ol{X}}))/(\cY,M_{\cY})$, 
 we have a quasi-isomorphism 
$$ \DR^{\d}_{\lam}(]\ol{X}[^{\log}_{\cP}/\cY, \cE) 
\os{=}{\lra} Rg_{K,*}\DR^{\d}_{\lam}(]\ol{X}[^{\log}_{\cQ},\cE).$$ 

\quad \\
First we prove claim 1. Because of the quasi-compactness, 
quasi-separetedness of $g_K$ (this can be checked after reducing to the 
case where $\iota$ is a homeomorphic exact closed immersion and it is 
clear in this case) and Lemma \ref{acyc1}, we have 
\begin{align*}
Rg_{K,*}(\varinjlim_{V} j_{V^{(\b)},\lam *}E_{V^{(\b)},\lam})
& = 
\varinjlim_{V} R(g_K \circ j_{V^{(\b)},\lam})_* 
E_{V^{(\b)},\lam} \\ 
& = 
\varinjlim_V Rj_{V, \lam *} R(g_K|_{V^{(\b)}_{\lam}})_* E_{V^{(\b)},\lam}. 
\end{align*}
So, if we have 
\begin{equation}\label{whitestar}
R(g_K|_{V^{(\b)}_{\lam}})_* E_{V^{(\b)},\lam} = E_{V,\lam}, 
\end{equation}
we have 
$$ 
\varinjlim_V Rj_{V, \lam *} R(g_K|_{V^{(\b)}_{\lam}})_* E_{V^{(\b)},\lam} 
= 
\varinjlim_V Rj_{V, \lam *} E_{V,\lam} = 
\varinjlim_V j_{V, \lam *} E_{V,\lam}
$$ 
and so we are done. Now note that the morphism 
$g_K|_{V^{(\b)}_{\lam}}: V^{(\b)}_{\lam} \lra V_{\lam}$ is obtained as 
the pull-back of the morphism $\cP^{(\b)}_K \lra \cP_K$ by 
$V_{\lam} \hra \cP^{\ex}_K \lra \cP_K$. So, for any affinoid admissible 
open set $U := \Spm (\Q \otimes_{\Z} A) \subseteq V_{\lam}$, 
the morphism $U \times_{V_{\lam}} V^{(\b)}_{\lam} \lra U$ is equal to the 
morphism 
$(\cP^{(\b)} \times_{\cP} \Spf A)_K \lra (\Spf A)_K$. So the 
quasi-isomorphism \eqref{whitestar} follows from \cite[3.9]{shiho3}. 
So we have proved claim 1. \par 
Next we prove claim 2. It is similar to the proof of claim 2 in 
the proof of Proposition \ref{logrig-wd}. 
Let $(\cQ^{\ex},M_{\cQ^{\ex}}), 
(\cP^{\ex},M_{\cP^{\ex}})$ be the exactification of $\iota_1, \iota_2$, 
respectively. Let $E$ be a coherent sheaf on a 
 strict neighborhood $V_0$ of $]X[^{\log}_{\cP}$ in 
 $]\ol{X}[^{\log}_{\cP}=\cP^{\ex}_K$ and for a strict neighborhood $V$ 
 contained in $V_0$, put $V_{\lam} := V \,\cap\, ]\ol{X}[^{\log}_{\cP,\lam}, 
 E_{V,\lam} := E|_{V_{\lam}}$ and denote the admissible open immersion 
$V_{\lam} \hra \cP^{\ex}_K$ by $j_{V,\lam}$. On the other hand, 
for a strict neighborhood $V'$ of $]\ol{X}[^{\log}_{\cQ}$ 
contained in $g_K^{-1}(V_0)$, 
put $V'_{\lam} := V' \,\cap\, ]\ol{X}[^{\log}_{\cQ,\lam}, 
(g_K^*E)_{V',\lam} := ((g_K|_{g_K^{-1}(V_0)})^*E)|_{V'_{\lam}}$ 
and denote the admissible open immersion 
$V'_{\lam} \hra \cP^{\ex}_K$ by $j'_{V,\lam}$. Then, 
by the same argument as 
the proof of Proposition \ref{logrig-wd}, we may assume that 
$\cP^{\ex} = \Spf A$ is affine, 
$\cQ^{\ex} = \Spf A[[t_1]]$ and it suffices to prove the quasi-isomorphism 
\begin{equation}\label{logrig-lam-qis}
\varinjlim_{V}j_{V,\lam *}E \os{=}{\lra} 
Rg_{K,*}(\varinjlim_{V'} j_{V',\lam *}(g_K^*E)_{V',\lam} \otimes \Omega^{\b}) 
\end{equation}
for any $E$ as above, where $\Omega^{\b}$ is the complex 
$[\cO_{\cQ^{\ex}} \ra \cO_{\cQ^{\ex}}dt_1]$. 
Let 
$\alpha_1, \cdots \alpha_c \in 
\Gamma(\cP^{\ex},\cO_{\cP^{\ex}})$ 
be a lift of generators of $\Ker(\Gamma(\ol{X},\cO_{\ol{X}}) \ra 
\Gamma(\ol{X}-X,\cO_{\ol{X}-X}))$ and 
for $\nu \in p^{\Q_{<0}}$, let 
$U_{\nu}$ be the admissible open subset $\bigcup_{i=1}^c\{|\alpha_i| 
\geq \nu\}$ 
in $\cP^{\ex}_K$. 
To prove the quasi-isomorphism \eqref{logrig-lam-qis}, 
it suffices to check on an affinoid admissible open subset $W$ 
of $\cP^{\ex}_K$ such that $W_{\nu} := W \cap U_{\nu}$ is also affinoid 
for any $\nu$ sufficiently close to $1$, by Lemma \ref{acyc0}. 
For 
$\lam'=p^{-b/a}, \nu \in p^{\Q_{<0}}$, let 
$W_{\lam',\nu}$ be the admissible open set 
$T_{a,b}(\cP^{\ex})_K \cap W_{\nu}$ of $W$, 
let $g_K^{-1}(W)_{\lam'}$ be the affinoid admissible open set 
$\{|t_1|\leq \lam'\} \cap g_K^{-1}(W_{\lam'})$ of $g_K^{-1}(W)$ 
 and 
let $g_K^{-1}(W)_{\lam',\nu}$ be the affinoid admissible open set 
$\{|t_1|\leq \lam'\} \cap g_K^{-1}(W_{\lam',\nu})$ of $g_K^{-1}(W)$. 
Then, when $\lam'$ is fixed, we have $W_{\lam',\nu} \subseteq V$, 
$g_K^{-1}(W)_{\lam',\nu} \subseteq g_K^{-1}(V)$ for $\nu$ sufficiently 
close to $1$. In this case, we denote the restriction of $E,g_K^*E$ to 
$W_{\lam',\nu}, g_K^{-1}(W)_{\lam',\nu}$ by 
$E_{\lam',\nu}, (g_K^*E)_{\lam',\nu}$, respectively. 
With this notation, we have 
\begin{equation*}
\Gamma(W,\varinjlim_V j_{V,\lam *}E_{V,\lam})= 
\varinjlim_V \Gamma(W \cap V_{\lam}, E_{V,\lam})
= 
\dvil_{\nu\to 1^-} \vpl_{\lam'\to\lam^-} 
\Gamma(W_{\lam',\nu}, E_{\lam',\nu}), 
\end{equation*}
\begin{align*}
& \phantom{=} 
R\Gamma(g_K^{-1}(W),\varinjlim_{V'}j_{V',\lam *}
(g_K^{*}E)_{V',\lam} \otimes \Omega^{\b}) \\ & = 
R\dvpl_{\mu\to1^-} 
R\Gamma(g_K^{-1}(W)_{\mu},\varinjlim_{V'}Rj'_{V',\lam *}
(g_K^*E)_{V',\lam} \otimes \Omega^{\b}) \\ 
& = 
R\dvpl_{\mu\to 1^-} 
\varinjlim_{V'}
R\Gamma(g_K^{-1}(W)_{\mu} \cap V'_{\lam}, 
(g_K^*E)_{V',\lam} \otimes \Omega^{\b}) \\ 
& = 
\varinjlim_{V'}
R\Gamma(g_K^{-1}(W)_{\lam} \cap V'_{\lam}, 
(g_K^*E)_{V',\lam} \otimes \Omega^{\b}) \\ 
& = 
\dvil_{\nu\to 1^-} 
R\Gamma(g_K^{-1}(W)_{\lam,\nu} \,\cap\, ]\ol{X}[^{\log}_{\cQ,\lam}, 
(g_K^*E)_{\lam,\nu} \otimes \Omega^{\b}) \\ 
& = 
\dvil_{\nu\to 1^-} 
\Gamma(g_K^{-1}(W)_{\lam,\nu} \cap ]\ol{X}[^{\log}_{\cQ,\lam}, 
(g_K^*E)_{\lam,\nu} \otimes \Omega^{\b}) \\ 
& = 
\varinjlim_{\nu\to 1^-} \varprojlim_{\lam'\to\lam^-}
\Gamma(g_K^{-1}(W)_{\lam',\nu}, (g_K^*E)_{\lam',\nu} \otimes \Omega^{\b}). 
\end{align*}
So, to prove the quasi-isomorphism \eqref{logrig-lam-qis}, it suffices 
to prove that the complex 
{\allowdisplaybreaks{
\begin{align*}
0 \ra 
\dvpl_{\lam'\to\lam^-} 
\Gamma(W_{\lam',\nu}, E_{\lam',\nu}) & \ra 
\dvpl_{\lam'\to\lam^-}
\Gamma(g_K^{-1}(W)_{\lam',\nu}, (g_K^*E)_{\lam',\nu}) \\ 
& \ra \dvpl_{\lam'\to\lam^-}
\Gamma(g_K^{-1}(W)_{\lam',\nu}, (g_K^*E)_{\lam',\nu} \otimes \Omega^1)
 \ra 0 
\end{align*}}}%
is exact. We can prove this assertion in the same way as the proof of 
claim 3 in the proof of Proposition \ref{logrig-wd}. So the proof of 
the proposition is finished. 
\end{pf*} 

Let us go back to the situation where we are given the 
diagram \eqref{logrig-diag} and 
log pairs $((X,M_X),(\ol{X},M_{\ol{X}})), ((Y,M_Y),(\ol{Y},M_{\ol{Y}}))$ 
satisfying $f(X) \subseteq Y$. Let $\cE$ be an overconvergent isocrystal 
on $((X,M_X),(\ol{X},M_{\ol{X}}))/\cY_K$, take the embedding system 
\eqref{logrig-emb2} and put $(X^{(\b)},M_{X^{(\b)}}) := 
(X,M_X) \times_{(\ol{X},M_{\ol{X}})} 
(\ol{X}^{(\b)},M_{\ol{X}^{(\b)}})$. 
Then, for $\lam \in (0,1)$, we have the canonical map of complexes 
\begin{equation}\label{lam-map0}
\DR^{\d}(]\ol{X}^{(\b)}[^{\log}_{\cP^{(\b)}}/\cY_K, \cE) \lra 
\DR^{\d}_{\lam}(]\ol{X}^{(\b)}[^{\log}_{\cP^{(\b)}}/\cY_K, \cE). 
\end{equation}
So we obtain the map 
\begin{equation}\label{lam-map1}
Rf_{(X,\ol{X})/\cY, \logrig *}\cE \lra 
Rf_{(X,\ol{X})/\cY, \logrig, \lam *}\cE 
\end{equation}  
from relative log rigid cohomology to its $\lam$-restriction. 
Then we have the following theorem: 

\begin{thm}\label{rig-lamrig}
Let the notation be as above and put $r := -\log_p\lam$. 
Then, if the morphism $f:(\ol{X},M_{\ol{X}}) \lra (\ol{Y},M_{\ol{Y}})$ is 
log smooth, $\iota$ is a homeomorphic exact closed immersion satisfying 
$r(\cO_{\cY},\Ker(\cO_{\cY} \ra \cO_{\ol{Y}})) > r$, 
the map \eqref{lam-map1} is a quasi-isomorphism. 
\end{thm} 

\begin{pf} 
It suffices to prove that the map \eqref{lam-map0} is a quasi-isomorphism, 
and to prove it, we may replace $\b$ by $n \in \N$. In this case, 
we may assume that $\ol{X}^{(n)}$ is affine. Moreover, 
the both hand sides are unchanged if we change 
the closed immersion $(\ol{X}^{(n)},M_{\ol{X}^{(n)}}) 
\hra (\cP^{(n)}, M_{\cP^{(n)}})$ by another closed immersion 
$(\ol{X}^{(n)},M_{\ol{X}^{(n)}}) 
\hra ({\cP'}^{(n)}, M_{{\cP'}^{(n)}})$ such that 
$({\cP'}^{(n)}, M_{{\cP'}^{(n)}})$ is formally log smooth over 
$(\cY,M_{\cY})$. So we may assume that the closed immersion 
$(\ol{X}^{(n)},M_{\ol{X}^{(n)}}) \hra (\cP^{(n)},M_{\cP^{(n)}})$ 
is exact and satisfies the isomorphism 
$(\cP^{(n)},M_{\cP^{(n)}}) \times_{(\cP,M_{\cP})} 
(\ol{X},M_{\ol{X}}) = (\ol{X}^{(n)},M_{\ol{X}^{(n)}})$. 
In this case, we have 
$r(\cO_{\cP^{(n)}}, \Ker(\cO_{\cP^{(n)}}\ra \cO_{\ol{X}^{(n)}})) 
\geq r(\cO_{\cY},\Ker(\cO_{\cY} \ra \cO_{\ol{Y}})) > r$. 
So we have $]\ol{X}^{(n)}[^{\log}_{\cP^{(n)},\lam} = 
]\ol{X}^{(n)}[^{\log}_{\cP^{(n)}}$. Then it is easy to see that the 
both hand sides of \eqref{lam-map0} ($\b$ replaced by $n$) coincide. 
So the theorem is proved. 
\end{pf} 

\begin{rem} 
In the case where all the log structures are trivial and 
$Y=\ol{Y}=\Spec k, \ol{Y}=\Spf V$ (where $V$ is a complete discrete 
valuation ring of mixed characteristic 
with residue field $k$), Crew defined in \cite{crew} 
the notion of 
rigid cohomology of radius $\lam$. There is a slight difference 
between his definition and our definition of $\lam$-restriction of 
rigid cohomology which is given here. (This is why we did not call 
our cohomology `the relative log rigid cohomology of radius $\lam$'.) 
However, these cohomologies are related: By Theorem \ref{rig-lamrig} and 
the analogue of it for Crew's cohomology, we can deduce that 
the both cohomology coincide for $\lam$ sufficiently close to $1$ 
(under the assumption $Y=\ol{Y}=\Spec k, \ol{Y}=\Spf V$) 
if $f:\ol{X} \lra \ol{Y}=\Spec k$ is smooth. 
\end{rem} 


\section{Relative rigid cohomology (I): log smooth case}

In this section, we give an affirmative answer to 
Conjecture \ref{bconj2} in the case 
where the given situation admits a `nice' log structure, by using 
the results in previous sections. In \cite{shiho3}, the coeffcient 
of the relative rigid cohomology (which is an overconvergent isocrystal) 
is assumed to extend to the boundary logarithmically, but 
here we do not impose such condition on the coefficient. So we 
prove the conjecture under less assumption. \par 
The first main theorem in this section is the following: 

\begin{thm}\label{main1}
Assume we are given a diagram 
\begin{equation}\label{maindiaglog}
(\ol{X},M_{\ol{X}}) \os{f}{\lra} (\ol{Y},M_{\ol{Y}}) 
\os{\iota}{\hra} (\cY,M_{\cY}), 
\end{equation}
where $f$ is a proper log smooth integral morphism in $\pLF$ having 
log smooth parameter and $\iota$ is a closed immersion. Assume moreover 
that we are given open immersions $j_X: X \hra \ol{X}, j_Y: Y \hra \ol{Y}$ 
satisfying $X \subseteq (\ol{X}, M_{\ol{X}})_{\triv}, 
Y \subseteq (\cY,M_{\cY})_{\triv}$ and $f^{-1}(Y)=X$. Then, 
for a locally free overconvergent 
isocrystal $\cE$ on $(X,\ol{X})/\cY_K$, the relative rigid cohomology 
$R^qf_{(X,\ol{X})/\cY, \rig *}\cE$ is a coherent 
$j_Y^{\d}\cO_{]\ol{Y}[_{\cY}}$-module for any $q \geq 0$. 
\end{thm} 

\begin{pf} 
By Theorem \ref{logrig-rig}, it suffices to prove that 
the relative log rigid cohomology $R^qf_{(X,\ol{X})/\cY, \logrig *}\cE$ 
is a coherent $j^{\d}_Y\cO_{]\ol{Y}[^{\log}_{\cY}}$-module for any $q\geq 0$. 
Let $(\ol{Y},M_{\ol{Y}}) \hra (\cY^{\ex},M_{\cY^{\ex}})$ 
be the exactification of $\iota$. We may assume that $\cY^{\ex}$ is affine. 
Put 
$\cI:=\Ker(\cO_{\cY^{\ex}} \ra \cO_{\ol{Y}})$, 
$\cJ:=\Ker(\cO_{\ol{Y}} \ra \cO_{\ol{Y}-Y})$ and 
let us take a lift $g_1, \cdots, g_c \in \Gamma(\cY^{\ex},\cO_{\cY^{\ex}})$ 
of generators of $\Gamma(\ol{Y},\cJ)$. \par 
Let us take 
an open covering 
$\ol{X} = \bigcup_{j \in J}\ol{X}_j$ by finite number of affine subschemes 
and exact closed immersions $i_j: 
(\ol{X}_j,M_{\ol{X}_j}) := 
(\ol{X}_j,M_{\ol{X}} |_{X_j}) \hra (\cP_j,M_{\cP_j}) \,(j \in J)$ 
into a fine log formal $\cB$-scheme $(\cP_j,M_{\cP_j})$ 
over $(\cY^{\ex},M_{\cY^{\ex}})$ satisfying 
$(\cP_j,M_{\cP_j}) 
\times_{(\cY^{\ex},}\allowbreak{}_{M_{\cY^{\ex}})} (\ol{Y},M_{\ol{Y}}) = 
(\ol{X}_j,M_{\ol{X}_j})$ 
such that each $(\cP_j,M_{\cP_j})$ is formally log smooth over 
 $(\cY^{\ex},M_{\cY^{\ex}})$ in the sense that 
$(\cP_j,M_{\cP_j})$ is $\cI\cO_{\cP_j}$-adically complete 
and that the morphism 
$$(\cP_j,M_{\cP_j}) \otimes_{\cO_{\cY^{\ex}}} \cO_{\cY^{\ex}}/\cI^n 
\lra (\cY^{\ex},M_{\cY^{\ex}}) \otimes_{\cO_{\cY^{\ex}}} 
\cO_{\cY^{\ex}}/\cI^n$$ is log smooth for any $n$. 
For any 
 non-empty subset $L \subseteq J$, let 
$(\ol{X}_L,M_{\ol{X}_L})$ (resp. $(\cP_L, M_{\cP_L})$) 
be the 
fiber product of $(\ol{X}_j,M_{\ol{X}_j})$'s (resp. 
$(\cP_j,M_{\cP_j})$'s) for $j \in L$ over 
$(\ol{X},M_{\ol{X}})$ (resp. $(\cY,M_{\cY})$) and 
for $m \in \N$, 
let $(\ol{X}^{(m)}, M_{\ol{X}^{(m)}})$ (resp. $(\cP^{(m)},M_{\cP^{(m)}})$) 
be the disjoint union of $(\ol{X}_L,M_{\ol{X}_L})$'s (resp. 
$(\cP_L,M_{\cP_L})$'s) for $L \subseteq J$ with $|L|=m+1$. 
Then $(\ol{X}^{(m)},M_{\ol{X}^{(m)}})$'s 
and $(\cP^{(m)},M_{\cP^{(m)}})$'s $(m \leq N:=|J|-1)$ 
naturally form diagrams 
$(\ol{X}^{(\b)},M_{\ol{X}^{(\b)}}), (\cP^{(\b)},M_{\cP^{(\b)}})$ 
indexed by the category $\Delta^+_N$ 
in the paragraph before Proposition 
\ref{logrig-fin} and we have a canonical diagram 
\begin{equation}\label{formalemb}
(\ol{X},M_{\ol{X}}) \os{g^{(\b)}}{\lla} (\ol{X}^{(\b)},M_{\ol{X}^{(\b)}}) 
\hra (\cP^{(\b)},M_{\cP^{(\b)}}) 
\end{equation} 
over $(\cY^{\ex},M_{\cY^{\ex}})$. \par 
Let $\{(T_{a,b}(\cY),M_{T_{a,b}(\cY)})\}_{(a,b) \in \A}$ be the system of 
universal enlargements of $((\cY,M_{\cY}), \allowbreak 
(\ol{Y}\allowbreak ,\allowbreak M_{\ol{Y}}))$. 
Let us take $\mu = p^{-b/a} \in p^{\Q_{<0}}$ and put 
$(\cY_{\mu},M_{\cY_{\mu}}):=(T_{a,b}(\cY),M_{T_{a,b}(\cY)})$.
Let us denote the canonical morphism 
$(\cY_{\mu},M_{\cY_{\mu}}) \lra (\cY,M_{\cY})$ by $\varphi_{\mu}$. 
Let us denote the base change of the diagram \ref{maindiaglog} by 
$\varphi_{\mu}$ by 
\begin{equation}\label{maindiaglogmu} 
(\ol{X}_{\mu},M_{\ol{X}_{\mu}}) \lra 
(\ol{Y}_{\mu},M_{\ol{Y}_{\mu}}) \hra (\cY_{\mu},M_{\cY_{\mu}})
\end{equation} 
and let us denote the base change of the diagram \eqref{formalemb} by 
$\varphi_{\mu}$ by 
\begin{equation}\label{embmu}
(\ol{X}_{\mu},M_{\ol{X}_{\mu}}) 
\lla (\ol{X}_{\mu}^{(\b)},M_{\ol{X}_{\mu}^{(\b)}}) 
\hra (\cP^{(\b)}_{\mu},M_{\cP^{(\b)}_{\mu}}). 
\end{equation} 
Put $X_{\mu}:=X \times_{\ol{X}} \ol{X}_{\mu}, 
Y_{\mu}:=Y \times_{\ol{Y}} \ol{Y}_{\mu}$ and  
$X^{(\b)}_{\mu} := X \times_{\ol{X}} \ol{X}^{(\b)}_{\mu}$. 
Then $]\ol{Y}[^{\log}_{\cY} = \bigcup_{\mu} 
]\ol{Y}[_{\cY_{\mu}} = \bigcup_{\mu} \cY_{\mu,K}$ is an 
admissible covering and we have 
$R^qf_{(X_{\mu},\ol{X}_{\mu})/\cY_{\mu}, \logrig *} \cE = 
(R^qf_{(X,\ol{X})/\cY, \logrig *}\cE) |_{\cY_{\mu,K}}$. 
So it suffices to prove that $R^qf_{(X_{\mu},}\allowbreak
{}_{\ol{X}_{\mu})/\cY_{\mu}, 
\logrig *} \cE$ is a coherent $j^{\d}_Y\cO_{\cY_{\mu,K}}$-module. \par 
Let $h$ be the morphism $]\ol{X}^{(\b)}_{\mu}[^{\log}_{\cP^{(\b)}_{\mu}} \lra 
]\ol{Y}_{\mu}[_{\cY_{\mu}} = \cY_{\mu,K}$ and put 
$\DR^{\d,(\b)} := 
\DR^{\d}(]\ol{X}^{(\b)}_{\mu}[^{\log}_{\cP^{(\b)}_{\mu}} \allowbreak 
/ \allowbreak \cY_{\mu,K},\cE)$. For any sufficiently small 
strict neighborhood $V^{(\b)}$ of $]X_{\mu}^{(\b)}[_{\cP_{\mu}^{(\b)}}$ 
in $]\ol{X}_{\mu}^{(\b)}[^{\log}_{\cP_{\mu}^{(\b)}}$ indexed by $\Delta^+_N$, 
we have the log de Rham complex $\DR_{V^{(\b)}} := \DR(V^{(\b)}/\cY_K,\cE)$ 
compatible with respect to $V^{(\b)}$ satisfying 
$\DR^{\d,(\b)} = \varinjlim_{V^{(\b)}} j_{V^{(\b)},*}\DR_{V^{(\b)}}$. 
(Here 
$j_{V^{(\b)}}$ denotes the admissible open immersion 
$V^{(\b)} \hra ]\ol{X}_{\mu}^{(\b)}[^{\log}_{\cP_{\mu}^{(\b)}}$.) 
For $V^{(\b)}$ as above and for $\lam \in (0,1)$, we put 
$V^{(\b)}_{\lam} := V^{(\b)} \,\cap\, 
]\ol{X}^{(\b)}_{\mu}[^{\log}_{\cP^{(\b)}_{\mu},\lam}$ and put 
$\DR_{V^{(\b)},\lam} := \DR_{V^{(\b)}} |_{V^{(\b)}_{\lam}}$. 
Denote the admissible 
open immersions $V^{(\b)}_{\lam} \hra V^{(\b)}$ by $i_{V^{(\b)},\lam}$ and 
denote the composite 
$V^{(\b)}_{\lam} \os{i_{V^{(\b)},\lam}}{\hra} V^{(\b)} 
\os{j_{V^{(\b)}}}{\hra} 
]\ol{X}^{(\b)}_{\mu}[^{\log}_{\cP^{(\b)}_{\mu}}$ by $j_{V^{(\b)},\lam}$. 
Then, as complexes of sheaves, we have 
$\DR_{V^{(\b)}} = \varprojlim_{\lam} i_{V^{(\b)},\lam *} 
\DR_{V^{(\b)},\lam}$. \par 
For $\nu = p^{-1/e} \in (0,1)$, 
let $U_{\nu}$ be the admissible open subset 
$\bigcup_{i=1}^c\{|g_i| \geq \nu\}$ of $\cY_{\mu,K}$ and denote the open 
immersion $U_{\nu} \hra \cY_{\mu,K}$ by $j_{Y,\nu}$. Let us denote 
the inverse image of $U_{\nu}$ in $]\ol{X}^{(\b)}[^{\log}_{\cP^{(\b)}}$ by 
$U^{(\b)}_{\nu}$ and let us put $U^{(\b)}_{\lam,\nu} := 
U^{(\b)}_{\nu} \,\cap\, ]\ol{X}^{(\b)}[^{\log}_{\cP^{(\b)},\lam}$. 
Then, for any fixed $\lam \in (0,1)$, the category of strict neighborhoods 
$V^{(\b)}$ satisfying the condition $V^{(\b)}_{\lam} = U^{(\b)}_{\lam,\nu}$ 
is cofinal with the category of all the strict neighborhoods. 
When this condition is satisfied, we denote the de Rham complex 
$\DR_{V^{(\b)},\lam}$ by $\DR^{(\b)}_{\lam,\nu}$ and denote the morphism 
$V^{(\b)}_{\lam} = U^{(\b)}_{\lam,\nu} \lra U_{\nu}$ by $h_{\lam,\nu}$. \par 
With the above notation, we have the following diagram: 
{\allowdisplaybreaks{
\begin{align*}
Rf_{(X_{\mu},\ol{X}_{\mu})/\cY_{\mu}, \logrig *}\cE 
& = 
Rh_*\DR^{\d,(\b)} = 
Rh_*\varinjlim_{V^{(\b)}}j_{V^{(\b)},*}\DR_{V^{(\b)}} \\ 
& = 
Rh_*\varinjlim_{V^{(\b)}}j_{V^{(\b)},*}\varprojlim_{\lam}i_{V^{(\b)},\lam *}
\DR_{V^{(\b)},\lam} \\ 
& = 
Rh_*\varinjlim_{V^{(\b)}}\varprojlim_{\lam}j_{V^{(\b)},\lam *}
\DR_{V^{(\b)},\lam} \\ 
& \os{\varphi_1}{\lra} 
Rh_*R\varprojlim_{\lam} \varinjlim_{V^{(\b)}} j_{V^{(\b)},\lam *}
\DR_{V^{(\b)},\lam} \\ 
& = 
R\varprojlim_{\lam}Rh_*\varinjlim_{V^{(\b)}} j_{V^{(\b)},\lam *}
\DR_{V^{(\b)},\lam} \\ 
& = 
R\varprojlim_{\lam}Rh_*\varinjlim_{V^{(\b)}} Rj_{V^{(\b)},\lam *}
\DR_{V^{(\b)},\lam} \qquad \text{(Lemma \ref{acyc1})} \\ 
& \os{\varphi_2}{\lla}
R\varprojlim_{\lam}\varinjlim_{V^{(\b)}} R(h \circ j_{V^{(\b)},\lam})_* 
\DR_{V^{(\b)},\lam} \\ 
& = 
R\varprojlim_{\lam}(\varinjlim_{\nu}Rj_{Y,\nu,*}Rh_{\lam,\nu *} 
\DR^{(\b)}_{\lam,\nu}). 
\end{align*}}}

By the above diagram, we see that, to prove the theorem, it suffices to 
prove the following claim: \\
\quad \\
{\bf claim.} \,\,\, With the notation above, \\
(1)\,\, The map $\varphi_1$ is a quasi-isomorphism. \\
(2)\,\, The map $\varphi_2$ is a quasi-isomorphism. \\
(3)\,\, If we fix $\lam \in (p^{-r(X,Y;\sqrt{\mu})},1)$, 
$Rh_{\lam,\nu *}\DR^{(\b)}_{\lam,\nu}$ has 
bounded and coherent cohomologies for $\nu$ sufficiently close to $1$, 
which are compatible with respect to $\nu$. \\
(4)\,\, $\varinjlim_{\nu}Rj_{Y,\nu,*}Rh_{\lam,\nu *} 
\DR^{(\b)}_{\lam,\nu}$ is independent of $\lam$ for 
$\lam \in (p^{-r(X,Y;\sqrt{\mu})},1).$ \\
\quad \\
Indeed, let us assume the claim is proved. Then we have 
$R^qf_{(X_{\mu},\ol{X}_{\mu})/\cY_{\mu}, \logrig *}\cE 
= \varinjlim_{\nu} \allowbreak j_{Y,\nu,*}R^qh_{\lam,\nu *} 
\DR^{(\b)}_{\lam,\nu}$ (because $j_{Y,\nu}$ is acyclic for coherent modules 
by the proof of Lemma \ref{acyc1}) and so 
$R^qf_{(X_{\mu},\ol{X}_{\mu})/\cY_{\mu}, \logrig *}\cE$ is a coherent 
$j_Y^{\d}\cO_{\cY_{\mu,K}}$-module. So, in the following, 
we will prove the above claim. \par 
First we prove the claim (1). Note that 
$R\varprojlim_{\lam}Rh_*\varinjlim_{V^{(\b)}} j_{V^{(\b)},\lam *}
\DR_{V^{(\b)},\lam}$, which is quasi-isomorphic to the target of the map 
$\varphi_1$, is nothing but 
$R\varprojlim_{\lam}Rf_{(X_{\mu},\ol{X}_{\mu})/\cY_{\mu}, \logrig,\lam *}\allowbreak\cE$. 
Since we have the canonical quasi-isomorphism 
$$Rf_{(X_{\mu},\ol{X}_{\mu})/\cY_{\mu}, \logrig *}\cE \os{=}{\lra}
Rf_{(X_{\mu},\ol{X}_{\mu})/\cY_{\mu}, \logrig,\lam *}\allowbreak \cE $$ 
for $\lam > \mu$ by Theorem \ref{rig-lamrig}, we have 
$R\varprojlim_{\lam}Rf_{(X_{\mu},\ol{X}_{\mu})/\cY_{\mu}, 
\logrig,\lam *}\cE =R \allowbreak f_{(X_{\mu},\ol{X}_{\mu})/\cY_{\mu},}
\allowbreak{}_{\logrig,\lam_0 *}\allowbreak \cE$ for 
any fixed $\lam_0 \in (\mu,1)$ and it is quasi-isomorphic to 
$Rf_{(X_{\mu},\ol{X}_{\mu})/\cY_{\mu}, \logrig *}\cE$ via $\varphi_1$. 
So claim (1) is proved. \par 
Next we prove (3). 
For $0 \leq n \leq N$ and $i=0,1,2$, 
let $(\cP^{(n)}_{\mu}(i),M_{\cP^{(n)}_{\mu}(i)})$ be the 
$(i+1)$-fold fiber product of $(\cP^{(n)}_{\mu}, M_{\cP^{(n)}_{\mu}})$ over 
$(\cY_{\mu}, M_{\cY_{\mu}})$. 
Then the restriction of $\cE$ to 
$I^{\d}((X_{\mu},\ol{X}_{\mu})/\cY_{\mu})=
I^{\d}(((X_{\mu},\ol{X}_{\mu})/\cY_{\mu})^{\log})$ 
induces a compatible family (with respect to $n=0,1,2$) of 
modules with stratifications on some strict neighborhood of 
$]X^{(n)}_{\mu}[_{\cP^{(n)}_{\mu}}$ in $]\ol{X}^{(n)}[^{\log}_{\cP^{(n)}}$: 
Precisely speaking, there exist strict neighborhoods 
$V^{n}(i) \, (0 \leq n \leq N, i=0,1,2)$ of 
$]X^{(n)}_{\mu}[_{\cP^{(n)}_{\mu}(i)}$ in 
$]\ol{X}^{(n)}_{\mu}[_{\cP^{(n)}_{\mu}(i)}$ 
with 
$p_k(V^{(n)}(1)) \subseteq V^{(n)}(0), p_{kl}(V^{(n)}(2)) \subseteq 
V^{(n)}(1), \Delta(V^{(n)}(0)) \subseteq V^{(n)}(1)$ 
(where $p_k, p_{kl}$ are projections 
\begin{align*}
& p_k:\,]\ol{X}^{(n)}_{\mu}[^{\log}_{\cP^{(n)}_{\mu}(1)} \lra 
]\ol{X}^{(n)}_{\mu}[^{\log}_{\cP^{(n)}_{\mu}} \,(k=1,2), \\ 
& p_{kl}:\,]\ol{X}^{(n)}_{\mu}[^{\log}_{\cP^{(n)}_{\mu}(2)} \lra 
]\ol{X}^{(n)}_{\mu}[^{\log}_{\cP^{(n)}_{\mu}(1)} \,(1 \leq k < l < 3) 
\end{align*}
and $\Delta$ is 
the diagonal $]\ol{X}^{(n)}_{\mu}[^{\log}_{\cP^{(n)}_{\mu}} \lra 
]\ol{X}^{(n)}_{\mu}[^{\log}_{\cP^{(n)}_{\mu}(1)}$) which is 
compatible with respsect to $n$ such that the restriction 
of $\cE$ to $I^{\d}((X_{\mu},\ol{X}_{\mu})/\cY_{\mu})$ is induced from 
the compatible family $(E^{(n)},\epsilon^{(n)})_{0\leq n \leq N}$, 
where $E^{(n)}$ is a locally free 
$\cO_{V^{(n)}(0)}$-module and $\epsilon^{(n)}$ is the isomorphism 
$p_2^*E^{(n)} \os{\cong}{\lra} p_1^*E^{(n)}$ on $V^{(n)}(1)$ satisfying 
$\Delta^*(\epsilon^{(n)})=\id$, 
$p_{12}^*(\epsilon^{(n)}) \circ p_{23}^*(\epsilon^{(n)}) \allowbreak = 
\allowbreak p_{13}^*(\epsilon^{(n)}).$ \par 
Now let us take $\lam \in (p^{-r(X,Y,\sqrt{\mu})},1)$ and fix it. 
For $\nu=p^{-1/e} \in (0,1) \, (e \in \N)$, 
let us denote the inverse image of 
$U_{\nu}$ in $]\ol{X}^{(n)}_{\mu}[^{\log}_{\cP^{(n)}_{\mu}(i)}$ by 
$U^{(n)}_{\nu}(i)$. Then, by taking $\nu$ suffiently close to $1$, 
we may assume that 
the inclusion $U_{\lam,\nu}^{(n)}(i) := U_{\nu}^{(n)}(i) \cap 
]\ol{X}^{(n)}_{\mu}[^{\log}_{\cP^{(n)}_{\mu}(i),\lam} 
\allowbreak \subseteq \allowbreak 
V^{(n)}(i)$ for all $0 \leq n \leq N,i=0,1,2$. 
Note that the restriction of $(E^{(n)},\epsilon^{(n)})_{0\leq n \leq N}$ 
to $U_{\lam,\nu}^{(n)}(i)$ induces the log de Rham complex 
$\DR_{\lam,\nu}^{(n)}$ above. \par 
Now, for $1 \leq j \leq c$ and 
$\nu=p^{-1/e}$ as above, we put 
$\cY_{\mu,\nu,j} := {\ul{\Spf}}_{\cO_{\cY_{\mu}}} 
(\cO_{\cY_{\mu}}[t]/(g_j^e t-p))^{\wedge}, 
M_{\cY_{\mu,\nu,j}} := M_{\cY_{\mu}}|_{\cY_{\mu,\nu,j}}$. 
Then we have $U_{\nu}=\bigcup_{j=1}^c \cY_{\mu,\nu,j,K}$. 
Let us denote the canonical morphism 
$(\cY_{\mu,\nu,j},M_{\cY_{\mu,\nu,j}}) \lra (\cY_{\mu},M_{\cY_{\mu}})$ 
by $\psi_{\nu,j}$, let us denote the base change of the diagrams 
\eqref{maindiaglogmu}, \eqref{embmu} by $\psi_{\nu,j}$ by 
\begin{equation}\label{maindiaglogmunu}
(\ol{X}_{\mu,\nu,j},M_{\ol{X}_{\mu,\nu,j}}) 
\lra (\ol{Y}_{\mu,\nu,j},M_{\ol{Y}_{\mu,\nu,j}})
\hra (\cY_{\mu,\nu,j},M_{\cY_{\mu,\nu,j}}), 
\end{equation} 
\begin{equation}\label{embmunu}
(\ol{X}_{\mu,\nu,j},M_{\ol{X}_{\mu,\nu,j}}) 
\lla (\ol{X}_{\mu,\nu,j}^{(\b)},M_{\ol{X}_{\mu,\nu,j}^{(\b)}}) 
\hra (\cP^{(\b)}_{\mu,\nu,j},M_{\cP^{(\b)}_{\mu,\nu,j}}), 
\end{equation} 
respectively and let us define $\cP^{(n)}_{\mu,\nu,j}(i)$ in the same way as 
$\cP^{(n)}_{\mu}(i)$. 
Then, for $0 \leq n \allowbreak \leq N$ and $i=0,1,2$, 
$]\ol{X}^{(n)}_{\mu,\nu,j}[^{\log}_{\cP^{(n)}_{\mu,\nu,j}(i)}$ is 
an admissible open subset of 
$]\ol{X}^{(n)}_{\mu}[^{\log}_{\cP^{(n)}_{\mu}(i)}$ and we have 
$U_{\lam,\nu}^{(n)}(i) \,\cap\, 
]\ol{X}^{(n)}_{\mu,\nu,j}[^{\log}_{\cP^{(n)}_{\mu,\nu,j}(i)} = 
]\ol{X}^{(n)}_{\mu,\nu,j}[^{\log}_{\cP^{(n)}_{\mu,\nu,j}(i),\lam}$. 
So $(E^{(n)},\epsilon^{(n)}) \,(0 \leq n \leq N)$ 
defines naturally an object in 
$\Strat''_{\lam}((\ol{X}^{(n)}_{\mu,\nu,j} \hra 
\cP^{(n)}_{\mu,\nu,j}/\cY_{\mu,\nu,j})^{\log})$ 
which is compatible with respect 
to $n$, that is, an object (which we denote by $\cE_{\lam,\nu,j}$) in 
$I_{\lamm\conv}((\ol{X}_{\mu,\nu,j}/\cY_{\mu,\nu,j})^{\log})$ such that 
the log de Rham complex associated to it is nothing but the complex 
$\DR^{(\b)}_{\lam,\nu}|_{]\ol{X}^{(\b)}_{\mu,\nu,j}[^{\log}_{\cP^{(\b)}_{\mu,\nu,j},\lam}}$. If we denote the morphism 
$]\ol{X}^{(\b)}_{\mu,\nu,j}[^{\log}_{\cP^{(\b)}_{\mu,\nu,j},\lam} \lra 
\cY_{\mu,\nu,j,K}$ by $h_{\lam,\nu,j}$, we have the isomorphism 
\begin{align*}
(R^qh_{\lam,\nu *}\DR^{(\b)}_{\lam,\nu})|_{\cY_{\mu,\nu,j}} & = 
R^qh_{\lam,\nu,j*}
(\DR^{(\b)}_{\lam,\nu}|_{]\ol{X}^{(\b)}_{\mu,\nu,j}[^{\log}_{\cP^{(\b)}_{\mu,\nu,j},\lam}}) \\ 
& = 
R^qf_{\ol{X}_{\mu,\nu,j}/\cY_{\mu,\nu,j}, \lamm\an *}\cE_{\lam,\nu,j}. 
\end{align*}
Since the morphisms $\varphi_{\mu}$ and $\psi_{\nu,j}$ are affine, 
we have $a(\ol{X}_{\mu,\nu,j}) \leq a(\ol{X})$. So we have 
$r(\ol{X}_{\mu,\nu,j},\cY_{\mu,\nu,j};\sqrt{\mu}) \geq 
r(\ol{X},\cY,\sqrt{\mu})$. Hence, by Theorem \ref{coherence0}, 
$R^qf_{\ol{X}_{\mu,\nu,j}/\cY_{\mu,\nu,j} \lamm\an *}\cE_{\mu,\nu,j}$
 is a coherent module on 
$]\ol{Y}_{\mu,\nu,j}[_{\cY_{\mu,\nu,j},\sqrt{\mu}} = 
\cY_{\mu,\nu,j,K}$. So we have proved the former part of claim (3). 
Moreover, the analytically flat base change implies that the above 
cohomology is compatible with respect to $\nu$. So the proof of claim (3) 
is finished. \par 
Next we prove claim (4). Take $\lam'>\lam > r(\ol{X},\ol{Y};\sqrt{\mu})$. 
Then, if $\nu$ is sufficiently close to $1$, the isocrystal 
$\cE_{\lam',\nu,j} \in 
I_{\lam'\text{-}\conv}((\ol{X}_{\mu,\nu,j}/\cY_{\mu,\nu,j})^{\log})$ 
analogous to $\cE_{\lam,\nu,j}$ above is defined and 
$R^qf_{\ol{X}_{\mu,\nu,j}/\cY_{\mu,\nu,j}, \lam'\text{-}\an *}
\cE_{\lam',\nu,j}$ is coherent. In this case, $\cE_{\lam,\nu,j}$ is 
the restriction of $\cE_{\lam',\nu,j}$ to 
$I_{\lamm\conv}((\ol{X}_{\mu,\nu,j}/\cY_{\mu,\nu,j})^{\log})$ and so 
we have the canonical isomorphism 
$R^qf_{\ol{X}_{\mu,\nu,j}/\cY_{\mu,\nu,j}, \lam'\text{-}\an *}
\cE_{\lam',\nu,j} \cong 
R^qf_{\ol{X}_{\mu,\nu,j}/\cY_{\mu,\nu,j}, \lamm\an *}
\cE_{\lam,\nu,j}$ by Proposition \ref{33}. So we have 
\begin{align*} 
\varinjlim_{\nu}Rj_{Y,\nu,*}Rh_{\lam',\nu *} 
\DR^{(\b)}_{\lam',\nu} & = 
\varinjlim_{\nu}j_{Y,\nu,*}Rh_{\lam',\nu *} 
\DR^{(\b)}_{\lam',\nu} \\ 
& = 
\varinjlim_{\nu}j_{Y,\nu,*}Rh_{\lam,\nu *} 
\DR^{(\b)}_{\lam,\nu} 
= \varinjlim_{\nu}Rj_{Y,\nu,*}Rh_{\lam,\nu *} 
\DR^{(\b)}_{\lam,\nu}. 
\end{align*}
So claim (4) is proved. \par 
Finally we prove claim (2). To prove claim (2), we may compare both hand sides 
before taking $R\vpl_{\lam}$ and we may replace $\b$ 
in the source and the target by $n \in \N$. 
Then, we see that the source and 
the target are unchanged even if we replace the closed immersion 
$(\ol{X}^{(n)}_{\mu},M_{\ol{X}^{(n)}_{\mu}}) 
\hra (\cP^{(n)}_{\mu},M_{\cP^{(n)}_{\mu}})$ by another closed immersion 
$(\ol{X}^{(n)}_{\mu},M_{\ol{X}^{(n)}_{\mu}}) 
\hra ({\cP'}^{(n)}_{\mu},M_{{\cP'}^{(n)}_{\mu}})$
such that $({\cP'}^{(n)}_{\mu},M_{{\cP'}^{(n)}_{\mu}})$ is formally log smooth 
over $(\cY_{\mu},M_{\cY_{\mu}})$: As for the source, it is true because 
the source is obtained as $\varinjlim_{\nu} Rj_{Y,\nu,*}$ of 
the relative log analytic cohomology of radius $\lam$ 
(as we saw above) and 
as for the target, it is true because 
the target is obtained as the $\lam$-restriction of 
the relative log rigid cohomology. So we may replace 
$(\cP^{(n)}_{\mu},M_{\cP^{(n)}_{\mu}})$ so that 
$(\cP^{(n)}_{\mu},M_{\cP^{(n)}_{\mu}}) 
\times_{(\cY_{\mu},M_{\cY_{\mu}})} (\ol{Y}_{\mu},M_{\ol{Y}_{\mu}}) = 
(\ol{X}^{(n)}_{\mu},M_{\ol{X}^{(n)}_{\mu}})$ holds. In this case, 
$Rh_*$ commutes with direct limit 
because $h:\,]\ol{X}^{(n)}_{\mu}[^{\log}_{\cP^{(n)}_{\mu}} = \cP^{(n)}_{\mu,K} 
\lra \cY_{\mu,K}$ is quasi-compact and quasi-separated. So the map 
$\varphi_2$ is a quasi-isomorphism. So claim (2) is proved. \par 
Since the claim is proved, the proof of the theorem is finished. 
\end{pf} 

Let us take a triple of the form $(S,S,\cS)$ 
and assume we are given a diagram 
$$ (\ol{X},M_{\ol{X}}) \os{f}{\lra} (\ol{Y},M_{\ol{Y}}) 
\os{g}{\lra} S $$
and open immersions $j_X:X \hra \ol{X}, j_Y:Y \hra \ol{Y}$ 
with $X \subseteq (\ol{X},M_{\ol{X}})_{\triv}, Y \subseteq 
(\ol{Y},M_{\ol{Y}})_{\triv}$ and $f^{-1}(Y)=X$, where 
$f$ is a proper log smooth integral 
morphism in $\LB$ having log smooth parameter 
and $g$ is a morphism in $\LB$. Let 
$\cE$ be a locally free overconvergent isocrystal on 
$(X,\ol{X})/\cS_K$. Then we have the following theorem, which is the main 
result in this section: 

\begin{thm}\label{thm00}
Let the notations be as above. Then there exists a diagram 
$$ \ol{Y} \os{g^{(0)}}{\lla} \ol{Y}^{(0)} \os{i^{(0)}}{\hra} \cY^{(0)}, $$ 
where $g^{(0)}$ is an etale surjective morphism, 
$i^{(0)}$ is a closed immersion and 
$\cY^{(0)}$ is formally smooth over $\cS$, which satisfies the following$:$ 
For $q \geq 0$, there exists uniquely an overconvergent isocrystal $\cF$ on 
$(Y,\ol{Y})/\cS_K$ such that, for any triple $(Z,\ol{Z},\cZ)$ over 
$(Y^{(0)},\ol{Y}^{(0)},\cY^{(0)})$ $($where $Y^{(0)}:=Y \times_{\ol{Y}} 
\ol{Y}^{(0)})$ satisfying $Z = Y^{(0)} \times_{\ol{Y}^{(0)}} \ol{Z}$ with 
$\cZ$ formally smooth over $\cY$, the restriction of $\cF$ to 
$I^{\d}((Z,\ol{Z})/\cS_K,\cZ)$ is given 
functorially by 
$(R^qf_{(X \times_Y Z,\ol{X} \times_{\ol{Y}} \ol{Z})/\cZ, \rig *}\cE, 
\epsilon)$, where $\epsilon$ is given by 
$$ 
p_2^*R^qf_{(X \times_Y Z,\ol{X} \times_{\ol{Y}} \ol{Z})/\cZ, \rig *}\cE 
\os{\simeq}{\rightarrow} 
R^qf_{(X \times_Y Z,\ol{X} \times_{\ol{Y}} \ol{Z})/\cZ \times_{\cS} \cZ, 
\rig *}\cE \os{\simeq}{\leftarrow} 
p_1^*R^qf_{(X \times_Y Z,\ol{X} \times_{\ol{Y}} \ol{Z})/\cZ, \rig *}\cE. 
$$ 
$($Here $p_i$ is the morphism 
$]\ol{Z}[_{\cZ \times_{\cS} \cZ} \lra \,]\ol{Z}[_{\cZ}$ induced by the 
$i$-th projection.$)$ 
\end{thm} 

\begin{pf} 
Let us take a diagram 
$$ (\ol{Y},M_{\ol{Y}}) \os{g^{(0)}}{\lla} (\ol{Y}^{(0)},M_{\ol{Y}^{(0)}})
\os{i^{(0)}}{\hra} (\cY^{(0)},M_{\cY^{(0)}}) $$ 
such that $g^{(0)}$ is strict etale surjective, $i^{(0)}$ is a closed 
immersion in $\pLF$, $(\cY^{(0)},M_{\cY^{(0)}})$ is formally log smooth over 
$\cS$ and $\cY^{(0)}$ is formally smooth over $\cS$. For $n \in \N$, let 
$(\ol{Y}^{(n)},M_{\ol{Y}^{(n)}})$ (resp. $(\cY^{(n)},M_{\cY^{(n)}})$) 
be the $(n+1)$-fold fiber product of $(\ol{Y}^{(0)},M_{\ol{Y}^{(0)}})$ 
(resp. $(\cY^{(0)},M_{\cY^{(0)}})$) over $(\ol{Y},M_{\ol{Y}})$ (resp. 
$\cS$) and denote the pull-back of $(\ol{X},M_{\ol{X}}), X$ by 
the morphism $(\ol{Y}^{(n)},M_{\ol{Y}^{(n)}}) \lra (\ol{Y},M_{\ol{Y}})$ by 
$(\ol{X}^{(n)},M_{\ol{X}^{(n)}}), X^{(n)}$, respectively. Then, by 
Theorem \ref{main1}, $R^qf_{(X^{(n)},\ol{X}^{(n)})/\cY^{(n)}, \rig *}\cE$ 
(resp. 
$R^qf_{(X^{(n)},\ol{X}^{(n)})/\cY^{(n)}\times_{\cS}\cY^{(n)},}\allowbreak
{}_{\rig *}\cE$) 
is a coherent $j_Y^{\d}\cO_{]\ol{Y}^{(n)}[_{\cY^{(n)}}}$-module 
(resp. a coherent $j_Y^{\d}\cO_{]\ol{Y}^{(n)}[_{\cY^{(n)} \times_{\cS} 
\cY^{(n)}}}$-module). 
Then, by the base change theorem of Tsuzuki (\cite[2.3.1]{tsuzuki3}), 
the arrows in the diagram 
$$ 
p_2^*R^qf_{(X^{(n)},\ol{X}^{(n)})/\cY^{(n)}, \rig *}\cE \ra 
R^qf_{(X^{(n)},\ol{X}^{(n)})/\cY^{(n)}\times_{\cS}\cY^{(n)}, \rig *}\cE 
\leftarrow 
p_1^*R^qf_{(X^{(n)},\ol{X}^{(n)})/\cY^{(n)}, \rig *}\cE $$ 
(where $p_i$ denotes the $i$-th projection 
$]\ol{Y}^{(n)}[_{\cY^{(n)} \times_{\cS} \cY^{(n)}} \lra 
]\ol{Y}^{(n)}[_{\cY^{(n)}}$) are isomorphisms. If we denote the composite 
of the above arrows by $\epsilon^{(n)}$, one can see (by using 
\cite[2.3.1]{tsuzuki3} again) that $\cF^{(n)}:=
(R^qf_{(X^{(n)},\ol{X}^{(n)})/\cY^{(n)}, \rig *}\cE,\epsilon^{(n)})$ 
defines an overconvergent isocrystal on $(Y^{(n)},\ol{Y}^{(n)})/\cS_K$ 
and that $\cF^{(n)}$ is compatible with respect to $n$. 
Then, by etale descent of overconvergent isocrystals (\cite[5.1]{shiho3}), 
we see that $\{\cF^{(n)}\}_{n=0,1,2}$ descents to an overconvergent 
isocrystal $\cF$ on $(Y,\ol{Y})/\cS_K$. \par 
Next we prove the uniqueness of $\cF$. 
Since the triple $(Y^{(n)},\ol{Y}^{(n)},\cY^{(n)})$ satisfies the 
condition required for $(Z,\ol{Z},\cZ)$ in the theorem, we see that the 
image of $\cF$ in $I^{\d}((Y^{(n)},\ol{Y}^{(n)})/\cS_K)$ should be 
functorially isomorphic to $\cF^{(n)}$ in the previous paragraph. Then, by 
\cite[5.1]{shiho3}, we see that this condition characterizes $\cF$ 
uniquely. \par 
Finally we prove that the overconvergent isocrystal $\cF$ satisfies the 
required condition. For a triple $(Z,\ol{Z},\cZ)$ as in the statement of 
the theorem, let us put $M_{\ol{Z}} := M_{\ol{Y}}|_{\ol{Z}}, 
M_{\cZ} := M_{\cY}|_{\cZ}$. Then the diagram 
$$ 
(\ol{Z},M_{\ol{Z}}) \times_{(\ol{Y},M_{\ol{Y}})} (\ol{X},M_{\ol{X}}) \lra 
(\ol{Z},M_{\ol{Z}}) \hra (\cZ,M_{\cZ}) 
$$ 
and the open immersions 
$Z \times_{Y} X \lra \ol{Z} \times_{\ol{Y}} \ol{X}, 
Z \hra \ol{Z}$ are as in the situation in Theorem \ref{main1}. 
So, by Theorem \ref{main1}, 
$R^qf_{(Z \times_Y X,\ol{Z} \times_{\ol{Y}} \ol{X})/\cZ, \rig *}\cE$ is 
a coherent $j^{\d}\cO_{]\ol{Z}[_{\cZ}}$-module. Then, by 
\cite[2.3.1]{tsuzuki3} again, we see that the restriction of $\cF$ to 
$I^{\d}((Z,\ol{Z})/\cS_K)$ is given by 
$(R^qf_{(X \times_Y Z,\ol{X} \times_{\ol{Y}} \ol{Z})/\cZ, \rig *}\cE, 
\epsilon)$ as in the statement of the theorem. 
Finally, the functoriality of the expression above is proved as the 
proof of \cite[4.8]{shiho3}. So we are done. 
\end{pf} 


\section{Preliminaries on alteration and hypercovering} 

In this section, we recall some preliminary facts on alteration 
proved by de Jong (\cite{dejong}, \cite{dejong2}) 
and using these, we prove the existence of certain 
diagram involving hypercovering. This diagram 
turns out to be useful to prove the overconvergence 
of relative rigid cohomology in general case, as is shown in the next 
section. \par 
{\it From now on in this paper, we assume that $k$ is perfect and 
$\cB = \Spf W(k)$ $($hence $B=\Spec k)$}. (So all the schemes are $k$-schemes separated of 
finite type and all the pairs are pairs separated of finite type 
over $k$.) \par 
First let us recall the notion of pluri nodal 
fibration, which was introduced in \cite{dejong2}: 

\begin{defn}
Let $S$ be a scheme. A pluri nodal fibration of relative dimension $d$ 
over $S$ is a system 
$(X_d \os{f_d}{\ra} X_{d-1} \os{f_{d-1}}{\ra} 
\cdots X_1 \os{f_1}{\ra} X_0=S, \{\sigma_{ij}\}_{1 \leq i \leq d, 
1 \leq j \leq n_i}, Z_0)$ satisfying the 
following conditions$:$ 
\begin{enumerate}
\item 
Each $f_i:X_i \lra X_{i-1}$ is a quasi-split semistable curve over $X_{i-1}$. 
\item 
$Z_0$ is a proper closed subset of $S$. 
\item 
For each $i$, $\sigma_{ij}:X_{i-1} \lra X_i \, (1 \leq j \leq n_i)$ are 
disjoint sections of $f_i$ into smooth locus. 
\item 
If we define $Z_i \subseteq X_i$ inductively by 
$Z_i := \bigcup_j \sigma_{ij}(X_{i-1}) \cup f_i^{-1}(Z_i)$, 
$f_{i+1}:X_{i+1} \lra X_i$ is smooth over $X_i-Z_i$. 
\end{enumerate}
\end{defn} 

Then de Jong proved in \cite[Thm 5.9]{dejong2} the following theorem: 

\begin{thm}\label{dejongthm}
Let $f:X \lra S$ be a proper morphism of integral schemes such that 
the generic fiber is geometrically irreducible of dimension $d \geq 1$. 
Then we have a commutative diagram 
\begin{equation*}
\begin{CD}
X @<{\psi}<< X_d \\ 
@VfVV @V{g}VV \\ 
S @<<< S', 
\end{CD}
\end{equation*}
where the horizontal arrows are alterations and $g$ is 
$($a part of a system of$)$ a pluri nodal fibration over $S'$. 
\end{thm} 

Let 
\begin{equation}\label{pluri}
X_d \os{g_d}{\ra} X_{d-1} \os{g_{d-1}}{\ra} \cdots X_1 
\os{g_1}{\ra} X_0=S', \qquad Z_i \subseteq X_i \,(0 \leq i \leq d) 
\end{equation} 
be a part of the data of pluri nodal fibration $g$ appearing in 
Theorem \ref{dejongthm}. Then, by \cite[Thm 4.1]{dejong}, 
there exists an alteration $\varphi:X'_0 \lra X_0$ such that 
$X'_0$ is regular and that $\varphi^{-1}(Z_0)$ is a simple normal 
crossing divisor in $X'_0$. So, by taking the base change by $\varphi$, 
we may assume that, in the diagram \eqref{pluri}, 
 $X_0=S'$ is regular and $Z_0$ is a simple normal crossing 
divisor in $X_0$ without changing the conclusion of 
Theorem \ref{dejongthm}. \par 
Next let us recall the following proposition, which is proved in 
\cite[Prop 5.11]{dejong}: 

\begin{prop}\label{dejongprop}
Let $S$ be a regular scheme, let $D \subseteq S$ be a normal crossing 
divisor and let $f:X \lra S$ be a quasi-split semi-stable curve which is 
smooth over $S-D$. Assume moreover that we are given 
disjoint sections $\sigma_1, \cdots, \sigma_n$ 
of $f$ into smooth locus. 
Then there exists a modification $X' \lra X$ 
satisfying the following conditions$:$ 
\begin{enumerate} 
\item 
$X'$ is regular and the center of the modification $X' \lra X$ 
is contained in non-regular locus of $X$. 
\item 
There exist lifts of sections 
$\ti{\sigma}_i:S \lra X'\,(1 \leq i \leq n)$ of $\sigma_i$ such that 
$D' \allowbreak := \bigcup_i \ti{\sigma}_i(S) \cup {f'}^{-1}(D)$ is 
a normal crossing divisor in $X'$, where $f'$ is the composite 
$X' \lra X \lra S$. 
\end{enumerate}
\end{prop} 

Moreover, by the proof of \cite[Prop 5.11]{dejong}, one has the 
following local expression of the morphism $f':X' \lra S$. 
If we take $x \in X'(\ol{k})$ and put $s:=f'(x) \in S(\ol{k})$ 
(where $\ol{k}$ is the algebraic closure of $k$), the homomorphism 
${f'}^*: \widehat{\cO}_{S,s} \lra \widehat{\cO}_{X,x}$ induced by $f'$ 
has one of the following two form: \\
(1) \,\, ${f'}^*$ has the form 
$$ \ol{k}[[t_1, \cdots, t_n]] \lra 
\ol{k}[[t_1, \cdots, t_n,u,v]]/(uv-t_1). $$
In this case, $D$ has the form $\bigcup_{i=1}^r\{t_i=0\}$ for some 
$r\geq 1$ and $D'$ has the form 
$\{u=0\} \cup \{v=0\} \cup (\bigcup_{i=2}^r\{t_i=0\})$. \\
(2) ${f'}^*$ has the form 
$$ \ol{k}[[t_1, \cdots, t_n]] \lra 
\ol{k}[[t_1, \cdots, t_n,y]]. $$
In this case, $D$ has the form $\bigcup_{i=1}^r\{t_i=0\}$ for some 
$r\geq 0$ and $D'$ has either the form 
$\bigcup_{i=1}^r\{t_i=0\}$ or the form 
$\bigcup_{i=1}^r\{t_i=0\} \cup \{y=0\}$. \\
\quad \\
Now, given a pluri nodal fibration 
\eqref{pluri} such that $X_0$ is regular and 
$Z_0$ is a normal crossing divisor, let us apply the following procedure 
for $i=1,\cdots,d$ inductively: 
Apply Proposition \ref{dejongprop} to the morphism 
$X_i \lra X_{i-1}$, the divisor $Z_{i-1} \subseteq X_{i-1}$ 
and the sections 
$\{\sigma_{ij}\}_j$ to obtain $Z'_i \subseteq X'_i \lra X_{i-1}$ and replace 
$Z_{i'}, X_{i'}, \{\sigma_{i'j}\} \,(i' \geq i)$ by 
$Z_{i'} \times_{X_i} X'_i, X_{i'} \times_{X_i} X'_i, 
\{\sigma_{i'j} \times \id \}$, respectively. After applying 
this procedure, we see that, for each $i$, the morphism 
$X_i \lra X_{i-1}$ satisfies the conditions in Proposition 
\ref{dejongprop} required for $f'$. So 
the morphism $X_d \os{g}{\lra} X_0$ and the divisors 
$Z_0 \subseteq X_0, Z_d \subseteq X_d$ 
have the following local expression: 
If we take $x \in X_d(\ol{k})$ and put $y:=g(x) \in X_0(\ol{k})$, 
the homomorphism 
$g^*: \widehat{\cO}_{X_0,y} \lra 
\widehat{\cO}_{X_d,x}$ induced by $g$ has the form 
\begin{align*}
& \ol{k}[[t_1, \cdots, t_n]] \lra \\
& \ol{k}[[x_{ij}\,(1 \leq i \leq r, 1 \leq j \leq m_i), 
t_{r+1}, \cdots t_n, 
y_1,\cdots, y_l]]/(\prod_{j=1}^{m_i}x_{ij}-t_i\,(1 \leq i \leq r)) 
\end{align*}
for some $n,r,l \geq 0$ with $r \leq n$, and $Z_0, Z_d$ have the form 
$Z_0 = \bigcup_{i=1}^r\{t_i=0\}$, 
$Z_d = \bigcup_{i=1}^r(\bigcup_{j=1}^{m_i}\{x_{ij}=0\}) \cup 
\bigcup_{i=1}^{l'}\{y_i=0\}$ (for some $l'\leq l$). In this case, 
$g^{-1}(Z_0) \subseteq Z_d$ corresponds to the subscheme 
$\bigcup_{i=1}^r(\bigcup_{j=1}^{m_i}\{x_{ij}=0\})$. So 
we see that $g^{-1}(Z_0) \subseteq X_d$ is also a normal crossing 
divisor. Moreover, if we denote the fine log structure on $X_0$ (resp. $X_d$) 
defined by $Z_0$ (resp. $g^{-1}(Z_0)$) by $M_{X_0}$ (resp. $M_{X_d}$), 
we see from the above local expression that the morphism 
$(X_d,M_{X_d}) \lra (X_0,M_{X_0})$ is proper log smooth integral. \par 
Summing up the above argument, we obtain the following proposition: 

\begin{prop}\label{dejongup}
Let $f:X \lra S$ be a proper morphism of integral schemes such that 
the generic fiber is geometrically irreducible of dimension $d \geq 1$. 
Then we have the commutative diagram 
\begin{equation}\label{djdiag}
\begin{CD}
X @<{\psi}<< X_d \\ 
@VfVV @V{g}VV \\ 
S @<<< X_0
\end{CD}
\end{equation}
and an open subset $U_0 \subseteq X_0$ such that 
the horizontal arrows are alterations, $Z_0:=X_0-U_0, 
Z_d:=X_d-g^{-1}(U_0)$ are normal crossing divisors of regular schemes 
$X_0,X_d$ respectively 
and that, if we denote the log structure on $X_0$ $($resp. $X_d)$ defined 
by $Z_0$ $($resp. $Z_d)$ by $M_{X_0}$ $($resp. $M_{X_d})$, the morphism 
$(X_d,M_{X_d}) \lra (X_0,M_{X_0})$ is proper log smooth integral. 
\end{prop} 

\begin{rem} 
Proposition \ref{dejongup} is true even for $d=0$: Indeed, 
in the case $d=0$, $f$ is birational. So, if we take an alteration 
$X_0 \lra X$ with $X_0$ regular and if we put $X_d:=X_0, g:=\id, U_0:=X_0$, 
we have the diagram \eqref{djdiag}. 
\end{rem} 

\begin{rem} 
It is stated (without proofs) in \cite{ak} that the statement of 
Proposition \ref{dejongup} is shown in \cite{dejong}, but it seems 
that we need results in \cite{dejong2}. We included the indication of 
the proof of the above proposition just because 
we could not find a proposition which states directly the content of 
the above proposition in the paper \cite{dejong2}. 
\end{rem} 

In the rest of this section, we prove, for a given proper morphism 
of schemes $X \lra Y$, the existence of certain diagram involving 
hypercovering by `nice' schemes, by using Proposition \ref{dejongup}. 
To this end, first we introduce some notions on morphism of 
pairs and schemes and prove 
an elementary lemma. Recall (\cite[2.3.3]{chts}) that a 
morphism of pairs $f:(U_X,X) \lra (U_Y,Y)$ is called strict if 
we have $f^{-1}(U_Y)=U_X$. We introduce the following terminology 
for certain covering of pairs (cf. \cite{tsuzuki4}): 

\begin{defn} 
A morphism of pairs $f:(U_X,X) \lra (U_Y,Y)$ is called a proper covering 
if $f$ is strict, $f:X \lra Y$ is proper and $f:U_X \lra U_Y$ is 
proper surjective. Moreover, if $f:X \lra Y$ is also surjective, $f$ is 
called a strongly proper covering. A hypercovering with respect to 
proper coverings 
is called a 
proper hypercovering.  
We can also define the notion of $n$-truncated proper hyperconvering 
in the same way. 
\end{defn} 

We also introduce the notion of a good proper surjective morphism 
and a good strongly proper covering as follows. 

\begin{defn} 
A proper surjective morphism of schemes $f:X \lra Y$ is called good 
if $f$ can be written as a composite of finite number of the morphisms 
of the following forms$:$ 
\begin{enumerate}
\item 
The morphism $\coprod_{i=1}^n Y_i \lra Y_{\red} \lra Y$ induced by the 
canonical inclusions, where $Y_{\red} = \bigcup_{i=1}^n Y_i$ is the 
decomposition of $Y_{\red}$ into connected components. 
\item 
The morphism $\coprod_{i=1}^n f_i: \coprod_{i=1}^n X_i \lra \coprod_{i=1}^n 
Y_i$, where each $f_i:X_i \lra Y_i$ 
is an alteration $($between integral schemes$)$. 
\end{enumerate} 
A strongly proper covering $f:(U_X,X) \lra (U_Y,Y)$ is called good if 
$f:X \lra Y$ is good. 
\end{defn} 

Note that 
a good proper surjective morphism $f:X \lra Y$ of integral schemes is 
necessarily an alteration. It is easy to see that a
good proper surjective morphism $f: X \lra Y$ has the following properties: 
\begin{enumerate} 
\item 
If $Z \subseteq X$ is a closed subset such that $X-Z \subseteq X$ is dense, 
$f(Z) \subseteq X$ is also a closed subset such that $Y-f(Z) \subseteq Y$ 
is dense. 
\item 
If $U \subseteq Y$ is an open dense subset, so is $f^{-1}(U) \subseteq X$. 
\end{enumerate} 

\begin{lem}\label{hc1}
Let $f:X \lra Y$ be a proper morphism of schemes and let 
$U \subseteq Y$ be a dense open subscheme. Then we have the diagram 
consisting of strict morphism of pairs 
\begin{equation}\label{*}
\begin{CD}
(U_X,X) @<{g_X}<< (U_{X'},X') \\ 
@VfVV @V{f'}VV \\ 
(U_Y,Y) @<{g_Y}<< (U_{Y'},Y')
\end{CD}
\end{equation}
with $(U_{Y'},Y') = \coprod_{j=1}^b (U_{Y'_j}, Y'_j)$, 
$(U_{X'},X') = \coprod_{j=1}^b(\coprod_{i=1}^{r_j} (U_{X'_{ji}}, X'_{ji}))$ 
$($decomposition into connected components with $f'(X'_{ji}) \subseteq Y'_j)$ 
satisfying the following conditions$:$ 
\begin{enumerate} 
\item $U_Y$ is contained in $U$ and dense in $Y$. 
\item $X'_{ji}$'s, $Y'_j$'s are integral. 
\item 
The map $\ol{g}_X: (U_{X'},X') \lra (U_X \times_{U_Y} U_{Y'}, X \times_Y Y')$ 
induced by $(g_X,f')$ 
is a proper covering and the map $g_Y$ is a good 
strongly proper covering. 
\item 
For each $i,j$, the generic fiber of 
$f'|_{X'_{ji}}: X'_{ji} \lra Y'_j$ is non-empty and geometrically 
irreducible. $($Attention$:$ It is possible that $r_j=0$ holds for some $j.)$
\end{enumerate} 
\end{lem} 

\begin{pf} 
First let $Y_{\red} = \bigcup_{j=1}^b Y_j$ 
be the decomposition of $Y_{\red}$ into irreducible components and 
let $f_j:X_j \lra Y_j$ be the base change of $f$ by $Y_j \hra 
Y_{\red} \hra Y$. Then each $Y_j$ is integral. Let 
$\eta_j$ be the generic point of $Y_j$ and Let $X_{\eta_j}$ be the generic 
fiber of $f_j$. Then, there exists a finite map $\eta'_j \lra \eta_j$ 
of spectra of fields such that any irreducible component of 
$X_{\eta'_j} := X_{\eta_j} \times_{\eta_j} \eta'_j$ is geometrically 
irreducible. Then, let $Y'_j$ be the normalization of $Y_j$ in 
$\eta'_j$, put $X''_j := X_j \times_{Y_j} Y'_j$ and denote the map 
$X''_j \lra Y'_j$ by $f''_j$. Then $\eta'_j$ is the generic point of 
$Y'_j$ and the generic fiber of $f''_j$ is $X_{\eta'_j}$. Note that 
the generic fiber 
$(X''_{j,\red})_{\eta'_j}$ of $f''_{j,\red}: 
X''_{j,\red} \hra X''_j \os{f''_j}{\lra} Y'_j$ 
is isomorphic to $(X_{\eta'_j})_{\red}$. 
Let $\xi_1, \cdots, \xi_{r_j}$ be the generic points of 
$(X''_{j,\red})_{\eta'_j}$ and let $X'_{ji}\,(1 \leq i \leq r_j)$ 
be the closure of $\xi_i$ in $X''_{j,\red}$. Then there exists a 
closed subscheme $X''_{j0}$ of $X''_{j,\red}$ such that  
$X''_{j,\red} = X''_{j0} \cup \bigcup_{i=1}^{r_j}X'_{ji}$ holds and that 
$f''_{j,\red}(X''_{j0})$ is a proper closed subset of $Y'_j$. 
(Since $Y'_j$ is integral, this implies that $Y'_j-f''_{j,\red}(X''_{j0})$ 
is dense in $Y'_j$.)
\par 
Now let us define $X', Y', f', g_X, g_Y, U_Y$ by 
$$X' := \coprod_{j=1}^b \coprod_{i=1}^{r_j} X'_{ji}, \qquad 
Y' := \coprod_{j=1}^b Y'_j,$$ 
$$f':\coprod_{j=1}^b \coprod_{i=1}^{r_j} X'_{ji} \hra 
\coprod_{j=1}^b X''_{j,\red} \os{\coprod_j f''_{j,\red}}{\lra} 
\coprod_{j=1}^b Y'_j,$$ 
$$g_X: \coprod_{j=1}^b \coprod_{i=1}^{r_j} X'_{ji} \lra 
\coprod_{j=1}^b X''_{j,\red} \hra 
\coprod_{j=1}^b X''_{j} \lra \coprod_j X_j \lra X, $$
$$g_Y: \coprod_{j=1}^b Y'_j \lra \coprod_{j=1}^b Y_j \lra Y, $$ 
$$ U_Y := U \cap (Y - g_Y(\bigcup_{j} f''_{j,\red}(X''_{j0}))) $$ 
and define the diagram \eqref{*} in order that the mophisms in 
\eqref{*} are all strict. Then one can check that 
all the required conditions are satisfied. So we are done. 
\end{pf} 

By using Proposition \ref{dejongup} and Lemma \ref{hc1}, we can prove 
the following proposition: 

\begin{prop}\label{hc2}
Let $f:X \lra Y$ be a proper morphism of schemes and let 
$U \subseteq Y$ be a dense open subscheme. Then we have the diagram 
consisting of strict morphism of pairs 
\begin{equation}\label{**}
\begin{CD}
(U_X,X) @<{g_X}<< (U_{\ti{X}},\ti{X}) \\ 
@VfVV @V{\ti{f}}VV \\ 
(U_Y,Y) @<{g_Y}<< (U_{\ti{Y}},\ti{Y})
\end{CD}
\end{equation}
satisfying the following conditions$:$ 
\begin{enumerate} 
\item $U_Y$ is contained in $U$ and dense in $Y$. 
\item $\ti{Y}$ is regular. 
\item 
The map $\ol{g}_X: (U_{\ti{X}},\ti{X}) 
\lra (U_X \times_{U_Y} U_{\ti{Y}}, X \times_Y \ti{Y})$ 
induced by $(g_X,\ti{f})$ 
is a proper covering and the map $g_Y$ is 
a good strongly proper covering. 
\item 
There exist fine log structures $M_{\ti{X}}, M_{\ti{Y}}$ on $\ti{X}, \ti{Y}$ 
respectively such that $U_{\ti{X}} \subseteq (\ti{X},M_{\ti{X}})_{\triv}, 
U_{\ti{Y}} \subseteq (\ti{Y},M_{\ti{Y}})_{\triv}$ holds and that 
there exists a proper log smooth integral morphism 
having log smooth parameter $(\ti{X},M_{\ti{X}}) \lra (\ti{Y},M_{\ti{Y}})$ 
whose underlying morphism of schemes is the same as $\ti{f}$. 
\end{enumerate}
\end{prop} 

\begin{pf} 
First let us take the diagram 
consisting of strict morphism of pairs 
\begin{equation*}
\begin{CD}
(V_X,X) @<{g_X}<< (U_{X'},X') \\ 
@VfVV @V{f'}VV \\ 
(V_Y,Y) @<{g_Y}<< (U_{Y'},Y')
\end{CD}
\end{equation*}
satisfying the conclusion of Lemma \ref{hc1}, and let 
$Y' = \coprod_{j=1}^b Y'_j$, $X' := 
\coprod_{j=1}^b(\coprod_{i=1}^{r_j} X'_{ji})$ the decomposition of $Y',X'$ 
into connected components with $f'(X'_{ji}) \subseteq Y'_j$. \par 
Let us fix $j$. We prove the following claim: \\
\quad \\
{\bf claim}\,\,\, 
For any $1 \leq s \leq r_j$, there exists a diagram 
consisting of strict morphism of pairs 
\begin{equation}\label{altdiag}
\begin{CD}
\coprod_{i=1}^{r_j}(X'_{ji},O'_{X,ji}) 
@<<< \coprod_{i=1}^{r_j} (X''_{ji},O''_{X,ji}) \\ 
@V{f'}VV @V{f''}VV \\
(Y'_j,O'_{Y,j}) @<<< (Y''_j,O''_{Y,j}) 
\end{CD}
\end{equation} 
satisfying the following conditions: 
\begin{enumerate}
\item 
$O'_{Y,j}$ is dense in $Y'_j$. 
\item 
$Y''_j$ is regular and connected (so it is integral). 
\item 
The upper horizontal arrow is a proper covering 
and the lower horizontal arrow is a good strongly proper covering. 
(This implies that the map $Y''_j \lra Y'_j$ is an alteration.) 
\item 
There exist 
fine log structures $M_{\coprod_{i=1}^sX''_{ji}}, 
M_{Y''_j}$ on $\coprod_{i=1}^sX''_{ji}, Y''_j$ respectively and 
a proper log smooth integral morphism 
having log smooth parameter $g:(\coprod_{i=1}^sX''_{ji}, 
M_{\coprod_{i=1}^sX''_{ji}}) \lra (Y''_j, M_{Y''_j})$ 
whose underlying morphism of schemes is $f''|_{\coprod_{i=1}^s X''_{ji}}$ 
such that $(Y''_j, M_{Y''_j})_{\triv}$ is dense open 
in $Y''_j$ and that 
$g^{-1}((Y''_j, M_{Y''_j})_{\triv})$ is contained in 
$(\coprod_{i=1}^sX''_{ji}, M_{\coprod_{i=1}^sX''_{ji}})_{\triv}$. 
\item 
For each $i \geq s+1$, $X''_{ji}$ is integral and the generic fiber of 
$f''|_{X''_{ji}}: X''_{ji} \lra Y''_j$ is 
non-empty and geometrically irreducible. 
\end{enumerate}

We prove the claim by induction on $s$: In the case $s=1$, 
we apply Proposition \ref{dejongup} to the morphism 
$X'_{j1} \lra Y'_j$ to obtain the diagram 
\begin{equation*}
\begin{CD}
X'_{j1} @<<< X''_{j1} \\ 
@VVV @VVV \\ 
Y'_j @<<< Y''_j 
\end{CD}
\end{equation*}
satisfying the conclusion of Proposition \ref{dejongup}. 
(In particular, $Y''_j$ is regular, connected and we have 
fine log structures $M_{X''_{j1}}, 
M_{Y''_j}$ on $X''_{j1}, Y''_j$ respectively and 
a proper log smooth integral morphism 
having log smooth parameter $g:(X''_{j1}, 
M_{X''_{j1}}) \lra (Y''_j, M_{Y''_j})$ 
whose underlying morphism of schemes is the same as the map 
$X''_{j1} \lra Y''_j$ in the above diagram 
such that $(Y''_j, M_{Y''_j})_{\triv}$ is dense open 
in $Y''_j$ and that 
$g^{-1}((Y''_j, M_{Y''_j})_{\triv})$ is contained in 
$(X''_{j1}, M_{X''_{j1}})_{\triv}$.) On the other hand, for $i \geq 2$, 
the generic fiber of $X'_{ji} \lra Y'_j$ is geometrically irreducible. 
So there exists a unique generic point of 
$(X'_{ji} \times_{Y'_j} Y''_j)_{\red}$ (which we denote by $\xi_{ji}$) 
which is sent to the generic point of $Y''_j$ by the morphism 
$(X'_{ji} \times_{Y'_j} Y''_j)_{\red} \lra Y''_j$. Let $X''_{ji}$ be the 
closure of $\xi_{ji}$ in $(X'_{ji} \times_{Y'_j} Y''_j)_{\red}$. 
Then there exists a closed subset $C_{ji}$ in 
$(X'_{ji} \times_{Y'_j} Y''_j)_{\red}$ with 
$(X'_{ji} \times_{Y'_j} Y''_j)_{\red} = X''_{ji} \cup C_{ji}$ such that 
the image of $C_{ji}$ in $Y''_j$ is a proper closed subset of $Y''_j$. 
Then, if we put 
$$ O'_{Y,j}:=Y'_j - \bigcup_{i=2}^{r_j}\text{(image of $C_{ji}$)} $$ 
and define the diagram \eqref{altdiag} in order that the morphisms are 
all strict, we can check that the diagram satisfies the required 
conditions. So we have proved the claim for $s=1$. \par 
If the claim is true for $s-1$, we have the diagram 
consisting of strict morphisms of pairs 
\begin{equation}
\begin{CD}
\coprod_{i=1}^{r_j} (X'_{ji},O'_{X,ji}) @<<< 
\coprod_{i=1}^{r_j} (X''_{ji},O''_{X,ji}) \\ 
@V{f'}VV @V{f''}VV \\
(Y'_j,O'_{Y',j}) @<<< (Y''_j,O''_{Y,j}), 
\end{CD}
\end{equation} 
fine log structures $M_{\coprod_{i=1}^{s-1}X''_{ji}}, 
M_{Y''_j}$ on $\coprod_{i=1}^{s-1}X''_{ji}, Y''_j$ respectively and 
a proper log smooth integral morphism 
having log smooth parameter $g:(\coprod_{i=1}^{s-1}X''_{ji}, 
M_{\coprod_{i=1}^{s-1}X''_{ji}}) \lra (Y''_j, M_{Y''_j})$ 
satisfying the properties stated in the claim. Then, by Proposition 
\ref{dejongup}, we have the diagram 
\begin{equation*}
\begin{CD}
X''_{js} @<<< X'''_{js} \\ 
@VVV @VVV \\ 
Y''_j @<<< Y'''_j. 
\end{CD}
\end{equation*}
and fine log structures $M^1_{X'''_{js}}, M^1_{Y'''_j}$ on 
$X'''_{js}, Y'''_j$ respectively, 
satisfying the conclusion of 
proposition \ref{dejongup}. (In particular, $Y'''_j$ is regular, connected 
and we have a proper log smooth integral morphism 
having log smooth parameter $g:(X'''_{js}, 
M_{X'''_{js}}) \lra (Y'''_j, M_{Y'''_j})$ 
whose underlying morphism of schemes is the same as the map 
$X'''_{js} \lra Y'''_j$ in the above diagram 
such that $(Y'''_j, M_{Y'''_j})_{\triv}$ is dense open 
in $Y'''_j$ and that 
$g^{-1}((Y'''_j, M_{Y'''_j})_{\triv})$ is contained in 
$(X'''_{js}, M_{X'''_{js}})_{\triv}$.) 
Let $\coprod_{i=1}^{s-1}X'''_{ji} \lra Y'''_j$ be the base change of 
$\coprod_{i=1}^{s-1}X''_{ji} \lra Y''_j$ by $Y'''_j \lra Y''_j$ and 
let $M^2_{Y'''_j}, M^2_{\coprod_{i=1}^{s-1}X'''_{ji}}$ 
be the pull-back of $M_{Y''_j}, M_{\coprod_{i=1}^{s-1}X''_{ji}}$ to 
$Y'''_j, \coprod_{i=1}^{s-1}X'''_{ji}$, respectively. Then we have 
proper log smooth integral morphisms having log smooth parameter 
$$ (X'''_{js},M^1_{X'''_{js}}) \lra (Y'''_j, M^1_{Y'''_j}), $$
$$ (\coprod_{i=1}^{s-1}X'''_{ji}, M^2_{\coprod_{i=1}^{s-1}X'''_{ji}}) 
\lra (Y'''_j, M^2_{Y'''_j}). $$
Moreover, we see that 
$(Y'''_j, M^1_{Y'''_j})_{\triv}, (Y'''_j, M^2_{Y'''_j})_{\triv}$ are 
dense in $Y'''_j$ and that the inverse images of them are contained in 
$(X'''_{js},M^1_{X'''_{js}})_{\triv}, 
(\coprod_{i=1}^{s-1}X'''_{ji}, M^2_{\coprod_{i=1}^{s-1}X'''_{ji}})_{\triv}$ 
respectively. Now let us put 
{\small{
$$ 
(Y'''_j,M_{Y'''_j}) := (Y''',M^1_{Y'''_j}) \times_{Y'''} (Y''',M^2_{Y'''_j}), 
\,\,\,\, 
(X'''_{js},M_{X'''_{js}}) := (X'''_{js},M^1_{X'''_{js}}) 
\times_{Y'''} (Y''',M^2_{Y'''_j}), $$ }}
$$ 
(\coprod_{i=1}^{s-1}X'''_{ji}, M_{\coprod_{i=1}^{s-1}X'''_{ji}}) := 
(Y''',M^1_{Y'''_j}) \times_{Y'''} 
(\coprod_{i=1}^{s-1}X'''_{ji}, M^2_{\coprod_{i=1}^{s-1}X'''_{ji}}). $$
Then we have the morphism 
\begin{equation}\label{logstr-s}
(\coprod_{i=1}^{s-1}X'''_{ji}, M_{\coprod_{i=1}^{s-1}X'''_{ji}})
\coprod (X'''_{js},M_{X'''_{js}}) \lra (Y'''_j, M_{Y'''_j}). 
\end{equation}
On the other hand, for $i \geq s+1$, define the closed subscheme 
$X'''_{ji}, C_{ji}$ with $(X''_{ji} \times_{Y''_j} Y'''_j)_{\red} = 
X'''_{ji} \cup C_{ji}$ as in the case $s=1$ (so that the image of 
$C_{ji}$ in $Y'''_j$ is a proper closed subset). Then we put 
$$ P'_{Y,j} := O'_{Y,j} \cap (Y'_j - \bigcup_{i=s+1}^{r_j}\text{(image of 
$C_{ji}$)}) $$ 
and define the diagram 
\begin{equation*}
\begin{CD}
\coprod_{i=1}^{r_j} (X'_{ji},P'_{X,ji}) @<<< 
\coprod_{i=1}^{r_j} (X'''_{ji},P'''_{X,ji}) \\ 
@V{f'}VV @V{f''}VV \\
(Y'_j,P'_{Y,j}) @<<< (Y'''_j,P'''_{Y,j}) 
\end{CD}
\end{equation*} 
in order that the morphisms of pairs are all strict. Then this diagram, 
together with the morphism \eqref{logstr-s} of log schemes, satisfies the 
required condition for $s$. So the claim is proved. \par 
By the claim for $s=r_j$ (for each $j$), we see that there exists a diagram 
consisting of strict morphism of pairs 
\begin{equation*}
\begin{CD}
(X',O'_X) @<<< (\ti{X},\ti{O}_{X}) \\ 
@V{f'}VV @V{\ti{f}}VV \\ 
(Y',O'_Y) @<<< (\ti{Y},\ti{O}_Y) 
\end{CD}
\end{equation*}
and fine log structures $M_{\ti{X}}, M_{\ti{Y}}$ on $\ti{X}, \ti{Y}$ 
respectively, satisfying the following conditions: 
\begin{enumerate}
\item 
$O'_Y \subseteq Y'$ is dense open. 
\item 
$\ti{Y}$ is regular (hence it is a disjoint union of integral schemes). 
\item 
The upper horizontal arrow is a proper covering 
and the lower horizontal arrow is a good strongly proper covering. 
(This implies that the morphism $\ti{Y} \lra Y'$ is a disjoint union of 
alterations.) 
\item  
There exists a proper log smooth integral morphism 
having log smooth parameter $(\ti{X},M_{\ti{X}}) \lra (\ti{Y},M_{\ti{Y}})$ 
whose underlying morphism of schemes is the same as $\ti{f}$ such that 
$(\ti{Y},M_{\ti{Y}})_{\triv}$ is dense in $\ti{Y}$ and 
$\ti{f}^{-1}((\ti{Y},M_{\ti{Y}})_{\triv}) \subseteq 
(\ti{X},M_{\ti{X}})_{\triv}$ holds. 
\end{enumerate}

Now let $\eta_j \,(1 \leq j \leq b)$ be the generic point of $Y'_j$, 
let $\xi_{j1}, \cdots, \xi_{jr_j}$ 
be the generic points of $X'$ over $\eta_j$ 
and let $\eta'_j$ be the generic point of $\ti{Y}$ over $\eta_j$. 
By the condition (3) in Lemma \ref{hc1}, 
$\xi_{ji} \times_{\eta_j} \eta'_j$ is irreducible. 
So $\{\xi_{ji} \times_{\eta_j} \eta'_j\}_{j,i}$ is 
bijective to the set of generic points of $X' \times_{Y'} \ti{Y}$ 
which are sent to one of the generic points of $\ti{Y}$ by second 
projection. Since $\ti{O}_{X} \lra O'_X$ is surjective and 
$\xi_{ji}$'s are contained in $O'_X$, 
there exists a point $\ti{\xi}_{ji}$ of $\ti{X}$ whose 
image in $X'$ is $\xi_{ji}$. Then $\ti{\xi}_{ji}$ is sent to 
$\eta_j$ by $\ti{X} \lra \ti{Y} \lra Y'$. By the property (3) above, 
$\eta'_j$ is the unique point of $\ti{Y}$ 
which is sent to $\eta_j \in Y'$. So $\ti{\xi}_{ji}$ is sent to 
$\eta'_j$ by $\ti{X} \lra \ti{Y}$. Hence $\ti{\xi}_{ji}$ is 
sent to $\xi_{ji} \times_{\eta_j} \eta'_j$ by the map 
$\ti{X} \lra X' \times_{Y'} \ti{Y}$. So 
the image of the map 
$\ti{X} \lra X' \times_{Y'} \ti{Y}$ contains the closure of 
$\{\xi_{ji} \times_{\eta_j} \eta'_j\}_{j,i}$. Therefore, we have 
a closed subset $C$ of $X' \times_{Y'} \ti{Y}$ such that 
$X' \times_{Y'} \ti{Y} = \text{(image of $\ti{X}$)} \cup C$ holds and that 
the image of $C$ in $\ti{Y}$ (which we denote by $Z_1$) is a 
closed subset of $\ti{Y}$ such that $\ti{Y}-Z_1$ is dense in $\ti{Y}$. \par 
Now let $Z_2 := \ti{Y}-(\ti{Y},M_{\ti{Y}})_{\triv}$, 
define $U_Y \subseteq Y$ by 
$$ U_Y := Y - ((Y-V_Y) \cup \text{(image of $Z_1 \cup Z_2$)} \cup 
\text{(image of $Y'-O'_{Y'}$)}), $$
and define the diagram \eqref{**} in order that the mophisms in 
\eqref{**} are all strict. Then one can check that 
all the required conditions are satisfied. So we are done. 
\end{pf} 

Now we prove our main theorem in this section: 

\begin{thm}\label{hc3}
Let $f:X \lra Y$ be a proper morphism of schemes and let $q \in \N$. 
Then we have the diagram 
consisting of strict morphism of pairs 
\begin{equation}\label{***}
\begin{CD}
(U_X,X) @<<< (U_{\ti{X}},\ti{X}) @<h<< 
(U_{\ti{X}^{(\b)}}, \ti{X}^{(\b)})\\ 
@VfVV @V{\ti{f}}VV @V{\ti{f}^{(\b)}}VV \\ 
(U_Y,Y) @<{g_Y}<< (U_{\ti{Y}},\ti{Y}) @= (U_{\ti{Y}},\ti{Y})
\end{CD}
\end{equation}
satisfying the following conditions$:$ 
\begin{enumerate} 
\item $U_Y$ is dense in $Y$.
\item $\ti{Y}$ is regular. 
\item 
The left square is Cartesian, $g_Y$ is a good strongly proper covering 
and $h$ is a $q$-truncated proper hypercovering by a 
$q$-truncated split simplicial pair. 
\item 
There exist fine log structures 
$M_{\ti{X}^{(n)}}, M_{\ti{Y}}$ on $\ti{X}^{(n)}, \ti{Y}$ 
respectively $(n \leq q)$ 
such that $U_{\ti{X}^{(n)}} \subseteq (\ti{X}^{(n)},M_{\ti{X}^{(n)}})_{\triv}, 
U_{\ti{Y}} \subseteq (\ti{Y},M_{\ti{Y}})_{\triv}$ holds and that 
for each $n \leq q$, there exists a proper log smooth integral morphism 
having log smooth parameter 
$(\ti{X}^{(n)},M_{\ti{X}^{(n)}}) \lra (\ti{Y},M_{\ti{Y}})$ 
whose underlying morphism of schemes is the same as $\ti{f}^{(n)}$. 
\end{enumerate}
\end{thm} 

\begin{pf} 
We prove the theorem by induction on $q$. When $q=0$ holds, let us take 
a diagram 
\begin{equation*}
\begin{CD}
(U_X,X) @<<< (U_{\ti{X}^{(0)}},\ti{X}^{(0)}) \\ 
@VfVV @V{\ti{f}^{(0)}}VV \\ 
(U_Y,Y) @<{g_Y}<< (U_{\ti{Y}},\ti{Y})
\end{CD}
\end{equation*}
satisfying the conclusion of Proposition \ref{hc2}. 
Then, if we put $(U_{\ti{X}},\ti{X}) := (U_X \times_{U_Y} U_{\ti{Y}}, 
X \times_Y \ti{Y})$, we obtain the required diagram. \par 
Now let us assume that the theorem is true for $q$, and take the 
diagram 
\begin{equation*}
\begin{CD}
(U_X,X) @<<< (U_{X'},X') @<h<< 
(U_{{X'}^{(\b)}}, {X'}^{(\b)})\\ 
@VfVV @V{f'}VV @V{{f'}^{(\b)}}VV \\ 
(U_Y,Y) @<{g_Y}<< (U_{Y'},Y') @= (U_{Y'},Y')
\end{CD}
\end{equation*}
satisfying the conclusion of the theorem for $q$. Then put 
$C := {\mathrm{cosk}}_q^{X'}({X'}^{(\b)})^{(q+1)}$ and 
starting from the morphism $C \lra Y'$, let us take a diagram 
\begin{equation*}
\begin{CD}
(V_C,C) @<<< (V_{\ti{X}^{(q+1)'}},\ti{X}^{(q+1)'}) \\ 
@VVV @VVV \\ 
(V_{Y'},Y') @<<< (V_{\ti{Y}},\ti{Y})
\end{CD}
\end{equation*}
satisfying the conclusion of Proposition \ref{hc2}. 
Now let us put 
$\ti{X}:= X' \times_{Y'} \ti{Y}, 
\ti{X}^{(\b)} := {X'}^{(\b)} \times_{Y'} \ti{Y}\,(\b\leq q)$, 
$\ti{C}:= C \times_{Y'} \ti{Y} = 
{\mathrm{cosk}}_q^{\ti{X}}({\ti{X}}^{(\b)})^{(q+1)}$. 
Let us replace $U_Y$ by 
$Y-(\text{(image of $Y'-V_{Y'}$)}\cup (Y-U_Y))$, 
define the diagram 
\begin{equation}\label{****}
\begin{CD}
(U_X,X) @<<< (U_{\ti{X}},\ti{X}) @<h<< 
(U_{\ti{X}^{(\b)}}, \ti{X}^{(\b)})\\ 
@VfVV @V{\ti{f}}VV @V{\ti{f}^{(\b)}}VV \\ 
(U_Y,Y) @<{g_Y}<< (U_{\ti{Y}},\ti{Y}) @= (U_{\ti{Y}},\ti{Y})
\end{CD}
\end{equation}
in order that the mophisms in 
the diagram are all strict, and define 
$U_{\ti{C}}, U_{\ti{X}^{(q+1)'}}$ by 
$U_{\ti{C}} := \ti{C} \times_{Y} U_Y, 
U_{\ti{X}^{(q+1)'}} := \ti{X}^{(q+1)'} \times_Y U_Y$. 
Then the diagram \eqref{****} satisfies the required condition for 
$q$ and we have a proper covering 
$(U_{\ti{X}^{(q+1)'}}, \ti{X}^{(q+1)'}) \lra (U_{\ti{C}},\ti{C})$ over 
$(U_{\ti{Y}},\ti{Y})$. Then, by using the recipe of 
\cite[5.1.3]{saintdonat}($=$\cite[7.2.5]{tsuzuki4}), 
we can form a $(q+1)$-truncated split simplicial pair 
$(U_{\ti{X}^{(\b)}}, \ti{X}^{(\b)})$ which fits into the diagram like 
\eqref{****} ($\b$ runs through $0 \leq \b \leq q+1$ now), and we see that 
it satisfies the required condition for $q+1$. So the proof of 
the theorem is finished. 
\end{pf} 


\section{Relative rigid cohomology (II): general case}

In this section, we prove (a version of) 
Berthelot's conjecture on the overconvergence of relative rigid 
cohomology of proper smooth morphisms (Conjecture \ref{bconj2}) 
under mild assumption. 
We also prove the generic overconvergence 
for proper morphisms which are not necessarily smooth. \par 
In this section, $B, \cB$ are as in the previous section: Namely, 
$B = \Spec k$ and $\cB = \Spf W(k)$, where $k$ is a perfect field of 
characteristic $p$. \par 
Before the proof, we need to prove the descent property of the category of 
overconvergent isocrystals for strongly proper \v{C}ech hypercovering. 
First we prove a lemma, which is essentially due to Ogus \cite[4.15]{ogus1}: 

\begin{lem}\label{proj1}
Let $N,a$ be positive integers and 
let us assume given the following diagram 
\begin{equation}\label{diag-proj}
\begin{CD}
\ol{Y}' @>>> (\Pr^N_{\cY})^a \\ 
@VfVV @VgVV \\ 
\ol{Y} @>>> \cY, 
\end{CD}
\end{equation}
where the horizontal arrows are closed immersions, $f$ is a 
projective surjective morphism between reduced $B$-schemes, 
$\cY$ is a flat formal $\cB$-scheme, $(\Pr^N_{\cY})^a$ is the $a$-fold 
fiber product of $\Pr^N_{\cY}$ over $\cY$ and $g$ is the natural morphism. 
Then, for any $m>0$, there exists $n>0$ such that, for any $n' \geq n$, 
the morphism $(g_K)^{-1}(T_{m,1}(\cY)_K) \cap T_{n',1}((\Pr^N_{\cY})^a)_K 
\lra T_{m,1}(\cY)_K$ induced by $g_K$ is surjective. 
\end{lem} 

\begin{pf} 
We prove the lemma by induction on $a$. In the case $a=1$, 
the lemma is proven in \cite[4.15]{ogus1}. In general case, we factor the 
diagram \eqref{diag-proj} as follows: 
\begin{equation*}
\begin{CD}
\ol{Y}' @>>> (\Pr^N_{\cY})^a \\ 
@V{f_1}VV @V{g_1}VV \\
\ol{Y}'' @>>> (\Pr^N_{\cY})^{a-1} \\
@V{f_2}VV @V{g_2}VV \\ 
\ol{Y} @>>> \cY, 
\end{CD}
\end{equation*}
where $\ol{Y}''$ denotes the image of $\ol{Y}' \ra (\Pr^N_{\cY})^a \ra 
(\Pr^N_{\cY})^{a-1}$. Then we can use the induction hypothesis: 
So, for any $m>0$, there exist $k$ and $n$ such that for any $n' \geq n$, 
we have surjections 
$$ 
(g_{2,K}^{-1})(T_{m,1}(\cY)_K) \cap T_{k,1}((\Pr^N_{\cY})^{a-1}_K) \lra 
T_{m,1}(\cY)_K, $$ 
$$ 
(g_{1,K}^{-1})(T_{k,1}(\Pr^N_{\cY})^{a-1}_K) 
\cap T_{n',1}((\Pr^N_{\cY})^{a}_K) \lra 
T_{k,1}((\Pr^N_{\cY})^{a-1})_K. $$ 
From these surjections, we obtain the assertion. 
\end{pf} 

\begin{cor}\label{oguscor}
Let $N,a$ be positive integers and 
let us assume given the following diagram 
\begin{equation*}
\begin{CD}
\ol{Y}' @>>> (\Pr^N_{\cY})^a \\ 
@VfVV @VgVV \\ 
\ol{Y} @>>> \cY, 
\end{CD}
\end{equation*}
where the upper horizontal arrow is a closed immersion, 
the lower horizontal arrow is a homeomorphic closed immersion, 
$f$ is a projective surjective morphism between $B$-schemes, 
$\cY$ is a formal $\cB$-scheme, $(\Pr^N_{\cY})^a$ is the $a$-fold 
fiber product of $\Pr^N_{\cY}$ over $\cY$ and $g$ is the natural morphism. 
Then there exists $n>0$ such that, for any $n' \geq n$, 
the morphism $T_{n',1}((\Pr^N_{\cY})^a)_K \lra 
\cY_K$ induced by $g_K$ is surjective. 
\end{cor} 

\begin{pf} 
By dividing $\cO_{\cY}$ by the ideal of $p$-torsions, we may assume that 
$\cY$ is flat over $\cB$. (This does not change the morphism $g_K$.) 
Moreover, we may replace $f$ by $f_{\red}:\ol{Y}'_{\red} \lra \ol{Y}_{\red}$ 
because the inductive system of 
admissible open sets $\{T_{n,1}((\Pr^N_{\cY})^a)_K\}_n$ are unchanged 
up to canonical isomorphism as inductive systems. Then the corollary 
is a special case of Lemma \ref{proj1}. 
\end{pf} 

Now we prove the descent property of the category of 
overconvergent isocrystals for strongly proper \v{C}ech hypercovering. 

\begin{prop}\label{properdescent}
Let $(Y,\ol{Y}) \os{g^{(\b)}}{\lla} (Y^{(\b)},\ol{Y}^{(\b)})$ 
be a strongly proper 
\v{C}ech hypercovering of pairs over a $B$-scheme $S$ and let 
$S \hra \cS$ be a closed immersion into a formal $\cB$-scheme $\cS$. 
Then, if we denote the category of descent data with respect to 
$I^{\d}((Y^{(n)},\ol{Y}^{(n)})/\cS_K) \, (n=0,1,2)$ by 
$I^{\d}((Y^{(\b)},\ol{Y}^{(\b)})/\cS_K)$, the restriction functor 
$$ I^{\d}((Y,\ol{Y})/\cS_K) \lra I^{\d}((Y^{(\b)},\ol{Y}^{(\b)})/\cS_K)$$ 
is an equivalence of categories. 
\end{prop} 

\begin{pf} 
Since we may work locally on $\ol{Y}$, we may assume that $\ol{Y}$ admits 
a closed immersion $i: \ol{Y} \hra \cY$ such that $\cY$ is formally smooth 
over $\cS$. By Chow's lemma, we can take a morphism 
$h:\ol{Z} \lra \ol{Y}^{(0)}$ such that both $h$ and $g^{(0)} \circ h$ 
are projective surjective morphisms. 
Since the proposition for the strongly proper \v{C}ech hypercoverings 
induced by $h$ and $g^{(0)} \circ h$ implies the proposition for $g^{(\b)}$, 
we may assume that $g^{(0)}$ is projective to prove the 
proposition. So we may assume that there exists a closed immersion 
$i^{(0)}:\ol{Y}^{(0)} \hra \Pr^N_{\cY}$ over $i$. Then, by taking 
fiber products, we obtain the following diagram: 
\begin{equation}\label{prodes1}
\begin{CD}
\ol{Y}^{(\b)} @>{i^{(\b)}}>> (\Pr^N_{\cY})^{\b+1} \\
@V{g^{(\b)}}VV @VVV \\ 
\ol{Y} @>i>> \cY. 
\end{CD}
\end{equation}
Let us denote the category of descent data for 
$I^{\d}((Y^{(n)},\ol{Y}^{(n)})/\cS, (\Pr^N_{\cY})^{n+1}) \, (n=0,1,2)$ by 
$I^{\d}((Y^{(\b)},\ol{Y}^{(\b)})/\cS, (\Pr^N_{\cY})^{\b+1})$. 
Then it suffices to prove that the restriction functor 
$$ I^{\d}((Y,\ol{Y})/\cS_K, \cY) \lra 
I^{\d}((Y^{(\b)},\ol{Y}^{(\b)})/\cS_K, (\Pr^N_{\cY})^{\b+1}) $$ 
is an equivalence of categories. \par 
For $m \in \N$, let $\cY^m$ be the $m$-fold fiber product of $\cY$ over 
$\cS$. Then, by taking fiber products, we obtain the following diagram 
from the diagram \eqref{prodes1}, where $\cZ = \cY^m$: 
\begin{equation}\label{prodes2}
\begin{CD}
\ol{Y}^{(\b)} @>{i^{(\b)}}>> (\Pr^N_{\cZ})^{m(\b+1)} \\
@V{g^{(\b)}}VV @VVV \\ 
\ol{Y} @>i>> \cZ. 
\end{CD}
\end{equation}
Since the category 
$I^{\d}((Y,\ol{Y})/\cS_K, \cY)$ (resp. 
$I^{\d}((Y^{(\b)},\ol{Y}^{(\b)})/\cS_K, (\Pr^N_{\cY})^{\b+1})$) 
is defined as the category of coherent 
$j^{\d}\cO_{]\ol{Y}[_{\cY}}$-modules (resp. coherent 
$j^{\d}\cO_{]\ol{Y}^{(\b)}[_{(\Pr^N_{\cY})^{\b+1}}}$-modules) 
endowed with an isomorphism between pull-backs on 
$]\ol{Y}[_{\cY^2}$ (resp. $]\ol{Y}^{(\b)}[_{(\Pr^N_{\cY^2})^{2(\b+1)}}$)
 satisfying the cocycle condition on 
$]\ol{Y}[_{\cY^3}$ (resp. $]\ol{Y}^{(\b)}[_{(\Pr^N_{\cY^3})^{3(\b+1)}}$). 
So, to prove the proposition, it suffices to prove the descent 
for coherent modules, that is, it suffices to prove that the restriction 
functor 
$$ \Coh(j^{\d}\cO_{]\ol{Y}[_{\cZ}}) \lra 
\Coh(j^{\d}\cO_{]\ol{Y}^{(\b)}[_{(\Pr^N_{\cZ})^{m(\b+1)}}}) $$ 
is an equivalence of categories for $\cZ = \cY^m$, 
where $\Coh(j^{\d}\cO_{]\ol{Y}^{(\b)}[_{(\Pr^N_{\cZ})^{m(\b+1)}}})$ denotes 
the category of descent data with respect to 
$\Coh(j^{\d}\cO_{]\ol{Y}^{(n)}[_{(\Pr^N_{\cZ})^{m(n+1)}}}) \, (n=0,1,2)$. 
To prove this, it suffices to prove the same assertion after we pull-back 
the diagram \eqref{prodes2} by $T_{l,1}(\cZ) \lra \cZ \, (l \in \N)$. 
So we may assume that $i$ is a homeomorphic closed immersion. 
(But we do not assume anymore that $\cZ$ is equal to $\cY^m$.) 
Moreover, we may assume that $\cZ$ is affine. \par 
By Corollary \ref{oguscor}, there exists a surjection 
$T_{k,1}((\Pr^N_{\cZ})^{m})_K \lra \cZ_K$ for some $k$. 
Let $T_{k,1}((\Pr^N_{\cZ})^{m})^{\b}$ be the $\b$-fold fiber 
product of $T_{k,1}((\Pr^N_{\cZ})^{m})$ over $\cZ$. Then 
$T_{k,1}((\Pr^N_{\cZ})^{m})^{\b+1}_K$ is an admissible open set of 
$]\ol{Y}^{(\b)}[_{(\Pr^N_{\cZ})^{m(\b+1)}}$. 
Let $g_1, \cdots, g_c \in \Gamma(\cZ,\cO_{\cZ})$ 
be a lift of generators of $\Gamma(\ol{Y},\Ker(\cO_{\ol{Y}} \ra 
\cO_{\ol{Y}-Y}))$, and for $\nu \in p^{\Q_{<0}}$, let us define the 
admissible open set $\cZ_{K,\nu}$ (resp. 
$T_{k,1}((\Pr^N_{\cZ})^{m})^{\b}_{K,\nu}$) of 
$\cZ_K$ (resp. $T_{k,1}((\Pr^N_{\cZ})^{m})^{\b}_{K}$) as the set 
$\bigcup_{j=1}^c \{|g_j| \geq \nu\}$. \par 
Now assume that we are given an object $E^{(\b)}$ in 
$\Coh(j^{\d}\cO_{]\ol{Y}^{(\b)}[_{(\Pr^N_{\cZ})^{m(\b+1)}}})$. 
Then, for $\nu$ sufficiently close to $1$, $E^{(\b)}$ defines 
a compatible family (with respect to $0 \leq \b \leq 2$) of 
coherent modules on $T_{k,1}((\Pr^N_{\cZ})^{m})^{\b+1}_{K,\nu}$. 
To show the assertion, it suffices to prove that it descents to 
a coherent module on $\cY_{\nu,K}$. That is, it suffices to prove the 
descent property of coherent modules for the map 
\begin{equation}\label{a}
T_{k,1}((\Pr^N_{\cZ})^{m})^{\b+1}_{K,\nu} \lra \cZ_{K,\nu}. 
\end{equation}
Note that this map is obtained as the pull-back 
from the \v{C}ech covering induced by the 
analytically flat, analytically 
surjective map of formal schemes 
$T_{k,1}((\Pr^N_{\cZ})^{m}) \lra \cZ$. So, by any affinoid admissible 
open set $\cQ_K$ of $\cZ_{K,\nu}$, the pull-back of the map \eqref{a} by 
$\cQ_K \allowbreak \lra \allowbreak \cZ_{K,\nu}$ is the \v{C}ech covering induced by an 
analytically flat, analytically surjective map of formal schemes. 
So the descent of coherent module for this map is true due to 
the rigid analytic faithfully flat descent of Gabber. So we have finished 
the proof of the proposition. 
\end{pf} 

Now we prove our main theorem on 
the generic overconvergence for proper morphisms. 

\begin{thm}\label{theorem1}
Take a triple of the form 
$(S,S,\cS)$ and let us assume given the diagram 
\begin{equation}\label{maindiag}
(X,\ol{X}) \os{f}{\lra} (Y,\ol{Y}) \os{g}{\lra} (S,S), 
\end{equation}
where $f$ is a strict proper morphism of pairs and 
$g$ is a morphism of pairs. 
Then, there exists 
a non-empty open set $U_Y$ of $Y$, a diagram 
$$ \ol{Y} \os{g_Y}{\lla} \ti{Y} \os{i}{\hra} \cY $$
$($where $g_Y$ is a composition of an etale surjective map and a proper 
surjective map and $i$ is a closed immersion into a $p$-adic formal scheme 
formally smooth over $\cS)$, such that, 
for any locally free overconvergent isocrystal $\cE$ on 
$(X,\ol{X})/\cS$ and $q \in \N$, there exists 
the unique overconvergent isocrystal $\cF$ on $(U_Y,\ol{Y})/\cS_K$ 
satisfying the following condition$:$ 
For any triple $(Z,\ol{Z},\cZ)$ over $(U_{\ti{Y}},\ti{Y},\cY)$ 
$($where $U_{\ti{Y}}:=U_Y \times_{\ol{Y}} \ti{Y})$ satisfying 
$Z = U_{\ti{Y}} \times_{\ti{Y}} \ol{Z}$ with $\cZ$ formally smooth over 
$\cY$, the restriction of $\cF$ to 
$I^{\d}((Z,\ol{Z})/\cS_K, \cZ)$ is given functorially by 
$(R^qf_{(X \times_{Y} Z, \ol{X} \times_{\ol{Y}} \ol{Z})/\cZ, \rig *}\cE, 
\epsilon)$, where $\epsilon$ is an isomorphism 
{\small{ $$ 
p_2^*R^qf_{(X \times_{Y} Z, 
\ol{X} \times_{\ol{Y}} \ol{Z})/\cZ, \rig *}\cE 
\os{\sim}{\ra} 
R^qf_{(X \times_{Y} Z, \ol{X} \times_{\ol{Y}} \ol{Z})/\cZ \times_{\cS} \cZ, 
\rig *}\cE \os{\sim}{\leftarrow}
p_1^*R^qf_{(X 
\times_{Y} Z, \ol{X} \times_{\ol{Y}} \ol{Z})/\cZ, \rig *}\cE
$$ }}
$(p_i$ denotes the $i$-th projection 
$]\ol{Z}[_{\cZ \times_{\cS} \cZ} \lra ]\ol{Z}[_{\cZ})$. 
\end{thm} 

\begin{pf} 
First we prove the following boundedness property 
(this is essentially the same as \cite[6.4.1]{tsuzuki4}): 
There exists a positive integer $q_0$ such that, for any overconvergent 
isocrystal $\cE$ on $(X,\ol{X})/\cS_K$ and any $(Y,\ol{Y})$-triple 
$(Z,\ol{Z},\cZ)$ over $(S,S,\cS)$, we have 
\begin{equation}\label{bound-rig}
R^qf_{(X \times_Y Z,\ol{X} \times_{\ol{Y}} \ol{Z})/\cZ, \rig *}\cE = 0
\end{equation} 
for any $q > q_0$. 
To prove this, it suffices to check the equality \eqref{bound-rig} 
(for $q$ sufficiently large) 
in the case where $\ol{Y}$ and $\cZ$ are affine, and by 
\cite[6.3.2]{tsuzuki4}, we may assume that $\ol{X}$ is affine. Then there 
exists a closed immersion $\ol{X} \times_{\ol{Y}} \ol{Z} \hra \cP$ such that 
$\cP$ is formally smooth over $\cZ$. Let us denote the morphism 
$]\ol{X} \times_{\ol{Y}} \ol{Z}[_{\cP} \,\lra\, ]\ol{Z}[_\cZ$ by $\varphi$. 
Then we have the spectral sequence 
$$ 
E_1^{s,t} = 
R^t\varphi_*(E \otimes_{j^{\d}\cO_{]\ol{X} \times_{\ol{Y}} \ol{Z}[_{\cP}}} 
j^{\d}\Omega^s_{]\ol{X} \times_{\ol{Y}} \ol{Z}[_{\cP}/]\ol{Z}[_{\cZ}}) 
\,\Longrightarrow\, 
R^{s+t}f_{(X \times_Y Z,\ol{X} \times_{\ol{Y}} \ol{Z})/\cZ, \rig *}\cE, 
$$ 
where $E$ denotes the coherent 
$j^{\d}\cO_{]\ol{X} \times_{\ol{Y}} \ol{Z}[_{\cP}}$-module associated to
 $\cE$. Then, by \cite[3.2.2 (i)(ii)]{tsuzuki4}, the $E_1^{s,t}$-term is 
 equal to zero for $t>1$ or sufficiently large $s$. 
So we have the claim. By this claim, 
it suffices to treat only the case $0 \leq q \leq q_0$ to prove the 
theorem. (So we assume this.) \par 
Put $N:= q_0(q_0+1)/2$ and 
let us take a diagram 
consisting of strict morphism of pairs 
\begin{equation*}
\begin{CD}
(U_X,\ol{X}) @<<< (U_{\ti{X}'},\ti{X}') @<{h'}<< 
(U_{\ti{X}^{(\b)'}}, \ti{X}^{(\b)'})\\ 
@VfVV @V{\ti{f}}VV @V{\ti{f}^{(\b)'}}VV \\ 
(U_Y,\ol{Y}) @<{g'_Y}<< (U_{\ti{Y}'},\ti{Y}') @= (U_{\ti{Y}'},\ti{Y}')
\end{CD}
\end{equation*}
satisfying the conclusion of Theorem \ref{hc3} (but $q$ is 
replaced by $N$): Then 
there exist fine log structures 
$M_{\ti{X}^{(n)'}}, M_{\ti{Y}'}$ on $\ti{X}^{(n)'}, \ti{Y}'$ 
respectively $(n \leq N$) 
such that $U_{\ti{X}^{(n)'}} \subseteq 
(\ti{X}^{(n)'},M_{\ti{X}^{(n)'}})_{\triv}, 
U_{\ti{Y}} \subseteq (\ti{Y}',M_{\ti{Y}'})_{\triv}$ holds and that 
for each $n \leq N$, there exists a proper log smooth integral morphism 
having log smooth parameter 
$(\ti{X}^{(n)'},M_{\ti{X}^{(n)'}}) \lra (\ti{Y}',M_{\ti{Y}'})$ 
whose underlying morphism of schemes is the same as $\ti{f}^{(n)'}$. 
Then let us take a strict etale covering 
$g''_Y: (\ti{Y},M_{\ti{Y}}) \lra (\ti{Y}',M_{\ti{Y}'})$ and a 
closed immersion $i: (\ti{Y},M_{\ti{Y}}) \hra (\cY,M_{\cY})$ such that 
the inverse image $U_{\ti{Y}}$ of $U_Y$ in $\ti{Y}$ is contained in 
$(\cY,M_{\cY})_{\triv}$, $(\cY,M_{\cY})$ is formally log smooth over 
$\cS$ and $\cY$ is formally smooth over $\cS$. 
Put $g_Y:=g''_Y \circ g'_Y$ and 
let us denote the base change of the diagram 
\begin{equation*}
\begin{CD}
(U_{\ti{X}'},\ti{X}') @<{h'}<< (U_{\ti{X}^{(\b)'}}, \ti{X}^{(\b)'})\\ 
@V{\ti{f}}VV @V{\ti{f}^{(\b)'}}VV \\ 
(U_{\ti{Y}'},\ti{Y}') @= (U_{\ti{Y}'},\ti{Y}')
\end{CD}
\end{equation*}
by $g''_Y: (U_{\ti{Y}},\ti{Y}) \lra (U_{\ti{Y}'},\ti{Y}')$ by 
\begin{equation*}
\begin{CD}
(U_{\ti{X}},\ti{X}) @<<< (U_{\ti{X}^{(\b)}}, \ti{X}^{(\b)})\\ 
@V{\ti{f}}VV @V{\ti{f}^{(\b)'}}VV \\ 
(U_{\ti{Y}},\ti{Y}) @= (U_{\ti{Y}},\ti{Y}). 
\end{CD}
\end{equation*}
Then we obtain the diagram 
consisting of strict morphism of pairs 
\begin{equation}\label{reduce}
\begin{CD}
(U_X,\ol{X}) @<<< (U_{\ti{X}},\ti{X}) @<<< 
(U_{\ti{X}^{(\b)}}, \ti{X}^{(\b)})\\ 
@VfVV @V{\ti{f}}VV @V{\ti{f}^{(\b)}}VV \\ 
(U_Y,\ol{Y}) @<{g_Y}<< (U_{\ti{Y}},\ti{Y}) @= (U_{\ti{Y}},\ti{Y})
\end{CD}
\end{equation}
satisfying the conditions in Theorem \ref{hc3} ($q$ is replaced by $N$) 
except that $g_Y$ is now a composition of 
an etale surjective map and a proper surjective map. 
Let us replace $U_Y$ by $Y \cap U_Y$ and make the diagram \eqref{reduce} 
in order that all the morphisms are strict. 
(Then $U_Y$ is no more dense in $\ti{Y}$ but only non-empty.) 
Let us define the simplicial scheme $\dti{X}^{(\b)}$ by 
$\dti{X}^{(\b)} := \cosk_N^{\ti{X}}(\ti{X}^{(\b)})$ 
and let $U_{\dti{X}^{(\b)}}$ be the inverse image of $U_Y$ in 
$\dti{X}^{(\b)}$. (Then the $N$-truncation of 
$(U_{\dti{X}^{(\b)}}, \dti{X}^{(\b)})$ is nothing but 
$(U_{\ti{X}^{(\b)}},\ti{X}^{(\b)})$.) Then, since 
$(U_{\dti{X}^{(\b)}}, \dti{X}^{(\b)}) \lra (U_{\ti{X}},\ti{X})$ is 
a proper hypercovering, we have the spectral sequence 
of relative rigid cohomology 
\begin{equation}\label{spseq}
E_1^{s,t} := R^tf_{(U_{\dti{X}^{(s)}},\dti{X}^{(s)})/\cY, \rig *}\cE 
\,\Longrightarrow\, 
R^{s+t}f_{(U_{\ti{X}},\ti{X})/\cY, \rig *}\cE 
\end{equation} 
established by Tsuzuki (see \cite{tsuzuki4}). 
Then, for $s \leq N$, the $E_1^{s,t}$-term 
$R^tf_{(U_{\dti{X}^{(s)}},\dti{X}^{(s)})/\cY, \rig *}\allowbreak \cE = 
R^tf_{(U_{\ti{X}^{(s)}},\ti{X}^{(s)})/\cY, \rig *}\cE$ is a coherent 
$j^{\d}\cO_{]\ti{Y}[_{\cY}}$-module (where $j$ is the admissible 
open immersion $]U_{\ti{Y}}[_{\cY} \hra ]\ti{Y}[_{\cY}$) by 
Theorem \ref{main1}. So, by the spectral sequence \eqref{spseq}, we see
 that $R^{q}f_{(U_{\ti{X}},\ti{X})/\cY, \rig *}\cE$ is also 
a coherent $j^{\d}\cO_{]\ti{Y}[_{\cY}}$-module. 
By the same reason, $R^{q}f_{(U_{\ti{X}},\ti{X})/\cY \times_{\cS} \cY, 
\rig *}\cE$ is also a coherent $j^{\d}\cO_{]\ti{Y}[_{\cY}}$-module. 
Then, by the base change theorem by Tsuzuki 
\cite[2.3.1]{tsuzuki3}, the morphisms 
$p_i^*R^qf_{(U_{\ti{X}}, \ti{X})/\cY, \rig *}\cE \allowbreak \lra 
R^qf_{(U_{\ti{X}}, \ti{X})/\cY \times_{\cS} \cY, \rig *}\cE$ 
for $i=1,2$ are isomorphisms, where $p_i$ denotes the $i$-th projection 
$]\ti{Y}[_{\cY \times_{\cS} \cY} \lra 
]\ti{Y}[_{\cY}$. if we define $\epsilon^{\lor}$ to be the 
isomorphism 
$$ 
p_2^*R^qf_{(U_{\ti{X}}, \ti{X})/\cY, \rig *}\cE 
\os{=}{\ra}
R^qf_{(U_{\ti{X}}, \ti{X})/\cY \times_{\cS} \cY, \rig *}\cE
\os{=}{\leftarrow}
p_1^*R^qf_{(U_{\ti{X}}, \ti{X})/\cY, \rig *}\cE, 
$$ 
then $\cF^{\lor}:= 
(R^qf_{(U_{\ti{X}}, \ti{X})/\cY, \rig *}\cE, \epsilon^{\lor})$ 
defines an overconvergent isocrystal on $(U_{\ti{Y}},\ti{Y})/\allowbreak 
\cS_K$. \par 
Now, for $m \in \N$, let $\ti{Y}^{\lmr}$ (resp. $\cY^{\lmr}$) 
be the $(m+1)$-fold fiber product of $\ti{Y}$ (resp. $\cY$) over 
$\ol{Y}$ (resp. $\cS$), let $U_{\ti{Y}^{\lmr}}$ be the inverse 
image of $U_Y$ by $\ti{Y}^{\lmr} \lra \ol{Y}$ and let 
$(U_{\ti{X}^{\lmr}},\ti{X}^{\lmr})$ be 
$(U_X \times_{U_Y} U_{\ti{Y}^{\lmr}}, \ol{X} \times_Y \ti{Y}^{\lmr})$. 
If we choose one of the projections $\ti{Y}^{\lmr} \lra \ti{Y}$ and if we 
put 
$(U_{\ti{X}^{\lmr,(\b)}},\ti{X}^{\lmr,(\b)}) :=  
(U_{\ti{X}^{(\b)}} \times_{U_{\ti{Y}}} U_{\ti{Y}^{\lmr}}, 
\ti{X}^{(\b)} \times_{\ti{Y}} \ti{Y}^{\lmr})$, the morphism 
$(U_{\ti{X}^{\lmr,(\b)}},\ti{X}^{\lmr,(\b)}) \lra (U_{\ti{X}^{\lmr}},
\ti{X}^{\lmr})$ is an $N$-truncated proper hypercovering and 
the morphism $\ti{X}^{\lmr,(n)} \lra \ti{Y}^{\lmr} \,(n \leq N)$ admits 
log structures (pull-back log structures from $\ti{X}^{(n)}, \ti{Y}$) 
which makes this morphism a proper log smooth integral morphism having 
log smooth parameter. So, by the argument using the spectral sequence 
like \eqref{spseq}, we see that 
$R^qf_{(U_{\ti{X}^{\lmr}}, \ti{X}^{\lmr})/\cY^{\lmr}, \rig *}\cE$ 
is a coherent $j^{\d}\cO_{]\ti{Y}^{\lmr}[_{\cY^{\lmr}}}$-module 
and $R^qf_{(U_{\ti{X}^{\lmr}}, \ti{X}^{\lmr})/\cY^{\lmr} \times_{\cS} 
\cY^{\lmr}, \rig *}\cE$ 
is a coherent $j^{\d}\cO_{]\ti{Y}^{\lmr}[_{\cY^{\lmr} \times_{\cS} 
\cY^{\lmr}}}$-module for 
each $m \in \N$. 
Then, by \cite[2.3.1]{tsuzuki3}, 
they naturally defines an overconvergent 
isocrystal $\cF^{\lmr}$ of the form 
$(R^qf_{(U_{\ti{X}^{\lmr}}, \ti{X}^{\lmr})/\cY^{\lmr}, \rig *}\cE,
\epsilon^{\lmr})$ on $(U_{\ti{Y}^{\lmr}},\ti{Y}^{\lmr})/\cS_K$. 
Again by \cite[2.3.1]{tsuzuki3}, $\cF^{\lmr}$ is compatible 
with respsect to $m$. So $\cF^{\lr}$ defines an object of 
$I^{\d}((U_{\ti{Y}^{\lr}}, \ti{Y}^{\lr})/\cS_K)$, where 
$I^{\d}((U_{\ti{Y}^{\lr}}, \ti{Y}^{\lr})/\cS_K)$ denotes the category 
of descent data with respect to 
$I^{\d}((U_{\ti{Y}^{\lmr}}, \ti{Y}^{\lmr})/\cS_K) \, (m=0,1,2)$. 
Now let us recall that the category of overconvergent isocrystals 
satisfies the descent property for strict etale hypercoverings 
(\cite[5.1]{shiho3}) and strict proper \v{C}ech hypercoverings (Proposition 
\ref{properdescent}). So 
it satisfies the descent property for \v{C}ech hypercoverings constructed 
from a morphism which is the composite of a strict etale covering 
and a strict proper covering. Hence we have an equivalence of categories 
\begin{equation}\label{des-et-pr}
I^{\d}((U_{\ti{Y}}, \ti{Y})/\cS_K) \os{=}{\lra}
I^{\d}((U_{\ti{Y}^{\lmr}}, \ti{Y}^{\lmr})/\cS_K), 
\end{equation} 
that is, $\cF^{\lr}$ descents to an overconvergent isocrystal 
on $(U_{\ti{Y}^{\lmr}}, \ti{Y}^{\lmr})/\cS$, which we denote by $\cF$. \par 
By definition, $\cF$ satisfies the required property when 
$(Z,\ol{Z},\cZ)$ is equal to $(U_{\ti{Y}^{\lmr}},
\ti{Y}\allowbreak{}^{\lmr},\allowbreak \cY^{\lmr}) \allowbreak \, (m \in \N)$. 
By the descent property \eqref{des-et-pr}, 
this property characterizes $\cF$. So we have 
the uniqueness of $\cF$. \par 
Finally we check that $\cF$ satisfies the required property. 
When a triple $(Z,\ol{Z},\cZ)$ as in the statement of the theorem is given, 
we have the diagram 
$$ 
(Z,\ol{Z}) \lla 
(X \times_{Y} Z, \ol{X} \times_{\ol{Y}} \ol{Z}) = 
(U_{\ti{X}} \times_{U_{\ti{Y}}} Z, \ti{X} \times_{\ti{Y}} \ol{Z}) \lla 
(U_{\ti{X}^{(\b)}} \times_{U_{\ti{Y}^{(\b)}}} Z, 
\ti{X}^{(\b)} \times_{\ti{Y}^{(\b)}} \ol{Z})
$$ 
and for each $n \leq N$, the morphism 
$(U_{\ti{X}^{(n)}} \times_{U_{\ti{Y}^{(n)}}} Z, 
\ti{X}^{(n)} \times_{\ti{Y}^{(n)}} \ol{Z}) \lra 
(Z,\ol{Z})$ admits log structure (the pull-back 
log structures from $\ti{X}^{(n)}, \ti{Y}$) 
which makes this morphism a proper log smooth integral morphism having 
log smooth parameter. So we have the coherence of 
$R^qf_{(X \times_{Y} Z, \ol{X} \times_{\ol{Y}} \ol{Z})/\cZ, \rig *}\cE$ and 
$R^qf_{(X \times_{Y} Z, \ol{X} \times_{\ol{Y}} \ol{Z})/\cZ \times_{\cS} \cZ, 
\rig *}\cE$, by the argument using spectral sequence like \eqref{spseq}. 
Then, by \cite[2.3.1]{tsuzuki3}, we see that the restriction 
of $\cF$ to $I^{\d}((Z,\ol{Z})/\cS_K)$ is given by 
$(R^qf_{(X \times_{Y} Z, \ol{X} \times_{\ol{Y}} \ol{Z})/\cZ, \rig *}\cE, 
\epsilon)$, which is as in the statement of the theorem. 
The functoriality of this expression can be shown in the same way as 
Theorem \ref{thm00} (see \cite[4.8]{shiho3}). So we are done. 
\end{pf} 

\begin{rem}\label{theorem1rem}
By looking at the proof carefully, we see the following: 
in Theorem \ref{theorem1}, we can take $U_Y$ to be open dense in $\ol{Y}$ 
when $Y \subseteq \ol{Y}$ is open dense. 
\end{rem} 

By using the above theorem, we can prove the following theorem, 
which gives the affirmative answer to Conjecture \ref{bconj2} 
(the version without Frobenius structure) in the case 
$\cS=\Spf V$ and $Y$ is smooth over $k$. 

\begin{thm}\label{theorem2}
Let $V$ be a complete discrete valuation ring of mixed characteristic 
with residue field $k$ and put $S := \Spec k, \cS := \Spf V$. 
Let us assume given the diagram 
\begin{equation}\label{diag-7-6}
(X,\ol{X}) \os{f}{\lra} (Y,\ol{Y}) \os{g}{\lra} (S,S), 
\end{equation}
where $f$ is a strict proper morphism of pairs 
such that $f|_X: X \lra Y$ is proper smooth, $Y$ is smooth over $k$ 
and $g$ is a morphism of pairs. 
Then, there exists 
a non-empty open set $U_Y$ of $Y$ and a diagram 
$$ \ol{Y} \os{g_Y}{\lla} \ti{Y} \os{i}{\hra} \cY $$
$($where $g_Y$ is a composition of an etale surjective map and a proper 
map and $i$ is a closed immersion into a $p$-adic formal scheme 
formally smooth over $\cS)$ such that, 
for an overconvergent isocrystal $\cE$ on $(X,\ol{X})/\cS_K$ 
$($which is automatically locally free$)$ and $q \in \N$, 
there exists 
the unique overconvergent isocrystal $\cF$ on $(Y,\ol{Y})/\cS_K$ 
satisfying the following condition$:$ For any triple $(Z,\ol{Z},\cZ)$ over 
$(Y \times_{\ol{Y}} \ti{Y}, \ti{Y}, \cY)$ satisfying either 
$\ol{Z} = Y \times_{\ol{Y}} \ol{Z}$ or $Z = U_Y \times_{\ol{Y}} \ol{Z}$ 
with $\cZ$ formally smooth over $\cY$, 
the restriction of $\cF$ to 
$I^{\d}((Z,\ol{Z})/\cS_K, \cZ)$ is given functorially by 
$(R^qf_{X  
\times_{Y} Z, \ol{X} \times_{\ol{Y}} \ol{Z})/\cZ, \rig *}\cE, 
\epsilon)$, where $\epsilon$ is an isomorphism 
{\small{ $$ 
p_2^*R^qf_{(X \times_{Y} Z, 
\ol{X} \times_{\ol{Y}} \ol{Z})/\cZ, \rig *}\cE 
\os{\sim}{\ra} 
R^qf_{(X\times_{Y} 
Z, \ol{X} \times_{\ol{Y}} \ol{Z})/\cZ \times_{\cS} \cZ, \rig *}\cE 
\os{\sim}{\leftarrow}
p_1^*R^qf_{(X 
\times_{Y} Z, \ol{X} \times_{\ol{Y}} \ol{Z})/\cZ, \rig *}\cE
$$ }}
$(p_i$ denotes the $i$-th projection 
$]\ol{Z}[_{\cZ \times_{\cS} \cZ} \lra ]\ol{Z}[_{\cZ})$. 
\end{thm} 

\begin{pf} 
First, let $\ol{Y}'$ be the closure of $Y$ in $\ol{Y}$ and put 
$\ol{X}':=\ol{Y}' \times_{\ol{Y}} \ol{X}$. Then we have the equivalences 
of categories 
$$ 
I^{\d}((X,\ol{X})/\cS_K) \cong I^{\d}((X,\ol{X}')/\cS_K), \,\,\,\,\, 
I^{\d}((Y,\ol{Y})/\cS_K) \cong I^{\d}((Y,\ol{Y}')/\cS_K). $$ 
So we may replace $\ol{X}, \ol{Y}$ by $\ol{X}', \ol{Y}'$ to prove the 
theorem, that is, we may assume that $Y$ is dense in $\ol{Y}$. \par 
By Theorem \ref{theorem1} and Remark \ref{theorem1rem}, 
there exists an open subset 
$U_Y \subseteq Y$ dense in $\ol{Y}$ and 
a diagram 
$$ \ol{Y} \os{g_Y}{\lla} \ti{Y} \os{i}{\hra} \cY $$
$($where $g_Y$ is a composition of an etale surjective map and a proper 
surjective map and $i$ is a closed immersion into a $p$-adic formal scheme 
formally smooth over $\cS$) and 
the unique overconvergent isocrystal $\cF_1$ on $(U_Y,\ol{Y})/\cS_K$ 
such that, for any triples $(Z,\ol{Z},\cZ)$ over 
$(Y \times_{\ol{Y}} \ti{Y}, \ti{Y},\cY)$ satisfying $Z = U_Y 
\times_{\ol{Y}} \ol{Z}$ with $\cZ$ formally smooth over $\cY$, 
the restriction of $\cF_1$ to 
$I^{\d}((Z,\ol{Z})/\cS_K, \cZ)$ is given functorially by 
$(R^qf_{(X 
\times_{Y} Z, \ol{X} \times_{\ol{Y}} \ol{Z})/\cZ, \rig *}\cE, 
\epsilon)$ as in the statement of the theorem. \par 
On the other hand, 
by \cite[5.14]{shiho3} and the descent 
property of $I^{\d}((Y,Y)/\allowbreak \cS_K)$ 
for etale and proper coverings, we see that there exists 
the unique overconvergent isocrystal $\cF_2$ on 
$(Y,Y)/\cS_K$ 
such that, for any triples $(Z,\ol{Z},\cZ)$ over 
$(Y \times_{\ol{Y}} \ti{Y}, \ti{Y},\allowbreak \cY)$ satisfying $\ol{Z} = Y 
\times_{\ol{Y}} \ol{Z}$ with $\cZ$ formally smooth over $\cY$, 
the restriction of $\cF_2$ to 
$I^{\d}((Z,\ol{Z})/\cS_K, \cZ)$ is given functorially by 
$(R^qf_{(X  
\times_{Y} Z, \ol{X} \times_{\ol{Y}} \ol{Z})/\cZ, \rig *}\cE, 
\epsilon)$ as in the statement of the theorem. 
(In fact, the unique existence of the restriction of $\cF_2$ 
to $I^{\d}((Y \times_{\ol{Y}} \ti{Y}, Y \times_{\ol{Y}} \ti{Y})/\cS_K)$ 
follows from \cite[5.14]{shiho3} and we can descent it to 
$\cF_2 \in I^{\d}((Y,Y)/K)$ by using the descent 
property.) \par 
Let us denote the restriction of $\cF_{i}\,(i=1,2)$ to 
$I^{\d}((U_Y,Y)/K)$ by $\cF'_i$. Then, for any 
triple $(Z,\ol{Z},\cZ)$ over 
$(Y \times_{\ol{Y}} \ti{Y}, \ti{Y},\cY)$ satisfying 
$Z = U_Y \times_{\ol{Y}} \ol{Z}$ and 
$\ol{Z} = Y \times_{\ol{Y}} \ol{Z}$ with $\cZ$ formally smooth over $\cY$, 
the restriction of $\cF'_i\,(i=1,2)$ to 
$I^{\d}((Z,\ol{Z})/\cS_K, \cZ)$ are both given functorially by 
$(R^qf_{(X 
\times_{Y} Z, \ol{X} \times_{\ol{Y}} \ol{Z})/\cZ, \rig *}\cE, 
\epsilon)$ as in the statement of the theorem. 
Moreover, the above condition 
charaterizes the overconvergent isocrystals $\cF'_i$ because of 
the descent property of $I^{\d}((U_Y,Y)/K)$ for 
\v{C}ech hypercoverings constructed 
from a morphism which is the composite of a strict etale covering 
and a strict strongly proper covering.
So we have the canonical isomorphism $\cF'_1 \cong \cF'_2$. Then, 
by \cite[5.3.7]{kedlaya2}, we see that there exists uniquely an overconvergent 
isocrystal $\cF$ on $(Y,\ol{Y})/\cS_K$ whose restriction to 
$I^{\d}((U,\ol{Y})/\cS_K)$ (resp. $I^{\d}((Y,Y)/\cS_K)$) is equal to 
$\cF_1$ (resp. $\cF_2$). One can see easily that $\cF$ satisfies the 
required property. The uniqueness of $\cF$ follows from that of 
$\cF_1$ and $\cF_2$. So we are done. 
\end{pf} 

Let us recall the following conjecture by Tsuzuki (\cite[1.2.1]{tsuzuki}): 

\begin{conj}[Tsuzuki]\label{tconj}
Let $X$ be a smooth separated $k$-scheme of finite type and 
let $X \hra \ol{X}$ be an open immersion such that $X$ is dense in 
$\ol{X}$. Then the restriction 
functor 
\begin{equation}\label{tconjfunctor}
I^{\d}((X,\ol{X})/\cS_K) \lra I^{\d}((X,X)/\cS_K) 
\end{equation}
is fully faithful. 
\end{conj} 

If we admit this conjecture, we can drop the condition `$Y$ is smooth' 
in Theorem \ref{theorem2}: 

\begin{thm}\label{theorem3}
Let us admit the validity of Conjecture \ref{tconj}. 
Let us put $S := \Spec k, \cS \allowbreak := \Spf V$ and 
let us assume given the diagram 
\begin{equation*}
(X,\ol{X}) \os{f}{\lra} (Y,\ol{Y}) \os{g}{\lra} (S,S), 
\end{equation*}
where $f$ is a strict proper morphism of pairs 
such that $f|_X: X \lra Y$ is proper smooth and $g$ is a morphism of pairs 
$($but $Y$ is not necessarily smooth$)$. 
Then the same conclusion as Theorem \ref{theorem2} holds. 
\end{thm} 

\begin{pf} 
We may assume that $Y$ is dense in $\ol{Y}$. 
By Theorem \ref{theorem1} and Remark \ref{theorem1rem}, 
we have an open subset 
$U \subseteq Y$ dense in $\ol{Y}$, 
proper surjective morphism 
$\ol{Y}^{\lor} \lra \ol{Y}$, an etale surjective morphism 
$\ti{Y} \lra \ol{Y}^{\lor}$, a closed immersion $\ti{Y} \hra \cY$ of 
$\ti{Y}$ into a formal $\cB$-scheme formally smooth over $\cS$ such that, 
for any overconvergent isocrystal $\cE$ on $(X,\ol{X})/\cS_K$ and $q \in \N$, 
there exists 
the unique overconvergent isocrystal $\cF_1$ on $(U,\ol{Y})/K$ satisfying 
the following condition: 
For any triple $(Z,\ol{Z},\cZ)$ over $(Y \times_{\ol{Y}} \ti{Y}, 
\ti{Y}, \cY)$ satisfying $Z = U \times_{\ol{Y}} \ol{Z}$ with 
$\cZ$ formally smooth over $\cY$, the restriction of $\cF_1$ to 
$I^{\d}((Z,\ol{Z})/\cS_K, \cZ)$ is given functorially by 
$(R^qf_{(X 
\times_{Y} Z, \ol{X} \times_{\ol{Y}} \ol{Z})/\cZ, \rig *}\cE, 
\epsilon)$ as in the statement of the theorem. 
By Theorem \ref{hc3} and the proof of Theorem \ref{theorem1}, 
we can take $\ol{Y}^{\lor}$ and $U$ in order that 
$\ol{Y}^{\lor}$ is smooth over $k$ and 
$U \times_{\ol{Y}} \ol{Y}^{\lor}$ is dense open in $\ol{Y}^{\lor}$. \par 
On the other hand, let $Y^{\lor}$ be the inverse image of 
$Y$ in $\ol{Y}^{\lor}$. Since $Y^{\lor}$ is smooth over $k$, 
there exists (by \cite[5.14]{shiho3} and Proposition \ref{properdescent}) 
the unique overconvergent isocrystal $\cF^{\lor}_2$ on 
$(Y^{\lor},Y^{\lor})/\cS_K$ 
such that, for any triple $(Z,\ol{Z},\cZ)$ over $(Y \times_{\ol{Y}} \ti{Y}, 
\ti{Y}, \cY)$ satisfying $\ol{Z} = Y \times_{\ol{Y}} \ol{Z}$ with 
$\cZ$ formally smooth over $\cY$, 
the restriction of $\cF^{\lor}_2$ to 
$I^{\d}((Z,\ol{Z})/\cS_K, \cZ)$ is given functorially by 
$(R^qf_{(X 
\times_{Y} Z, \ol{X} \times_{\ol{Y}} \ol{Z})/\cZ, \rig *}\cE, 
\epsilon)$ as in the statement of the theorem. 
For $m \in \N$, let $Y^{\lor,m}$ be the $m$-fold fiber product of 
$Y^{\lor}$ over $Y$. Then, since the morphism $\ol{X} \times_{\ol{Y}} 
Y^{\lor,m} \lra Y^{\lor,m}$ has log smooth parameter, we can prove 
(again by \cite[5.14]{shiho3} and Proposition \ref{properdescent}) 
the unique existence of overconvergent 
isocrystals $\cF_2^{\lor,m}$ on $(Y^{\lor,m},Y^{\lor,m})/\cS_K \, (m \in \N)$ 
which is compatible with respect to $m$. Hence, by descent property of 
overconvergent isocrystals with respect to the proper covering 
$Y^{\lor,m} \lra Y$, we see that the family $\{\cF^{\lor,m}_2\}_{m=1,2,3}$ 
descents to the overconvergent isocrystal on $(Y,Y)/\cS_K$, which we denote 
by $\cF_2$. \par 
Let us denote the restriction of $\cF_{i}\,(i=1,2)$ to 
$I^{\d}((U,Y)/\cS_K)$ by $\cF'_i$. Then, we can see the following in the same 
way as the proof of Theorem \ref{theorem2}: For any triple 
$(Z,\ol{Z},\cZ)$ over $(Y \times_{\ol{Y}} \ti{Y}, 
\ti{Y}, \cY)$ satisfying $Z = U \times_{\ol{Y}} \ol{Z}$ and 
$\ol{Z} = Y \times_{\ol{Y}} \ol{Z}$ with 
$\cZ$ formally smooth over $\cY$, 
the restriction of $\cF'_i\,(i=1,2)$ to 
$I^{\d}((Z,\ol{Z})/\cS_K, \cZ)$ are both given functorially by 
$(R^qf_{(X 
\times_{Y} Z, \ol{X} \times_{\ol{Y}} \ol{Z})/\cZ, \rig *}\cE, 
\epsilon)$ as in the statement of Theorem \ref{theorem2}. 
This implies that there exists the canonical isomorphism 
$\cF'_1 \cong \cF'_2$. So it suffices to prove that the functor 
\begin{equation}\label{fib}
\phi: 
I^{\d}((Y,\ol{Y})/\cS_K) \lra I^{\d}((U,\ol{Y})/\cS_K) 
\times_{I^{\d}((U,Y)/K)} I^{\d}((Y,Y)/\cS_K) 
\end{equation} 
induced by the restriction functors is an equivalence of categories. 
Note that we can take (by using the results in \cite{dejong}) 
a split proper hypercovering 
$\ol{Y} \lla \ol{Y}^{\lr}$ such that $\ol{Y}^{\lor}$ is given as above and 
that each $\ol{Y}^{\lmr}$ is smooth over $k$. Let $Y^{\lr}, U^{\lr}$ 
be the inverse image of $Y, U$ in $\ol{Y}^{\lr}$, respectively. Then, 
by proper descent 
(Proposition \ref{properdescent}) of overconvergent isocrystals and 
\cite[(3.3.4.2)]{saintdonat}, the restriction functor 
\begin{equation}\label{eqcat-pq}
I^{\d}((P,Q)/\cS_K) \lra I^{\d}((P^{\lr},Q^{\lr})/\cS_K)
\end{equation}
for $(P,Q)=(Y,\ol{Y}), (U,\ol{Y}), (Y,Y), (U,Y)$ are equivelences of 
categories, where $I^{\d}((P^{\lr}, \allowbreak Q^{\lr})/\cS_K)$ denotes the 
category of descent data with respect to $I^{\d}((P^{\lmr},Q^{\lmr})/\cS_K) \, 
\allowbreak (m\allowbreak =\allowbreak 0,1,2)$. 
For $m \in \N$, let us denote the functor 
$$
I^{\d}((Y^{\lmr},\ol{Y}^{\lmr})/K) \lra I^{\d}((U^{\lmr},\ol{Y}^{\lmr})/K) 
\times_{I^{\d}((U^{\lmr},Y^{\lmr})/K)} I^{\d}((Y^{\lmr},Y^{\lmr})/K) 
$$
by $\phi^{\lmr}$. 
Then, by the equivalence of categories \eqref{eqcat-pq}, it suffices to prove 
that $\phi^{\lor}$ is an equivalence of categories and 
that $\phi^{\lmr}$ is fully faithful for $m \in \N$ to prove the 
equivalence of $\phi$. \par 
First, $\phi^{\lor}$ is an equivalence of categories by 
\cite[5.3.7]{kedlaya2}, 
since $\ol{Y}^{\lor}$ is smooth over $k$ and $U^{\lor}$ is dense open 
in $\ol{Y}^{\lor}$. Next, let us prove the full-faithfulness of 
$m\geq 1$. Since $\ol{Y}^{\lmr}$ is smooth over $k$, we may consider 
connected componentwise. So we may assume that $\ol{Y}^{\lmr}$ is connected. 
Then, if $U^{\lmr}$ is non-empty, it is dense in $\ol{Y}^{\lmr}$. In this 
case, $\phi^{\lmr}$ is an equivalence of categories by 
\cite[5.3.7]{kedlaya2}. On the other hand, if $U^{\lmr}$ is empty, 
then $\phi^{\lmr}$ is fully-faithful by Conjecture \ref{tconj}. 
So $\phi^{\lmr}$ is fully faithful for $m \geq 1$, and 
the proof of the theorem is now finished. 
\end{pf} 

Next we give a result on Frobenius structure on the overconvergent 
isocrystal constructed in Theorems \ref{thm00}, 
\ref{theorem1}, \ref{theorem2}. \par 
Let $V$ be a complete discrete valuation ring of mixed characteristic 
with perfect residue field $k$, let $\pi$ be a uniformizer of $V$ and put 
$\cS := \Spf V$, $S := \Spec V/\pi V = B$. 
Let us fix an integer $q$ which is a power of $p$. For a scheme 
$X$ over $\F_p$, let $F_X$ be the $q$-th power Frobenius 
(the morphism induced by $q$-th power endomorphism of structure sheaf). 
Assume that we have an endomorphism $\sigma: \cS \lra \cS$ which 
lifts $F_S:S \lra S$ and fix it. Then, if we have a morphism of pairs 
$(X,\ol{X}) \os{f}{\lra} (S,S),$ 
we have the canonical functor
$$
F^*_X: I^{\d}((X,\ol{X})/\cS_K) \lra I^{\d}((X,\ol{X})/\cS_K) $$
induced by the morphisms $F_{\ol{X}}$ and $\sigma$. 
An overconvergent $F$-isocrystal on 
$(X,\ol{X})/\cS_K$) is defined to be a pair $(\cE,\alpha)$, where 
$\cE$ is 
an overconvergent isocrystal on $(X,\ol{X})/\cS_K$ and $\alpha$ is 
an isomorphism $F^*_X\cE \os{\sim}{\lra} \cE$ (which we call a Frobenius 
structure on $\cE$). 
Then we have the following: 

\begin{thm}\label{theorem4}
Let the notations be as in Theorem \ref{thm00}, Theorem \ref{theorem1} or 
Theorem \ref{theorem2}. 
Assume moreover that $\cS=\Spf V, S=\Spec k$ $(V$ is as above$)$ 
and that we have a fixed endomorphism 
$\sigma: \cS \lra \cS$ which lifts $F_S:S \lra S$. 
Then, if $\cE$ has a Frobenius structure, 
$\cF$ in the theorem has a canonical Frobenius structure induced from 
that on $\cE$. 
\end{thm} 

\begin{pf} 
Since the proof is quite similar to that of \cite[5.16]{shiho3}, we only 
give a sketch: we are reduced to proving that, for any $(Z,\ol{Z},\cZ)$ 
as in the statement of Theorem \ref{thm00}, 
\ref{theorem1} or \ref{theorem2} endowed with 
an endomorphism $\sigma_{\cZ}: \cZ \lra \cZ$ compatible with $F_{\ol{Z}}$ and 
$\sigma$, there exists a functorial Frobenius structure on 
the restriction $\cF_{\cZ}$ of $\cF$ to $I^{\d}((Z,\ol{Z})/\cS_K,\cZ)$. 
To simplify the notation, we replace $(Z,\ol{Z},\cZ), 
(X \times_Y Z, \ol{X} \times_{\ol{Y}} \ol{Z}), \sigma_{\cZ}$ by 
$(Y,\ol{Y},\cY), (X,\ol{X}), \sigma_{\cY}$, respectively. Then we have 
the canonical homomorphism 
$$ 
\sigma^*_{\cY}R^qf_{(X,\ol{X})/\cY, \rig *}\cE 
\os{F^*}{\lra} R^qf_{(X,\ol{X})/\cY, \rig *}F_X^*\cE
 \os{\alpha}{\lra} R^qf_{(X,\ol{X})/\cY, \rig *}\cE. $$
(where $F^*$ is induced by $F_{\ol{X}}^*$ and $\sigma_{\cY}$). 
We can prove that it is an isomorphism in exactly the same way as 
\cite[5.16]{shiho3}. 
\end{pf} 

Using Theorem \ref{theorem4}, we can prove the following theorem, 
which gives 
the affirmative answer to Conjecture \ref{bconj2} 
(the version with Frobenius structure) in the case 
$\cS=\Spf V$. (This is the `Frobenius version' of 
Theorem \ref{theorem3}, but we do not assume 
Conjecture \ref{tconj} here): 

\begin{thm}\label{theorem5}
Let $S := \Spec k, \cS := \Spf V$ and 
let us fix an endomorphism 
$\sigma: \cS \lra \cS$ which lifts $F_S:S \lra S$.
Let us assume given the diagram 
\begin{equation*}
(X,\ol{X}) \os{f}{\lra} (Y,\ol{Y}) \os{g}{\lra} (S,S), 
\end{equation*}
where $f$ is a strict proper morphism of pairs 
such that $f|_X: X \lra Y$ is proper smooth and $g$ is a morphism of pairs 
Then, for an overconvergent $F$-isocrystal $\cE$ on 
$(X,\ol{X})/\cS_K$ and 
$q \in \N$, the same conclusion as Theorem \ref{theorem2} 
holds and $\cF$ has a canonical Frobenius structure induced by that on $\cE$. 
\end{thm} 

\begin{pf} 
We may assume that $Y$ is dense in $\ol{Y}$. 
Let the notations be as in the proof of Theorem \ref{theorem3}. 
Then the argument in the proof of Theorems \ref{theorem3}, \ref{theorem4} 
shows that there exists a canonical Frobenius structure on $\cF_1$ and 
$\cF_2$ with the isomorphism $\cF'_1 \cong \cF'_2$ as overconvergent 
$F$-isocrystals. Let us denote the category of overconvergent 
$F$-isocrystals on a pair 
$(P,Q)$ by $\FI^{\d}((P,Q)/\cS_K)$. Then it suffices to 
prove that the restriction functor 
$$
\phi: \FI^{\d}((Y,\ol{Y})/\cS_K) \lra  
\FI^{\d}((U,\ol{Y})/\cS) \times_{\FI^{\d}((U,Y)/\cS_K)} 
\FI^{\d}((Y,Y)/\cS_K) 
$$
is an equivalence of categories. Note that, for a smooth separated $k$-scheme 
$Z$ and a dense 
open immersion $Z \hra \ol{Z}$, the restriction functor 
$$ \FI^{\d}((Z,\ol{Z})/\cS_K) \lra \FI^{\d}((Z,Z)/\cS_K) $$ 
is known to be fully faithful \cite{kedlaya3}, \cite[4.2.1]{kedlaya-val}. 
So we see that the functor 
\begin{align*}
\phi^{\lmr}: & \FI^{\d}((Y^{\lmr},\ol{Y}^{\lmr})/\cS_K) \\ 
& \lra 
\FI^{\d}((U^{\lmr},\ol{Y}^{\lmr})/\cS) \times_{\FI^{\d}((U^{\lmr},
Y^{\lmr})/\cS_K)} \FI^{\d}((Y^{\lmr},Y^{\lmr})/\cS_K) 
\end{align*}
is an equivalence of categories for $m=0$ by \cite[5.3.7]{kedlaya2} and 
it is fully-faithful in general. (See the proof of Theorem \ref{theorem3}. 
Here we do not use Conjecture \ref{tconj} because we use the above-mentioned 
result of Kedlaya instead of Conjecture \ref{tconj}.) 
From this fact, we can deduce that $\phi$ is an 
equivalence. So we are done. 
\end{pf} 


\end{document}